\documentclass[11pt]{article}

\usepackage{amssymb}
\usepackage{amstext}
\usepackage{amsmath}
\usepackage{amscd}
\usepackage{latexsym}
\usepackage[all]{xy}
\newfont{\gothique}{eufm10 scaled 1100}  %% gotic
\newcommand{\goth}[1]{\mbox{\gothique{#1}}}

\newcommand{\PP}{{\mathbb{P}}}
\newcommand{\CC}{{\bf{C}}}
\newcommand{\ZD}{{\cal{Z}}}
\newcommand{\JJ}{{\cal{J}}} 
\newcommand{\JG}{{\bf{J}}_{\Gamma}}

\newcommand{\OO}{{\cal O}}
\newcommand{\GS}{{ \cal G}}
\newcommand{\GT}{{\tilde{ \cal G}}}

\newcommand{\EZ}{Ext_{Z}^{1}}
\newcommand{\SE}{{\cal E}}
\newcommand{\SEE}{{\SE},[e]}
\newcommand{\EEF}{{\cal E}{\it{xt}}}
\newcommand{\ENDO}{{\cal E}{\it{nd}}}

\newcommand{\JA}{{\bf{J}}( X;L,d)} 
\newcommand{\JAA}{{\bf{J}}}

\newcommand{\XD}{X^{[d]}}
\newcommand{\HT}{{\bf \tilde{H}}}
\newcommand{\HH}{{\bf{H}}}
\newcommand{\ID}{{\cal{I}}}

\newcommand{\FIT}{{\bf{\tilde{ F}}}}
\newcommand{\FI}{{\bf F}}
\newcommand{\GA}{\Gamma}
\newcommand{\GAC}{\Gamma^{(0)}_{conf}}
\newcommand{\GAO}{{\stackrel{\circ}{\Gamma^r_d}}}
\newcommand{\GG}{{\Gamma}^r_d}

\newcommand{\JAB}{{\bf \breve{J}}} 
\newcommand{\HOM}{{\cal H}\it{om}}
\newcommand{\HB}{\bf H}
\newcommand{\QB}{{\bf q}}
\newcommand{\QT}{{\bf {\tilde{q}}}}

\newcommand{\HO}{H^0(\OO_}
\newcommand{\HOZE}{H^0(\OO_{Z_e})}

\newcommand{\SLA}{\tilde{\GS}}
\newcommand{\ZA}{([Z],[\alpha])}
\newcommand{\FT}{\tilde{\cal F}}
\newcommand{\FF}{{\cal F}}
\newcommand{\CS}{C^r (L,d)}
\newcommand{\TE}{{\bf \Theta}} 
\newcommand{\RI}{\cite{[R1]}}
\newcommand{\CO}{Conf_d (X)}
\newcommand{\CSA}{C^r_{adm} (L,d)}
\newcommand{\LG}{l_{\GA}}
\newcommand{\LAG}{\GS_{\mbox{\unboldmath$\GA$}} }
\newcommand{\LAGT}{\GT_{\mbox{\unboldmath$\GA$}} }
\newcommand{\CG}{{\cal C}_{\mbox{\unboldmath$\GA$}} }
\newcommand{\BM}{\boldmath}
\newcommand{\UB}{\unboldmath}
\newcommand{\BEN}{\begin{equation}}
\newcommand{\EEN}{\end{equation}}
\newcommand{\GAB}{\breve{\GA}}
\newcommand{\JABG}{\JAB_{\GA}}
\newcommand{\PF}{{}^{\prime}{\FI}} %august 3%
\newcommand{\FL}{{\cal{FL}}}
\newcommand{\FLA}{\mbox{\BM${\FL}$}_{\overrightarrow{h^{\prime}}_{\GA}}}
\newcommand{\FLAO}{\mbox{\BM${\FL}$}_{\overleftarrow{h^{\prime}}_{\GA}}}
\newcommand{\TPI}{T_{\pi}}

\newcommand{\FTB}{\FT^{\prime \bullet}}
\newcommand{\FTP}{\FT^{\prime}}
\newcommand{\PG}{p_{\GA}}
\newcommand{\OPG}{{}^{op}p_{\GA}}
\newcommand{\NI}{{\cal N}^{\prime}}
\newcommand{\DE}{d^{\prime}_{\GA}}
\newcommand{\LIG}{\mbox{\BM$G^{\prime}_{\mbox{\UB$\GA$}}$}}
\newcommand{\BK}{{\bf k}}
\newcommand{\BO}{{\bf o}}
\newcommand{\IG}{\mbox{\BM$Gr_{\mbox{\UB$\GA$}}$}}
\newcommand{\LGR}{\mbox{\BM$G^{\prime \mbox{rel}}_{\mbox{\UB$\GA$}}$}}
\newtheorem{thm}{Theorem}[section]
\newtheorem{lem}[thm]{Lemma}
\newtheorem{pro}[thm]{Proposition}
\newtheorem{cor}[thm]{Corollary}
\newtheorem{rem}[thm]{Remark}

\newtheorem{defi}[thm]{Definition}  %[section]
\newtheorem{cl}[thm]{Claim}
\newtheorem{pro-defi}[thm]{Proposition-Definition}
\newlength{\myskip}
\newenvironment{pf}{
     \addvspace{\myskip}  

     \noindent {\it Proof.$\, $}}
     {$\Box$

     \addvspace{\myskip}
     }

\makeatletter
\renewcommand{\@seccntformat}[1]{\S \/ {\csname the#1\endcsname}\hspace{0.5em}}
\makeatother

\title{NONABELIAN JACOBIAN OF SMOOTH PROJECTIVE SURFACES AND REPRESENTATION THEORY  \\
          }
\author{Igor Reider}

\begin{document}
\bibliographystyle{amsplain}

\maketitle

\setcounter{section}{-1}
\numberwithin{equation}{section}

\begin{abstract}
The paper studies representation theoretic aspects of a nonabelian version of the Jacobian for
a smooth complex projective surface $X$ introduced in \cite{[R1]}.
The sheaf of reductive Lie algebras 
$\mbox{\BM$\GT$}$ associated to the nonabelian Jacobian is determined and its Lie algebraic properties are explicitly related to the geometry of configurations of points on $X$.
 In particular, it is shown
that the subsheaf of centers of $\mbox{\BM$\GT$}$ determines a distinguished decomposition of configurations into the disjoint union of 
subconfigurations. 
Furthermore, it is shown how to use
${\bf sl_2}$-subalgebras associated to certain nilpotent elements of 
$\mbox{\BM$\GT$}$  to write equations defining configurations of $X$
in appropriate projective spaces.

The same nilpotent elements are used to establish a relation of the nonabelian Jacobian with such fundamental objects in the representation theory 
as nilpotent orbits, Springer resolution and Springer fibres of simple Lie algebras of type ${\bf A_n}$, for appropriate values of $n$.
This leads to a construction of distinguished collections of objects in the category of representations of symmetric groups as well as 
in the category of perverse sheaves on the appropriate Hilbert schemes of points of $X$. 
Hence two ways of categorifying the second Chern class of vector bundles of rank 2 on smooth projective surfaces.

We also give a `loop' version of the above construction by relating the nonabelian Jacobian to the Infinite Grassmannians of 
simple Lie groups of type $ {\bf SL_n ( C )}$, for appropriate values of $n$. 
This gives, via the geometric version of the Satake isomorphism, a distinguished collection of irreducible representations of the Langlands dual groups
thus indicating  a relation of the nonabelian Jacobian to the Langlands duality on smooth projective surfaces.      
\end{abstract}
\tableofcontents

\section{Introduction}
It is hard to overestimate the role of the Jacobian in the theory of smooth complex projective curves. 
The celebrated theorem of Torelli says that a curve of genus $\geq2$ is determined, up to isomorphism, by its Jacobian and its theta-divisor.
Virtually all projective geometric features of a curve can be extracted from its Jacobian. But the Jacobian of a curve has its intrinsic importance and beauty. 
It is enough to recall that it is a principally polarized abelian variety with an incredibly rich and beautiful theory of theta-functions.

When one turns to higher dimensional projective varieties one quickly discovers that a comparable theory does not exist.
However, in the end of 1960's Griffiths initiated a far reaching theory of Variation of Hodge structure (abbreviated in the sequel by VHS).
Some of its goals include generalizations of the Theorem of Torelli and a study of algebraic cycles. From Griffiths' theory
emerges a substitute for the Jacobian - Griffiths period domain. This is an open subset of a certain flag variety (factored out by the action of a certain discrete group). In particular, the theory of VHS continues to have strong ties with the theory of Lie groups.
Furthermore, a VHS comes with the period map and Griffiths suggested to view its derivative as a substitute for the principal polarization of the classical Jacobian (see \cite{[G]} for an overview and references therein).

In \cite{[R1]} we proposed a new version of Jacobian for a smooth complex projective surface $X$. We suggested to call it nonabelian
Jacobian for the simple reason that it parametrizes a distinguished family of rank 2 bundles on $X$. More precisely, similar to
  its classical counterpart, our nonabelian Jacobian is, on the one hand, related to the moduli stack of torsion free sheaves\footnote{sheaves are of rank $2$, 
contrary to the classical situation of line bundles.} on $X$, and, on the other hand, to
the Hilbert scheme of points on $X$. 
It also carries a distinguished divisor which can be viewed as a nonabelian analogue of the classical theta-divisor.
 But a new feature of our Jacobian is that it is also related to the Griffiths' ideas of the VHS and period maps. One of the consequences of this is an appearance
  of a sheaf of reductive 
Lie algebras canonically attached to our Jacobian. 
This can be viewed as an analogue of the Lie algebraic structure of the classical Jacobian.

In this paper we undertake a study of this sheaf of Lie algebras. 
Our considerations are naturally divided into two parts:
\\
\\
1) establish a dictionary between the properties of the sheaf of reductive Lie algebras attached to our Jacobian and geometric properties of 
$X$.
\\
\\
2) Use the representation theory to define interesting objects (e.g. sheaves, complexes of sheaves) which can serve as new invariants of vector
bundles on $X$ as well as invariants of the surface itself.
\\
\\
\indent
For the first part we are able to uncover: 
\\
\\
a) a precise relationship between the center of the reductive Lie algebras in question and canonical decompositions
of configurations of points on $X$ into disjoint union of subconfigurations,
\\
b) how to use particular ${\bf sl_2}$-subalgebras of our reductive Lie algebras to gain an insight into the geometry of configurations of points on $X$.
\\
\\
\indent
For the second part we show how to use the sheaf of reductive Lie algebras associated to our nonabelian Jacobian to attach to $X$: 
\\
\\
1) a distinguished
collection of objects in the category of representations of symmetric groups,
\\
2) a distinguished
collection of objects in the category of perverse sheaves on the appropriate
Hilbert schemes of points on $X$,
\\
3) a distinguished
collection of irreducible representations of the Langlands dual group
${}^L {\bf  SL}_n (\CC)\linebreak ={\bf PGL}_n (\CC)$, for appropriate values of $n$.   
\\
\\
\indent
These results come from the fact that our Jacobian connects in a natural way to such fundamental objects in geometric representation theory as the
Springer resolution of the nilpotent cone of simple Lie algebras (of type ${\bf A_n}$), Springer fibres, loop algebras and Infinite Grassmannians.

In the rest of this introduction, following a brief summary of \cite{[R1]},
we give a more detailed account of the results of this paper.
\subsection {Nonabelian Jacobian $\JA$ (a summary of \cite{[R1]}).}
A new version of the Jacobian for smooth projective surfaces was proposed in \cite{[R1]}. Our construction is based on viewing the 
Jacobian of a smooth projective curve as the parameter space for line bundles with a fixed Chern class. We suggested that for a smooth
projective variety $X$ of dimension $n\geq2$, the Jacobian could be the parameter space of a distinguished family of vector bundles of rank 
$n= dim_{\bf C} X$ with fixed Chern invariants. Using this analogy for a smooth projective surface $X$, we have constructed the scheme
$\JA$, whose closed points are pairs $({\SE},[e])$, where $\SE$ is a torsion free sheaf of rank $2$ on $X$ with Chern invariants 
$c_1 ({\SE})=L$ and $c_2({\SE}) = d$, where $L$ is a suitably fixed divisor on $X$ and $d$ is a fixed positive integer, and where
$[e]$ is the homothety class of a global section $e$ of $\SE$, whose zero-locus
$Z_e = (e=0)$ is a subscheme of codimension $2$ (equivalently, dimension $0$) of $X$. We suggested to call $\JA$ a nonabelian Jacobian of
$X$ (of type $(L,d)$).

By definition $\JA$ is a scheme over the Hilbert scheme $X^{[d]}$, the scheme parametrizing the subschemes $Z$ of $X$ having dimension zero and
length $d$. The natural morphism
\begin{equation}\label{pi}
\pi: \JA \longrightarrow X^{[d]} 
\end{equation}
sends a pair $(\SEE)$ to the point $[Z_e] \in X^{[d]}$ corresponding to the subscheme $Z_e = (e=0)$ of $X$.

As in the classical case, $\JA$, over a suitable subscheme of $X^{[d]}$, comes with a distinguished Cartier divisor
$\TE(X;L,d)$, whose closed points parametrize pairs $(\SEE)$, where the sheaf $\SE$ is not locally free.
But there is also a new phenomenon: $\JA$ carries a natural structure resembling a VHS \`a la Griffiths. More precisely,
for every point $(\SEE) \in \JA$, one has a distinguished filtration on $\HOZE$
\begin{equation}\label{filtH}
0= \tilde{\cal H}_{0} (\SEE) \subset \tilde{\cal H}_{-1} (\SEE) \subset \ldots \subset \tilde{\cal H}_{-l_{Z_e} -1} (\SEE) = \HOZE\,,
\end{equation}
where the integer $l_{Z_e}$ is intrinsically associated to $Z_e$.

Furthermore, if $(\SEE)$ is in a certain constructible subset $\JAB$ of $\JA$, the filtration (\ref{filtH}) splits.
By this we mean that $\HOZE$ admits a distinguished direct sum decomposition
 \begin{equation}\label{odH}
\HOZE = \bigoplus_{p=0}^{l_{Z_e}} {\HB}^p (\SEE)
\end{equation} 
with a natural identification 
$$
{\HB}^p (\SEE) \cong \tilde{\cal H}_{-(p+1)} (\SEE) / \tilde{\cal H}_{-p}(\SEE)\,,
$$
for $p=0, \ldots, l_{Z_e} $. This direct sum decomposition could be thought of as some kind of periods for the points in 
$\JAB$. Thus our nonabelian Jacobian possesses features of the classical Jacobian as well as a period map in the spirit of 
Griffiths theory of VHS.

The decomposition (\ref{odH}) together with the obvious ring structure on $\HOZE$ gives rise to a reductive Lie subalgebra
{\boldmath${\tilde{\GS}}$}$(\SEE)$ of ${\bf{gl}} (\HOZE)$. By varying $(\SEE)$ in $\JAB$ we obtain the sheaf
 {\boldmath${\tilde{\GS}}$}$(X;L,d)$ of reductive Lie algebras naturally associated to $\JA$. 
This could be viewed as a generalization of the Lie algebraic nature of the classical Jacobian.

One of the features of the sheaf {\boldmath$\tilde{\GS}$}$(X;L,d)$ is that it gives rise to a natural family of Higgs structures
in the sense of Simpson, \cite{[S]}. The parameter space $ H$ of this family turns out to be a toric (singular) Fano variety whose hyperplane sections are, 
in general, singular Calabi-Yau varieties. This $H$ could be viewed as a nonabelian Albanese of $\JA$.

It should be pointed out that $H$ depends only on the properties of the sheaf {\boldmath$\tilde{\GS}$}$(X;L,d)$ of reductive Lie algebras
and the decomposition (\ref{odH}). All this can be encapsulated in the the following trivalent graph
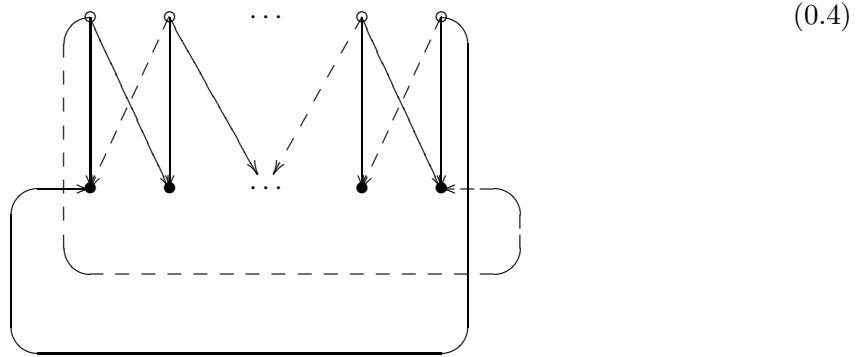
\begin{equation}\label{graph}
\xymatrix{
&*=0{\circ} \ar[2,0] \ar[2,1] \ar@{--{>}}`l[d]                               
                               `[3,1]                               
                               `[2,5]                               
                                `[2,4]
                                [2,4]
                              &
                               *=0{\circ}\ar[2,0] \ar[rdd] \ar@{--{>}}[ddl]&
\cdots&*=0{\circ} \ar[2,0] \ar@{--{>}}[ldd] \ar[2,1]&
*=0{\circ} \ar[2,0] \ar@{--{>}}[2,-1] \ar`r[d] 
                              `[4,-2]
                              `[2,-5]                              
                              `[2,-4]
                              [2,-4]&\\
& & & & & & \\
&*=0{\bullet}&*=0{ \bullet}&\cdots&*=0{\bullet}&*=0{\bullet}&\\
& & & & & &\\
& & & & & &}
\end{equation}
where the vertical levels represent the first $l_{Z_e}$ summands of the decomposition (\ref{odH}) and the slanted arrows represent certain degree
$\pm 1$ operators which are among the generators of {\boldmath${\tilde{\GS}}$}$(X;L,d)$.
\\
\indent
The features of $\JA$ enumerated above show that our Jacobian relates in a natural way to 
\\
- Lie algebras and their representations (the sheaf of reductive Lie algebras {\boldmath{$\tilde{\GS}$}}$(X;L,d)$)
\\
-toric geometry and Calabi-Yau varieties (the nonabelian Albanese $H$)
\\
-low dimensional topology (trivalent graph (\ref{graph})).
\\
Being such a multifaceted object it seems to us that $\JA$ is worthy of a serious study.
         \\
\\
\indent
In this paper we undertake a study of the Lie algebraic aspect of our Jacobian. In the following subsections of the introduction we summarize
the key results of this paper.
\subsection{The center of the Lie algebra {\boldmath{$\SLA$}}$(\SEE)$ and geometry of $Z_e$.}
We determine the reductive algebras {\boldmath$\SLA $}$(\SEE)$ attached to points of the Jacobian $\JAB$. It turns out  that the center of these algebras completely determines 
the Lie algebra {\boldmath{$\SLA$}}$(\SEE)$  and is related to the geometry of the zero-locus $Z_e =(e=0)$ associated to $(\SEE )\in \JAB$. More precisely, we show
\begin{thm} \label{tcd}
The zero locus $Z_e =(e=0)$ decomposes into the disjoint union
 \begin{equation}\label{cd}
Z_e = \bigcup_{i=1}^{\nu} Z_{e}^{(i)}\,,
\end {equation}
where $\nu$ is the dimension of the center of the Lie algebra {\boldmath$\SLA$}$(\SEE)$ attached to $(\SEE) \in \JAB$.
Furthermore, the Lie algebra {\boldmath$\SLA$}$(\SEE)$ and hence, its center act on the subspace $\HT_{-\LG} (\SEE)$ of the filtration
of $\HOZE$ in (\ref{filtH}). This action of the center determines the weight decomposition
$$
\HT_{-\LG} (\SEE) =\bigoplus^{\nu}_{i=1} V_i (\SEE)
$$
which possesses the following properties:
\begin{enumerate}
\item[1)]
$\HO{Z_e^{(i)}}) \cong V_i (\SEE) \cdot \HO{Z_e})$,
\item[2)]
 one has a natural isomorphism
\begin{equation}\label{Lid}
\mbox{\boldmath$\tilde{\GS} $}(\SEE) \cong \bigoplus^{\nu}_{i=1} {\bf gl}(V_i (\SEE))\,.
\end {equation}
\end{enumerate}
\end{thm}

This result establishes a precise dictionary between the decomposition of the Lie algebra
{\boldmath$\SLA$}$(\SEE)$ into the direct sum of matrix algebras and the geometric decomposition of $Z$ into the disjoint union of subschemes
in (\ref{cd}).
\\
\\
\indent
 It turns out that the Lie algebra {\boldmath{$\SLA$}}$ (\SEE)$ also controls the properties of the derivative of the period map associated to
$\JAB$.
\begin{thm}\label{infTor}
The derivative of the period map attached to $\JAB$ is injective precisely at the points $(\SEE)$ for which
{\boldmath{$\SLA$}}$ (\SEE) \cong {\bf gl_{d^{\prime}} (C)}$, where $d^{\prime}= dim (\HT_{-\LG} (\SEE))$ and where $\HT_{-\LG} (\SEE)$ is as in the filtration in (\ref{filtH}).
\end{thm}

This is a version of the Infinitesimal Torelli Theorem for $\JA$. Thus in our story the Infinitesimal Torelli property, i.e.
the injectivity of the differential of the period map, has a precise geometric meaning:
it fails exactly when the decomposition (\ref{cd}) is non-trivial. 

These results constitute a semisimple aspect of the representation theory of {\boldmath{$\SLA$}}$(\SEE)$ in a sense that it takes into account the action on the space $\HOZE$ of
 the center of {\boldmath{$\SLA$}}$ (\SEE)$, which is
composed of semisimple elements.
There is also a nilpotent aspect  which is much more involved.
\subsection{Nilpotent aspect of {\boldmath{$\GS$}}$ (\SEE)$}
Let {\boldmath{$\GS$}}$ (\SEE)$ be the semisimple part of {\boldmath{$\SLA$}}$(\SEE)$.
From the construction of the Lie algebra 
{\boldmath{$\GS$}}$ (\SEE)$
it follows that we can attach a nilpotent element $D^{+}(v)$ of {\boldmath{$\GS$}}$ (\SEE)$ with every vertical\footnote{throughout the paper `vertical' means in 
the direction of the fibres of the projection
$\pi$ in (\ref{pi}).}\label{footnote:vv} tangent vector $v$ of
$\JAB$ at a point $(\SEE) \in \JAB$. On the diagrammatic representation (\ref{graph}) the elements $D^{+}(v)$ are depicted by the right-handed arrows.
As $v$ runs through the space $T_{\pi} (\SEE)$ of the vertical tangent vectors of $\JAB$ at $(\SEE)$ we obtain the linear map
\begin{equation}\label{D+}
D_{(\SEE)}^{+} : T_{\pi} (\SEE) \longrightarrow \mbox{\boldmath${\cal N} ({\GS}$}(\SEE))
\end{equation}
into the nilpotent cone {\boldmath${\cal N} ({\GS}$}$(\SEE))$ of {\boldmath{$\GS$}}$ (\SEE)$.

From the well-known fact that
{\boldmath${\cal N} ({\GS}$}$(\SEE))$ is partitioned into a finite set of nilpotent orbits we deduce that
the map $D_{(\SEE)}^{+} $ assigns to $(\SEE)$ a finite collection of nilpotent orbits of {\boldmath${\cal N} ({\GS}$}$(\SEE))$.
These are the orbits intersecting the image of $D_{(\SEE)}^{+} $. Varying $(\SEE)$ in the suitable subvarieties of $\JAB$ we deduce the following.
\begin{thm}\label{nilo}
The Jacobian $\JA$ gives rise to a finite collection
${\cal V}$ of quasi-projective subvarieties of $\XD$ such that every $\Gamma \in {\cal V}$ determines a finite collection $O(\Gamma )$ of nilpotent orbits in 
${\bf sl_{d^{\prime}_{\Gamma}} (C)}$, where  $d_{\Gamma}^{\prime} \leq d$ is an integer intrinsically associated to $\Gamma$.
\end{thm}

 Recalling that nilpotent orbits in $\bf sl_n (C)$ are parametrized by the set of partitions $P_n$ of $n$, the above  result can be rephrased
by saying that every $\Gamma$ in ${\cal V}$ distinguishes a finite collection $P(\Gamma)$ of partitions of $d_{\Gamma}^{\prime}$.
Since partitions of $n$ also parametrize isomorphism classes of irreducible representations of the symmetric group
$S_n$ we obtain the following equivalent version of Theorem \ref{nilo}.
\begin{thm}\label{irrSn}
The Jacobian $\JA$ gives rise to a finite collection
${\cal V}$ of quasi-projective subvarieties of $\XD$ such that every $\Gamma \in {\cal V}$
 determines a finite collection $R_{d^{\prime}_{\Gamma}}(\Gamma )$ of irreducible representations of the
symmetric group $S_{d^{\prime}_{\Gamma}}$, where  $d^{\prime}_{\Gamma} \leq d$ is an integer intrinsically associated to $\Gamma$.
\end{thm}

One way to express this result is by saying that the Jacobian $\JA$ elevates a single topological invariant $d$, the degree of the second Chern
class of sheaves parametrized by certain subvarieties of $\JA$, to the level of modules of symmetric groups.
Thus our Jacobian gives rise to new invariants with values in the categories of modules of symmetric groups.

But there is more to it. The partitions distinguished by $\JA$ contain a great deal of geometry of subschemes parametrized by $ \Gamma$'s in 
Theorem \ref{nilo}. In down to earth terms one can say that the partitions picked out by points $(\SEE)$ of $\JAB$ yield equations defining
the image of $Z_e$ under certain morphisms into appropriate projective spaces. 

The process of obtaining these equations is somewhat evocative
of the classical method of Petri (see \cite{[Mu]} for an overview). However, the essential ingredient in our approach is representation theoretic. It turns on the use of ${\bf sl_2}$-subalgebras of {\boldmath{$\GS$}}$ (\SEE)$ associated to the nilpotent elements
$D_{(\SEE)}^{+} (v)$, the values of the map $D_{(\SEE)}^{+}$ in (\ref{D+}). The operator $D_{(\SEE)}^{+} (v)$ in our considerations plays the role of the operator $L$ in the Lefschetz decomposition in the Hodge theory. Completing it to an ${\bf sl_2}$-subalgebra of 
{\boldmath{$\GS$}}$ (\SEE)$ in an appropriate way and considering  its representation on 
$\HOZE$, gives a sort of Lefschetz decomposition of $\HOZE$. This combined with the orthogonal decomposition 
in (\ref{odH}) yields a bigrading of $\HOZE$ thus revealing a much finer structure than the initial grading (\ref{odH}).

Once this bigrading is in place, writing down the equations defining $Z_e$ in a certain projective space is rather straightforward. 
This is discussed in details in \S\ref{sec-equations}. The equations themselves can be complicated and, in general, not very illuminating. 
What is essential in our approach is that this complicated set of equations is encoded in an appropriate 
${\bf sl_2}$-decomposition of $\HOZE$. This in turn can be neatly ``packaged" in the properties of the partitions singled out by the points
$(\SEE)$ of $\JAB$ ``polarized" by operators $D_{(\SEE)}^{+} (v)$, with $v$ varying in $T_{\pi} (\SEE)$ as in (\ref{D+}).  

To summarize, one can say that the nilpotent aspect of the representation theory of {\boldmath{$\GS$}}$(X;L,d)$ provides new geometric insights as well as new invariants of the representation theoretic nature.

 This turns out to be only a part of the story. In fact, we can go further by relating $\JA$ to the category of perverse sheaves 
on $\XD$.
\begin{thm}\label{ps}
The Jacobian $\JA$ determines a finite collection 
{\boldmath${\cal P}$}$(X;L,d)$ of perverse sheaves on $\XD$.
These perverse sheaves are parametrized by pairs
$ (\Gamma,\lambda)$, where $\Gamma$ is a subvariety in ${\cal V}$ as in Theorem \ref{nilo} and $\lambda$ is a partition in $P(\Gamma)$.
\end{thm}

This result subsumes two previous theorems since  the perverse sheaves 
${\cal C}(\Gamma,\lambda)$ in {\boldmath${\cal P}$}$(X;L,d)$ have the following properties:
\\
a) ${\cal C}(\Gamma,\lambda)$ is the Intersection Cohomology complex
$IC(\Gamma, {\cal L}_{\lambda})$ associated to the local system ${\cal L}_{\lambda}$ on $\Gamma$.
\\
b) The local system ${\cal L}_{\lambda}$ corresponds to a representation 
\begin{equation}\label{repgf}
\rho_{\Gamma,\lambda} : \pi_1 ( \Gamma,[Z]) \longrightarrow Aut(H^{\bullet}(B_{\lambda}, {\bf C}))
\end{equation}
of the fundamental group $\pi_1 ( \Gamma,[Z])$ of $\Gamma$ based at a point $[Z]\in \Gamma$
and where 
$H^{\bullet}(B_{\lambda}, {\bf C})$ is the cohomology ring (with coefficients in $\bf C$) of a Springer fibre\footnote{
a Springer fibre $B_{\lambda}$ is a fibre
of the Springer resolution
\begin{equation*}\label{Spres}
\sigma: \mbox{\boldmath$\tilde{\cal N}$} \longrightarrow \mbox{\boldmath${\cal N}$}({\bf sl}_{d_{\Gamma}^{\prime}} ({\bf C}))
\end{equation*}
of the nilpotent cone $\mbox{\boldmath${\cal N}$}({\bf sl}_{d_{\Gamma}^{\prime}} ({\bf C}))$ of ${\bf sl}_{d_{\Gamma}^{\prime}} (\bf C)$
and where a fibre $B_{\lambda}$ is taken over the nilpotent orbit {\boldmath$O_{\mbox{\unboldmath$\lambda$}}$}
in $\mbox{\boldmath${\cal N}$}({\bf sl}_{d_{\Gamma}^{\prime}} ({\bf C}))$ corresponding to a partition $\lambda$ of $d_{\Gamma}^{\prime}$.}
$B_{\lambda}$ over the nilpotent orbit {\boldmath$O_{\mbox{\unboldmath$\lambda$}}$} of ${\bf sl}_{d_{\Gamma}^{\prime}} (\bf C)$
corresponding to the partition $\lambda$.
\\
c) The representation $\rho_{\Gamma,\lambda}$ admits the following factorization
\begin{equation}\label{Sprep}
\rho_{\Gamma,\lambda} : \pi_1 ( \Gamma,[Z]) \stackrel{\rho^{\prime}}{\longrightarrow} S_{d_{\Gamma}^{\prime}} \stackrel{sp_{\lambda}}{\longrightarrow} 
Aut(H^{\bullet}(B_{\lambda}, {\bf C}))\,,
\end{equation}
where 
$S_{d_{\Gamma}^{\prime}} \stackrel{sp_{\lambda}}{\longrightarrow} Aut(H^{\bullet}(B_{\lambda}, {\bf C}))$ 
is the Springer representation of the Weyl group
$W=S_{d_{\Gamma}^{\prime}}$ of ${\bf sl}_{d_{\Gamma}^{\prime}} (\bf C)$ on the cohomology of a Springer fibre $B_{\lambda}$.
\\
\\
\indent
Using the fact that the category of perverse sheaves is semisimple, the collection {\boldmath${\cal P}$}$(X;L,d)$ gives rise to a distinguished
collection, denoted $\mbox{\BM$C$}(X;L,d)$, of {\it irreducible} perverse sheaves on $\XD$. This in turn defines the abelian category
${\cal A} (X;L,d)$ whose objects are isomorphic to finite direct sums of complexes of the form ${\cal C}[n]$, where 
${\cal C} \in \mbox{\BM$C$}(X;L,d)$ and $n\in {\bf Z}$. 

This construction parallels the construction of local systems on the classical Jacobian.
Recall that if $J(C)$ is the Jacobian of a smooth projective curve $C$, then isomorphism classes of irreducible local systems on $J(C)$ are
parametrized by the group of characters $Hom (H_1 (J(C)), \CC^{\times})$. So we suggest to view the collection of irreducible perverse sheaves
$\mbox{\BM$C$}(X;L,d)$ as a nonabelian analogue of the group of characters of the classical Jacobian, while
the abelian category ${\cal A} (X;L,d)$ could be envisaged as an analogue of the group-ring of $Hom (H_1 (J(C)), \CC^{\times})$.

Though objects of  ${\cal A} (X;L,d)$ are complexes of sheaves on the Hilbert scheme $\XD$, they really descend from $\JA$ and one of the ways to
remember this is the following
\begin{thm}\label{exp-int}
Let $\stackrel{\circ}{\JAA} (X;L,d) = \JA \setminus \mbox{\BM$\Theta$} (X;L,d)$ be the complement of the theta-divisor
$\mbox{\BM$\Theta$} (X;L,d)$ in $\JA$ and let
${\cal T}^{\ast}_{\stackrel{\circ}{\JAA} (X;L,d)/ {\XD}}$ be the sheaf of relative differentials of 
$\stackrel{\circ}{\JAA} (X;L,d)$ over $\XD$. 
Then there is a natural map
$$
exp\left(\int\right): H^0 ({\cal T}^{\ast}_{\stackrel{\circ}{\JAA} (X;L,d)/ {\XD}} ) \longrightarrow  {\cal A} (X;L,d)\,.
$$
\end{thm}

The map in the above theorem could be viewed as a reincarnation of the classical map
$$
H^0 ({\cal T}^{\ast}_{J(C)} ) \longrightarrow Hom (H_1 (J(C)), \CC^{\times})\,,
$$
where ${\cal T}^{\ast}_{J(C)}$ is the cotangent bundle of $J(C)$. This map sends a holomorphic $1$-form $\omega$ on $J(C)$ to the exponential of the linear functional
$$
\int (\omega) : H_1 (J(C)) \longrightarrow \CC
$$
given by integrating $\omega$ over $1$-cycles on $J(C)$ (the notation `$exp(\int)$' in Theorem \ref{exp-int} is an allusion to this classical
map).
\\
\\
\indent
Relations of the Hilbert schemes of points of surfaces to partitions is not new. Notably, Haiman's work on the Macdonald positivity conjecture,
\cite{[Hai]}, makes an essential use of such a relation. The same goes for an appearance of perverse sheaves on $\XD$: the work of 
G\"{o}ttsche and Soergel, \cite{[Go-So]}, uses the decomposition theorem of \cite{[BBD]} for the direct image of the Intersection
cohomology complex $IC(\XD)$ under the Hilbert-Chow morphism to compute  the cohomology of Hilbert schemes. In both of these works the partitions appear from the outset because the authors exploit the points of the Hilbert scheme corresponding to the zero-dimensional subschemes $Z$ of $X$, where the points in $Z$ are allowed to collide according to the pattern determined by partitions.
 In our constructions it is essential to work over the open
part $Conf_d (X)$ of $\XD$, parametrizing configurations of $d$ distinct points of $X$. So there are no partitions seen on the level of the Hilbert scheme.
The partitions become visible only on the Jacobian $\JA$ via the Lie algebraic invariants attached to it. One can say that our constructions turn a configuration of distinct points with no interesting structure on it into a dynamical object. The dynamics is given by certain
linear operators acting on the space of complex valued functions on a configuration. In particular, the operators 
$D^{+} (v)$ obtained as values of the morphism $D^{+}$ in (\ref{D+}) give rise to the ``propagations" and ``collisions" in the direct sum decomposition (\ref{odH}). This is not an actual, physical, collision of points in a configuration but rather algebro-geometric constraints for
a configuration to lie on hypersurfaces in the appropriate projective spaces. The partitions attached to the nilpotent operators $D^{+} (v)$
can be viewed as a combinatorial (or representation theoretic) measure of this phenomenon, while the perverse sheaves in Theorem \ref{ps} could be envisaged as its categorical
manifestation.
\subsection{From $\JA$ to Affine Lie algebras.}\label{AffLie}
One of the major developments of the last 15 years about the Hilbert schemes of points of complex projective surfaces is the discovery of
Grojnowski and Nakajima of the action of affine Lie algebras on the direct sum of the cohomology rings (with rational coefficients) of the
Hilbert schemes $X^{[n]} (n\in {\bf Z_{+}}) $ (see \cite{[N]} and the references therein for more details). However, as Nakajima points out
 in the Introduction of \cite{[N]}, until now one has no good explanation of this phenomenon. In this subsection we explain
how our Jacobian can be used to address this problem.

It is clear that formally we can replace the Lie algebra {\boldmath$\GS$}$(\SEE)$ attached to a point $(\SEE) \in \JA$ by its loop Lie algebra
{\boldmath$\GS$}$(\SEE)[z^{-1},z]$, where $z$ is a formal variable. However, there is a more natural and explicit reason for appearance of loop
Lie algebras in our story. To explain this we recall that the Lie algebra
{\boldmath$\tilde{\GS}$}$(\SEE)$ is obtained as follows.

For every $h$ in the summand $\HH^0(\SEE)$ of the decomposition (\ref{odH}), we consider the operator $D(h)$ of multiplication by $h$ in the ring
$\HOZE$. Decomposing this operator according to the direct sum in (\ref{odH}) yields a triangular decomposition
\begin{equation}\label{tri-d}
D(h)= D^{-}(h) +D^0(h) + D^{+}(h)\,,
\end{equation}
where $D^{\pm}(h)$ are linear operators of degree $\pm 1$ with respect to the grading in (\ref{odH}) and $D^0(h)$ is a grading preserving operator. 
In particular, the operators $D^{+}(h)$, for $h\in \HH^0(\SEE)$, are essentially the same as the values of the morphism in
(\ref{D+}), due to the canonical identification of the relative tangent space $T_{\pi} (\SEE)$ with a codimension one subspace of
$\HH^0(\SEE)$.

 It is quite natural and immediate to turn (\ref{tri-d}) into a loop
\begin{equation}\label{l-tri-d}
D(h,z)= z^{-1} D^{-}(h) +D^0(h) + zD^{+}(h)\,,
\end{equation} 
where $z$ is a formal parameter. Morally, this natural one-parameter deformation of the multiplication in $\HOZE$ is behind the following
loop version
of the map (\ref{D+}):
\begin{equation}\label{lD+}
LD^{+}_{(\SEE)} :\, \stackrel{\circ}{T}_{\pi} (\SEE) \longrightarrow Gr(\mbox{\boldmath$\GS$}(\SEE))\,,
\end{equation}
where $Gr(\mbox{\boldmath$\GS$}(\SEE))$ is the loop or Infinite Grassmannian of the semisimple Lie algebra $\mbox{\boldmath$\GS$}(\SEE)$
and $\stackrel{\circ}{T}_{\pi} (\SEE)$ is an appropriate Zariski open subset of the vertical tangent space ${T}_{\pi} (\SEE)$ of $\JA$ at 
$(\SEE)$ .
This gives the following `loop' version of Theorem \ref{nilo}
\begin{thm}\label{lo}
The Jacobian $\JA$ gives rise to a finite collection $\cal V$ (the same as in Theorem \ref{nilo}) of subvarieties $\Gamma$ of
$\XD$.
Every such $\Gamma$ determines a finite collection 
$LO(\Gamma)$ of orbits of the Infinite Grassmannian
$Gr({\bf SL}_{d^{\prime}_{\Gamma}} ({\bf C}))$
of ${\bf SL}_{d^{\prime}_{\Gamma}}({\bf C})$, 
where $d^{\prime}_{\Gamma}$ is the same as in Theorem \ref{nilo}.
\end{thm}

Taking the Intersection Cohomology complexes
$IC({O})$ of the orbits $O$ in $LO(\Gamma)$, for every $\Gamma$ in $\cal V$, we pass to the category of perverse sheaves on 
$Gr({\bf SL}_{d^{\prime}_{\Gamma}} ({\bf C}))$. A beautiful and profound result of Ginzburg, \cite{[Gi]}, and Mirkovi\v{c} and Vilonen, \cite{[M-V]},
which establishes an equivalence between the category of perverse sheaves (subject to a certain equivariance condition) on the Infinite Grassmannian $Gr({\bf G})$ of a semisimple
Lie group ${\bf G}$ and the category of finite dimensional representations of the Langlands dual group $\bf {}^L G$ of $\bf G$, gives a Langlands dual
version of Theorem \ref{nilo}.
\begin{thm}\label{LD}
For every subvariety $\Gamma$ in $\cal V$ in Theorem \ref{lo} the Jacobian
$\JA$ determines a finite collection ${}^L R(\Gamma)$ of irreducible representations of the Langlands dual
group
${}^L {\bf SL}_{d^{\prime}_{\Gamma}}({\bf C}) = {\bf PGL}_{d^{\prime}_{\Gamma}}({\bf C})$.
\end{thm}

In retrospect a connection of our Jacobian with the Langlands duality could have been foreseen. After all, the nature of
$\JA$ as the moduli space of pairs $(\SEE)$ resembles the moduli space of pairs of Drinfeld in \cite{[Dr]}.
The fundamental difference is that the groups
${\bf SL}_{d^{\prime}_{\Gamma}}({\bf C})$ and their Langlands duals in our story have nothing to do with the structure group
($\bf GL_2 (C)$) of bundles parametrized by $\JA$. These groups rather reflect the geometric underpinnings of our construction
related to the Hilbert scheme $\XD$. Noting this difference, we also point out one of the key features of $\JA$:
\\
\\
{\it
it transforms the vertical vector fields of $\JA$ (i.e. sections of the relative tangent sheaf
${\cal T}_{\pi} ={\cal T}_{ \JA / \XD }$) to perverse sheaves on $\XD$. }
\\
\\
This feature is essentially the map in Theorem \ref{exp-int} and it can be viewed as a ``tangent" version of Grothendieck's ``functions-faisceaux dictionnaire",
 which plays an important role in a reformulation of the classical, number theoretic, Langlands correspondence into the geometric one (see \cite{[Fr]},
 for an excellent introduction to the subject of the geometric Langlands program).

\subsection{Concluding remarks and speculations}

The results of the paper show that the Lie algebraic aspects of our Jacobian are useful in addressing various issues related to algebro-geometric properties of
configurations of points on surfaces. It also enables us to attach to the degree of the second Chern class of vector bundles such objects as irreducible
representations of symmetric groups and perverse sheaves of the representation theoretic origin. In fact, we believe that the tools developed in the paper allow one to transfer virtually any object/invariant
of the geometric representation theory to the realm of smooth projective surfaces. For example, one should be able to have a version of Theorem \ref{irrSn}, where
the representations of the symmetric groups are replaced by the representations of the corresponding Hecke algebras as well as Affine Hecke algebras.

To our mind all these invariants fit into a sort of `secondary' type invariants for vector bundles in the sense of Bott and Chern 
in \cite{[B-C]}.
 Indeed, our construction begins by replacing the second Chern class of a bundle 
$\SE$ (of rank 2)
 by its geometric
realization, i.e. the zero-locus $Z$ of a suitable global section $e$ of $\SE$. This is followed by a distinguished orthogonal decomposition (\ref{odH}) of
the space of functions $\HO{Z})$ on $Z$. The decomposition gives rise to the Lie subalgebra {\boldmath{$\GT$}}$(\SEE)$ of ${\bf gl}(\HO{Z}))$ which is intrinsically associated to the pair
$(\SEE)$. This Lie subalgebra could be viewed as the `secondary' structure Lie algebra associated to $\SE$. While the structure group 
(${\bf GL_2 (C)}$) with
its Lie algebra provide the topological invariants of $\SE$, i.e. its Chern classes, the secondary structure Lie algebra detects various algebro-geometric properties of the subscheme $Z$.
 For example, Theorem \ref{tcd} can be interpreted as a statement of reduction of the secondary structure Lie algebra to a proper Lie subalgebra of
${\bf gl}(\HO{Z}))$ (see (\ref{Lid})).  A geometric significance of such a reduction is the decomposition of $Z$ in (\ref{cd}). 
Furthermore, if the structure group and its Lie algebra
yield the Chern invariants of $\SE$ by evaluating the basic structure group-invariant polynomials on a curvature form of $\SE$, it is plausible to expect that our
 secondary Lie algebra
should provide many more representation theoretic invariants of $(\SEE)$, which would reflect  properties of geometric representatives of the Chern invariants of $\SE$.
 Other theorems stated in the introduction could be viewed as a confirmation of this heuristic reasoning.

Theorem \ref{irrSn} and Theorem \ref{ps} could also be viewed as two kinds of categorifications of the second Chern class of rank 2 vector bundles on projective surfaces.
The latter result and the tools developed to obtain it suggest that there might be a categorification of the representation of affine Lie algebras
on the direct sum of the cohomology rings of the Hilbert schemes
discovered by I.Grojnowski and H.Nakajima (see the discussion in \S \ref{AffLie})

The results of \S \ref{AffLie} indicate a relation of our Jacobian to the Langlands duality. On the other hand it is conceptually sound to suggest that a
formulation of the geometric Langlands program for higher dimensional varieties could involve correspondences in the middle dimension.\footnote{
what we have in mind here is that
correspondences in the middle dimension could be taken as a geometric substitute for the Galois side of the Langlands correspondence.}
Now
the very idea of the Jacobian as a tool to study correspondences goes back to A.Weil (see \cite{[W]}).
In fact, one of our main motivations for introducing and studying $\JA$ was to study correspondences in the case of projective surfaces. 
Thus what emerges from our considerations is
the following triangular relation
\begin{equation}\label{tri}
\xymatrix{
&\JA \ar[dl] \ar[dr]&\\
*\txt{Correspondences\\ of X} \ar[rr] & &*\txt{Langlands Duality}} 
\end{equation}
A precise discussion of these interrelations will appear elsewhere but we hope that the results and tools developed in this paper will convince the reader
that the nonabelian Jacobian $\JA$ exhibits strong ties with the base of the above triangle.
\subsection{Organization of the paper}
There is a number of different topics discussed in the paper and we would like to summarize here how they fit together in our exposition.

To begin with it is not very realistic to make this paper self-contained since it draws heavily on the results of 
\RI. However, for the convenience of the reader \S1 is devoted to a concise summary of the main properties of our nonabelian Jacobian obtained in that paper. This is also a place
 to introduce the main notation and conventions used throughout the paper. Thus it is with \S2 that this paper truly begins.
The essential results here are Lemma \ref{HTl=triv} and its geometric realization in Corollary \ref{cor-Fpr}. These results are of technical nature and are in preparation
 for the determination of the Lie algebras attached to points of $\JA$.

In \S3 these Lie algebras are explicitly determined. This is done in two stages: 
\\
-in \S\ref{Center} we consider the center of the Lie algebras in question; the geometric consequences of this study are given in 
Corollary \ref{Z-c-dec};
\\
-in \S\ref{sec-ses} we determine the semisimple part of the Lie algebras attached to points of $\JA$: the main technical result here is
Proposition \ref{ga-dec}. 
\\
A combination of these two stages constitutes the results of Theorem \ref{tcd} of the Introduction.

In \S\ref{sec-periods} we switch to a more geometric point of view on our constructions by defining the period maps for our Jacobian.
We show that the period maps satisfy Griffiths transversality condition (Proposition \ref{pro-Grif-trans}) and compute their differentials  in terms of the operators $D^{\pm} (h)$ of the triangular decomposition in (\ref{tri-d}). This gives a purely algebraic formulas to compute the derivatives of our period maps (Lemma \ref{der=d+-}, Proposition \ref{pro-der=c(d)}) and links the geometry of the periods maps with the Lie algebraic considerations of the previous sections.

In \S\ref{sec-Torelli} we define Torelli property for our period maps and show that it is entirely controlled by the center of the Lie algebras
attached to points of $\JA$ (Corollary \ref{InfTorelli=s}, Theorem \ref{Torelli=InfTor}).

Next three sections are devoted to ${\bf sl_2}$-subalgebras associated to the operators $D^{\pm} (h)$ of the triangular decomposition in
(\ref{tri-d}).

In \S\ref{const-fib} we consider ${\bf sl_2}$-subalgebras associated to the operators  $D^{+} (h)$. This gives rise to bigraded structures on 
$\HOZE $ in (\ref{odH}). The main properties of these bigradings and the action of $D^{+} (h)$ are given in
Proposition \ref{pro-dgr}. In \S\ref{const-sh} we give a sheaf version of the above structures.

In \S\ref{sec-sl2bis} we consider the adjoint action of the ${\bf sl_2}$-subalgebras in \S\ref{sec-sl2} on the sheaf of Lie algebras attached to $\JA$. This results in a bigraded structure of the Lie algebras attached to points of $\JA$. The properties of this bigrading
can be found in Lemma \ref{range-p,q} and in Proposition-Definition \ref{pro-gr-g-bigr}.

In \S\ref{sec-inv} we change from operator $D^{+} (h)$ to $D^{-} (h)$ and consider ${\bf sl_2}$-subalgebras associated to
$D^{-} (h)$. The formalism is of course the same and the main issue here is the interaction of the two structures.
 In Proposition \ref{pro-two-filt} and Corollary \ref{cor-two-bigr} it is shown how the two ${\bf sl_2}$-structures are related.
The result is reminiscent of the Hodge-Riemann bilinear relations  in Hodge theory.

In \S\ref{sec-strat} we return to geometric considerations. In particular, we show how to use  ${\bf sl_2}$-subalgebras studied
in previous sections to define a stratification of the relative tangent sheaf of $\JA$. The resulting strata are indexed by
certain upper triangular, integer-valued matrices which we call multiplicity matrices (Definition \ref{mult-mat-def}, 
Proposition \ref{pro-strat}) or, equivalently, by partitions associated to the nilpotent operators $D^{+} (h)$ (Proposition \ref{pro-strat1}).

\S\ref{sec-equations} is devoted to applications of the theory built so far to various algebro-geometric questions concerning configurations of points on $X$.

 In \S\S\ref{sec-sl2-basis}- \ref{sec-conf-eqns} we present a general method of
using ${\bf sl_2}$-subalgebras considered in \S\ref{sec-sl2} to obtain equations of hypersurfaces cutting out configurations in an appropriate projective space.
In \S\ref{sec-K3} the general method is applied to a particular case : complete intersections on a $K3$-surface. In this case everything can be computed quite explicitly. In particular, one obtains a complete list of very simple quadratic hypersurfaces (of rank $\leq 4$) cutting out
complete intersections (see Proposition \ref{quadrics}). This gives a {\it hyperplane section} version of Mark Green's theorem
on quadrics of rank 4 in the ideal of a canonical curve in \cite{[Gr]}.

In \S9.6 the ${\bf sl_2}$-subalgebras considered in \S\ref{sec-inv} are put to use to study geometry of configurations of points on $X$
with respect to the adjoint linear system
 $\left| L+ K_X \right|$. Our considerations show how the partition associated to the nilpotent
operator $D^{-} (h)$ in (\ref{tri-d}) determines a special subvariety in $\PP(H^0(\OO_X (L+ K_X))^{\ast})$, passing through the image of a configuration under the morphism
defined by $\left| L+ K_X \right|$. This is Theorem \ref{italian-gen} which generalizes a well-known classical result saying that
$d$ points ($d\geq 4$) in general position in the projective space $\PP^{d-3}$ lie on a rational normal curve. 

In \S\ref{sec-RT-const} we return to general considerations with the intention to use nilpotent elements
$D^{+} (h)$ in a more conceptual way. This leads to a relation of $\JA$ to the nilpotent cone and the Springer resolution of simple Lie algebras of type
${\bf sl_n}$. The main results in \S10.2 are Proposition \ref{pro-T-orbits} and Theorem \ref{nilo1} (which is equivalent to Theorem \ref{nilo}
of the Introduction).

In \S\ref{sec-ps} the Springer resolution and Springer fibres are used to construct perverse sheaves on the Hilbert scheme $\XD$ 
(Theorem \ref{thm-IC}). This yields the collection $\mbox{\BM${\cal P}$}(X;L,d)$ of perverse sheaves on $\XD$ as in Theorem \ref{ps}
of the Introduction.

In \S\ref{A(XLd)} the collection $\mbox{\BM${\cal P}$}(X;L,d)$ is put to use to construct the abelian category ${\cal A}(X;L,d)$
appearing in Theorem \ref{exp-int}. The relation of relative differentials of $\JA$ with objects of ${\cal A}(X;L,d)$
(the map $exp(\int)$ in Theorem \ref{exp-int}) is given in Theorem \ref{thm-cot-ps-map} (see also Proposition \ref{tan-ps} and
Remark \ref{analogy}).

In \S\ref{sec-LD} a relation of $\JA$ and the Infinite Grassmannian of type ${\bf SL_n (\CC)}$ is established (Proposition \ref{d+loop}).
This leads to Theorem \ref{lo} and Theorem \ref{LD} of the Introduction (stated respectively as Proposition \ref{LO(Gam)} and 
Proposition \ref{LD1}).
\\
\\
{\bf Acknowledgments.} It is a pleasure to thank Vladimir Roubtsov for his unflagging interest to this work. Our thanks go
to the referee of \RI\, who also suggested in his report a possible connection of our Jacobian with perverse sheaves.

\section{Nonabelian Jacobian $\JA$: main properties}

In this section we introduce the main objects of our study and recall the main results of \cite{[R1]}.
\subsection{Construction of $\JA$.}

By analogy with the classical Jacobian of a smooth projective curve our Jacobian is supposed to be the parameter space of a 
certain distinguished family of torsion free sheaves of rank 2 over $X$, having the Chern invariants $(L,d)$. Its formal definition
is as follows.

One starts with the Hilbert scheme $\XD$ of closed zero-dimensional subschemes of $X$ having length $d$.
Over $\XD$ there is the universal scheme $\ZD$ of such subschemes
\begin{equation}\label{uc}
\xymatrix{
&\ZD \ar[dl]_{p_1} \ar[dr]^{p_2} \ar@{^{(}{-}{>}}[r]& X \times \XD\\
X& &\XD &}
\end{equation}
where $p_i (i=1,2)$ is the restriction to $\ZD$ of the projections $pr_i,\,\,i=1,2,$ of the Cartesian product $X \times \XD$
onto the corresponding factor. For a point $\xi \in \XD$, the fibre $p^{\ast}_2 (\xi)$ is isomorphic via $p_1$ with the subscheme $Z_{\xi}$ of $X$ corresponding to $\xi$, i.e.
\begin{equation}\label{fc}
 Z_{\xi} = p_{1 \ast} (p^{\ast}_2 (\xi) ).
\end{equation}
In the sequel we often make no distinction between $Z_{\xi}$  
and the fibre  
$
p^{\ast}_2 (\xi) 
$
itself. If $Z$ is a closed subscheme of dimension zero and length $d$, then $[Z]$ will denote the corresponding point in the Hilbert scheme
$
\XD
$.

The next step is to fix a line bundle
$\OO_X (L)$ corresponding to a divisor $L$ on $X$.
It will be assumed throughout the paper that $\OO_X (L)$
satisfies the following condition
\begin{equation}\label{vc}
H^0 (\OO_X (-L)) =H^1 (\OO_X (-L)) =0.
\end{equation}
We are aiming at geometric applications, where the divisor $L$ will
be sufficiently positive (e.g. $L$ is ample), so the above condition is quite natural.

Once a divisor $L$ and a positive integer $d$ are fixed we 
 consider the following morphism of sheaves on $\XD$
\begin{equation}\label{morph}
\xymatrix{
\HO{X} (L+K_X)) \otimes \OO_{\XD} \ar[r]^(.55){\rho}& p_{2 \ast}\big{(} p^{\ast}_1 \OO_X (L+K_X) \big{)}. }
\end{equation}

We define
\begin{equation}\label{dJac}
\JA := {\bf Proj}(S^{\bullet} coker \rho),
\end{equation}
where $S^{\bullet} coker \rho$ is the symmetric algebra of $coker \rho$.
By definition
$\JA$
comes with the natural projection
\begin{equation}\label{pi1}
\pi: \JA \longrightarrow \XD
\end{equation}
and the  invertible sheaf 
$
\OO_{\JA} (1)
$
such that the direct image
\begin{equation}\label{O1}
\pi_{\ast}  \OO_{\JA} (1) = coker \rho
\end{equation}
(when $X,L$ and $d$ are fixed and no ambiguity is likely, we will omit these parameters in the notation for the Jacobian and simply write 
$\JAA$ instead of $\JA$).

Observe that the set of closed points of the fibre of  $\pi$ over a point
$[Z]$ in $\XD$ is naturally homeomorphic to the projective space
$\PP(H^1 (\ID_{Z} (L+K_X))^{\ast})$.
By Serre duality on $X$
\begin{equation}\label{SD}
H^1 (\ID_{Z} (L+K_X))^{\ast} = Ext^1(\ID_{Z}(L+K_X), \OO(K_X))=Ext^1(\ID_{Z}(L), \OO_X).
\end{equation}
To simplify the notations the last space will be denoted by $\EZ$ throughout the paper. Thus the set of closed points of 
$\JA$ is in one to one correspondence with the set of pairs
$([Z],[\alpha])$, where 
$[Z] \in \XD$ 
and
$[\alpha] \in \PP(\EZ)$.
Alternatively, a pair $\ZA$ can be thought of as the pair
$(\SEE)$,
where
$\SE$ is the torsion free sheaf sitting in the middle of the extension sequence defined by the class 
$\alpha$ in $\EZ$
\begin{equation}\label{exs}
\xymatrix@1{
0\ar[r]& {\OO_X} \ar[r]& {\SE}\ar[r]&{\ID_Z} (L) \ar[r] & 0 }
\end{equation}
and where $[e]$ is the point in the projective space 
$\PP(H^0 (\SE))$
corresponding to the image of 
$H^0 (\OO_X)$
under the monomorphism in (\ref{exs}). Thus closed points of $\JA$
parametrize the set of pairs
$(\SEE)$,
where 
$\SE$ is a torsion free sheaf on $X$ having rank $2$ and the Chern invariants
$(L,d)$,
and $[e]$ is a one-dimensional subspace of 
$H^0 (\SE)$, whose generator $e$ is a global section of 
$\SE$
with the scheme of zeros
$Z_e =(e=0)$ having dimension $0$.
As it was explained in \cite{[R1]}, p.439, our Jacobian $\JA$
is the moduli stack of such pairs. From this description it follows that we have a morphism of stacks
\begin{equation}\label{morst}
\xymatrix@1{
{h:\,\,\JA} \ar[r]& {\bf M}_X (2,L,d), }
\end{equation}
where $ {\bf M}_X (2,L,d) $ is the moduli stack of  torsion free sheaves on $X$ having rank $2$, fixed determinant $\OO_X (L)$ and the second Chern class
of degree $d$. This morphism sends a pair
$(\SEE)$ in $\JA$ to the point $[\SE] \in {\bf M}_X (2,L,d) $, corresponding to the sheaf $\SE$.
In particular, the fibre of $h$ over a point 
$[\SE]$ in ${\bf M}_X (2,L,d) $ is the Zariski open subset
$U_{\SE}$ of 
$\PP(H^0 (\SE))$
parametrizing sections of 
$\SE$ (up to a non-zero scalar multiple) having $0$-dimensional locus of zeros.

Putting together (\ref{pi1}) and (\ref{morst}) we obtain the following diagram
\begin{equation}\label{HJM}
\xymatrix{
&{\JA} \ar[dl]_{\pi} \ar[dr]^h& \\
{\XD}& & {\bf M}_X (2,L,d)  }
\end{equation}

Thus $\JA$ can be thought of as a kind of thickening of $\XD$ and ${\bf M}_X (2,L,d) $.
In both cases the thickening is obtained by inserting over the points, of either the Hilbert scheme
$\XD$ or the moduli stack ${\bf M}_X (2,L,d) $, of rational fibres:
the projective space
$\PP(\EZ)$ over a point $[Z] \in \XD$
and the Zariski open subset 
$U_{\SE}$ of $\PP(H^0 (\SE))$ over a point $[\SE]$ in 
${\bf M}_X (2,L,d) $.

Another, equivalent, way of saying this is that through every point
$(\SEE)$ of $\JA$ pass two rational subvarieties, 
$U_{\SE}$ and $\PP(Ext^1_{Z_e})$, which are respectively the fibre of $h$ over $[\SE]$
and the fibre of $\pi$ over $[Z_e]$. These subvarieties have the following geometric meaning:
\\
- the subvariety $U_{\SE}$ is the space of rationally equivalent geometric realizations of the second Chern class of $\SE$,
\\
-the subvariety $\PP(Ext^1_{Z_e})$ is the space of natural deformations of the pair $(\SEE)$.

\subsection{A stratification of $\JA$.}
The Hilbert scheme $\XD$ acquires a distinguished stratification by the degeneracy loci of the morphism
$\rho$ in (\ref{morph}). Namely, set 
$D^{r} (L,d)$ to be the subscheme of zeros of the exterior power
$\wedge^{d-r} \rho$ of $\rho$ of degree $(d-r)$. Set-theoretically
$D^{r} (L,d)$ is the subset
\begin{equation}\label{grd}
\{\xi \in \XD \mid {\rm{dim}} (coker(\rho (\xi))) \geq r+1 \}.
\end{equation}
Denote by 
$\GG (L)$ the subscheme $D^{r} (L,d)$ taken with its reduced structure, i.e. set-theoretically 
$\GG (L)$ is the same as the set in (\ref{grd}), but as a subscheme of $\XD$ it is defined by the radical of the ideal sheaf of
$D^{r} (L,d)$. In particular, all $\GG(L)$'s are reduced closed subschemes of $\XD$ and they give a stratification of $\XD$:
\begin{equation}\label{strH}
\XD \supset \GA^0_d (L) \supset \GA^1_d (L) \supset \ldots \supset \GG(L) \supset \GA^{r+1}_d (L) \supset \GA^{d-1}_d (L)
\supset \GA^{d}_d =\emptyset.
\end{equation}
Denote by $\stackrel{\circ}{\Gamma^r_d} (L)$ the complement
$\GG(L) \setminus \GA^{r+1}_d (L)$.
This is a Zariski open subset of $\GG (L)$.
If 
$\GAO (L) \neq \emptyset$,
then it is well-known that $\GA^{r+1}_d (L)$ is contained in the singular locus of $\GG$ (see e.g. \cite{[ACGH]}, Ch2).
In particular, the smooth part 
$reg(\GG(L))$ of $\GG(L)$ is equal to 
$reg (\GAO (L))$, the smooth part of 
$\GAO (L)$.

The stratification (\ref{strH}) can be lifted via the projection $\pi$ in (\ref{pi1}) to define the stratification of $\JA$:
\begin{equation}\label{strJ}
\JA=\JAA^0 \supset \JAA^1 \supset \ldots \supset \JAA^r \supset \JAA^{r+1} \supset \ldots
\end{equation}
where 
$\JAA^r = {\bf Proj}(S^{\bullet} coker(\rho) \otimes \OO_{\GG(L)}  )$.
In particular,
the stratum
\begin{equation}\label{JO}
{\stackrel{\circ}{\JAA^r}} = \JAA^r \setminus \JAA^{r+1}
\end{equation}
is a $\PP^r$-bundle over $\GAO (L)$. 

\subsection{A nonabelian theta-divisor $\TE(X;L,d)$.}
For a closed subscheme $Z \subset X$ with $[Z] \in \GAO (L)$, the integer
\begin{equation}\label{is}
\delta (L,Z):= h^1(\ID_Z (L+K_X)) =r+1
\end{equation}
is called the index of $L$-speciality of $Z$. Following Tyurin, \cite{[Ty]} (see also \cite{[R1]}, Definition 1.1), we define
the notion of $L$-stability.
\begin{defi}\label{L-s}
A zero-dimensional subscheme $Z$ of $X$ is called $L$-stable iff
$$
\delta (L,Z^{\prime}) < \delta (L,Z)\,,
$$
for any proper subscheme $Z^{\prime}$ of $Z$.
\end{defi}

Denote by 
${{}^s}\! {\XD_{lci}}$
the subscheme of $\XD$ parametrizing $L$-stable subschemes which are local complete intersections (lci).
The two conditions are open so it is a Zariski open subset of $\XD$.

Let
$\JAA_{{{}^s} \!{\XD_{lci}}} = \pi^{-1} ({{}^s}\! {\XD_{lci}})$
be the part of $\JA$ lying over ${{}^s}\!{ \XD_{lci}}$.
It was shown in \cite{[R1]}, \S 1.2, that 
$\JAA_{{{}^s}\!{ \XD_{lci}}}$
carries a distinguished Cartier divisor denoted
 $\TE(X;L,d)$
and called the theta-divisor of 
$\JA$.
The closed points of 
 $\TE(X;L,d)$
parametrize pairs
$(\SEE)$, 
where the sheaf $\SE$ is not locally free.
Scheme-theoretically,
$\TE(X;L,d)$
is the subscheme of zeros of a distinguished section of the invertible sheaf
$\OO_{\JA} (d) \otimes \pi^{\ast} {\cal L}$,
for some line bundle ${\cal L}$ on $\XD$  (see \cite{[R1]}, (1.19)).
In particular,
the fibre
$\TE_{Z}$
of $\TE(X;L,d)$
over 
$[Z] \in {{}^s}\!{ \XD_{lci}}$
is a hypersurface of degree $d$ in
$\PP(\EZ)$.
Furthermore, one can show that set-theoretically
$\TE_{Z}$
is the union of hyperplanes
$H_z$ in $\PP(\EZ)$,
where $z$ runs through the set of closed points in $Z$.
Thus the divisor
$\TE(X;L,d)$
captures geometry of zero-dimensional subschemes of $X$ parametrized by the underlying points of the Hilbert scheme $\XD$.
\begin{rem}\label{morst1}
Let
$\JAA^{\prime}_{{{}^s}\!{ \XD_{lci}}}$
be the complement of 
$\TE(X;L,d)$ in $\JAA_{{{}^s}\!{ \XD_{lci}}}$.
Then the restriction of the morphism $h$ in (\ref{morst}) to
$\JAA^{\prime}_{{{}^s}\!{ \XD_{lci}}}$ gives the morphism
\begin{equation}\label{h-pr}
\xymatrix@1{
h^{\prime} : \JAA^{\prime}_{{{}^s}\!{ \XD_{lci}}} \ar[r]& {\bf B}_X (2,L,d), }
\end{equation}
where 
${\bf B}_X (2,L,d) $
is the moduli stack of locally free sheaves on $X$ having rank $2$ and the Chern invariants
$(L,d)$.
\end{rem}

From now on we will be working over 
${{}^s}\!{ \XD_{lci}}$.
By Definition \ref{L-s} we have an inclusion
${{}^s}\!{ \XD_{lci}} \subset \GA^0_d (L)$, provided $d\geq 2$.
Setting
\begin{equation}\label{s-str}
\xymatrix@1{
{}^{\prime} \GG(L) = {{}^s}\!{ \XD_{lci}}  \bigcap \GG(L)&*\txt{and}&
 {}^{\prime} \GAO (L) = {}^{\prime} \GG(L) \setminus  {}^{\prime} \GA^{r+1}_d (L)\,, }
\end{equation}
we obtain the stratification
\begin{equation}\label{s-str1}
{{}^s}\!{ \XD_{lci}} = \bigcup_{r\geq 0}  {}^{\prime} \GAO (L)
\end{equation}
of ${{}^s}\!{ \XD_{lci}}$ by locally closed sets 
${}^{\prime} \GAO (L)$.
Taking the inverse image by
$\pi$ of this stratification, we obtain the stratification
\begin{equation}\label{s-strJ}
\JAA_{{{}^s}\!{ \XD_{lci}}} = \pi^{-1} ({{}^s}\!{ \XD_{lci}}) = \bigcup_{r\geq 0} {}^{\prime} {\stackrel{\circ}{\JAA^r}}\,,
\end{equation}
where
$ {}^{\prime} {\stackrel{\circ}{\JAA^r}} =  \pi^{-1} ( {}^{\prime} \GAO (L) )$
is a $\PP^r$-bundle over 
${}^{\prime} \GAO (L) $.

\subsection{The sheaf $\tilde{\cal F}$ on $\JA$.} \label{secFT}
Recall the universal scheme $\ZD$ in (\ref{uc}) and set
\begin{equation}\label{F}
\FF= p_{2\ast} \OO_{\ZD}
\end{equation}
to be the direct image of the structure sheaf $\OO_{\ZD}$ of $\ZD$ with respect to the projection
$p_2$ in (\ref{uc}). This is a locally free sheaf of rank $d$ on $\XD$.

Let
\begin{equation}\label{FT}
\FT= \pi^{\ast} \FF
\end{equation}
be the pullback of $\FF$ under the projection $\pi$ in (\ref{pi1}).
This is also a locally free sheaf of rank $d$ on $\JA$.

One of the main points of the constructions in \cite{[R1]} is a distinguished subsheaf
$\HT$ of $\FT \otimes \OO_{\JAA^r}$ defined for every stratum
$\JAA^r$ in (\ref{strJ}) (see \cite{[R1]},\S1.3, for details).
Since $\FT$ is a sheaf of rings, the multiplication in
$\FT$
gives rise to a distinguished filtration of
$\FT \otimes \OO_{\JAA^r}$
\begin{equation}\label{filtHT}
0=\HT_0 \subset \HT_{-1} \subset \ldots \subset \HT_{-i} \subset \ldots \subset \FT \otimes \OO_{\JAA^r}\,,
\end{equation}
where
\begin{equation}\label{HT-i}
\xymatrix@1{
{\HT_{-1} = \HT}&*\txt{and}& \HT_{-i}= im\left(S^i {\HT}\right. \ar[r]&\left.\FT \otimes \OO_{\JAA^r} \right) }
\end{equation}
is the image of the morphism
\begin{equation}\label{m-i}
\xymatrix@1{
m_i: S^i {\HT} \ar[r]& {\FT \otimes \OO_{\JAA^r}} }
\end{equation}
induced by the multiplication in $\FT$, where $S^i {\HT}$ is the $i$-th symmetric power of $\HT$.

Fix a stratum $\GG(L)$ in (\ref{strH}) such that the open part
${}^{\prime}\GAO (L)$,
 defined by (\ref{s-str}),
is non-empty and consider the smooth part
$reg({}^{\prime}\GG (L))$
of ${}^{\prime}\GG (L)$ 
(recall from the discussion following (\ref{strH}) that 
$reg({}^{\prime}\GG (L))$
is contained in ${}^{\prime}\GAO (L) $).
Denote by
\begin{equation}\label{pi-0}
C^r (L,d) := \pi_0 (reg({}^{\prime}\GG (L)))
\end{equation}
the set of connected components of
$reg({}^{\prime}\GG (L))$. This is also the set of connected components of
$reg({}^{\prime}\JAA^r)$, the smooth locus of 
${}^{\prime}\JAA^r = \pi^{-1} ({}^{\prime}\GG (L))$.
For every connected component
$\GA \in \CS$,
denote by
\begin{equation}\label{JG}
\JG := \pi^{-1} (\GA)
\end{equation}
the corresponding connected component of
$reg({}^{\prime}\JAA^r)$.

Fix $\GA \in \CS$ and consider the restriction of the filtration (\ref{filtHT}) to
$\JG$.
The sheaves 
$\HT_{-i} \otimes \OO_{\JG}$
are non-zero and torsion free, for every $i\geq 1$. So their ranks are well-defined.
Set
\begin{equation}\label{h-i}
h^{i-1}_{\GA} = rk(\HT_{-i} \otimes \OO_{\JG}) - rk( \HT_{-i+1} \otimes \OO_{\JG})\,,
\end{equation}
for every $i \geq 1$.

 Denote by $l_{\GA}$ the largest index $i$ for which
$h^{i-1}_{\GA} \neq 0$
and call it {\it the length of the filtration} of
$\HT_{-\bullet} \otimes \OO_{\JG}$.
Thus on 
$\JG$
the filtration (\ref{filtHT}) stabilizes at 
$\HT_{-l_{\GA}} \otimes \OO_{\JG}$.
We also agree to assign to
$\FT\otimes \OO_{\JG}$
the index 
$-(l_{\GA} +1)$
and use the notation
\begin{equation}\label{lterm}
\FT\otimes \OO_{\JG} = \HT_{-l_{\GA} -1}\,.
\end{equation}
So the filtration (\ref{filtHT}) restricted to $\JG$ has the following form
\begin{equation}\label{filtHT-JG}
0=\HT_0 \otimes \OO_{\JG} \subset \HT_{-1} \otimes \OO_{\JG} \subset \ldots \subset \HT_{-l_{\GA} +1} \otimes \OO_{\JG}
 \subset  \HT_{-l_{\GA}} \otimes \OO_{\JG} \subset   \HT_{-l_{\GA}-1} = \FT \otimes \OO_{\JG}\,.
\end{equation}
\begin{lem}\label{hg=c}
Set
$$
h_{\GA} =(h^0_{\GA},\ldots, h^{l_{\GA} -1} _{\GA},h^{l_{\GA} } _{\GA})
$$
and call it the Hilbert vector of $\GA$.

The Hilbert vector $h_{\GA}$ and its components
$h^i_{\GA},\,\,(i=0,\ldots, l_{\GA}),$
have the following properties:
\begin{enumerate}
\item[1)]
$h_{\GA}$ is a composition of $d$, i.e.
$$
\sum^{l_{\GA}}_{i=0} h^i_{\GA} = d\,,
$$
\item[2)]
$h^0_{\GA}= rk(\HT_{-1} \otimes \OO_{\JG}) =rk(\HT \otimes \OO_{\JG}) = r+1$,
\item[3)]
$h^i_{\GA} >0$, for $i=0,\ldots, l_{\GA} -1$, and
$h^{l_{\GA}}_{\GA} \geq 0$.
\end{enumerate}
\end{lem}
\begin{pf}
From the defining equations (\ref{h-i}) it follows
$$
\sum^{l_{\GA}}_{i=0} h^i_{\GA} = rk(\FT \otimes \OO_{\JG}) =  d\,.
$$

The second assertion follows from \cite{[R1]}, Proposition 1.4, and 3) is obvious.
\end{pf}

Let $C(d)$ be the set of compositions of $d$. From Lemma \ref{hg=c} it follows
that the assignment of the Hilbert vector
$h_{\GA}$
to the components $\GA$ in $\CS$ gives a map
\begin{equation}\label{CS-Hv}
\xymatrix@1{
{h(L,d,r): \CS} \ar[r]& C(d)}\,.
\end{equation} 
\begin{rem}\label{kap}
The sheaves 
$\HT_{-i}$
in (\ref{filtHT}) and hence their ranks are related to the geometry of the underlying points of the Hilbert scheme.
This is the content of \cite{[R1]}, Remark 1.5, which we reproduce here for convenience of the reader.

Let $ \ZA$ be a point of 
$\JG$.
The fibre 
$\HT \ZA$
of
$\HT$
at
$\ZA$
can be viewed as a linear system on $Z$. It is always base point free, since
$\HT \ZA$
contains the constant functions on $Z$ (see \cite{[R1]}, Remark 1.3).
Thus
$\HT \ZA$
defines a morphism
\begin{equation}\label{kappa}
\xymatrix@1{
{\kappa \ZA }: Z \ar[r]& {\PP( \HT \ZA^{\ast})}\,. }
\end{equation}
Hence the Hilbert function of the image 
$Z^{\prime}(\alpha)$
of 
$\kappa \ZA$
is given by the ranks of
$\HT_{-i}$'s at $\ZA$, for $i=1,\ldots, l_{\GA}$.
\end{rem}

From Remark \ref{kap} it follows that the Hilbert vector
$h_{\GA}$
encodes the Hilbert function of the image of
$\kappa \ZA$,
for all 
$\ZA$ varying in the complement of the singularity loci of the sheaves
$\HT_{-i}$,
 for $i=1,\ldots, l_{\GA}$.
This is a non-empty Zariski open subset of
$\JG$ which we denote by
$\JG^{\prime}$.
Set
\begin{equation}\label{G(0)}
\GA^{(0)} = \pi (\JG^{\prime})
\end{equation}
to be its image under the projection $\pi$ in (\ref{pi1}). This is a Zariski open subset of
$\GA$. The following lemma relates $\JG^{\prime}$ and the complement of the theta-divisor in $\JG$.
\begin{lem}\label{lem-rk}
Let 
$ \TE_{\GA^{(0)}} =\TE(X;L,d) \cap \pi^{-1} ( \GA^{(0)})$
be the theta-divisor over 
$\GA^{(0)}$
and let
$$
\JG^{(0)} = \pi^{-1} ( \GA^{(0)}) \setminus  \TE_{\GA^{(0)}}
$$
be its complement in
$\pi^{-1} ( \GA^{(0)})$.
Then
$
\JG^{(0)} \subset \JG^{\prime}
$.
\end{lem}
\begin{pf}
Let 
$[Z] \in \GA^{(0)}$
and let
$\JAA_{Z}$ (resp. $\TE_{Z}$) be the fibre of
$\JG$ (resp. $\TE_{\GA^{(0)}}$) over $[Z]$.
By definition
$\JG^{\prime}$ intersects 
$\JAA_{Z}$
along a non-empty Zariski open set. In particular, there is
$[\alpha]$ in $\JAA_{Z} \setminus \TE_{Z}$
such that
$\ZA \in \JG^{\prime}$.
Hence the ranks of the sheaves
$\HT_{-i}$ at $\ZA$
are given by
\begin{equation}\label{rk-f}
rk( \HT_{-i} \ZA )= \sum ^{i-1}_{k=0} h^k_{\GA}\,,
\end{equation}
for $i=1, \ldots, l_{\GA}$. We claim that the left hand side of
(\ref{rk-f}) stays constant for all
$[\beta] \in (\JAA_{Z} \setminus \TE_{Z} )$.
This will give the assertion of the lemma.
In fact, we claim that the following holds.
\begin{cl}\label{cl-ab}
For any 
$[\alpha], [\beta] \in (\JAA_{Z} \setminus \TE_{Z} )$
one has an isomorphism
$$
\xymatrix@1{
{\phi_{[\alpha], [\beta]} : \HT \ZA} \ar[r]& {\HT ([Z],[\beta])} }
$$
which induces isomorphisms
$$
\xymatrix@1{
{\phi^{i}_{[\alpha], [\beta]} : \HT_{-i} \ZA} \ar[r]& {\HT_{-i} ([Z],[\beta])}\,, }
$$
for every 
$i=1, \ldots, l_{\GA}$.
\end{cl}
{\it {Proof of }} {\bf Claim \ref{cl-ab}}.
From \cite{[R1]}, Proposition 1.4, $\alpha$ induces an isomorphism
\begin{equation}\label{HT=Ext}
\xymatrix@1{
{\HT \ZA} \ar[r]^(0.55){\alpha}& {\EZ} }
\end{equation}
This can be seen explicitly by recalling that
$$
\EZ = Ext^1 (\ID_{Z} (L), \OO_X)
$$
can be identified as a subspace of
$
H^0 ( \EEF^2 (\OO_Z (L), \OO_X))
$
(see \RI, (1.16)). The latter space is
$
H^0 ( \omega_Z  \otimes \OO_X (-L -K_X))
$,
where
$\omega_Z$ 
is the dualizing sheaf of $Z$.
Since $Z$ is a local complete intersection,
$\omega_Z$
is invertible and we think of
$\EZ$ as a linear subspace of sections of an invertible sheaf on $Z$.
Furthermore,
$\alpha$ is an extension class corresponding to a locally free sheaf
$\SE_{\alpha}$ sitting in the middle of the short exact sequence
$$
\xymatrix@1{
0\ar[r]& {\OO_X} \ar[r]& {\SE_{\alpha}}\ar[r]&{\ID_Z} (L) \ar[r]& 0 }
$$
defined by $\alpha$. This is equivalent to
$\alpha$ being nowhere vanishing\footnote{as a section of 
$\omega_Z  \otimes \OO_X (-L -K_X)$.}
 on $Z$. Hence (\ref{HT=Ext}) implies the following identification
\begin{equation}\label{HT=fr}
\HT \ZA = \left\{\left. \frac{\gamma}{\alpha} \right| \gamma \in \EZ \right\}\,.
\end{equation}

Define $\phi_{[\alpha], [\beta]}$ to be the multiplication by
$\displaystyle{\frac{\alpha}{\beta}}$
to obtain the isomorphism
\begin{equation}\label{a/b}
\xymatrix@1{
 {\frac{\alpha}{\beta}: \HT \ZA} \ar[r]& {\HT ([Z],[\beta]) }\,. }
\end{equation}
This proves the first assertion of the claim. 

From (\ref{a/b}) it also follows
\begin{equation}\label{a/b1}
\HT ([Z],[\beta]) = \frac{\alpha}{\beta} \HT \ZA\,.
\end{equation}
This implies
\begin{equation}\label{a/b2}
\HT_{-i} ([Z],[\beta]) = \left(\frac{\alpha}{\beta} \right)^{i} \HT_{-i} \ZA\,,
\end{equation}
for every 
$i\geq 1$. Thus taking
$\phi^{i}_{[\alpha], [\beta]}$
to be the multiplication by
$\displaystyle{\left(\frac{\alpha}{\beta} \right)^{i}}$, 
yields the second assertion of the claim.
\end{pf}

As a consequence we deduce the following.
\begin{cor}\label{cor-lf}
For every 
$\GA \in \CS$, there exists a non-empty Zariski open subset
$\GA^{(0)} \subset \GA$ (defined in (\ref{G(0)})) such that
on the open part
$$
\JG^{(0)} = \pi^{-1} (\GA^{(0)} ) \setminus \TE_{\GA^{(0)} }
$$
of $\JG$ all non-zero sheaves of the filtration
$\HT_{-\bullet}$ in (\ref{filtHT}) are locally free and their ranks are determined by the Hilbert vector
$$
h_{\GA} =(h^0_{\GA},\ldots, h^{l_{\GA} -1} _{\GA},h^{l_{\GA} } _{\GA})
$$
from Lemma \ref{hg=c} by the formula
$$
rk (\HT_{-i} \otimes \OO_{\JG^{(0)}} ) = \sum^{i-1}_{k=0} h^k_{\GA}\,.
$$
\end{cor}

From Remark \ref{kap} it follows that our considerations are non-trivial provided 
\begin{equation}\label{r}
r \geq 1\,.
\end{equation}
This will be assumed for the rest of the paper.

\subsection{Orthogonal decomposition of $\FT$.}\label{sec-ord}
The filtration $\HT_{-\bullet}$ in (\ref{filtHT}) acquires more structure over the points of
$\GG(L)$ corresponding to the reduced subschemes of $X$.
Let 
$Conf_d (X)$
be the locus of the Hilbert scheme 
$\XD$
parametrizing the subschemes of $d$ distinct points of $X$. This is a Zariski open subset of
$\XD$, since it can be described as the complement of the branching divisor of the ramified covering
$$
\xymatrix@1{
{p_2 :\ZD} \ar[r]& {\XD} }
$$
in (\ref{uc}). The subschemes $Z$ of $X$, with
$[Z] \in \CO$, will be called configurations (of $d$ points) on $X$.

We are interested in the connected components 
$\GA \in \CS$
having a non-empty intersection with
$\CO$.
\begin{defi}\label{csa}
A component 
$\GA \in \CS$
is called admissible if
$$
\GA \cap \CO \neq \emptyset\,.
$$
The set of admissible components in $\CS$ will be denoted by
$C^r_{adm} (L,d)$.
\end{defi}

For a subset $Y$ in
$\XD$
we denote by
$Y_{conf}$
the intersection
$Y \cap \CO$.
In particular, for  
$\GA \in \CS$,
the subset $\GA_{conf}$ is Zariski open in $\GA$ and it is non-empty if and only if 
$\GA \in \CSA$ .
The subset $\GA_{conf}$ will be called the configuration subset of $\GA$.

We will now recall why configurations are important in our constructions (see \RI, \S2, for details).
The sheaf $\FF$ (resp. $\FT$) admits the trace morphism
\begin{equation}\label{tr}
\xymatrix@1{
{Tr: \FF=p_{2\ast} \OO_{\ZD}} \ar[r]& {\OO_{\XD}}& *\txt{(resp.}& {\tilde{Tr}: \FT} \ar[r]& {\OO_{\JA}}) }
\end{equation}
It can be used to define the bilinear, symmetric pairing
$\QB$ (resp. ${\bf \tilde{q}}$) on $\FF$ (resp. $\FT$)
\begin{equation}\label{q}
\QB(f,g) =Tr(fg)\,\,({\rm{resp.}}\,\, {\bf \tilde{q}}(f,g) = \tilde{Tr} (fg)),
\end{equation}
for every pair $(f,g)$ of local sections of $\FF$ (resp. $\FT$).

This pairing is non-degenerate precisely over
$\CO$. Using it we obtain a natural splitting of the filtration
$\HT_{-\bullet} \otimes \OO_{\JG}$ in (\ref{filtHT-JG})
on a certain Zariski open subset of
$\JAA_{\GA^{(0)}_{conf}} =\pi^{-1} (\GA^{(0)}_{conf})$,
for every admissible component
$\GA \in \CSA$, where $\GA^{(0)}$ is a Zariski open subset of $\GA$ defined in (\ref{G(0)}).
 More precisely, set
\begin{equation}\label{ort-c}
\FI^i =\big{(} \HT_{-i} \otimes \OO_{\GA^{(0)}_{conf}} \big{)}^{\perp}
\end{equation}
to be the subsheaf of
$\FT  \otimes \OO_{\GA^{(0)}_{conf}}$
orthogonal to 
$ \HT_{-i} \otimes \OO_{\GA^{(0)}_{conf}} $
with respect to the quadratic form
$\QT$ in (\ref{q}). It was shown in \RI, Corollary 2.4, that there exists a non-empty Zariski open subset
$\JAB_{\GA}$ of $\JG$ subject to the following properties:
\begin{enumerate}
\item[(a)]
the open set 
$\JAB_{\GA}$ lies over 
$ \GA^{(0)}_{conf}$, i.e. the morphism $\pi$ in (\ref{pi1}) restricted to 
$\JAB_{\GA}$ gives the surjective morphism
\begin{equation}\label{pi-JAB}
\xymatrix@1{
{\pi: \JAB_{\GA}} \ar[r]& {\GA^{(0)}_{conf}}, }
\end{equation}
\item[(b)]
$\JAB_{\GA}$ lies in the complement of the theta-divisor, i.e.
\begin{equation}\label{JAB-inc}
\JAB_{\GA} \subset \JG^{(0)},
\end{equation}
 where $ \JG^{(0)} $ is as in Lemma \ref{lem-rk},  
\item[(c)]
$\FT $, restricted to $\JAB_{\GA}$, admits the orthogonal direct sum decomposition
\begin{equation}\label{FT-ors}
\FT \otimes \OO_{\JAB_{\GA}} =\HT_{-i} \otimes \OO_{\JAB_{\GA}} \oplus \FI^i  \otimes \OO_{\JAB_{\GA}},
\end{equation}
for every $i=0, 1, \ldots, l_{\GA} +1$.
\end{enumerate}
This gives rise to the filtration
\begin{equation}\label{filtF}
\FT \otimes \OO_{\JAB_{\GA}} =\FI^0  \otimes \OO_{\JAB_{\GA}} \supset \FI^1  \otimes \OO_{\JAB_{\GA}} \supset \ldots \supset 
\FI^{ l_{\GA}} \otimes \OO_{\JAB_{\GA}} \supset \FI^{ l_{\GA} +1} =0.
\end{equation}
Putting together the filtrations $\HT_{-\bullet} \otimes \OO_{\JAB_{\GA}}$ and $\FI^{\bullet}  \otimes \OO_{\JAB_{\GA}}$,
we define the subsheaves
\begin{equation}\label{H}
\HH^{i-1} =\big{(} \HT_{-i} \otimes \OO_{\JAB_{\GA}} \big{)}\cap \big{(} \FI^{i-1}  \otimes \OO_{\JAB_{\GA}}\big{)},
\,\,{\rm{for}} \,\,i=1,\ldots,  l_{\GA} +1.
\end{equation}
This definition together with (\ref{FT-ors}) yield the following decomposition of 
$\FT  \otimes \OO_{\JAB_{\GA}}$
into the orthogonal sum
\begin{equation}\label{ordFT}
\FT  \otimes \OO_{\JAB_{\GA}} = \bigoplus^{ l_{\GA}}_{p=0} \HH^p.
\end{equation}
\begin{rem}\label{rkH}
\begin{enumerate}
\item[1)]
Observe that the ranks of the summands
$\HH^p$'s in (\ref{ordFT}) form the Hilbert vector
$h_{\GA}$ as defined in Lemma \ref{hg=c}, i.e.
\begin{equation}\label{rkHp}
rk(\HH^p)= h^p_{\GA}.
\end{equation}
This follows from the definition of 
$ h^p_{\GA}$ in (\ref{h-i}), the inclusion (\ref{JAB-inc})
and the orthogonal decomposition 
\begin{equation}\label{ordH-i}
\HT_{-i} \otimes \OO_{\JAB_{\GA}} = \bigoplus^{ i-1}_{p=0} \HH^p.
\end{equation}
In particular, for $i=1$ one obtains
\BEN\label{ordH-1}
\HT =\HT_{-1} =\HH^0.
\EEN
\item[2)]
From the orthogonal decomposition (\ref{ordFT}) it follows that the subsheaves 
$\FI^i$ of the filtration $\FI^{\bullet}$ in (\ref{filtF}) admit the following orthogonal decomposition
\BEN\label{ordF-i}
\FI^i = \bigoplus^{\LG}_{p=i} \HH^p.
\EEN
\end{enumerate}
\end{rem}

The decomposition (\ref{ordFT}) together with the multiplicative structure of 
$\FT$ play the crucial role in our considerations. In particular, they give rise to
the sheaf of Lie algebras 
{\boldmath$\GT$} whose definition will be recalled in  the next subsection.
We close this one by giving the dual version of the filtration
$\HT_{-\bullet}$ in (\ref{filtHT}) which coincides with the filtration
$\FI^{\bullet}$ once restricted to $\JAB_{\GA}$. However, it has a virtue of being more geometric.

The starting point of the dual construction is another natural sheaf on $\XD$
which takes account of the divisor $L$. Namely, we consider the sheaf
\begin{equation}\label{def-FL}
\FF(L) = p_{2\ast} \left( p^{\ast}_1 \OO_X (L+K_X) \right)\,,
\end{equation}
where $p_i (i=1,2)$ are as in (\ref{uc}).  Taking its pullback via $\pi$ in 
(\ref{pi1}) we obtain the sheaf
\begin{equation}\label{FTL}
\FT(L) = \pi^{\ast} \left( \FF(L) \right)\,.
\end{equation}

In \RI, \S1.3, it was shown that there is a natural morphism
\begin{equation}\label{mor-R}
\xymatrix@1{
{{\bf R^r}: \FT \otimes \OO_{\JAA^r}} \ar[r]& H^0(L+K_X)^{\ast} \otimes \OO_{\JAA^r} (1), }
\end{equation}
where 
$\OO_{\JAA^r} (1) $
is the restriction to $\JAA^r$ of the tautological invertible sheaf
$ \OO_{\JA} (1)$  (see (\ref{O1}) for notation).
In particular, the subsheaf 
$\HT$, encountered in \S\ref{secFT}, is defined in \RI, (1.21), as the kernel of 
${\bf R^r}$. Furthermore, we have morphisms
\begin{equation}\label{mor-RT-i}
\xymatrix@1{
{{\bf  \tilde{R}^r_i}: S^i \HT} \ar[r]^{m_i}&{\FT \otimes \OO_{\JAA^r}} \ar[r]^(0.35){\bf R^r}& H^0(L+K_X)^{\ast} \otimes \OO_{\JAA^r} (1), }
\end{equation}
where $m_i$ is as in (\ref{m-i}). 

Dualizing and tensoring with 
$ \OO_{\JAA^r} (1)$
yields
$$
\xymatrix@1{
H^0(L+K_X) \otimes \OO_{\JAA^r}  \ar[r] & {\left( S^i \HT \right)^{\ast} \otimes \OO_{\JAA^r} (1) } }
$$
Setting $\FIT_{i}$ to be the kernel of this morphism, we obtain the following filtration
$$
H^0(L+K_X) \otimes \OO_{\JAA^r}= \FIT_1 \supset \FIT_2 \supset  \ldots \supset \FIT_i \supset \FIT_{i+1} \supset \ldots
$$

Each 
$\FIT_i $
contains the sheaf 
$
\tilde{\JJ}_{\ZD} = \pi^{\ast} \big{(} pr_{2\ast} \big{(} \JJ_{\ZD} \otimes pr^{\ast}_1 \OO_X (L+K_X) \big{)} \big{)}
$,
where 
$ \JJ_{\ZD}$
is the sheaf of ideals of the universal subscheme $\ZD$ in $X \times \XD$ (see (\ref{uc}) for notation) 
and $pr_j (j=1,2)$ are the projections of
$X \times \XD$
onto the corresponding factor.\footnote{
the inclusion
$\pi^{\ast} \big{(} pr_{2\ast} \big{(} \JJ_{\ZD} \otimes pr^{\ast}_1 \OO_X (L+K_X) \big{)} \big{)} \subset \FIT_i $
 is proved in \RI, Proposition 1.6.} Factoring out by 
$\displaystyle{\tilde{\JJ}_{\ZD}}$,
 one obtains the following filtration of 
$\FT(L)$:
\begin{equation}\label{filtFd}
\FT(L) \otimes \OO_{\JAA^r} =\FI_0   \supset \FI_1   \supset \ldots \supset \FI_i \supset \FI_{ i +1} \supset \ldots
\end{equation}
where
\begin{equation}\label{F-i}
 \FI_i  =\FIT_i  \big{/} \tilde{ \JJ}_{\ZD}.
\end{equation}
To relate this filtration to the one in (\ref{filtF}) one observes that there is a natural morphism
\begin{equation}\label{mor-RT}
\xymatrix@1{
{{\bf \tilde{R}}: \FT  \otimes \OO_{\JAA^r}} \ar[r]& {(\FT (L))^{\ast} \otimes \OO_{\JAA^r} (1)} }
\end{equation}
(see \RI, (1.27) and (1.19) for details). Furthermore, this morphism is an isomorphism precisely on the complement of the theta-divisor
$\TE(X;L,d)$,  since the latter is defined by the vanishing of the determinant of 
${\bf \tilde{R}}$ 
(see the formula for $\TE(X;L,d)$ below \RI, (1.19)).
Taking the dual of
${\bf \tilde{R}}$, 
we obtain a natural identification of
$\FT (L) \otimes \OO_{\JAA^r} (-1) $
with
$ \FT^{\ast}  \otimes \OO_{\JAA^r}$
on the complement of the theta divisor in $\JAA^r$.
Restricting further to the configurations and using the self-duality\footnote{the self-duality of $\FT$
over $\CO$ is provided by the quadratic form $\QT$ in (\ref{q}).}
 of
$\FT$
over $\CO$, we obtain a natural identification of
$\FT (L) \otimes \OO_{\JAA^r} (-1) $
and
$ \FT  \otimes \OO_{\JAA^r}$
over the complement
$$
\pi^{-1} (\GAC) \setminus \TE(X;L,d),
$$
for every admissible component
$\GA \in \CSA$. In particular, this identification holds on 
$\JAB_{\GA}$ in view of the inclusion in (\ref{JAB-inc}).

With the above identification in hand, we can transfer the filtration
$\FI_{\bullet} \otimes \OO_{\JAB_{\GA}} (-1)$ of 
$\FT(L) \otimes \OO_{\JAB_{\GA}} (-1)$
in (\ref{filtFd}) (twisted by  $\OO_{\JAB_{\GA}} (-1)$)
to a filtration of
$\FT\otimes \OO_{\JAB_{\GA}}$.
 The point is that the resulting filtration is the filtration
$\FI^{\bullet}$, defined previously in (\ref{filtF}) via orthogonality
(see more detailed discussion in \RI, \S2).

By definition, the filtration 
$\HT_{-\bullet}$
is related to the geometry of the morphisms
$\kappa \ZA$ in (\ref{kappa})
associated to the linear systems
$\left| \HT\ZA \right|$
on $Z$, as $[Z]$ varies through the points of the admissible components
$\GA$.
On the other hand the filtration $\FI^{\bullet}$ in (\ref{filtF}), in view
of its identification with
$\FI_{\bullet} \otimes \OO_{\JAB_{\GA}} (-1)$,
reflects geometric properties of the subschemes $Z$
(parametrized by $\GA$) with respect to the adjoint linear system
$\left| L+K_X \right|$ on $X$.
Thus the orthogonal decomposition
(\ref{ordFT}) contains information about the geometry of the subschemes
$Z$ with respect to {\it both} linear systems.

\subsection{The sheaf of the Lie algebras {\boldmath$\GT_{\mbox{\unboldmath${\GA}$}}$}. }\label{Lie}
To attach Lie algebras to points of $\JA$ we view local sections of the sheaf
$\HT =\HT_{-1}$
in (\ref{filtHT}) as operators of the multiplication in  the sheaf of rings $\FT$, i.e. we consider the inclusion
\begin{equation}\label{m-D}
\xymatrix@1{
{D: \HT} \ar[r]& {\ENDO(\FT)}  }
\end{equation}
which sends a local section $t$ of $\HT$ to the operator
$D(t)$ of the multiplication by $t$ in $\FT$.

Over the configuration subset
$\GAC$ of an admissible component $\GA \in \CSA$
we have defined the subscheme
$\JAB_{\GA}$ (see (\ref{pi-JAB}) for notation),
where the orthogonal decomposition
(\ref{ordFT}) holds. Using this decomposition we write
\begin{equation}\label{d-Dt}
D(t) = D^{-} (t) + D^{0}(t)  + D^{+}(t)\,,
\end{equation}
where the components 
$D^{\pm} (t)$
have degree $\pm 1$, while 
$ D^{0}(t) $
is of degree $0$, with respect to the grading in (\ref{ordFT})
(see \RI, Remark 3.8, for more details).
Thus on $\JAB_{\GA}$ the morphism $D$ in (\ref{m-D}) admits the triangular decomposition
\begin{equation}\label{d-D}
 D = D^{-}  + D^{0}  + D^{+}
\end{equation}
and we define 
{\boldmath$\GT_{\mbox{\unboldmath${\GA}$}}$}
to be the subsheaf of Lie subalgebras of
$\ENDO(\FT)$
generated by the subsheaves
$ D^{\pm} (\HT)$ and $D^{0} (\HT)$.

It was observed in \RI, \S7, that 
{\boldmath$\GT_{\mbox{\unboldmath${\GA}$}}$}
is a sheaf of reductive Lie algebras and its representation theory is intimately related to the geometry of the subschemes of $X$
parametrized by $\GA$.

The main objective of this paper is to pursue the investigation of this relation.
Our considerations logically fall into two parts.
The first one seeks to use the representation theory of 
{\boldmath$\GT_{\mbox{\unboldmath${\GA}$}}$}
to gain an insight into geometry of the subschemes of $X$ parametrized by
$\GA$.
The second uses various representation theoretic constructions related to
{\boldmath$\GT_{\mbox{\unboldmath${\GA}$}}$}
to obtain interesting objects (e.g. sheaves, complexes of sheaves) on 
$\JA$ or on the underlying Hilbert scheme $\XD$. 

\subsection{Conventions and notation}\label{c-n}
In this section we summarize all notation and conventions introduced so far and which will be used throughout the rest of the paper.
\begin{enumerate}
\item[{\bf \S}]
{\bf \ref{c-n}.0.}
The Chern datum $(L,d)$, consisting of a divisor $L$ (up to the rational equivalence) on $X$ and a positive integer 
$d$, is
fixed once and for all. The divisor $L$ is always subject to the vanishing assumptions in (\ref{vc}).
\item[{\bf \S}]
{\bf \ref{c-n}.1.}\label{c-n1}
The Hilbert scheme (resp. Jacobian)  $\XD$ (resp. $\JA$) is equipped with the stratification defined in (\ref{strH})
(resp. (\ref{strJ})).

We always consider the strata $\GG$ with $r \geq 1$ and
${}^{\prime} \GAO$ non-empty (see (\ref{s-str}) for notation).
For such a stratum $\GG$ we denote by
$reg( {}^{\prime} \GAO)$
its smooth part  
and we let 
$\CS$ 
to be the set of its connected components.
Furthermore, we denote by 
$\CSA$ the set of {\it admissible} components (see Definition\ref{csa})
of
$reg( {}^{\prime} \GAO)$.
\item[{\bf \S}] 
{\bf \ref{c-n}.2.}\label{ssec-G(0)}
For a component
$\GA \in \CSA$, we set
\begin{equation}\label{JG-TG}
\JG = \pi^{-1} (\GA)\,\,{\rm{and}}\,\, \TE_{\GA} =\TE(X;L,d) \cap \JG,
\end{equation}
where $\pi$ is the morphism defined in (\ref{pi1}).
Set 
$\JG^{\prime}$ to be the largest Zariski open subset of 
$\JG$ over which all non-zero sheaves of the filtration
in (\ref{filtHT-JG}) are locally free and let
$$
\GA  ^{(0)} = \pi (\JG^{\prime})
$$
be the corresponding Zariski open subset of 
$\GA$.
\\
\indent
Denote by
$\JG^{(0)}$
the complement of the theta-divisor in 
$\pi^{-1} (\GA^{(0)})$. From Corollary \ref{cor-lf} it follows that all
non-zero sheaves in the filtration
in (\ref{filtHT-JG}) are locally free  over 
$\JG^{(0)}$
and their ranks are determined by the Hilbert vector (see Lemma \ref{hg=c})
\begin{equation}\label{HvG}
h_{\GA} =(h^0_{\GA},\ldots, h^{l_{\GA} -1} _{\GA},h^{l_{\GA} } _{\GA}) 
\end{equation}
by the formula in Corollary \ref{cor-lf}.
\item[{\bf \S}] 
{\bf \ref{c-n}.3.}\label{breve}
For a component
$\GA \in \CSA$, denote by 
$\JAB_{\GA}$
the largest Zariski  open subset of 
$\JG^{(0)}$
over which the orthogonal decomposition
(\ref{ordFT}) holds. It is known that it projects onto $\GA^{(0)}_{conf}$
(see (\ref{pi-JAB})), i.e. the morphism
\begin{equation}\label{pi-JAB1}
\xymatrix@1{
{\pi: \JAB_{\GA}} \ar[r]& {\GA^{(0)}_{conf}} }
\end{equation}
is surjective. In the sequel, to simplify the notation, we set
\begin{equation}\label{GAb}
\breve{\GA} := \GA^{(0)}_{conf}.
\end{equation}
\item[{\bf \S}] 
{\bf \ref{c-n}.4.}
For an admissible component $\GA$ in $\CS$, we always denote by
{\boldmath$\GT_{\mbox{\unboldmath${\GA}$}}$}
the sheaf of reductive Lie algebras defined in
\S\ref{Lie}. Its semisimple part
{\boldmath$[ \GT_{\mbox{\unboldmath$\GA$}}, \GT_{\mbox{\unboldmath$\GA$}}]$}
will be denoted by
{\boldmath$ \GS_{\mbox{\unboldmath$\GA$}}$}.
\item[{\bf \S}] 
{\bf \ref{c-n}.5.}
Working over an admissible component
$\GA$
most of the time we consider sheaves over
$\JAB_{\GA}$. So, to simplify the writing,
we often omit tensoring with
$\OO_{\JAB_{\GA}}$
in the notation of restriction to
$\JAB_{\GA}$
of a sheaf defined on a larger space.
Thus, for example, the filtration
$\HT_{-\bullet}$ in (\ref{filtHT}) will be considered over 
$\JAB_{\GA}$, unless said otherwise, and we will continue to write
$\HT_{-\bullet}$ instead of the more cumbersome
$\HT_{-\bullet} \otimes \OO_{\JAB_{\GA}}$.
\item[{\bf \S}]\label{sh-bdl} 
{\bf \ref{c-n}.6.}
At certain parts of the paper we will need to distinguished between locally free sheaf on a variety and the corresponding vector bundle.
To do this we use capital calligraphic letter for the former and the {\it same} capital letter, but in Roman type, for the latter. Thus, if we have a variety
 $Y$ with a locally free sheaf $\cal{A}$ on it, then we denote by $A$ the corresponding vector bundle over $Y$.
\end{enumerate}

\section{Some properties of the filtration ${\HT_{-\bullet}}$}
We fix a stratum
$\GG$ in (\ref{strH}) according to the conventions in \S \ref{c-n1}
and consider an admissible component
$\GA$ in $\CSA$.
\begin{lem}\label{HTl=triv}
Let 
$[Z] \in \GA$
and
let 
$\JAA_{Z} = \pi ^{-1} ([Z])$
be the fibre of $\pi$ over $[Z]$. Denote by
$\JAA^{(0)}_{Z}$
the complement of the theta-divisor
$$
\TE_{Z} = \TE (X;L,d) \cap \JAA_{Z} 
$$
in
$\JAA_{Z}$.

Then
\begin{enumerate}
\item[1)]
The rank $ rk(\HT_{-i})$ of the sheaf $\HT_{-i}$ in (\ref{filtHT-JG}) is constant along
$\JAA^{(0)}_{Z}$,
for all $i \geq 1$.
\item[2)]
$\HT_{-l_{\GA}} \otimes \OO_{\JAA^{(0)}_{Z}}$
is a trivial subbundle of
$\FT \otimes \OO_{\JAA^{(0)}_{Z}} = \HO{Z}) \otimes \OO_{\JAA^{(0)}_{Z}} $.
\end{enumerate}
\end{lem}
\begin{pf}
The first assertion is a restatement of Claim \ref{cl-ab}.
To see the second assertion we take two distinct points
$[\alpha]$ and $[\beta]$ in $\JAA^{(0)}_{Z}$
and we go back to the identity
$$
\HT ([Z],[\beta]) = \frac{\alpha}{\beta} \HT \ZA
$$
in (\ref{a/b1}). Write
$$
 \frac{\alpha}{\beta} =\frac{1}{\frac{\beta}{\alpha}} =\frac{1}{1+t}\,,
$$
where
$t= \frac{\beta}{\alpha} -1$ is in 
$\HT \ZA$. This gives the identity
$$
\HT ([Z],[\beta]) = \frac{1}{1+t} \HT \ZA\,.
$$
Hence every element
$h \in \HT ([Z],[\beta])$
can be written in the form
\begin{equation}\label{h=fr}
h = \frac{1}{1+t} s\,,
\end{equation}
for some $s \in \HT \ZA$. In particular, for
$\beta$ in a small neighborhood\footnote{in the complex topology of $\EZ$.} of $\alpha$,
we can expand (\ref{h=fr}) in a convergent power series
$$
h = \sum^{\infty}_{n=0} t^n s\,,
$$
where the terms of the series are in
$\HT_{-\LG} \ZA$, for all $n\geq 0$.
This implies that
$
\HT ([Z],[\beta]) \subset \HT_{-\LG} \ZA
$.
Since 
$\HT_{-\LG} \ZA$
is closed under the multiplication in
$\HO Z)$,
we obtain an inclusion
$$
\HT_{-\LG} ([Z],[\beta]) \subset \HT_{-\LG}  \ZA\,.
$$
By the first part of the lemma the dimensions of the two vector spaces are equal.
This yields an equality
\begin{equation}\label{eqHTl}
\HT_{-\LG} ([Z],[\beta]) = \HT_{-\LG}  \ZA\,, 
\end{equation}
for all $[\beta]$ in a small open neighborhood of $[\alpha]$.
Since $\JAA^{(0)}_{Z}$ is path connected, it follows that the equality
(\ref{eqHTl}) holds for all
$[\beta] \in \JAA^{(0)}_{Z}$.
\end{pf}
\begin{rem}\label{r-Zpr}
The subring
$\HT_{-\LG} \ZA$ of $\HO Z)$ has the following geometric meaning.

Recall the morphism
$$
\xymatrix@1{
{\kappa \ZA }: Z \ar[r]& {\PP( \HT \ZA^{\ast})} } 
$$
in (\ref{kappa}) and let
$Z^{\prime} (\alpha)$ be the image of
$\kappa \ZA$.
Then the space 
$\HT_{-\LG} \ZA$
is isomorphic to 
$\HO{Z^{\prime} (\alpha)})$ with the isomorphism given by the pullback by 
$\kappa \ZA $. More precisely, we have
\begin{equation}\label{Zpr}
Z^{\prime} (\alpha) = Spec (\HT_{-\LG} \ZA)
\end{equation}
and
$$
\xymatrix@1{
{(\kappa \ZA)^{\ast}: \HO{Z^{\prime} (\alpha)})} \ar[r]& {\HT_{-\LG} \ZA} }
$$
is an isomorphism.

From (\ref{Zpr}) and Lemma \ref{HTl=triv},2), it also follows that the scheme
$Z^{\prime} (\alpha)$ is independent of 
$[\alpha] \in \JAA^{(0)}_{Z}$. In the sequel it will be denoted by $Z^{\prime}$.
\end{rem}
\begin{cor}\label{cor-Fpr}
Let 
$\GA^{(0)}$
be as in \S\ref{ssec-G(0)} and let
$\ZD_{\GA^{(0)}}$
be the universal subscheme over 
$\GA^{(0)}$.

There exists a subsheaf
$\FF^{\prime}$
of 
$\FF \otimes \OO_{\GA^{(0)}}$
such that
\begin{equation}\label{pi-Fpr}
\pi^{\ast} \FF^{\prime} = \HT_{-\LG} \otimes \OO_{\JG^{(0)}}.
\end{equation}
Furthermore, 
$\FF^{\prime}$
is a subsheaf of subrings of $\FF \otimes \OO_{\GA^{(0)}}$
and one has the following factorization
\begin{equation}\label{Z-Zpr}
\xymatrix{
{\ZD_{\GA^{(0)}}} \ar[rr]^f  \ar[dr]_{p_2}& &{\ZD^{\prime}_{\GA^{(0)}}} \ar[dl]^{p^{\prime}_2} \\
&{\GA^{(0)}}& }
\end{equation}
where
$\ZD^{\prime}_{\GA^{(0)}} =Spec (\FF^{\prime} )$ and
$$
\xymatrix{
{f: \ZD_{\GA^{(0)}}} \ar[r] &{\ZD^{\prime}_{\GA^{(0)}}} }
$$
is the morphism corresponding to the inclusion of sheaves of rings
$ 
\FF^{\prime} \hookrightarrow \FF\otimes \OO_{\GA^{(0)}}
$. In particular, one has a canonical identification
$$
\FF^{\prime} = p^{\prime}_{2\ast} \OO_{\ZD^{\prime}_{\GA^{(0)}}}.
$$
\end{cor}
\begin{pf}
Set
\begin{equation}\label{d-pr-G}
d^{\prime}_{\GA} = rk(\HT _{-\LG} \otimes \OO_{\JG^{(0)}})
\end{equation}
and let
\begin{equation}\label{Gr-rel}
{\bf Gr}_{\GA^{(0)}} = Gr ( d^{\prime}_{\GA} , \FF\otimes \OO_{\GA^{(0)}})
\end{equation}
be the relative Grassmannian of
$d^{\prime}_{\GA}$-planes in 
$\FF\otimes \OO_{\GA^{(0)}}$.
We have the diagram
\begin{equation}\label{JG-Gr}
\xymatrix{
{\JG^{(0)}} \ar[rr]^{\tilde{\gamma}}  \ar[dr]_{\pi}& &{\bf Gr}_{\GA^{(0)}} \ar[dl]^{pr_{\GA^{(0)}}} \\
&{\GA^{(0)}}& }
\end{equation}
where the morphism $\tilde{\gamma}$ corresponds to the inclusion
$$
\HT _{-\LG} \otimes \OO_{\JG^{(0)}} \hookrightarrow \FT\otimes \OO_{\JG^{(0)}} = \pi^{\ast} ( \FF\otimes \OO_{\GA^{(0)}})\,.
$$
The morphism $\tilde{\gamma}$ sends a closed point
$\ZA$ of $\JG^{(0)}$ to the
$d^{\prime}_{\GA}$-plane
$\HT _{-\LG} \ZA$
of $\HO Z) = \FF ([Z])$, the fibre of 
$\FF$ at $[Z]$.
From Lemma \ref{HTl=triv}, 2), it follows that the morphism
$\tilde{\gamma}$
is constant along the fibres of $\pi$. Hence it factors through
$\GA^{(0)}$ yielding a section 
$$
\gamma: \GA^{(0)} \longrightarrow {\bf Gr}_{\GA^{(0)}}
$$
of the natural projection $pr_{\GA^{(0)}}$ in (\ref{JG-Gr}), i.e. we have
\begin{equation}\label{g-gt-pr}
\tilde{\gamma} = \gamma \circ \pi\,\,{\rm{and}}\,\, 
pr_{\GA^{(0)}} \circ \gamma = id_{\GA^{(0)}}.
\end{equation}

Let
${\cal U}$ be the universal subbundle of
$pr^{\ast}_{\GA^{(0)}} (\FF\otimes \OO_{\GA^{(0)}})$
on 
the relative Grassmannian 
${\bf Gr}_{\GA^{(0)}}$.
Applying $\gamma^{\ast}$ to the inclusion
$$
{\cal U} \hookrightarrow pr^{\ast}_{\GA^{(0)}} (\FF\otimes \OO_{\GA^{(0)}})
$$
gives the subbundle
\begin{equation}\label{def-Fpr}
\FF^{\prime} =\gamma^{\ast} {\cal U} \hookrightarrow \gamma^{\ast} \big{(} pr^{\ast}_{\GA^{(0)}} (\FF\otimes \OO_{\GA^{(0)}}) \big{)} =
(pr_{\GA^{(0)}} \circ \gamma)^{\ast} \big{(} \FF\otimes \OO_{\GA^{(0)}} \big{)} = \FF\otimes \OO_{\GA^{(0)}},
\end{equation}
where the last equality comes from the second identity in (\ref{g-gt-pr}).

Next we check that the sheaf 
$\FF^{\prime}$ is subject to the identity (\ref{pi-Fpr}) of the corollary.
For this we apply $\pi^{\ast}$ to the equalities in (\ref{def-Fpr}) to obtain
$$
\pi^{\ast} \FF^{\prime} =\pi^{\ast} (\gamma^{\ast} {\cal U}) = (\gamma \circ \pi)^{\ast} {\cal U} = {\tilde{\gamma}}^{\ast} {\cal U} =
\HT _{-\LG} \otimes \OO_{\JG^{(0)}},
$$
where the third equality comes from the first identity in (\ref{g-gt-pr}), while the last one follows from the definition of the morphism
$\tilde{\gamma}$ in (\ref{JG-Gr}).

From the equality 
$\pi^{\ast} \FF^{\prime} =\HT _{-\LG} \otimes \OO_{\JG^{(0)}}$
and the fact that the latter sheaf is a subsheaf of subrings of
$\FT \otimes \OO_{\JG^{(0)}}$ it follows that
$ \FF^{\prime}$ is a sheaf of subrings of
$ \FF\otimes \OO_{\GA^{(0)}}$.
Hence the monomorphism 
$$
\FF^{\prime} \hookrightarrow  \FF\otimes \OO_{\GA^{(0)}}
$$
of sheaves of rings defines a surjective morphism of schemes
$$
f: \ZD_{\GA^{(0)}}= Spec (\FF\otimes \OO_{\GA^{(0)}})  \longrightarrow  Spec (\FF^{\prime}) = \ZD^{\prime}_{\GA^{(0)}}
$$
over $\GA^{(0)}$. This yields the commutative diagram asserted in (\ref{Z-Zpr}).
\end{pf}

Recall the notation of 
$\GAB$ and $\JABG$ in \S\ref{breve}. By definition over 
$\JABG$ the orthogonal decomposition (\ref{ordFT}) holds.
In particular, we have
\BEN\label{l-ord}
\FT \otimes \OO_{\JABG} = \HT_{-\LG} \otimes \OO_{\JABG} \oplus \HH^{\LG}.
\EEN
From Corollary \ref{cor-Fpr} it follows easily that this decomposition is the pullback by $\pi$
of the orthogonal decomposition
\BEN\label{l-ord1}
\FF\otimes \OO_{\GAB} = \FF^{\prime}\otimes \OO_{\GAB} \oplus  \big{(}\FF^{\prime}\otimes \OO_{\GAB} \big{)}^{\perp}\,,
\EEN
where 
$\big{(}\FF^{\prime}\otimes \OO_{\GAB} \big{)}^{\perp}$
is the subsheaf of 
$\FF\otimes \OO_{\GAB}$
orthogonal to
$\FF^{\prime}\otimes \OO_{\GAB} $
with respect to the quadratic form $\QB$ in (\ref{q}), i.e. the following holds
\begin{eqnarray}
\HH^{\LG} &=& \pi^{\ast} \left(\big{(}\FF^{\prime}\otimes \OO_{\GAB} \big{)}^{\perp} \right)\,,\\ \label{Hl=pb}
\FT \otimes \OO_{\JABG} &=&\pi^{\ast} (\FF^{\prime}\otimes \OO_{\GAB} ) \oplus 
\pi^{\ast} \left(\big{(}\FF^{\prime}\otimes \OO_{\GAB} \big{)}^{\perp} \right)\,. \label{ordl=pb}
\end{eqnarray}
\begin{rem}\label{geom-ordl}
From the algebro-geometric perspective of the factorization in (\ref{Z-Zpr}) the above decomposition and its constituents can be described 
as follows.

Let $[Z]$ be a point in $\GAB$ and consider the diagram (\ref{Z-Zpr}) over it. This gives
the morphism
\BEN\label{f-Z}
f: Z=p^{-1}_2 ([Z]) \longrightarrow Z^{\prime} =p^{{\prime} -1}_2 ([Z]) 
\EEN
which on the level of the rings of functions translates into the injective homomorphism of rings
\BEN\label{f-Z-rings}
\FF^{\prime} ([Z]) =\HO{Z^{\prime}}) \longrightarrow \HO {Z} )  = \FF ([Z])
\EEN
given by the pullback $f^{\ast}$ of functions. In particular, in the decomposition 
\BEN\label{HoZ-dec}
  \HO {Z} )  = \HT_{-\LG} \ZA \oplus \HH^{\LG} \ZA
\EEN
the summand
$\HT_{-\LG} ([Z], [\alpha])$ is identified with 
$\HO{Z^{\prime}}) $ via the pullback in (\ref{f-Z-rings}), for every
$[\alpha]$ in $\JABG$ lying over $[Z]$. With this condition on $[\alpha]$, assumed for the rest of this discussion, the space
space $\HT_{-\LG} ([Z], [\alpha])$ (resp. $\HH^{\LG} ([Z], [\alpha])$) is independent of $[\alpha]$ and we denote it by
$\HT_{-\LG} ([Z])$ (resp. $\HH^{\LG} ([Z])$) (this is proved in Lemma \ref{HTl=triv},2)).

To describe the orthogonal complement
$\HH^{\LG} ([Z])$ of 
$\HT_{-\LG} ([Z])$ in
$\HO {Z} )$ we set
$Z_{z^{\prime}} =f^{-1}(z^{\prime})$,
for every closed point of 
$Z^{\prime}$. Then $Z$ admits the following decomposition
\BEN\label{Z-f-dec}
Z =\sum_{z^{\prime} \in Z^{\prime}} Z_{z^{\prime}}\,.
\EEN
Let 
$$
\delta_{Z_{z^{\prime}}} =\sum_{z \in Z_{z^{\prime}}} \delta_z = f^{\ast}(\delta_{z^{\prime}} )
$$
be the pullback of the delta-function on $Z^{\prime}$ supported at $z^{\prime}$.
These functions form a basis of
$\HT_{-\LG} ([Z], [\alpha])$ as $z^{\prime}$ runs through the closed points of 
$Z^{\prime}$. Furthermore, every $h \in \HO {Z} )$ can be written uniquely
$$
h= \sum_{z^{\prime} \in Z^{\prime}} h_{Z_{z^{\prime}}}\,,
$$
where for each $z^{\prime} \in Z^{\prime}$ the component
$ h_{Z_{z^{\prime}}} =h \delta_{Z_{z^{\prime}}}$ is supported on $Z_{z^{\prime}}$.
With these preliminaries in mind, we can now describe the space
$\HH^{\LG} ([Z])$ as follows
\BEN \label{Hl-Z}
\HH^{\LG} ([Z]) = \left\{ \left.h \in  \HO {Z} ) \right| Tr(h_{Z_{z^{\prime}}} ) = 0,\,\,\forall z^{\prime} \in Z^{\prime} \right\}\,,
\EEN
where $Tr$ stands for the trace morphism in (\ref{tr}).
\end{rem}
 
\section{The sheaf of Lie algebras {\boldmath$\GT_{\mbox{\unboldmath$\GA$}}$}} \label{sec-Lie1}
In this section we establish the basic properties of the sheaf
{\boldmath
$\LAGT$}
(see \S\ref{Lie} for its definition).
From \RI, Proposition 7.2, we know that it is a sheaf of reductive Lie algebras. By the structure theorem of reductive Lie algebras 
(see e.g. \cite{[Bour]}) one obtains
the following decomposition
{\boldmath
\begin{equation}\label{str-d}
\LAGT = \CG \oplus \LAG,
\end{equation}
where
$\CG$ is the center of
$\LAGT$ and
$\LAG =[\LAGT, \LAGT]$
}
is a sheaf of semisimple Lie algebras.

By definition,
{\boldmath
$\LAGT$}
comes together with a faithful representation on
$\HT_{-\LG}$,
i.e.
{\boldmath
$\LAGT$}
is defined as a subsheaf of
$\ENDO(\HT_{-\LG})$.

One of the features of
{\BM
$\LAGT$}
is that it comes along with a distinguished subsheaf of Cartan subalgebras determined by the image
of the morphism
$D$ defined in (\ref{m-D}). Namely, define the subsheaf
{\BM
${\cal C} (\HT)$}
of
{\BM
$\LAGT$}
to be the centralizer of
$D(\HT)$, i.e.
a local section $x$ of 
{\BM
$\LAGT$}
belongs to 
{\BM
${\cal C} (\HT)$}
if and only if
\begin{equation}\label{c-CHT}
[x, D(t)] =0, \,\,{\rm{for\,\, any\,\, local\,\,section\,\,}} t \,\,{\rm{of}}\,\, \HT.
\end{equation}
\begin{pro}\label{Cartan}
{\BM
${\cal C} (\HT)$}
is a subsheaf of Cartan subalgebras of 
{\BM
$\LAGT$}.
\end{pro}
\begin{pf}
This is the result of \RI, Lemma 7.5.
\end{pf}
\begin{rem}\label{act-m}
The proof of Lemma 7.5 in \RI \,\, implies that
{\BM
${\cal C} (\HT)$}
can be naturally identified with a subsheaf of
$\pi^{\ast} \FF^{\prime}$,
where $\FF^{\prime}$ is as in Corollary \ref{cor-Fpr}.
Furthermore, it acts on
$$
\pi^{\ast} \FF^{\prime} = \HT_{-\LG}
$$
via the multiplication in $ \HT_{-\LG}$.
\end{rem}
The sheaf
{\BM
${\cal C} (\HT)$}
decomposes according to the structure decomposition
(\ref{str-d})
{\BM
\begin{equation}\label{C-d}
{\cal C} (\HT) =\CG \oplus {\cal H}_{\mbox{\UB$\GA$}},
\end{equation}
where
\begin{equation}\label{CA-def}
{\cal H}_{\mbox{\UB$\GA$}} = {\cal C} (\HT) \cap \LAG
\end{equation}
is a subsheaf of Cartan subalgebras of
$
\LAG
$.

\subsection{The center $\CG$ of  $\LAGT$}\label{Center}
We will investigate the action of the center
$\CG$
}
on the sheaf
$\HT_{-\LG}$. 
To begin with consider the situation fibrewise.

Fix
$\ZA \in \JAB_{\GA}$ and let
\begin{equation}\label{cen-fb}
{\bf c} = \mbox{\BM${\cal C}$} \ZA
\end{equation}
be the fibre of 
{\BM
$\CG$}
at $\ZA$.

Consider the affine version of the morphism
$\kappa\ZA$ in (\ref{kappa}):
$$
Z \longrightarrow  {\HT\ZA}^{\ast}
$$
for which we use the same notation.
Denote by 
\begin{equation}\label{def-Zpr}
Z^{\prime} = Spec(\HT_{-\LG} \ZA) = Spec( \FF^{\prime} ([Z])),
\end{equation}
where
$\FF^{\prime} ([Z])$ denotes the fibre of 
$\FF^{\prime} $ at $[Z]$ and where 
the last equality in (\ref{def-Zpr}) comes from the identity
(\ref{pi-Fpr}) in Corollary \ref{cor-Fpr}. In particular, we have a natural identification
\begin{equation}\label{f-Zpr=HTl}
\HT_{-\LG} \ZA \cong \HO{Z^{\prime} }).
\end{equation}
In the sequel we use freely this identification, by switching frequently from one space to another,
without explicitly invoking this isomorphism. Thus, for example, the fact that the action of 
{\BM
$\LAGT$}$\ZA$
 on 
$\HT_{-\LG} \ZA$ 
implies the action on 
$\HO{Z^{\prime} })$ will be taken for granted.

Consider the action of the center
${\bf c}$
on
$\HO{Z^{\prime} })$. It is known that this action is semisimple.
Thus we obtain the following weight decomposition
\begin{equation}\label{Zpr-wd}
\HO{Z^{\prime} }) = \bigoplus_{\lambda \in {\bf c}^{\ast}} V_{\lambda} \ZA,
\end{equation} 
where 
$V_{\lambda} \ZA$ is the weight space corresponding to a weight $\lambda$ and the direct sum is taken over the weights
of this action.
\begin{pro}\label{pro-wd}
\begin{enumerate}
\item[1)]
The weight spaces
$V_{\lambda} \ZA$
are ideals in 
$\HO{Z^{\prime} }) $,
\item[2)]
$V_{\lambda} \ZA  \cdot V_{\mu} \ZA =0$,
for any two weights $\lambda \neq \mu$
occurring in (\ref{Zpr-wd}).
\end{enumerate}
\end{pro}
\begin{pf}
1) By Remark \ref{act-m} the center ${\bf c}$ can be viewed as a subspace of
$\HO{Z^{\prime} }) $
and its action is identified with the multiplication in
$\HO{Z^{\prime} }) $. Hence we have
\begin{equation}\label{act-com}
h(fv)=f(hv),
\end{equation}
for any
$f,v \in \HO{Z^{\prime} }) $ and any
$h \in {\bf c}$, where the operation in (\ref{act-com}) is
the multiplication in 
$\HO{Z^{\prime} }) $.
In particular,
if 
$v \in V_{\lambda} \ZA$
we obtain
$$
h(fv)=\lambda (h) (fv),
$$
for all $h \in {\bf c}$. Hence $fv \in V_{\lambda} \ZA$, for all 
$f \in \HO{Z^{\prime} }) $. This shows that
$V_{\lambda} \ZA$ is an ideal in 
$\HO{Z^{\prime} }) $.

2) Take 
$v \in V_{\lambda} \ZA$ and 
$v^{\prime} \in V_{\mu} \ZA$ 
and consider the product
$vv^{\prime}$ in
$\HO{Z^{\prime} }) $.
Applying 
$h \in {\bf c}$ to
$vv^{\prime}$ gives
$$
h(vv^{\prime}) =( hv)v^{\prime}= \lambda (h) vv^{\prime}
$$
as well as
$$
h(vv^{\prime}) = v(hv^{\prime})= \mu (h) vv^{\prime}.
$$
These two equalities yield
$$
(\mu (h) -\lambda (h)) vv^{\prime} =0,\,\,\forall h\in{\bf c}.
$$
If $\lambda \neq \mu$, then the above implies
$$
v v^{\prime} =0,\,\,\forall v\in V_{\lambda} \ZA ,\,\, \forall v^{\prime} \in V_{\mu} \ZA.
$$
This proves the second assertion of the proposition.
\end{pf}
\begin{lem}\label{lem-Vlam}
Let
$Z^{\prime (\lambda, \alpha)}$
be the subscheme of 
$Z^{\prime }$
corresponding to the ideal
$V_{\lambda} \ZA$
and let
$Z^{\prime}_{ (\lambda, \alpha)}$
be its complement in 
$Z^{\prime}$.
Then
\begin{equation}\label{Vlam}
V_{\lambda} \ZA = \bigoplus_{p^{\prime}} \CC \delta_{p^{\prime}},
\end{equation}
where the sum is taken over the distinct closed points of 
$Z^{\prime}_{ (\lambda, \alpha)}$ and 
$ \delta_{p^{\prime}}$
is the delta-function supported at 
$p^{\prime}$.

Furthermore, the elements of ${\bf c}$, viewed as functions on
$Z^{\prime}$, are constant on $Z^{\prime}_{ (\lambda, \alpha)}$.
More precisely, if $h\in {\bf c}$, then
$$
h(p^{\prime}) = \lambda (h),
$$
for all closed points $p^{\prime}$ in 
$Z^{\prime}_{ (\lambda, \alpha)}$.
\end{lem}
\begin{pf}
The functions in
$V_{\lambda} \ZA$ must vanish on 
$Z^{\prime (\lambda, \alpha)}$.
Hence their support is in
$Z^{\prime}_{ (\lambda, \alpha)}$.
This implies an inclusion
$$
V_{\lambda} \ZA \subset  \bigoplus_{p^{\prime}} \CC \delta_{p^{\prime}},
$$
where the sum is taken over the distinct closed points of 
$Z^{\prime}_{ (\lambda, \alpha)}$.
On the other hand every
$ \delta_{p^{\prime}}$, the delta function having support at a closed point
$p^{\prime}$ of 
$Z^{\prime}_{ (\lambda, \alpha)}$, vanishes on
$Z^{{\prime}{ (\lambda, \alpha)}}$ and hence belongs to
$V_{\lambda} \ZA$.
This proves the equality (\ref{Vlam}).

To prove the second assertion consider the action of
$h \in {\bf c}$ on
$\delta_{p^{\prime}}$'s in (\ref{Vlam})
\begin{equation}\label{act-del}
h(\delta_{p^{\prime}}) = \lambda(h) \delta_{p^{\prime}}.
\EEN
On the other hand since $h$ acts on 
$\HO{Z^{\prime}})$
by multiplication the left hand side in (\ref{act-del}) can be written as follows
$$
h(\delta_{p^{\prime}}) = h \delta_{p^{\prime}} = h(p^{\prime}) \delta_{p^{\prime}}.
$$
This and (\ref{act-del}) imply
$$
\lambda(h) = h (p^{\prime}), 
$$
for every closed point of
$Z^{\prime}_{ (\lambda, \alpha)}$.
\end{pf}

The ring
$\HT_{-\LG} \ZA$,
viewed as a subring of
$\HO Z)$ (see (\ref{filtHT-JG})), is generated by the
the subspace
$\HT \ZA$.
So it is natural to ask for a relation of the weight spaces
$V_{\lambda} \ZA$
in (\ref{Zpr-wd}) and
$\HT \ZA$.
The following proposition answers this question.
\begin{pro} \label{HT-wd}
Set
$$
 \HT_{\lambda}  \ZA  = V_{\lambda} \ZA \bigcap \HT \ZA.
$$
Then
$$
\HT \ZA = \bigoplus_{\lambda \in {\bf c}^{\ast}} \HT_{\lambda}  \ZA 
$$
and
$$
 V_{\lambda} \ZA = im \left( S^{\bullet} \big{(} \HT_{\lambda}  \ZA \big{)} \longrightarrow  \HT_{-\LG} \ZA =\HO {Z^{\prime}}) \right)\,, 
$$
for every weight $\lambda$ occurring in the decomposition (\ref{Zpr-wd}).
\end{pro}
\begin{pf}
Let 
$t \in \HT \ZA$. Decompose it according to the weight decomposition in (\ref{Zpr-wd}):
\BEN \label{t-wd}
t= \sum_{\lambda \in {\bf c}^{\ast}} t_{\lambda},
\EEN
where 
$t_{\lambda}$ is the component of $t$ in
$V_{\lambda} \ZA$, for every weight 
$\lambda$ occurring in (\ref{Zpr-wd}).
We claim that each 
$ t_{\lambda}$
is in 
$\HT_{\lambda}  \ZA $.
To see this recall that the elements of
$\HT$ are characterized by the property of being annihilated by the operators
$D^{-} (t^{\prime})$ in the triangular decomposition (\ref{d-Dt}), for any local section $t^{\prime}$ of $\HT$
(see \RI, Remark 7.8).
Applying
$D^{-} (t^{\prime})$
to the both sides in (\ref{t-wd}) yields
\BEN \label{D-t-wd}
0=D^{-} (t^{\prime}) (t) = \sum_{\lambda \in {\bf c}^{\ast}}  D^{-} (t^{\prime}) (t_{\lambda}).
\EEN
But by definition
$D^{-} (t^{\prime}) $ are in
{\BM $\LAGT$}$\ZA$, the fibre
of
 {\BM $\LAGT$} at $\ZA$, for all 
$t^{\prime}$ in
$\HT \ZA$.
Hence 
$D^{-} (t^{\prime}) $ commute with
${\bf c}$,  for all 
$t^{\prime}$ in
$\HT \ZA$. This implies
that
$D^{-} (t^{\prime}) (t_{\lambda}) \in  V_{\lambda} \ZA$,
for every $\lambda$ in the sum of (\ref{D-t-wd}).
Combining this with the equation (\ref{D-t-wd}) yields
$$
D^{-} (t^{\prime}) (t_{\lambda}) =0,
$$
for every 
$t^{\prime} \in \HT\ZA$.
Using \RI, Remark 7.8, once again, we obtain that
$t_{\lambda} \in \HT\ZA$,
for every $\lambda$.

Turning to the second assertion we use the fact that
$V_{\lambda} \ZA$
is an ideal in
$\HT_{-\LG} \ZA$ and hence closed under the multiplication.
This gives an inclusion
\BEN \label{inc-Vlam}
im \left( S^{\bullet} \left( \HT_{\lambda}  \ZA \right) \longrightarrow  \HT_{-\LG} \ZA  \right)  \subset V_{\lambda} \ZA.
\EEN
On the other hand
$$
\HT_{-\LG} \ZA =im \left( S^{\bullet} \left( \HT  \ZA \right) \longrightarrow  \HT_{-\LG} \ZA  \right) 
$$
and from the weight decomposition (\ref{Zpr-wd}) and Proposition \ref{pro-wd}, 2), it follows
\begin{eqnarray*}
(i)\,\,S^{\bullet} \big{(} \HT\ZA \big{)} = &\bigotimes_{\lambda}  S^{\bullet} \big{(} \HT_{\lambda}  \ZA \big{)} \\
(ii)\,\,\,\HT_{-\LG} \ZA = & 
im \left( \bigotimes_{\lambda}  S^{\bullet} \big{(} \HT_{\lambda}  \ZA \big{)} \longrightarrow  \HT_{-\LG} \ZA  \right) \\
& =\bigoplus_{\lambda} im \left( S^{\bullet} \big{(} \HT_{\lambda}  \ZA \big{)} \longrightarrow  \HT_{-\LG} \ZA  \right)\,.
\end{eqnarray*}
This together with (\ref{inc-Vlam}) yield the asserted equality
$$
V_{\lambda} \ZA = im \left( S^{\bullet} \left( \HT_{\lambda}  \ZA \right) \longrightarrow  \HT_{-\LG} \ZA  \right)\,. 
$$
\end{pf}

By Remark \ref{act-m} we can identify 
${\bf c}$ with a subspace of 
$\HT_{-\LG} \ZA = \HO{Z^{\prime}})$.
Our next task will be to locate the center
${\bf c}$ as a subspace of 
$\HO{Z^{\prime}})$. 
\begin{pro}\label{pro-c-basis}
Let
$$
\delta_{Z^{\prime}_{(\lambda, \alpha)}} = \sum_{p^{\prime} \in Z^{\prime}_{(\lambda, \alpha)}} \delta_{p^{\prime}},
$$
where  the sum is taken over the closed points of the subscheme
$Z^{\prime}_{(\lambda, \alpha)}$ of $Z^{\prime}$
defined in Lemma \ref{lem-Vlam}, and let
${\bf {\tilde c}}$ be the span of these functions
in
$\HO{Z^{\prime}}) =\HT_{-\LG} \ZA $.
Then
\begin{enumerate}
\item[1)]
${\bf {\tilde c}} \subset \HT\ZA $.
\item[2)]
$
{\bf {\tilde c}} = \bigoplus_{\lambda} \CC  \delta_{Z^{\prime}_{(\lambda, \alpha)}}.
$
\item[3)]
$ {\bf {\tilde c}}$
is a subring of  $\HT_{-\LG} \ZA $.
\item[4)]
The morphism $D$ in (\ref{m-D})
identifies  
${\bf {\tilde c}}$ with
${\bf c}$.
\end{enumerate}
\end{pro}

\begin{pf}

To prove the the first assertion it is enough to show that
$ \delta_{Z^{\prime}_{(\lambda, \alpha)}}$
belongs to  $ \HT\ZA $,
for every weight $\lambda$ occurring in the decomposition (\ref{Zpr-wd}).
For this observe that Proposition \ref{pro-wd}, 2), implies
$$
Z^{\prime}_{(\lambda, \alpha)} \bigcap  Z^{\prime}_{(\mu, \alpha)} = \emptyset,
$$
for all $\lambda \neq \mu$. This and the direct sum
(\ref{Zpr-wd}) give the decomposition of 
 $Z^{\prime}$ into the disjoint union
$$
Z^{\prime} = \bigcup_{\lambda}  Z^{\prime}_{(\lambda, \alpha)} 
$$
over the weights in (\ref{Zpr-wd}). In particular, the constant function
$1 \in \HO{Z^{\prime}})$ can be written as follows
$$
1 = \sum _{\lambda}   \delta_{Z^{\prime}_{(\lambda, \alpha)}}.
$$
This is the weight decomposition of $1$ because
$ \delta_{Z^{\prime}_{(\lambda, \alpha)}} \in V_{\lambda} \ZA$,
for every
$\lambda$. Furthermore, the constant
$1$ lies in $\HT \ZA$ ( see \RI, Remark 1.3) 
and by Proposition \ref{HT-wd} its $\lambda$-components
$\delta_{Z^{\prime}_{(\lambda, \alpha)}} \in \HT_{\lambda} \ZA$,
for all $\lambda$.

The argument above also shows that the family
of functions
$\delta_{Z^{\prime}_{(\lambda, \alpha)}}$, as 
$\lambda$
runs through the distinct weights in (\ref{Zpr-wd}), is linear independent
in $\HT \ZA$. This yields the second assertion.

From  the above it also follows that $1$ is contained in ${\bf {\tilde c}}$, while the identities
$$
 \delta_{Z^{\prime}_{(\lambda, \alpha)}} \delta_{Z^{\prime}_{(\mu, \alpha)}} =0,\,\,{\rm {for}} \,\, \lambda \neq \mu,
\,\,\delta^2_{Z^{\prime}_{(\lambda, \alpha)}} =\delta_{Z^{\prime}_{(\lambda, \alpha)}}
$$
assure that
${\bf {\tilde c}}$
is closed under the multiplication. This proves the third assertion.

Turning to the last assertion we observe that by the second assertion of
Lemma \ref {lem-Vlam} the center
${\bf c}$ is identified with a subspace of 
the space
${\bf {\tilde c}}$.
So it will be enough to check that
the morphism
$D$ takes the functions 
$\delta_{Z^{\prime}_{(\lambda, \alpha)}}$
to
${\bf c}$.
Indeed,
 the operator
$D(\delta_{Z^{\prime}_{(\lambda, \alpha)}})$ is the operator of multiplication
by $\delta_{Z^{\prime}_{(\lambda, \alpha)}}$ in the ring
$\HO Z^{\prime})$
and it acts as the identity on
$V_{\lambda} \ZA$ and by zero on all other
weight spaces in (\ref{Zpr-wd}). Hence 
$D(\delta_{Z^{\prime}_{(\lambda, \alpha)}})$ belongs to the center
${\bf c}$,
for every $\lambda$ in (\ref{Zpr-wd}).
\end{pf}
Set
\BEN \label{J-Z}
\JAB_{Z} = \pi^{-1} ([Z])
\EEN
to be the fibre of $\JAB_{\GA}$ over $[Z] \in \breve{\GA}$,
where
$\pi$ and $\breve{\GA}$ are as defined in (\ref{GAb}).
We will now show that the weight decomposition (\ref{Zpr-wd}) does not depend on
$[\alpha] \in \JAB_{Z}$.
\begin{lem}\label{c-Z}
The center
${\bf c} =\mbox{\BM${\cal{C}}$}_{\GA} \ZA$
and the weight decomposition (\ref{Zpr-wd})
do not depend on 
$[\alpha] \in \JAB_{Z}$.
\end{lem}
\begin{pf}
First observe that the independence of the center on a point in 
$\JAB_{Z}$ implies the same for the weight decomposition (\ref{Zpr-wd}).
So only the first assertion needs to be proved. In order to do this we take
 $[\alpha]$ and $[\beta]$ to be two distinct points of  $ \JAB_{Z}$
and let
${\bf c}_{[\alpha]}$ and ${\bf c}_{[\beta]}$
be the centers of the Lie algebras
{\BM
$\LAGT$}$\ZA$
and
{\BM
$\LAGT$}$([Z],[\beta])$,
respectively.
To relate the two centers we use their explicit description as subspaces
of 
$\HT\ZA$ and $\HT([Z],[\beta])$, respectively, obtained in
Proposition \ref{pro-c-basis}. So our first step will be to relate
the weight decomposition of
$\HO {Z^{\prime}})$ under  
the ${\bf c}_{[\alpha]}$-action with the subspace
$\HT([Z],[\beta])$
(recall that by Corollary \ref{cor-Fpr} we have
$\HO {Z^{\prime}})= \HT_{-\LG} \ZA = \HT_{-\LG} ([Z],[\beta])$).
For this we use the identity
$$
\HT ([Z],[\beta]) = \frac{\alpha}{\beta} \HT \ZA
$$
in (\ref{a/b1}). Putting it together with the ${\bf c}_{[\alpha]}$-weight decomposition of
$\HT \ZA$ in Proposition \ref{HT-wd}, we obtain
\BEN \label{HTb-Va}
\HT ([Z],[\beta]) = \bigoplus_{\lambda} \frac{\alpha}{\beta} \HT_{\lambda} \ZA.
\EEN
Since $V_{\lambda} \ZA$ is an ideal ( Proposition \ref{pro-wd}, 1)), it follows that
$$
\frac{\alpha}{\beta} \HT_{\lambda} \ZA \subset V_{\lambda} \ZA.
$$
This and (\ref{HTb-Va}) give an inclusion
$$
\HT ([Z],[\beta]) \bigcap  V_{\lambda} \ZA \supset \frac{\alpha}{\beta} \HT_{\lambda} \ZA,
$$
for every weight $\lambda$ occurring in (\ref{HTb-Va}). However, the equality in
(\ref{HTb-Va}) implies that the above inclusion is actually an equality
\BEN \label{HTb-int-Va}
\HT ([Z],[\beta]) \bigcap  V_{\lambda} \ZA  = \frac{\alpha}{\beta} \HT_{\lambda} \ZA,
\EEN
for every weight $\lambda$ in the decomposition (\ref{HTb-Va}) .

The equality in (\ref{HTb-int-Va}) generalizes to the entire filtration
$\HT_{-\bullet} ([Z],[\beta]) $. In fact, set
$$
\big{(} \HT_{\lambda} \ZA \big{)}_{-i} = im \big{(}  S^i \big{(} \HT_{\lambda} \ZA \big{)} \longrightarrow \HO {Z^{\prime}}) \big{)}. 
$$
From Proposition \ref{HT-wd} it follows that
$$
\big{(} \HT_{\lambda} \ZA \big{)}_{-i} =  V_{\lambda} \ZA \bigcap \HT_{-i} \ZA.
$$
This together with the  weight decomposition (\ref{Zpr-wd}) yield
\BEN \label{HT-i-wd}
 \HT_{-i} \ZA = \bigoplus_{\lambda} \big{(} \HT_{\lambda} \ZA \big{)}_{-i}\,\,.
\EEN
Using the identity
$$
\HT_{-i} ([Z],[\beta]) = \left(\frac{\alpha}{\beta} \right)^{i} \HT_{-i} \ZA
$$
in (\ref{a/b2}) together with the decomposition in (\ref{HT-i-wd}), we deduce the identity analogous to the one in (\ref{HTb-int-Va}):
\BEN \label{HTbi-int-Va}
\HT_{-i} ([Z],[\beta]) \bigcap V_{\lambda} \ZA =\left(\frac{\alpha}{\beta} \right)^{i} \big{(} \HT_{\lambda} \ZA \big{)}_{-i}\,\,,
\EEN
for every $i\geq 1$. This yields the decomposition
\BEN \label{HTbi-int-Va1}
\HT_{-i} ([Z],[\beta]) = \bigoplus_{\lambda} \HT_{-i} ([Z],[\beta]) \bigcap V_{\lambda} \ZA.
\EEN
With these considerations accomplished, we proceed to relating the centers
${\bf c}_{[\alpha]}$ and ${\bf c}_{[\beta]}$. 

From Proposition \ref{pro-c-basis} we know that 
${\bf c}_{[\alpha]}$
is spanned by the operators
$D(\delta_{Z^{\prime}_{(\lambda,\alpha)}})$ as 
$\lambda$ runs through the weights of the decomposition in (\ref{Zpr-wd}).
So to show the inclusion
${\bf c}_{[\alpha]} \subset {\bf c}_{[\beta]}$
it would be enough to prove 
\BEN \label{bas-a-in-cb}
D(\delta_{Z^{\prime}_{(\lambda,\alpha)}}) \in  {\bf c}_{[\beta]},\,\, \forall \lambda.
\EEN
 Exchanging the roles of $\alpha$ and $\beta$ will give the opposite inclusion and hence the
equality
${\bf c}_{[\alpha]} = {\bf c}_{[\beta]}$.
Thus our argument will be completed once we show (\ref{bas-a-in-cb}).

{\it Proof of  (\ref{bas-a-in-cb}) :}
take the element
$\frac{\beta}{\alpha} \in \HT\ZA$ and write it according to the weight decomposition in (\ref{HT-wd}):
\BEN \label{b/a-wd}
\frac{\beta}{\alpha} = \sum_{\lambda} \left(\frac{\beta}{\alpha} \right)_{\lambda},
\EEN
where $\displaystyle{\left(\frac{\beta}{\alpha} \right)_{\lambda}} \in \HT_{\lambda} \ZA$ is the $\lambda$-component of
$\displaystyle{\frac{\beta}{\alpha}}$, for every weight 
$\lambda$ occurring in the decomposition (\ref{Zpr-wd}).

Multiplying (\ref{b/a-wd}) by 
$\frac{\alpha}{\beta}$ we obtain
$$
1 =   \sum_{\lambda}  \frac{\alpha}{\beta} \left(\frac{\beta}{\alpha} \right)_{\lambda}\,.
$$
From the identity in (\ref{HTb-int-Va}) it follows that 
$\displaystyle{
\frac{\alpha}{\beta} \left(\frac{\beta}{\alpha} \right)_{\lambda} \in \HT ([Z],[\beta]) \bigcap V_{\lambda} \ZA
}$.
On the other hand, from the proof of Proposition \ref{pro-c-basis} we know that the 
$\lambda$-components of $1$ are
$\delta_{Z^{\prime}_{(\lambda,\alpha)}}$'s. Hence  the equality
$$
\delta_{Z^{\prime}_{(\lambda,\alpha)}} = \frac{\alpha}{\beta} \left(\frac{\beta}{\alpha} \right)_{\lambda},
$$
for every weight $\lambda$ in (\ref{Zpr-wd}). This implies that the elements 
$\delta_{Z^{\prime}_{(\lambda,\alpha)}}$ are all in
$\HT ([Z],[\beta]) $. Hence the operators
$D(\delta_{Z^{\prime}_{(\lambda,\alpha)}})$ are
in the Lie algebra
{\BM
$\LAGT$}$ ([Z],[\beta])$.

It remains to see that
$D(\delta_{Z^{\prime}_{(\lambda,\alpha)}})$ are central
in the Lie algebra
{\BM
$\LAGT$}$ ([Z],[\beta])$.
This is done by examining their action on the filtration
$\HT_{-\bullet} ([Z],[\beta])$.

Recall that for every $\lambda$ the operator of multiplication
$D(\delta_{Z^{\prime}_{(\lambda,\alpha)}})$
acts on 
$V_{\lambda} \ZA$ as the identity
$id_{V_{\lambda} \ZA}$ and by zero on all other weight spaces in (\ref{Zpr-wd}).
In view of the decomposition in 
(\ref{HTbi-int-Va1}) this implies that the multiplication by 
$ \delta_{Z^{\prime}_{(\lambda,\alpha)}}$
preserves
$\HT_{-i} ([Z],[\beta])$, for every $i\geq 1$.
Hence 
$D^{+} (\delta_{Z^{\prime}_{(\lambda,\alpha)}}) =0$
in
{\BM
$\LAGT$}$ ([Z],[\beta])$.
Since $D^{-} (\cdot)$ and $D^{+} (\cdot)$ are adjoint to each other with respect to the quadratic form
$ \QT$ defined in (\ref{q}) (see \RI, Lemma 3.7), it follows that
$D^{-} (\delta_{Z^{\prime}_{(\lambda,\alpha)}}) =0$
in
{\BM
$\LAGT$}$ ([Z],[\beta])$ as well. Hence we obtain
\BEN \label{D-D0}
D(\delta_{Z^{\prime}_{(\lambda,\alpha)}}) = D^{0} (\delta_{Z^{\prime}_{(\lambda,\alpha)}})
\EEN
in
{\BM
$\LAGT$}$ ([Z],[\beta])$. From this identity it follows that
$D(\delta_{Z^{\prime}_{(\lambda,\alpha)}})$
is central in
{\BM
$\LAGT$}$ ([Z],[\beta])$.
Indeed, by definition the Lie algebra
{\BM
$\LAGT$}$ ([Z],[\beta])$
is generated by elements
$D^{\pm} (t),\,D^{0} (t)$ as $t$ varies in the vector space
$\HT ([Z],[\beta])$ and we have the triangular decomposition
$$
D(t) =D^{-} (t) + D^{0} (t) + D^{+} (t)
$$
as recalled in \S\ref{Lie}. From the commutativity of the multiplication in
$\HO{Z^{\prime}})$ it follows
\begin{eqnarray*}
\lefteqn{0 = [D(\delta_{Z^{\prime}_{(\lambda,\alpha)}}) , D(t) ] =[D^{0} (\delta_{Z^{\prime}_{(\lambda,\alpha)}}), D(t)]=}\\
&[D^{0} (\delta_{Z^{\prime}_{(\lambda,\alpha)}}), D^{-}(t)] + [D^{0} (\delta_{Z^{\prime}_{(\lambda,\alpha)}}), D^0(t)] +
[D^{0} (\delta_{Z^{\prime}_{(\lambda,\alpha)}}), D^{+}(t)],
\end{eqnarray*}
where on the right hand side we have the sum of operators of degree $-1$, $0$, $1$ with respect to the grading 
provided by the orthogonal decomposition in (\ref{ordFT}). This implies that the component of each degree
on the right hand side vanishes.
Thus 
$ D(\delta_{Z^{\prime}_{(\lambda,\alpha)}}) = D^{0} (\delta_{Z^{\prime}_{(\lambda,\alpha)}})$ commutes
with
$D^{\pm} (t), D^{0} (t)$, for every 
$t \in \HT ([Z],[\beta])$.
Hence 
$ D(\delta_{Z^{\prime}_{(\lambda,\alpha)}}) $ is in the center
${\bf c}_{[\beta]}$.
This completes the proof of (\ref{bas-a-in-cb}) and of the proposition.
\end{pf}
\begin{rem}\label{Zpr-lam}
From Lemma \ref{c-Z} it follows that the subschemes 
$Z^{\prime (\lambda,\alpha)}$ (resp. $Z^{\prime}_ {(\lambda,\alpha)}$)
introduced in Lemma \ref{lem-Vlam} are independent of 
$[\alpha]$ varying in the fibre $\JAB_Z$ of $\JABG$ over $[Z]$. So from now on we denote them
by
$Z^{{\prime} {\lambda}}$ (resp. $Z^{\prime}_ {\lambda}$).The same goes for the weight spaces
$V_{\lambda} \ZA$ in (\ref{Zpr-wd}) - they will be denoted $V_{\lambda} ([Z])$.
With this notation the weight decomposition (\ref{Zpr-wd}) takes the form
\BEN\label{Zpr-wd2}
\HO{Z^{\prime}}) = \bigoplus_{\lambda \in {\bf c}^{\ast}} V_{\lambda} ([Z]).
\EEN 

Analogously, the functions
$\delta_{Z^{\prime}_ {(\lambda,\alpha)}}$ defined in Proposition \ref{pro-c-basis}
will be denoted by 
$\delta_{Z^{\prime}_ {\lambda}}$.
From Proposition \ref{pro-c-basis}
the operators
$D(\delta_{Z^{\prime}_ {\lambda}})$ form a basis of the center ${\bf c}$. Furthermore,
$D(\delta_{Z^{\prime}_ {\lambda}})$ acts on the decomposition in (\ref{Zpr-wd2}) as
the identity $id_{V_{\lambda} ([Z])}$ on 
$V_{\lambda} ([Z])$ and by zero on all other summands.
\end{rem}

We will now sheafify the above results.
In what follows we use an additional hypothesis that
the rank of the sheaf of centers
{\BM
$\CG$}
is constant\footnote{the rank of {\BM
$\CG$}
is {\it a priori} constant on some non-empty Zariski open subset of $\JABG$.}
 on
$\JABG$.
After a complete description of the sheaf
{\BM
$\LAGT$}
(see Corollary \ref{C-Lie=lf}) this assumption will be superfluous.
\begin{pro}\label{sh-c}
Assume the rank of
{\BM
$\CG$}
to be constant on 
$\JABG$.
Then there exists a subsheaf
$\FF^{\prime}_c $
of the sheaf 
$\FF^{\prime} \otimes \OO_{\GAB}$ ( the sheaf $\FF^{\prime}$ is defined in Corollary \ref{cor-Fpr})
such that its pullback
$\pi^{\ast} \FF^{\prime}_c$
is a subsheaf of 
$\HT$ which is identified with the center
{\BM
$\CG$}
via the morphism
$D$ in (\ref{m-D}), i.e.
$$
D(\pi^{\ast} \FF^{\prime}_c) =\mbox{\BM${\cal C}$}_{\GA}.
$$
Furthermore, $\FF^{\prime}_c $ is a subsheaf of subrings of
$\FF^{\prime} \otimes \OO_{\GAB}$.
\end{pro}
\begin{pf}
From Proposition \ref{pro-c-basis}, 4) it follows that there is a subsheaf
$\FT^{\prime}_c$ of
$\HT$ such that
\BEN \label{FT-c}
D( \FT^{\prime}_c) = \mbox{\BM${\cal C}$}_{\GA}.
\EEN
Furthermore, by Lemma \ref{c-Z}
the restriction of $\FT^{\prime}_c$ to
fibres of the morphism
$$
\pi :\JABG \longrightarrow \GAB
$$
is the trivial bundle. So we are in the situation analogous to the one in Corollary \ref{cor-Fpr}.
Arguing as in the proof there, we obtain that there is a subsheaf
$\FF^{\prime}_c $ of $\FF^{\prime} \otimes \OO_{\GAB}$ such that
$$
\pi^{\ast} (\FF^{\prime}_c ) =\FT^{\prime}_c\,.
$$
This together with (\ref{FT-c}) give the equality
$$
D(\pi^{\ast} \FF^{\prime}_c) =\mbox{\BM${\cal C}$}_{\GA}.
$$
From Proposition \ref{pro-c-basis} it also follows that
$\FT^{\prime}_c= \pi^{\ast} (\FF^{\prime}_c )$ is a subsheaf of subrings
of $\FT^{\prime} =\pi^{\ast} (\FF^{\prime} \otimes \OO_{\GAB})$. This yields that
$\FF^{\prime}_c$ is a subsheaf of subrings of 
$\FF^{\prime} \otimes \OO_{\GAB}$ as well.
\end{pf}
\begin{cor}\label{scZpr-c}
The sheaf $\FF^{\prime}_c $ of Proposition \ref{sh-c} determines the scheme
$$
\ZD^{\prime}_c =Spec (\FF^{\prime}_c)
$$ 
over $\GAB$ with the structure morphism
$$
p^{\prime}_c :  \ZD^{\prime}_c \longrightarrow \GAB
$$
and a surjective morphism of $\GAB$-schemes
$$
f^{\prime}_c : \ZD^{\prime}_{\GAB} \longrightarrow  \ZD^{\prime}_c\,,
$$
where 
$ \ZD^{\prime}_{\GAB} = Spec( \FF^{\prime} \otimes \OO_{\GAB})$, i.e.
one has the commutative diagram
$$
\xymatrix{
{\ZD^{\prime}_{\GAB}} \ar[rr]^{f^{\prime}_c} \ar[dr]_{p^{\prime}_2} & &{\ZD^{\prime}_c } \ar[dl]^{p^{\prime}_c} \\
& {\GAB} }
$$
where all morphisms are surjective.
\end{cor}
\begin{pf}
The inclusion of the sheaves of rings
$$
\FF^{\prime}_c \hookrightarrow \FF^{\prime} \otimes \OO_{\GAB}
$$
determines a surjective morphism
$$
f^{\prime}_c : \ZD^{\prime}_{\GAB} = Spec( \FF^{\prime} \otimes \OO_{\GAB}) \longrightarrow   Spec(\FF^{\prime}_c) =\ZD^{\prime}_c
$$
which by definition commutes with the structure morphisms (which are respectively $p^{\prime}_2$ and $p^{\prime}_c$) onto
$\GAB$.
\end{pf}
\begin{defi}\label{sc-cw}
The scheme 
$\ZD^{\prime}_c = Spec (\FF^{\prime}_c) $
defined in Corollary \ref{scZpr-c} will be called the scheme of
central weights of $\GA$. It is a finite scheme over 
$\GAB$ with the structure morphism
$$
p^{\prime}_c : \ZD^{\prime}_c \longrightarrow \GAB
$$
such that 
$$
p^{\prime}_{c \ast} \OO_{\ZD^{\prime}_c} = \FF^{\prime}_c  \cong \mbox{\BM${\cal C}$}_{\GA}.
$$
In particular, $deg(p^{\prime}_{c }) = rk (\mbox{\BM${\cal C}$}_{\GA})$.
\end{defi}

Denote by $\ZD_{\GAB}$ the part of the universal scheme
$\ZD$ in (\ref{uc}) lying over $\GAB$. This means that
$\ZD_{\GAB} = Spec (\FF \otimes \OO_{\GAB})$ and it comes with
the structure morphism
\BEN \label{def-Z-GAb}
p_2 : \ZD_{\GAB} \longrightarrow \GAB
\EEN
which is an {\it unramified} covering of degree $d$.

In Corollary \ref{cor-Fpr} we found that this covering factors through the scheme
$ \ZD^{\prime}_{\GAB} = Spec (\FF^{\prime})$, while the discussion above gives further factorization imposed by the center
{\BM
$\CG$}
of
{\BM
$\LAGT$}.
This is summarized in the following commutative diagram of various morphisms introduced so far.
\BEN\label{diag-morph}
\xymatrix{
{\ZD_{\GAB}} \ar[r]^f \ar[dr]_{p_2} & { \ZD^{\prime}_{\GAB}} \ar[r]^{f^{\prime}_c} \ar[d]^{p^{\prime}_2} &
 {\ZD^{\prime}_c} \ar[dl]^{p^{\prime}_c} \\
          &{\GAB}&                }
\EEN
\begin{rem}\label{et}
By construction all morphisms in the above diagram are surjective, finite and flat. Since the morphism 
$p_2$ in (\ref{diag-morph}) is unramified it follows that all other morphisms are unramified coverings as well.
\end{rem}
 
Set $f_c = f^{\prime}_c \circ f$ to be the composition of the horizontal arrows in (\ref{diag-morph}) and consider the resulting 
factorization of $\ZD_{\GAB}$ through the scheme of central weights $\ZD^{\prime}_c$
\BEN\label{diag-c}
\xymatrix{
{\ZD_{\GAB}} \ar[rr]^{f_c} \ar[dr]_{p_2} &  &{\ZD^{\prime}_c} \ar[dl]^{p^{\prime}_c} \\
          &{\GAB}&                }
\EEN
The scheme of central weights $\ZD^{\prime}_c$ parametrizes the weights of the action of the center
{\BM
$\CG$}
on the sheaf $\FF^{\prime} = p^{\prime}_{2\ast} \OO_{ \ZD^{\prime}_{\GAB}}$, where 
$p^{\prime}_2$ is as in (\ref{diag-morph}). The meaning of the factorization in (\ref{diag-c}) is that it decomposes
the configurations on $X$ (the fibres of $p_2$) parametrized by $\GAB$ into the disjoint union of subconfigurations.
The following statement summarizes some of the basic properties of this decomposition.
\begin{cor}\label{Z-c-dec}
For every $[Z] \in \GAB$ the factorization in (\ref{diag-c}) gives a decomposition of the configuration
$Z$ into the disjoint union of subschemes (subconfigurations) given by the formula
$$
Z= \sum_{\lambda \in p^{\prime -1}_c ([Z])} f^{-1} (Z^{\prime}_{\lambda}),
$$
where
$Z^{\prime}_{\lambda} = f^{{\prime}-1}_c (\{ \lambda \})$, for 
$\lambda \in p^{\prime -1}_c ([Z])$.

For a weight $\lambda$, let  
$Z_{\lambda} = f^{-1}_c (\{ \lambda \})  = f^{-1} (Z^{\prime}_{\lambda})$
be the subscheme of $Z$ corresponding to this weight and let
$$
\delta_{Z_{\lambda}} = f^{\ast} ( \delta_{Z^{\prime}_{\lambda}})
$$
be the pullback by $f$ of the delta-functions
defined in Remark \ref{Zpr-lam} . Then the subshemes 
$Z_{\lambda}$ have the following properties.
\begin{enumerate}
\item[(i)]
$\HO{Z_{\lambda}}) \cong  V_{\lambda} ([Z]) \oplus \delta_{Z_{\lambda}} \cdot \HH^{\LG} \ZA$,
where 
$V_{\lambda} ([Z])$ is the weight space corresponding to $\lambda$ in the weight decomposition 
of $\HO{Z^{\prime}})$ in (\ref{Zpr-wd2}), where
$Z^{\prime} =f(Z) = p^{{\prime} -1}_2 ([Z])$, the fibre of $p^{\prime}_2$ over $[Z]$
in the diagram (\ref{diag-morph}).
\item[(ii)]
$Ext^1_{Z_{\lambda}} := Ext^1 (\ID_{Z_{\lambda}} (L), \OO_X) \cong \HT_{\lambda} \ZA$,
for some
$[\alpha] \in \JAB_Z$, where 
$\ID_{Z_{\lambda}}$ is the sheaf of ideals of $Z_{\lambda}$ on $X$, 
$\HT_{\lambda} \ZA$ is as in Proposition \ref{HT-wd} and 
$\JAB_Z$ is as in (\ref{J-Z}).
\item[(iii)]
There is a natural morphism
$$
p_{\lambda} : \JAB_Z \longrightarrow \PP (Ext^1_{Z_{\lambda}})
$$
which sends a point $[\alpha] \in \JAB_Z $ to the extension class 
$[\alpha_{\lambda}]$ (up to a non-zero scalar)
$$
\xymatrix@1{
0\ar[r]&{\OO_X} \ar[r] & {\SE}_{\alpha_{\lambda}} \ar[r]&{\ID_{Z_{\lambda}}} (L) \ar[r]& 0 }
$$
where the sheaf $\SE_{\alpha_{\lambda}}$ is locally free of rank 2 and having Chern invariants
$(L,d_{\lambda})$ with $d_{\lambda} = degZ_{\lambda}$.
\item[(iv)]
 The filtration
$\HT_{-\bullet} ([Z_{\lambda}],[\alpha_{\lambda}])$
as well as the orthogonal decomposition
$$
\HO{Z_{\lambda}}) =\bigoplus_{p\geq 0} \HH^p ([Z_{\lambda}],[\alpha_{\lambda}])
$$
are obtained from the one's for $\ZA$ using the isomorphism in (i). More precisely,
\begin{eqnarray*}
\HT_{-i} ([Z_{\lambda}],[\alpha_{\lambda}]) &\cong & \HT_{-i} \ZA \bigcap V_{\lambda} ([Z]),\,\,\forall i\leq \LG\\
\HH^p ([Z_{\lambda}],[\alpha_{\lambda}]) &\cong & \HH^{p} \ZA \bigcap V_{\lambda} ([Z]),\,\,\forall p\leq l_{Z_{\lambda}} -1,
\end{eqnarray*}
where $l_{Z_{\lambda}} (\leq \LG)$ is the length of the filtration
$\HT_{-\bullet} ([Z_{\lambda}],[\alpha_{\lambda}])$.
Furthermore, 
$$
\HH^{l_{Z_{\lambda}}} \cong \delta_{Z_{\lambda}} \cdot \HH^{\LG} \ZA.
$$
\end{enumerate}
\end{cor}
\begin{pf}
Set $Z^{\lambda}$ to be the subscheme of $Z$ complementary to 
$Z_{\lambda}$.  The two subschemes give rise to the following commutative diagram
of sheaves on $Z$
\BEN\label{Zlam-Zlamc}
\xymatrix{
 & &0 \ar[d] & &\\
 & &{\JJ_{Z_{\lambda}}} \ar[d] \ar[dr]& & \\ 
0\ar[r]&{ \JJ_{Z^{\lambda}}} \ar[r] \ar[dr] &{\OO_Z} \ar[r]  \ar[d] &{ \OO_{Z^{\lambda}}} \ar[r]& 0 \\
& &{\OO_{Z_{\lambda}}} \ar[d]& & \\
& &0& & }
\EEN
where $\JJ_{Z_{\lambda}}$ (resp. $ \JJ_{Z^{\lambda}}$) is the sheaf of ideals of 
$Z_{\lambda}$ (resp. $Z^{\lambda}$) in $Z$. This implies the identification
\BEN\label{id=o}
 \JJ_{Z^{\lambda}} \cong \OO_{Z_{\lambda}}
\EEN
which in turn gives an isomorphism
\BEN\label{id=o1}
\HO{Z_{\lambda}}) \cong  H^0 (\JJ_{Z^{\lambda}} ).
\EEN
By definition 
$Z^{\lambda} = f^{-1} (Z^{\prime \lambda})$,
where 
$Z^{\prime \lambda}$ is as in Remark \ref{Zpr-lam}.
From the proof of Lemma \ref{lem-Vlam} the latter subscheme is defined by the ideal
$V_{\lambda} ([Z])$. This implies 
\BEN\label{id=o2}
H^0 (\JJ_{Z^{\lambda}} ) = V_{\lambda} ([Z]) \cdot \HO{Z}).
\EEN
We know (see Remark \ref{geom-ordl}) that
$$
\HO{Z}) = \HT_{-\LG} \ZA \oplus \HH^{\LG} \ZA \cong \HO{Z^{\prime}}) \oplus \HH^{\LG} \ZA. 
$$
Substituting this in (\ref{id=o2}) we obtain
\BEN\label{id=o3}
H^0 (\JJ_{Z^{\lambda}} ) = V_{\lambda} ([Z])  \oplus V_{\lambda} ([Z]) \cdot \HH^{\LG} \ZA. 
\EEN
On the other hand the multiplication of $\HH^{\LG} \ZA $ by elements of 
$\HT_{-\LG} \ZA$ preserves $\HH^{\LG} \ZA $. Furthermore, the functions in
$V_{\lambda} ([Z])$ lie in $\HT_{-\LG} \ZA$ and have support in 
$Z_{\lambda}$. This gives the inclusion
$$
V_{\lambda} ([Z]) \cdot \HH^{\LG} \ZA \subset \delta_{Z_{\lambda}} \cdot \HH^{\LG} \ZA. 
$$
The inclusion in the opposite direction is obvious since
$\delta_{Z_{\lambda}} \in V_{\lambda} ([Z])$. Thus we obtain the equality
$$
V_{\lambda} ([Z]) \cdot \HH^{\LG} \ZA  =\delta_{Z_{\lambda}} \cdot \HH^{\LG} \ZA.
$$
Combining it with (\ref{id=o3}) and (\ref{id=o1}) we deduce the first assertion.

To see the second assertion recall that an extension class
$\alpha \in \EZ$ defines the cup-product
\BEN\label{cp-a}
\xymatrix@1{
{\HO{Z}) }\ar[r]^(.2){\alpha} &{ Ext^2 (\OO_Z , \OO_X (-L)) = \HO{Z} (K_X +L))^{\ast}}, }
\EEN
where the equality is the Serre Duality on $X$ (see \RI, \S1.2, for details).
Furthermore, if $[\alpha]$ belongs to the complement of the theta-divisor the homomorphism above is an isomorphism.
Composing it with the dual of the restriction map
$$
\rho (Z) :  \HO{X} (K_X +L)) \longrightarrow \HO{Z} (K_X +L))
$$
we obtain the homomorphism
\BEN\label{R-Za}
{\bf R^r}\ZA : \HO{Z})  \longrightarrow \HO{X} (K_X +L))^{\ast}
\EEN
which is the value of the morphism ${\bf R^r}$ in (\ref{mor-R}) at 
$\ZA$ in $\JABG$. By definition the kernel of this homomorphism is
$\HT \ZA$ (see \RI, (1.21), for details) and $\alpha$ in (\ref{cp-a})
restricted to $\HT \ZA$ induces the isomorphism
\BEN\label{cp-a1}
\xymatrix@1{
{\HT\ZA }\ar[r]^(.6){\alpha} &{ \EZ} }
\EEN
which we already encountered in (\ref{HT=Ext}).

To calculate $Ext^1_{Z_{\lambda}}$ we use the direct sum decomposition
$$
{\HO{Z} ) } =H^0 (\JJ_{Z^{\lambda}} ) \oplus H^0 (\JJ_{Z_{\lambda}} ) 
$$
and its twisted version 
$$
{\HO{Z}(K_X +L) ) } =H^0 (\JJ_{Z^{\lambda}} (K_X +L)) \oplus H^0 (\JJ_{Z_{\lambda}}(K_X +L) )
$$
coming from the diagram (\ref{Zlam-Zlamc}) (resp. (\ref{Zlam-Zlamc}) tensored with $\OO_X  (K_X +L)$).
Substituting these equalities in (\ref{cp-a}) yields the isomorphism
$$
\xymatrix@1{
H^0 (\JJ_{Z^{\lambda}} ) \ar[r]^(.35){\alpha} & H^0 (\JJ_{Z^{\lambda}} (K_X +L))^{\ast}. }
$$
Using the isomorphism in (\ref{id=o1}) we deduce the isomorphism
\BEN\label{cp-alam}
\HO{Z_{\lambda}}) \longrightarrow \HO{Z_{\lambda}} (K_X+L))^{\ast}.
\EEN
This homomorphism will be denoted by $\alpha_{\lambda}$. Combining it with the dual of the restriction homomorphism
$$
\rho (Z_{\lambda}) : \HO{X} (K_X+L) ) \longrightarrow \HO{Z_{\lambda}} (K_X+L))
$$
yields the following diagram
\BEN\label{diag-ExtZlam}
\xymatrix{
 & 0 \ar[d] \\
 & Ext^1_{Z_{\lambda}} \ar[d] \\
{\HO{Z_{\lambda}})} \ar[r]^(.35){\alpha_{\lambda}} \ar[rd] & {\HO{Z_{\lambda}} (K_X+L))^{\ast}} \ar[d]^{\rho (Z_{\lambda})^{\ast} } \\
      & {\HO{X} (K_X+L ))^{\ast}} }
\EEN
Set 
$$
{\bf R^r}_{\lambda} \ZA = \rho (Z_{\lambda})^{\ast} \circ \alpha_{\lambda}
$$
and observe that (\ref{diag-ExtZlam}) implies the isomorphism
\BEN\label{ExtZlam}
\alpha_{\lambda} : ker ({\bf R^r}_{\lambda} \ZA) \longrightarrow Ext^1_{Z_{\lambda}}.
\EEN
On the other hand 
using the identification (\ref{id=o1}) we can identify the homomorphism
${\bf R^r}_{\lambda} \ZA$ with the restriction of 
${\bf R^r} \ZA$ in (\ref{R-Za})
to the subspace 
$H^0 (\JJ_{Z^{\lambda}} )$. This identification allows us to calculate
$ker ({\bf R^r}_{\lambda} \ZA)$ as follows:
\begin{eqnarray*}
\lefteqn{
 ker ({\bf R^r}_{\lambda} \ZA )= H^0 (\JJ_{Z^{\lambda}} ) \bigcap \HT\ZA =}\\
&\big{(} V_{\lambda} ([Z]) \oplus \delta_{Z_{\lambda}} \cdot \HH^{\LG} \ZA \big{)} \bigcap \HT\ZA =
 \big{(} V_{\lambda} ([Z]) \big{)} \bigcap \HT\ZA =\HT_{\lambda} \ZA, 
\end{eqnarray*}
where the second equality follows from the first assertion of the corollary, while the last one is the defining
identity in Proposition \ref{HT-wd}.
Thus the isomorphism in (\ref{ExtZlam}) takes the form
\BEN\label{ExtZlam1}
\alpha_{\lambda} : \HT_{\lambda} \ZA) \longrightarrow Ext^1_{Z_{\lambda}}
\EEN
as stated in (ii) of the corollary.

For part (iii) we observe that the function
$\delta_{Z_{\lambda}}$ belongs to 
$\HT_{\lambda} \ZA$ (Proposition \ref{pro-c-basis}). The isomorphism
in (\ref{id=o1}) identifies it with the unit $1_{Z_{\lambda}}$ of the ring 
$\HO{Z_{\lambda}})$. Applying to it the isomorphism in (\ref{ExtZlam1})
we obtain an extension class in 
$Ext^1_{Z_{\lambda}}$ which we denote by $\alpha_{\lambda}$
(the notation is justified since with this notation the homomorphism
in (\ref{ExtZlam1}) becomes the multiplication by this extension class).
Since $\alpha_{\lambda} \neq 0$ we obtain the map
$$
\JAB_Z \longrightarrow \PP(Ext^1_{Z_{\lambda}})
$$
which sends $[\alpha] \in \JAB_Z $ to the point 
$[\alpha_{\lambda}] \in \PP(Ext^1_{Z_{\lambda}})$. Furthermore, the extension class $\alpha_{\lambda}$ is nowhere vanishing on $Z_{\lambda}$.
This implies by a lemma of Serre, \cite{[O-S-S]}, Lemma 5.1.2, that the sheaf ${\cal E}_{\alpha_{\lambda}}$ sitting in the middle of the exact sequence in (iii) is locally free. 

Turning to (iv) we recall from \S\ref{secFT}, (\ref{HT-i}), that the filtration
$\HT_{-\bullet} ([Z_{\lambda}],[\alpha_{\lambda}])$ is defined as follows:
\begin{equation}\label{HT-i-lam}
\HT_{-1}  ([Z_{\lambda}],[\alpha_{\lambda}]) \cong \HT_{\lambda} \ZA\,\,{\rm{and}}\,\,
 \HT_{-i}  ([Z_{\lambda}],[\alpha_{\lambda}]) = im(S^i (\HT_{\lambda} \ZA)\longrightarrow \HO{ Z_{\lambda}}) ), 
\end{equation}
where the first isomorphism is provided by (ii). From the proof of Lemma \ref{c-Z}, (\ref{HT-i-wd}), it follows that
$$
 im(S^i (\HT_{\lambda} \ZA) \longrightarrow \HO{ Z_{\lambda}}) )  \cong \HT_{-i} \ZA \bigcap V_{\lambda} ([Z]).
$$
This and (\ref{HT-i-lam}) imply the identification
$$
\HT_{-i}  ([Z_{\lambda}],[\alpha_{\lambda}]) \cong \HT_{-i} \ZA \bigcap V_{\lambda} ([Z]).
$$
In particular, one sees that the length $l_{Z_{\lambda}}$ of the filtration $\HT_{-\bullet}  ([Z_{\lambda}],[\alpha_{\lambda}])$ is the smallest index $i$ for which $V_{\lambda} ([Z]) \subset \HT_{-i} \ZA$.

From the proof of Lemma \ref{c-Z} it also follows that the center 
of 
{\BM
$\LAGT$}$\ZA$
acts by endomorphisms of degree $0$ on the orthogonal decomposition
$$
\HO{Z}) = \bigoplus^{\LG}_{p=0} \HH^p \ZA
$$
(see (\ref{D-D0})). This implies that the summands $ \HH^p \ZA$ admit the weight decomposition
$$
 \HH^p \ZA = \bigoplus_{\lambda}  V_{\lambda} ([Z]) \bigcap  \HH^p \ZA.
$$
This implies the identification
$$
\HH^p  ([Z_{\lambda}],[\alpha_{\lambda}]) \cong  V_{\lambda} ([Z]) \bigcap  \HH^p \ZA,
$$
for all $p\leq l_{Z_{\lambda}}-1$ as asserted in (iv). 

Finally, combining this with the isomorphism in 
(i) yields the orthogonal decomposition
\begin{eqnarray}\label{ortd-Zlam}
\HO{Z_{\lambda}}) \cong & V_{\lambda} ([Z]) \oplus \delta_{Z_{\lambda}} \cdot \HH^{\LG} \ZA =& \\ \nonumber
  &\left( \bigoplus^{\LG-1}_{p=0} V_{\lambda} ([Z]) \bigcap \HH^p \ZA \right) \oplus \delta_{Z_{\lambda}} \cdot \HH^{\LG} \ZA =&  \\ \nonumber
&\left( \bigoplus^{l_{Z_{\lambda}}-1}_{p=0} V_{\lambda} ([Z]) \bigcap \HH^p \ZA \right)\oplus \delta_{Z_{\lambda}} \cdot \HH^{\LG} \ZA \cong & \\ \nonumber
& \left( \bigoplus^{l_{Z_{\lambda}}-1}_{p=0} \HH^p  ([Z_{\lambda}],[\alpha_{\lambda}]) \right) 
\oplus \delta_{Z_{\lambda}} \cdot \HH^{\LG} \ZA &
\end{eqnarray}
which implies the last assertion in (iv).
\end{pf}
\subsection{The sheaf of semisimple Lie algebras {\BM$\LAG$}}\label{sec-ses}
In this subsection we determine the semisimple part
{\BM$\LAG$} of the sheaf of Lie algebras 
{\BM$\LAGT$}.
The essential part of the argument is a consideration of the restriction of the sheaf
{\BM$\LAG$}
to the fibres of the projection
\BEN\label{pi-JAB2}
\pi: \JABG \longrightarrow \GAB.
\EEN
So we fix the fibre
$\JAB_Z$ of the morphism $\pi$ in (\ref{pi-JAB2}) over a point
$[Z] \in \GAB$ and denote by
\BEN\label{LAa}
{\bf \tilde{g}}([\alpha]) =\mbox{\BM${\GT}$}\ZA\,\,
({\rm{resp.}}\,\,
{\bf g}([\alpha]) =\mbox{\BM${\GS}$}\ZA),
\EEN
 the fibre of 
{\BM
$\LAGT$ (resp. $\LAG$)}
at $\ZA \in \JAB_Z$. In particular, we have the structure decomposition
$$
{\bf \tilde{g}}([\alpha]) = {\bf c \oplus g}([\alpha]),
$$
where 
${\bf c}$ is the fibre of the center
{\BM
$\CG$}
at $\ZA  \in \JAB_Z$. We have seen in \S\ref{Center},  Lemma \ref{c-Z}, that
${\bf c}$ depends only on $[Z]$. We show that the same holds for
${\bf  {g}}([\alpha])$ and hence for ${\bf \tilde{g}}([\alpha])$.

Let 
$Z^{\prime} =Spec( \HT_{-\LG}( [Z]))$ as defined\footnote{we use the fact, proved in Lemma \ref{HTl=triv}, 2), that $\HT_{-\LG}$ is constant along $\JAB_Z$.} in (\ref{def-Zpr})
and let
\BEN\label{Zpr-wd3}
\HO{Z^{\prime}}) = \bigoplus_{\lambda} V_{\lambda} ([Z])
\EEN
be the weight decomposition of $\HO{Z^{\prime}}) $ under
the ${\bf c}$-action as in (\ref{Zpr-wd2}). In particular, we have
\begin{enumerate}
\item[(a)]
${\bf \tilde{g}}([\alpha]) \subset \bigoplus_{\lambda} {\bf gl} (V_{\lambda} ([Z]))$,
\item[(b)]
${\bf c} \cong \bigoplus_{\lambda} \CC id_{V_{\lambda} ([Z])} \cong Center \big{(}\bigoplus_{\lambda} {\bf gl} (V_{\lambda} ([Z])) {)}$,
\end{enumerate}
where the first inclusion is obvious, while the first isomorphism in (b) is
Remark \ref{Zpr-lam}.

\begin{pro}\label{ga-dec}
If 
${\bf  {g}}([\alpha]) \neq 0$, for some 
$[\alpha] \in \JAB_Z$, then
$$
{\bf  {g}}([\alpha]) = \bigoplus_{\lambda} {\bf sl} (V_{\lambda} ([Z])).
$$
In particular, 
${\bf  {g}}([\alpha])$ is independent of $[\alpha] \in \JAB_Z$.
\end{pro}
\begin{pf}
Let 
${\bf \tilde{h}}([\alpha])$ be the centralizer of
$D(\HT\ZA)$, where $D$ is as in (\ref{m-D}). By Proposition \ref{Cartan}
this is a Cartan subalgebra of 
${\bf \tilde{g}}([\alpha])$ and we have the decomposition
$$
{\bf \tilde{h}}([\alpha]) = {\bf c} \oplus {\bf h}([\alpha]),
$$
where ${\bf h}([\alpha])$ is a Cartan subalgebra of 
${\bf  {g}}([\alpha])$.

From Remark \ref{act-m} we know that 
${\bf \tilde{h}}([\alpha])$ can be identified with a subspace of the ring of functions
$\HO{Z^{\prime}})$ and its action on 
$\HO{Z^{\prime}})$ is identified with the multiplication in
$\HO{Z^{\prime}})$. In particular, the vectors
$\delta_{p^{\prime}} (p^{\prime} \in Z^{\prime})$, are the weight vectors of the action of 
${\bf h}([\alpha])$ on 
$\HO{Z^{\prime}})$.

Let 
${\bf h}([\alpha])^{\ast}$ be the vector space dual to
${\bf h}([\alpha])$
and let
$R([\alpha])$
be the set of roots of
${\bf  {g}}([\alpha])$
with respect to the Cartan subalgebra ${\bf h}([\alpha])$.

Fix the points in $Z^{\prime}$ in some order.  This fixes the order
on the basis
\BEN\label{basis}
\{ \delta_{p^{\prime}} \mid p^{\prime} \in Z^{\prime} \}
\EEN
and from now on we can identify
 ${\bf gl}(\HO{Z^{\prime}}))$
with the Lie algebra
${\bf gl}_{d^{\prime}} (\CC)$, where we put
\BEN\label{d-pr}
d^{\prime} = deg Z^{\prime}.
\EEN
Then the set of coroots
${\check{R}}([\alpha])$ will be identified with a subset of integer-valued matrices
in
${\bf gl}_{d^{\prime}} (\CC)$.
Hence
the set of coroots
${\check{R}}([\alpha])$ 
is independent of the continuous parameter $[\alpha] \in \JAB_Z$.
Since 
${\check{R}}([\alpha])$ spans the Cartan subalgebra
$ {\bf h}([\alpha])$ we deduce that
${\bf h}([\alpha]) ={\bf h}$
is independent of $[\alpha] \in \JAB_Z$. This combined with Lemma \ref{c-Z} implies
\BEN\label{ht-cnst}
{\bf \tilde{h}}([\alpha] )={\bf \tilde{h}} = {\bf c \oplus h}
\EEN
is independent of $[\alpha] \in \JAB_Z$.

Observe that the basis of 
$\HO{Z^{\prime}})$ fixed in (\ref{basis})
allows us to identify the space of diagonal matrices of
${\bf gl}_{d^{\prime}} (\CC)$ with
the space of functions
$\HO{Z^{\prime}})$ on $Z^{\prime}$. This identification
implies an inclusion
\BEN\label{Cartan-sub-f}
{\bf \tilde{h}} \subset \HO{Z^{\prime}}).
\EEN
We claim that the equality holds.
\begin{cl}\label{Cartan=f}
${\bf \tilde{h}} ={\bf c \oplus h} =\HO{Z^{\prime}})$.
\end{cl}

Let us assume this and complete the proof of the proposition.
For this
consider the Cartan decomposition of 
${\bf g}([\alpha])$ with respect to ${\bf {h}}$
$$
{\bf g}([\alpha]) = {\bf {h}} \oplus \left(\bigoplus_{\xi \in R([\alpha])} {\bf g}([\alpha])_{\xi} \right)\,,
$$
where
${\bf g}([\alpha])_{\xi}$ is the root space of ${\bf g}([\alpha])$ corresponding to a root
$\xi \in R([\alpha])$. Choose a root vector
$E_{\xi} ([\alpha])$, a generator of the root space ${\bf g}([\alpha])_{\xi}$, for every
$\xi \in R([\alpha])$.

Let $E_{p^{\prime}, q^{\prime}}$ be the endomorphism of
$\HO{Z^{\prime}})$
which takes $\delta_{q^{\prime}}$ to $\delta_{p^{\prime}}$
and kills all other vectors of our basis in (\ref{basis}). Then the set
$\{ E_{p^{\prime}, q^{\prime}} \}_{p^{\prime}, q^{\prime} \in  Z^{\prime}}$
forms the standard basis of 
${\bf gl}(\HO{Z^{\prime}}))$.
We will show that the root vectors
$E_{\xi} ([\alpha])$ can be chosen to be in the standard basis. For this write
\BEN\label{root}
E_{\xi} ([\alpha]) = \sum_{p^{\prime}, q^{\prime}} c_{p^{\prime}, q^{\prime}} E_{p^{\prime}, q^{\prime}}\,\,.
\EEN

Extend $\xi$ by zero on the center
${\bf c}$ and view the roots of 
${\bf g}([\alpha])$
as linear functions on ${\bf  \tilde{h} =c \oplus h }= \HO{Z^{\prime}})$, where the second equality is Claim \ref{Cartan=f}.
Applying
$ad(h)$ to (\ref{root}), for $h \in {\bf  \tilde{h}}$, we obtain
\BEN\label{ad-h}
\xi (h) E_{\xi} ([\alpha])  = \sum_{p^{\prime}, q^{\prime}} c_{p^{\prime}, q^{\prime}} [h,E_{p^{\prime}, q^{\prime}}]=
\sum_{p^{\prime}, q^{\prime}} c_{p^{\prime}, q^{\prime}} (h(p^{\prime}) -h(q^{\prime}))E_{p^{\prime}, q^{\prime}} \,.
\EEN
This yields
\BEN\label{rel1}
 c_{p^{\prime}, q^{\prime}} (h(p^{\prime}) -h(q^{\prime}) -\xi (h)) =0,\,\,\forall p^{\prime}\neq q^{\prime} \in  Z^{\prime}
\,\,{\rm{and}}\,\,\forall h\in \HO{Z^{\prime}})\,.
\EEN
Observe that in this relation we use the identification
${\bf  \tilde{h} }= \HO{Z^{\prime}})$ provided by Claim \ref{Cartan=f} and view $h$ as a function on $ Z^{\prime}$.

Set
$$
 Z^{\prime}_{\xi} = \left\{ ( \left.q^{\prime}, p^{\prime}) \in Z^{\prime} \times Z^{\prime} \right| c_{p^{\prime}, q^{\prime}} \neq 0 \right\}.
$$
For every $( q^{\prime}, p^{\prime}) \in  Z^{\prime}_{\xi} $ the equality (\ref{rel1}) yields
\BEN\label{rel2}
h(p^{\prime}) -h(q^{\prime}) -\xi (h) =0,\,\,\forall h\in \HO{Z^{\prime}}).
\EEN
This implies that the restrictions of the projections
$$
\epsilon_j: Z^{\prime}_{\xi} \longrightarrow Z^{\prime}
$$
on the $j$-th factor of $Z^{\prime} \times Z^{\prime}$, for $j=1,2$, is injective.

Set $Z^{j}_{\xi} = \epsilon_j ( Z^{\prime}_{\xi})$, for $j=1,2$, and let
\BEN\label{Z1-Z2}
l_{\xi} = \epsilon_2 \circ \epsilon^{-1}_1 : Z^1_{\xi} \longrightarrow  Z^2_{\xi}
\EEN
be the corresponding bijection. With this notation in mind we rewrite (\ref{root}) as follows
\BEN\label{root1}
E_{\xi} ([\alpha]) = \sum_{ q^{\prime} \in Z^1_{\xi}} c_{l_{\xi} (q^{\prime}), q^{\prime}} E_{l_{\xi} (q^{\prime}), q^{\prime}}\,.
\EEN
We want to show that $Z^1_{\xi}$ is a single point. For this examine the positions of the subsets
$Z^1_{\xi}$ and $Z^2_{\xi}$ in $ Z^{\prime}$ relative to each other.

{\bf Case 1:} $Z^1_{\xi} \neq  Z^2_{\xi}$. In this situation we can find a point
$p^{\prime}$ in $ Z^{\prime}$ which belongs exactly to one of these sets, say
$p^{\prime} \in Z^1_{\xi}$ and $p^{\prime} \not\in Z^2_{\xi}$.
Take the function $\delta_{p^{\prime}}$ and observe that the multiplication by this function in
$\HO{Z^{\prime}})$ corresponds to the endomorphism
$H_{p^{\prime}} :=E_{p^{\prime},p^{\prime}} \in {\bf \tilde{h}}$.
Applying
$ad(H_{p^{\prime}})$ to (\ref{root1}) we obtain
$$
\xi(H_{p^{\prime}}) E_{\xi} ([\alpha]) =\sum_{ q^{\prime} \in Z^1_{\xi}} c_{l_{\xi} (q^{\prime}), q^{\prime}} 
(\delta_{p^{\prime}} (l_{\xi} (q^{\prime})) - \delta_{p^{\prime}} (q^{\prime}) ) E_{l_{\xi} (q^{\prime}), q^{\prime}} =
-  c_{l_{\xi} (p^{\prime}), p^{\prime}} E_{l_{\xi} (p^{\prime}), p^{\prime}}\,.
$$
This implies that the root vector $E_{\xi} ([\alpha])$ can be chosen to be
$E_{l_{\xi} (p^{\prime}), p^{\prime}}$, for some $ p^{\prime} \in  Z^{\prime}$.

{\bf Case 2:} $Z^1_{\xi} =  Z^2_{\xi}$. We show that this situation is impossible.
Indeed, take $p^{\prime} \in Z^1_{\xi}$ and go through the above calculation to obtain the following
\BEN\label{rel3}
\xi(H_{p^{\prime}}) E_{\xi} ([\alpha]) =-  c_{l_{\xi} (p^{\prime}), p^{\prime}} E_{l_{\xi} (p^{\prime}), p^{\prime}} +
c_{p^{\prime}, l^{-1}_{\xi} (p^{\prime})} E_{p^{\prime}, l^{-1}_{\xi} (p^{\prime})}\,.
\EEN
This implies 
$$
 Z^1_{\xi} =\{ p^{\prime}, \,l^{-1}_{\xi} (p^{\prime}) \} =  Z^2_{\xi} = \{ p^{\prime}, \,l_{\xi} (p^{\prime}) \}. 
$$
Hence $Z^1_{\xi} =  Z^2_{\xi} = \{ p^{\prime}, q^{\prime} \}$, for some two distinct points
$ p^{\prime}, q^{\prime} \in Z^{\prime}$ and $l_{\xi}$ is the transposition of these points.
With this notation the equation (\ref{rel3}) takes the following form
$$
\xi(H_{p^{\prime}}) E_{\xi} ([\alpha]) =-  c_{q^{\prime} , p^{\prime}} E_{ q^{\prime}, p^{\prime}} +
c_{p^{\prime},  q^{\prime}} E_{p^{\prime}, q^{\prime}}\,.
$$
This yields the system of equations
$$
\left\{
\begin{array}{lr}
\xi(H_{p^{\prime}}) c_{q^{\prime} , p^{\prime}} =& -  c_{q^{\prime} , p^{\prime}} \\
\xi(H_{p^{\prime}}) c_{p^{\prime} , q^{\prime}} =& c_{p^{\prime},  q^{\prime}}
\end{array} \right.
$$
which is clearly impossible.

Thus we have shown that root vectors of
${\bf g}([\alpha])$ can be chosen to lie in the set
$\{ E_{p^{\prime},  q^{\prime}} \}_{p^{\prime} \neq  q^{\prime} \in Z^{\prime} }$.
Furthermore, 
${\bf g}([\alpha])$ must preserve the subspaces
$$
V_{\lambda} ([Z]) = \bigoplus_{p^{\prime} \in Z^{\prime}_{\lambda}} \CC \delta_{p^{\prime}}\,,
$$
for every  weight $\lambda$ occurring in (\ref{Zpr-wd3}) and where the description of $V_{\lambda} ([Z])$ is given in Lemma \ref{lem-Vlam}. Hence the root vectors
of
${\bf g}([\alpha])$ can be chosen to be in the subset
$$
\bigcup_{\lambda} \{  E_{p^{\prime},  q^{\prime}} \}_{p^{\prime} \neq  q^{\prime} \in Z^{\prime}_{\lambda}}\,,
$$
where the union is taken over the weights $\lambda$ with
$d^{\prime}_{\lambda} = deg Z^{\prime}_{\lambda} \geq 2$.

Fix such $Z^{\prime}_{\lambda}$. It remains to be proved that
every $ E_{p^{\prime},  q^{\prime}}$, for 
$p^{\prime} \neq  q^{\prime} \in Z^{\prime}_{\lambda}$
is a root vector of ${\bf g}([\alpha])$.

First we show that for every $p^{\prime} \in Z^{\prime}_{\lambda}$
an element 
$ E_{p^{\prime},  q^{\prime}}$ is a root vector of ${\bf g}([\alpha])$, for some 
$ q^{\prime} \in  Z^{\prime}_{\lambda} \setminus \{p^{\prime} \}$.
Suppose this is not the case. Then $H_{p^{\prime}}$ commutes with
${\bf g}([\alpha])$ and hence, belongs to the center ${\bf c}$ of 
${\bf \tilde{g}}([\alpha])$. By Lemma \ref{lem-Vlam},
$H_{p^{\prime}}$, viewed as a function on $Z^{\prime}$, must be constant on
$Z^{\prime}_{\lambda}$.
However, $H_{p^{\prime}}$ corresponds to the function
$\delta_ {p^{\prime}}$. This means that
$ Z^{\prime}_{\lambda} =\{ p^{\prime} \}$ 
contrary to the assumption 
$ deg Z^{\prime}_{\lambda} \geq 2$.

Let 
$\{ p^{\prime}_1, \ldots, p^{\prime}_{d^{\prime}_{\lambda}} \}$
be an ordering of the points in
$Z^{\prime}_{\lambda}$ such that the root vectors
$$
\{ E_{1,2}, \ldots, E_{j-1,j} \}
$$
form a longest uninterrupted string of root vectors in the set
$\{E_{k,m}\}_{1 \leq k,m \leq d^{\prime}_{\lambda}}$, where
we write
$E_{k,m}$ instead of $E_{p^{\prime}_k, p^{\prime}_m}$.
If $j=d^{\prime}_{\lambda}$, then we are done.
Assume $j < d^{\prime}_{\lambda}$.
Then the elements 
$E_{k,m}$  with $1 \leq k \neq m \leq j$ are all in ${\bf g}([\alpha])$
and they generate the Lie subalgebra of
${\bf g}([\alpha])$
which can be identified with
${\bf sl} (V_j)$,
where
$V_j = \bigoplus^j_{k=1} \CC \delta_ {p^{\prime}_k}$.
Furthermore, the root vectors of ${\bf g}([\alpha])$ which are
not in 
${\bf sl} (V_j)$ are of the form
$E_{s,t}$ with 
$s,t >j$, since otherwise we could make our string longer.
Set
$V^j = \bigoplus_{k>j}  \CC \delta_ {p^{\prime}_k}$.
Then the considerations above imply an inclusion
$$
{\bf g}([\alpha]) \subset {\bf sl} (V_j) \oplus {\bf sl} (V^j).
$$
From this it follows that the endomorphism
$$
H_{\leq j} = \sum^j_{k=1} H_k
$$
is in the center of 
${\bf \tilde{g}}([\alpha])$. Hence the corresponding function
$$
\delta_{\leq j} = \sum^j_{k=1} \delta_{p^{\prime}_k}
$$ 
must be constant on $Z^{\prime}_{\lambda}$ (see Lemma \ref{lem-Vlam}).
But this contradicts the assumption that
$j<d^{\prime}_{\lambda} = degZ^{\prime}_{\lambda}$.
\\
\\
\indent
We now turn to  Claim \ref{Cartan=f}. Identify the Cartan subalgebra
${\bf \tilde{h}} ={\bf c \oplus h}$ as a subspace of 
$H^0(\OO_{Z^{\prime}})$. We already know from the first part of the proof that this subspace does not depend on
$[\alpha]$ in $\JAB_Z$ (see (\ref{ht-cnst})). On the other hand we have an inclusion
$$
\HT \ZA \subset {\bf \tilde{h}},
$$
for all $[\alpha] \in \JAB_Z$. In particular, consider a path
$\alpha(\epsilon) = \alpha + \epsilon \beta$ in 
$\EZ$ passing through $ \alpha$ in the direction of $\beta$.
For all $ \epsilon \in \CC$ with $\left|  \epsilon \right|$ sufficiently small, the
points $[\alpha(\epsilon)]$ lie in $\JAB_Z$ and we have an inclusion
\BEN\label{inc-e}
\HT ([Z],[\alpha(\epsilon)]) \subset {\bf \tilde{h}},
\EEN
for all $\epsilon$ in a small disk 
$B_{\epsilon}$ around $0 \in \CC$.

As in the proof of Claim \ref{cl-ab}, (\ref{a/b1}), we write
$$
\HT ([Z],[\alpha(\epsilon)]) = \frac{\alpha}{\alpha(\epsilon)} \HT \ZA =
\frac{\alpha}{\alpha + \epsilon \beta} \HT \ZA = 
\frac{1}{1 + \epsilon t} \HT \ZA, 
$$
where $t=\frac{\beta}{\alpha} \in \HT\ZA$.
This and (\ref{inc-e}) imply
\BEN\label{inc-e1}
\frac{1}{1 + \epsilon t} \HT \ZA \subset  {\bf \tilde{h}},
\EEN
for all $\epsilon \in B_{\epsilon}$. Taking the $\epsilon$-expansion we
deduce that
$t^k \HT\ZA \subset {\bf \tilde{h}}$,
for all $k \in {\bf Z}_{\geq 0}$ and for all
$t\in \HT\ZA$.
Hence
$\HT_{-\LG} \ZA = \HO{Z^{\prime}} ) \subset {\bf \tilde{h}}$. This together with 
(\ref{Cartan-sub-f}) yield the equality
$$
{\bf \tilde{h}} =\HO{Z^{\prime}})
$$
asserted in Claim \ref{Cartan=f}.
\end{pf}

At this stage we have a complete description of the Lie algebras
{\BM
$\LAGT$}$\ZA$, for $\ZA \in \JABG$, and their relation to the geometry of the
underlying configurations $Z \subset X$ with $[Z] \in \GAB$. This is summarized in the
following statement.
\begin{thm}\label{th-Lie-dec}
Let $Z$ be a configuration on $X$ with $[Z] \in \GAB$, where $\GA$ is an admissible component
 in $\CS$ subject to the conventions in \S\ref{c-n1}. Let 
$\JAB_Z$ be the fibre of the projection
$$
\pi: \JABG \longrightarrow \GAB
$$
over $[Z] \in \GAB$.
Then the following holds.
\begin{enumerate}
\item[1)]
The Lie algebra
{\BM
$\LAGT$}$\ZA$
 and its center
{\BM
$\CG$}$\ZA$, the fibres, respectively of
{\BM
$\LAGT$} and 
{\BM
$\CG$}, at $\ZA$,
are independent of $[\alpha] \in \JAB_Z$. These Lie algebras will be denoted
{\BM
$\LAGT$}$([Z])$
and
{\BM
$\CG$}$([Z])$ 
respectively.
\item[2)]
The subspace
$\HT_{-\LG} ([Z])=\HT_{-\LG} \ZA$ of the filtration
$\HT_{-\bullet}$ in (\ref{filtHT-JG}) at $\ZA$ is independent of $[\alpha] \in \JAB_Z$
and it decomposes into the direct sum 
$$
\HT_{-\LG} ([Z]) = \bigoplus_{\lambda} V_{\lambda} ([Z])
$$
of the weight spaces under the action of the center 
{\BM
$\CG$}$[Z]$,
where
$V_{\lambda} ([Z])$ is the weight space  
 corresponding to a weight $\lambda$.
Furthermore, the weights occurring in the above decomposition form a
basis of 
$(\mbox{\BM${\cal C}$}_{\GA} ([Z]))^{\ast}$,
the space dual to
{\BM
$\CG$}$[Z]$.
\item[3)]
The Lie algebra
{\BM
$\LAGT$}$([Z])$
has the following form
$$
\mbox{\BM${\GT}$}_{\GA} ([Z]) = \bigoplus_{\lambda} {\bf gl} (V_{\lambda} ([Z])).
$$
\item[4)]
The configuration $Z$ admits a decomposition into the disjoint union
$$
Z= \bigcup_{\lambda} Z_{\lambda}
$$
of subconfigurations $Z_{\lambda}$ indexed by the weights occurring in the weight decomposition
of $\HT_{-\LG} ([Z])$ in 2). 
Each subconfiguration $Z_{\lambda}$ has the following properties
\begin{enumerate}
\item[(i)]
$\HO{Z_{\lambda}}) \cong V_{\lambda} ([Z]) \cdot \HO{Z})$
where ``$\cdot$'' stands for the multiplication in the ring
$\HO{Z})$.
\item[(ii)]
$Z_{\lambda}$ is $L$-special whose index of $L$-speciality  (defined in (\ref{is})) is given by the formula
$$
\delta(L,Z_{\lambda}) = dim\left( V_{\lambda} ([Z]) \bigcap \HT_{-1} \ZA\right)\,,
$$
where $\HT_{-1} (\ZA) =\HT\ZA$ is the fibre of the sheaf
$\HT_{-1} =\HT$ at $\ZA \in \JABG$.
\end{enumerate}
\end{enumerate}
\end{thm}
\begin{pf}
Everything has been already proved. For the convenience of the reader we give the list
of references, where the proofs could be found.

Part 1) is proved in Lemma \ref{c-Z}, for the center 
{\BM
$\CG$}$([Z])$, and in Proposition \ref{ga-dec}, for the Lie algebra
 {\BM
$\LAGT$}$([Z])$.

Part 2) is proved in Lemma \ref{HTl=triv}, Lemma \ref{c-Z}.
The last assertion in 2) follows from 
Proposition \ref{pro-c-basis}, 2).

Parts 3) and 4) are Proposition \ref{ga-dec} and Corollary \ref{Z-c-dec},
respectively.
\end{pf}
With the fibrewise study of 
{\BM
$\LAGT$}
completed we turn now to its global properties.
\begin{pro}\label{Lie=lf}
\begin{enumerate}
\item[1)]
If the sheaf of semisimple algebras 
{\BM
$\LAG$}$=0$,
then
{\BM
$$
\LAGT = \CG \cong \HT \mbox{\UB$\otimes \OO_{\JABG} = \HT_{-\LG}$},
$$}
where the second identification is given by the morphism
$$
D: \HT \longrightarrow \mbox{\BM$\LAGT$} \subset \ENDO(\FT )
$$
in (\ref{m-D}).
\item[2)]
If the sheaf
{\BM
$\LAG$}$\neq 0$,
then it is locally free.
\end{enumerate}
\end{pro}
\begin{pf}
If 
{\BM
$\LAG$}$=0$,
then for every local section $t$ of $\HT$ the components
$D^{\pm} (t) =0$ (\RI, Lemma 7.6). This implies that 
{\BM
$\LAGT$}
is abelian. Hence the first equality and the second isomorphism in 1). Furthermore, the multiplication by
$t$ preserves $\HT$, i.e.
$\HT$ is a sheaf of subrings of 
$\HT_{-\LG}$. Since the latter, by definition, is the subsheaf of rings in
$\FT$ generated by $\HT$, we deduce the second equality in 1).

The assertion 2) follows from a well-known fact that local deformations of a semisimple Lie algebras are trivial.
\end{pf} 
\begin{cor}\label{C-Lie=lf}
The subsheaf of centers
{\BM
$\CG$}
and the sheaf of Lie algebras 
{\BM
$\LAGT$}
are locally free.
\end{cor}
\begin{pf}
The first assertion together with the structure decomposition in (\ref{str-d})
and Proposition \ref{Lie=lf},2) imply that the sheaf
{\BM
$\LAGT$}
is locally free. So it is enough to prove that
{\BM
$\CG$} is locally free. For this we argue according to two cases in Proposition \ref{Lie=lf}.

If 
{\BM
$\LAG$}$=0$, then by Proposition \ref{Lie=lf},1), the center
{\BM
$\CG$}
is isomorphic to the sheaf
$\HT$ which is locally free on $\JABG$.

If
{\BM
$\LAG$}$\neq 0$, 
then we consider the sheaf  
{\BM
${\cal{C}}(\HT)$} which is, according to Proposition \ref{Cartan},
 a subsheaf of Cartan subalgebras of
{\BM
$\LAGT$} and has the direct sum decomposition
{\BM
\BEN\label{C-d1}
{\cal{C}}(\HT) = \CG \oplus {\cal{H}}_{\mbox{\UB$\GA$}},
\EEN}
where 
{\BM
$ {\cal{H}}_{\mbox{\UB$\GA$}}$}
is a subsheaf of Cartan subalgebras of 
{\BM
$\LAG$} (see (\ref{C-d}) and (\ref{CA-def})).
From Claim \ref{Cartan=f} it follows 
$$
\mbox{\BM${\cal{C}}$} (\HT) \cong \pi^{\ast} (\FF^{\prime}) = \HT_{-\LG}\,,
$$
 where the equality comes from Corollary \ref{cor-Fpr}. In particular,
{\BM
${\cal{C}}(\HT)$}
is locally free.  This together with the direct sum decomposition in (\ref{C-d1})
imply that 
{\BM
$\CG$}
is locally free.
\end{pf}
\begin{rem}\label{assum-unnec}
From Corollary \ref{C-Lie=lf} it follows that the assumption on the rank of
{\BM
$\CG$}
made in Proposition \ref{sh-c} is unnecessary.
\end{rem}

Our study distinguishes two types of components $\GA$ in 
$\CSA$ according to whether 
the sheaf of Lie algebras
{\BM
$\LAGT$} is abelian or not.
This motivates the following definition.
\begin{defi}\label{qa-c}
\begin{enumerate}
\item[1)]
 A component  $\GA$ in $\CSA$ is called quasi-abelian (we will often abbreviate - q-a) if 
the sheaf of Lie algebras
{\BM
$\LAGT$}
is abelian or, equivalently, the sheaf 
{\BM
$\LAG$}$=0$. 
\item[2)]
A configuration $Z \subset X$ is called quasi-abelian if there
exists a q-a component $\GA \in \CSA$ such that
$[Z] \in \GAB$.
\end{enumerate}
\end{defi}

The following statement characterizes quasi-abelian components in terms of properties of schemes and morphisms
appearing in Corollary \ref{scZpr-c}.
\begin{cor}\label{cor-qa}
If $\GA \in \CSA$ is q-a, then the scheme of central weights
$\ZD^{\prime}_c$ (Definition \ref{sc-cw}) coincides with the scheme
$\ZD^{\prime}_{\GA}$, i.e. the morphism $f^{\prime}_c$ in
Corollary \ref{scZpr-c} is the identity. In particular,
$$
p^{\prime}_2 : \ZD^{\prime}_{\GA} \longrightarrow \GAB
$$
is an unramified covering of degree $(r+1)$.
\end{cor}
\begin{pf}
 From Proposition \ref{Lie=lf}, 1), it follows that the
sheaf $\FF^{\prime}_c$ defining the scheme 
$\ZD^{\prime}_c$ (see Corollary \ref{scZpr-c}) pulled back by 
$\pi$ coincides with $\HT$, i.e.
$$
\pi^{\ast}\FF^{\prime}_c =\HT.
$$
  The second equality 
$\HT =\HT_{-\LG}$ in Proposition \ref{Lie=lf}, 1),
and the definition of the scheme
$\ZD^{\prime}_{\GA}$ in Corollary \ref{cor-Fpr} imply the
equality
 $\ZD^{\prime}_c =\ZD^{\prime}_{\GA}$.
\end{pf}

On the opposite extreme of quasi-abelian components one has components
$\GA \in \CSA$ with the center
{\BM
$\CG$}
being trivial.
\begin{defi}\label{s-c}
\begin{enumerate}
\item[1)]
 A component  $\GA$ in $\CSA$ is called simple if the rank
$rk(\mbox{\BM${\cal C}$}_{\GA})$ of the center
{\BM
$\CG$}
is equal to 1.
\item[2)]
A configuration $Z \subset X$ is called simple if there
exists a simple component $\GA \in \CSA$ such that
$[Z] \in \GAB$.
\end{enumerate}
\end{defi}

The above terminology is justified in view of the following result.
\begin{cor}\label{cor-s}
Let
$\GA$ be a simple component in $\CSA$.
Then
$$
\mbox{\BM$\LAGT$} = \pi^{\ast} ({\bf gl}( \FF^{\prime})) = \OO_{\JABG} \oplus \pi^{\ast} \big{(} {\bf sl}( \FF^{\prime})\big{)}\,,
$$
where ${\bf gl}({\cal A})$ (resp. ${\bf sl}( {\cal A})$) stands for the sheaf of endomorphisms (resp. traceless endomorphisms) of a locally free
sheaf
${\cal A}$.
In particular, an admissible component $\GA$ is simple if and only if the sheaf 
{\BM
$\LAG$} is a sheaf of simple Lie algebras of type $\bf A_{d^{\prime}_{\GA}}$,
where 
$d^{\prime}_{\GA} = rk ( \FF^{\prime})$.
\end{cor}
\begin{pf}
Observe that the identity endomorphism $id_{\FT}$ of $\FT$ is always included in
$H^0( \mbox{\BM$\LAGT$} )$. So it generates the subsheaf of   
{\BM
$\CG$} 
isomorphic to $\OO_{\JABG}$. Thus by Definition \ref{s-c} a component
$\GA \in \CSA$ is simple if and only if
\BEN\label{c=oo}
\mbox{\BM${\cal C}$}_{\GA} \cong \OO_{\JABG}\,.
\EEN
This implies that we must be in the 
 case 2) of Proposition \ref{Lie=lf}, i.e.
{\BM
$\LAG$}$\neq 0$ (this is because $rk(\HT) =r+1 \geq 2$, where the
inequality is our convention of $r\geq1$ in \S\ref{c-n1}).  We now apply
Proposition \ref{ga-dec} to deduce the equality
$$
\mbox{\BM$\LAG$} ={\bf sl}(\pi^{\ast} \FF^{\prime}) = \pi^{\ast}({\bf sl} (\FF^{\prime}) )\,.
$$
This together with (\ref{c=oo}) imply that the structure decomposition of
{\BM
$\LAGT$} 
is as follows
$$
 \mbox{\BM$\GT$}_{\GA} = \OO_{\JABG} \oplus \pi^{\ast}({\bf sl} (\FF^{\prime}) ) = \pi^{\ast}({\bf gl}(\FF^{\prime}) )\,.
$$
\end{pf}

A supply of simple configurations is given by the classical algebro-geometric notion of points in general position
(see e.g. \cite{[G-H]}).
\begin{cor}\label{s-gp}
Let $Z$ be a configuration of $d$ points on $X$ such that the index of $L$-speciality
$\delta(L,Z) =r+1 \geq 2$ and $d \geq r+2$. Assume $Z$ to be in general position with respect to the adjoint linear system
$\left| K_X + L \right|$.
Then $Z$ is simple.
\end{cor}
\begin{pf}
Let $\GA$ be an admissible component in $\CSA$ containing $[Z]$.
Then by \RI, Corollary 7.13,
$\HO{Z})$ is an irreducible 
$\mbox{\BM$\LAGT$}\ZA$-module, for any
$\alpha \in \EZ$ such that $\ZA \in \JABG$.
This implies that the center
$\mbox{\BM${\cal C}$}\ZA$ is one dimensional.
Hence the rank of the center
$\mbox{\BM${\cal C}$}_{\GA}$ is equal to 1. By Definition \ref{s-c} the component
$\GA$ is simple.
\end{pf}

The two kinds of components - quasi-abelian (Definition \ref{qa-c}) and simple (Definition \ref{s-c}) -
are prototypical in a sense that a general situation can be reduced to these two types.
This is explained in the discussion below.

Let us go back to the diagram (\ref{diag-morph}). To begin with, we explain how the reduction alluded to above works on a fibre
of the morphism $p_2$.

Let $[Z]$ be a point in $\GA$ and consider the decomposition
\BEN\label{Z-c-dec1}
Z = \sum_{\lambda \in p^{\prime -1}_c ([Z])} Z_{\lambda}
\EEN
as in Corollary \ref{Z-c-dec}. We know that this decomposition is determined by the weight decomposition
$$
\HO{Z^{\prime}}) = \bigoplus_{\lambda \in p^{\prime -1}_c ([Z])} V_{\lambda} ([Z])
$$
in (\ref{Zpr-wd2}), where $Z^{\prime} = f(Z)$. Denote by 
$\Lambda_{[Z]} =p^{\prime -1}_c ([Z])$ the set of weights of this decomposition and divide it into two subsets
according to the dimension of the corresponding weight spaces
\BEN\label{Lam1-2}
\Lambda^1_{[Z]} =\{\lambda \in \Lambda_{[Z]} \mid dim(V_{\lambda} ([Z])) =1 \}
\,\,\mbox{and}\,\,\,
\Lambda^{\geq 2}_{[Z]} =\{\lambda \in \Lambda_{[Z]} \mid dim(V_{\lambda} ([Z])) \geq 2 \}.
\EEN
This separates the terms in (\ref{Z-c-dec1}) into two parts
\BEN\label{Z-ty-dec}
Z=Z^1 + \sum_{\lambda \in \Lambda^{\geq 2}_{[Z]}} Z_{\lambda}\,,
\EEN
where 
$\displaystyle{Z^1 =\sum_{\lambda \in \Lambda^{1}_{[Z]}} Z_{\lambda} }$.

On the side of the Lie algebra 
{\BM
$\LAGT$}$\ZA$ we have the following decomposition
\BEN\label{Lie-ty-dec}
\mbox{\BM$\GT$}_{\GA} \ZA = \mbox{\BM${\cal C}$}_{\GA} \ZA \oplus \mbox{\BM$\GS$}_{\GA} \ZA
= \mbox{\BM${\cal C}$}_{\GA} \ZA \oplus \left( \bigoplus_{\lambda \in \Lambda^{\geq 2}_{[Z]}} {\bf sl}(V_{\lambda} ([Z])) \right)\,,
\EEN
where the second equality comes from Proposition \ref{ga-dec}. Furthermore, we can write the center in the following way
$$
\mbox{\BM${\cal C}$}_{\GA} \ZA  =\left( \bigoplus_{\lambda \in \Lambda^1_{[Z]}} {\bf gl}(V_{\lambda} ([Z])) \right) \oplus
\left( \bigoplus_{\lambda \in \Lambda^{\geq 2}_{[Z]} } \CC id_{V_{\lambda} ([Z])} \right)\,.
$$

Comparing the geometric decomposition in (\ref{Z-ty-dec}) with the Lie algebraic decomposition in
(\ref{Lie-ty-dec}), we see that the subconfigurations
$Z_{\lambda}$ with $\lambda \in \Lambda^{\geq 2}_{[Z]}$ are precisely the ones which contribute simple factors
into the decomposition of 
$\mbox{\BM$\GT$}_{\GA} \ZA$, while the subconfiguration
$Z^1$ contributes to the center of
$\mbox{\BM$\GT$}_{\GA} \ZA$ only. Thus on the basis of this matching between the subconfigurations in
 (\ref{Z-c-dec1}) and the summands in 
(\ref{Lie-ty-dec}) we deduce the following.
\begin{thm}\label{ty-dec}
Let $\GA$ be a component in $\CSA$ and let $[Z] \in \GAB$.
Then the decomposition
$$
Z= Z^1 + \sum_{\lambda \in \Lambda^{\geq 2}_{[Z]}} Z_{\lambda}
$$
described in (\ref{Z-ty-dec}) provides the decomposition of $Z$ into the disjoint union of subconfigurations of two types
\begin{enumerate}
\item[1)]
$Z_{\lambda}$ is simple, for all $\lambda \in \Lambda^{\geq 2}_{[Z]}$,
\item[2)]
$Z^1$ is quasi-abelian, provided\footnote{the assumption $Card(\Lambda^{1}_{[Z]}) \geq 2$
 is needed to insure that the index of
$L$-speciality 
$\delta (L, Z^1) \geq 2$ which is our convention in \S\ref{c-n1}.} 
$Card(\Lambda^{1}_{[Z]}) \geq 2$.
\end{enumerate}
Furthermore,
$Z$ is quasi-abelian (resp. simple) if and only if
$\Lambda^{\geq 2}_{[Z]} =\emptyset$
(resp. $\Lambda^{1}_{[Z]}=\emptyset$ and $Card(\Lambda^{\geq 2}_{[Z]}) =1$).
\end{thm}

This result generalizes readily to the whole family of configurations
$p_2 : \ZD_{\GAB} \longrightarrow \GAB$. To do this consider
the morphism
$$
p^{\prime}_c : \ZD^{\prime}_c \longrightarrow \GAB
$$
in (\ref{diag-morph}). Take our sheaf 
$\FF^{\prime}$ 
on $\GAB$ and consider its pullback
$(p^{\prime}_c )^{\ast}(\FF^{\prime})$ to $\ZD^{\prime}_c$.
By definition $\ZD^{\prime}_c$ is the variety parametrizing the weights of the action of the center
{\BM
$\CG$}
on $\FF^{\prime}$. So we can think of points of  $\ZD^{\prime}_c$ as pairs
$([Z],\lambda)$, where $[Z] \in \GAB$ and $\lambda$ is a weight occurring in the decomposition
of
$\HO{Z^{\prime}})$ in (\ref{Zpr-wd2}). In this way we see that the fibre
$(p^{\prime}_c )^{\ast}(\FF^{\prime})_{([Z],\lambda)} =\HO{Z^{\prime}})$
at a point $([Z],\lambda) \in \ZD^{\prime}_c$ comes along with a distinguished subspace
$V_{\lambda} ( [Z])$, the weight subspace of $\HO{Z^{\prime}})$ corresponding to $\lambda$.
As $([Z],\lambda)$ varies in $\ZD^{\prime}_c$ the subspaces $V_{\lambda} ( [Z])$ fit together to form a distinguished subsheaf
${\cal{V}}$ of $(p^{\prime}_c )^{\ast}(\FF^{\prime})$. Furthermore, the dimension
of the fibre $V_{\lambda} ( [Z])$ of ${\cal{V}}$ at $([Z],\lambda) \in \ZD^{\prime}_c$
is equal to the degree of the fibre $Z^{\prime}_{\lambda}$ of the morphism\footnote{the equality
$dim(V_{\lambda} ( [Z])) =degZ^{\prime}_{\lambda}$ follows from the fact that $V_{\lambda} ( [Z])$ can be canonically identified with
$\HO{Z^{\prime}_{\lambda}})$ as it was done in the proof of Corollary \ref{Z-c-dec}.}
 $f^{\prime}_c$ in (\ref{diag-morph}) over 
$([Z],\lambda)$. In particular, the fibre dimension  of ${\cal{V}}$ is constant over every connected component of $\ZD^{\prime}_c$. Hence 
the restriction of
${\cal{V}}$ to every connected component of $\ZD^{\prime}_c$ is locally free.

We can now define a continuous version of the sets $\Lambda_1$ and $\Lambda_2$ in Theorem \ref{ty-dec}.
Namely, set
$\pi_{0} (\ZD^{\prime}_c)$ to be the set of connected\footnote{since $\ZD^{\prime}_c$ is smooth, this is the same as the set of irreducible components of $\ZD^{\prime}_c$.} components of
$\ZD^{\prime}_c$. For a connected component $W \in \pi_{0} (\ZD^{\prime}_c)$ denote by 
${\cal{V}}_W$ the restriction of the sheaf ${\cal{V}}$ to $W$. By analogy with (\ref{Lam1-2}) we divide the set of components
$\pi_{0} (\ZD^{\prime}_c)$ into two disjoint subsets
\BEN\label{pi01-pi02}
\pi^{1}_{0} (\ZD^{\prime}_c) = \left\{ \left.W \in \pi_{0} (\ZD^{\prime}_c) \right| rk({\cal{V}}_W) =1 \right\}\,\, 
{\mbox{and}}\,\,
\pi^{\geq 2}_{0} (\ZD^{\prime}_c) = \left\{ \left.W \in \pi_{0} (\ZD^{\prime}_c) \right| rk({\cal{V}}_W)  \geq 2 \right\}.
\EEN

 For every connected component $W \in \pi_{0} (\ZD^{\prime}_c)$, set
$$
\ZD^{\prime}_W =f^{\prime -1}_c (W)\,\,\,\, \ZD_W = f^{-1} (\ZD^{\prime}_W) =f^{-1}_c (W)
$$
to be the inverse image of $W$ by 
$f^{\prime }_c$ and $f_c =f^{\prime }_c \circ f$ respectively. This gives us the diagram analogous to the one in (\ref{diag-morph})
\BEN\label{diag-morph-W}
\xymatrix{
{\ZD_{W}} \ar[r]^{{}^W \!f} \ar[dr]_{{}^W \!p_2} & { \ZD^{\prime}_{W}} \ar[r]^{{}^W\! f^{\prime}_c} \ar[d]^{{}^W \!p^{\prime}_2} &
 {W} \ar[dl]^{{}^W\! p^{\prime}_c} \\
          &{\GAB}&                }
\EEN
where
${}^W\! f$ and ${}^W\! p_2$ (resp. ${}^W\! f^{\prime}_c$ and ${}^W\! p^{\prime}_2$) denote the restrictions to
$\ZD_{W}$ (resp. $\ZD^{\prime}_{W}$) of the morphisms $f$ and $p_2$ (resp. $f^{\prime }_c$ and $p^{\prime}_2$)
in (\ref{diag-morph}). 
Composing the horizontal arrows in (\ref{diag-morph-W}) gives the diagram
\BEN\label{diag-c-W}
\xymatrix{
{\ZD_{W} } \ar[rr]^{{}^W\!f_c} \ar[dr]_{{}^W\! p_2} &  &W \ar[dl]^{{}^W\! p^{\prime}_c} \\
          &{\GAB}&                }
\EEN
where ${}^W\! f_c = {}^W\! f^{\prime}_c  \circ {}^W\! f$. With this notation in mind we have the following
decomposition of $\ZD_{\GAB}$ into the disjoint union
$$
\ZD_{\GAB} = \bigsqcup_{W \in  \pi_{0} (\ZD^{\prime}_c)} \ZD_{W}\,.
$$
Setting
\BEN\label{Z1-cont}
\ZD^1_{\GAB} = \bigsqcup_{W \in  \pi^1_{0} (\ZD^{\prime}_c)} \ZD_{W}\,\,
\EEN
gives the following decomposition of $\ZD_{\GAB} $ into the disjoint union
\BEN\label{Z-c-dec-cont}
\ZD_{\GAB} =\ZD^1_{\GAB} \sqcup \left( \bigsqcup_{\scriptscriptstyle{W \in  \pi^{\geq 2}_{0} (\ZD^{\prime}_c)}} \ZD_{W} \right)\,.
\EEN
Thus one obtains the following 
``continuous'' analogue of Theorem \ref{ty-dec}.
\begin{thm}\label{ty-dec-cont}
Let $\GA$ be a component in $\CSA$ and let
$\ZD^{\prime}_c$ be its scheme of central weights (Definition \ref{sc-cw}).
Then the decomposition
$$
\ZD_{\GAB} =\ZD^1_{\GAB} \sqcup \left( \bigsqcup_{\scriptscriptstyle{W \in  \pi^{\geq 2}_{0} (\ZD^{\prime}_c)}} \ZD_{W} \right)
$$
described in (\ref{Z-c-dec-cont}) provides the decomposition of $\ZD_{\GAB}$ into the disjoint union of  families of
subconfigurations of two types
\begin{enumerate}
\item[1)]
${}^W\!f_c : \ZD_W \longrightarrow W$ is a family of simple configurations, for all $W \in \pi^{\geq 2}_{0} (\ZD^{\prime}_c)$,
\item[2)]
$\ZD^1 \longrightarrow \GAB$ is a family of  quasi-abelian configurations, provided the covering
$$
p^{\prime}_c :  \bigsqcup_{W \in  \pi^{1}_{0} (\ZD^{\prime}_c)} W \longrightarrow \GAB
$$
has degree\footnote{the degree assumption is needed for the same reason as in the footnote in Theorem \ref{ty-dec}, 2).}
 $\geq 2$, i.e. $\sum_{W \in  \pi^{1}_{0} (\ZD^{\prime}_c)} deg({}^{W}\!{p^{\prime}_c }) \geq 2$.
\end{enumerate}

Furthermore,
$\GA$ is quasi-abelian (resp. simple) if and only if
$ \pi^{\geq 2}_{0} (\ZD^{\prime}_c) =\emptyset$
(resp. $\ZD^{\prime}_c = \GAB$).
\end{thm}

This result shows that the study of configurations on $X$ can be reduced to either
quasi-abelian or simple ones. It should be also clear that the quasi-abelian configurations are quite special and to our mind are akin to hyperelliptic
divisors on curves.\footnote{a study of the quasi-abelian configurations will appear elsewhere.}
On the other hand, if $\GA$ is neither quasi-abelian nor simple, then the set
$\pi^{\geq 2}_0 (\ZD^{\prime}_c)$ is not empty.
So replacing the original family 
$p_2 : \ZD_{\GAB} \longrightarrow \GAB$
by the family
$$
{}^{W}\!f_c :\ZD_W \longrightarrow W
$$
in (\ref{diag-c-W}), corresponding to a connected component 
$W \in  \pi^{\geq 2}_0 (\ZD^{\prime}_c)$,
we  obtain a reduction to a family of simple configurations. Thus in studying the components of $\CS$, the assumption that
the set $\CSA$ contains simple components is not essential.

\subsection{A natural grading of {\BM$\LAG$}}\label{sec-gr}
Let $\GA$ be a component in $\CSA$ and 
let 
{\BM
$\LAGT$ }
be the corresponding sheaf of Lie algebras
on $\JABG$. As it was recalled in \S\ref{Lie} this sheaf is generated by certain local sections
of $\ENDO(\HT_{-\LG})$ of degree $\pm1$ and $0$ with respect to the grading on
$\HT_{-\LG}$ given by the orthogonal decomposition 
\BEN\label{ordH-l}
\HT_{-\LG } =\bigoplus^{\LG-1}_{p=0} \HH^p
\EEN
 as in (\ref{ordH-i}), for $i=\LG$. Thus the sheaf
{\BM
$\LAGT$}
comes along with a natural grading
{\BM
\BEN\label{Lie-grad}
\LAGT =\bigoplus^{\LG - 1}_{i=-(\LG-1)} \GT^i_{\mbox{\UB$\GA$}}\,,
\EEN
where the subsheaf
$\GT^i_{\mbox{\UB$\GA$}}$ is formed by local sections
$\phi$ of $\LAGT$ of degree $i$ with respect to the grading in (\ref{ordH-l}), i.e.
the restriction of $\phi$ to a summand
$\HH^p$ is a local section of
$\HOM(\HH^p, \HH^{p+i})$}, 
for every $p\in \{0,1,\ldots, \LG-1 \}$. With this gradation the sheaf
{\BM
$\LAGT$}
becomes a sheaf of {\it graded} Lie algebras (see \RI, (7.7), for details).

The same holds for the subsheaf 
{\BM 
$\LAG =[\LAGT,\LAGT]$:
\BEN\label{sLie-grad}
\LAG = \bigoplus^{\LG - 1}_{i=-(\LG-1)} \GS^i_{\mbox{\UB$\GA$}}\,,
\EEN
while from the study of the center
$\CG$ in \S\ref{Center} we know that its local sections are grading preserving, i.e.
\BEN\label{C-gr0}
\CG \subset \GT^0_{\mbox{\UB$\GA$}}.
\EEN
Observe that 
$ \GT^0_{\mbox{\UB$\GA$}}$ is the subsheaf of Lie {\it subalgebras} of
$\LAGT$. Then the structure decomposition in (\ref{str-d}) together with  the inclusion in (\ref{C-gr0}) give the following
\BEN\label{Lie0-dec}
\GT^0_{\mbox{\UB$\GA$}} = \CG \oplus \GS^0_{\mbox{\UB$\GA$}}.
\EEN
Furthermore, by \RI, Proposition 7.17, 
$\GS^0_{\mbox{\UB$\GA$}}$ is a subsheaf of reductive Lie subalgebras of 
$\LAG$  and the structure decomposition for it yields
\BEN\label{str-dec-0}
\GS^0_{\mbox{\UB$\GA$}} ={\cal C}^0_{\mbox{\UB$\GA$}}\oplus  {}^s{\GS^0_{\mbox{\UB$\GA$}}}\,,
\EEN
where 
${\cal C}^0_{\mbox{\UB$\GA$}}$  
and
$ {}^s{\GS^0_{\mbox{\UB$\GA$}}} = [\GS^0_{\mbox{\UB$\GA$}}, \GS^0_{\mbox{\UB$\GA$}}]$
are, respectively, the center and the semisimple part of
$\GS^0_{\mbox{\UB$\GA$}}$.
\\
\\
\indent
We will now  give a more detailed description of the gradation of 
$\LAG$
}
in (\ref{sLie-grad}) in the case of $\GA$ being a simple component in $\CSA$. 

By Corollary \ref{cor-s} the sheaf 
{\BM
$\LAGT$ (resp. $\LAG$) is the pullback of
}
${\bf gl}(\FF^{\prime})$ (resp. ${\bf sl}(\FF^{\prime})$) by the natural projection
$
\pi: \JABG \longrightarrow \GAB
$. Thus we obtain
\begin{eqnarray}
\mbox{\BM$\LAGT$} = \pi^{\ast}{\bf gl}(\FF^{\prime})& =&{\bf gl} \big{(} \pi^{\ast} (\FF^{\prime}) \big{)}= 
{\bf gl} (\FT^{\prime}) \\ \label{Lie-s1}
\mbox{\BM$\LAG$} = \pi^{\ast}{\bf sl}(\FF^{\prime})& =&{\bf sl} \big{(} \pi^{\ast} (\FF^{\prime}) \big{)}= 
{\bf sl} (\FT^{\prime})\,, \label{Lie-s2}
\end{eqnarray}
where we set
\BEN\label{FTpr-def}
\FT^{\prime} =\pi^{\ast} (\FF^{\prime}).
\EEN

Our first step in understanding the grading (\ref{sLie-grad})  in the case of $\GA$ being simple is  to calculate
{\BM
$\GS^0_{\mbox{\UB$\GA$}}$ (resp.  ${}^s{\GS^0_{\mbox{\UB$\GA$}}}$).
}
\begin{pro}\label{pro-Lie0-s}
Let $\GA$ be a simple component in $\CSA$. Then
{\BM
$$
{}^s{\GS^0}_{\mbox{\UB$\GA$}} = \bigoplus^{\LG-1}_{p=0}  sl(\HH^p)
$$
and the center
${\cal C}^0_{\mbox{\UB$\GA$}}$ is the subsheaf of $\GS^0_{\mbox{\UB$\GA$}}$
},
whose local sections $\phi$ have the following form
$$
\phi= \sum^{\LG-1}_{p=0} c_p id_{\HH^p}\,,
$$
where $c_p$'s are local sections of 
$\OO_{\JABG}$ such that
$$
\sum^{\LG-1}_{p=0} c_p =0
$$
\end{pro}
\begin{pf}
The result follows immediately from (\ref{Lie-s2}) and the orthogonal decomposition of
$\FT^{\prime}$ in (\ref{ordH-l}).
\end{pf}

 To  describe other graded pieces
{\BM
$\GS^i_{\mbox{\UB$\GA$}}$ of
$\LAG$}
in the decomposition  (\ref{sLie-grad}) it is useful to make a general observation:
\begin{center}
{\it each
{\BM 
$\GS^i_{\mbox{\UB$\GA$}}$
is a 
$\GS^0_{\mbox{\UB$\GA$}}$-module and, in particular, it is
${\cal C}^0_{\mbox{\UB$\GA$}}$-module.}}
\end{center}

We return now to the case of $\GA$ being simple and describe
{\BM
$\GS^i_{\mbox{\UB$\GA$}}$, for $i\neq 0$, together with 
its weight decomposition under the action of
${\cal C}^0_{\mbox{\UB$\GA$}}$.
}
 
Set $e_p =id_{\HH^p}$ to be the identity endomorphism of $\HH^p$ and let
\BEN\label{C0}
C^0=\CC \{e_0, \ldots, e_{\LG-1}  \}
\EEN
be the complex vector space spanned by
$e_0, \ldots, e_{\LG-1} $. Denote by
$(C^0)^{\ast}$ the dual of
$C^0$, equipped with the basis
$\mu_0, \ldots, \mu_{\LG-1}$ dual to
$e_0, \ldots, e_{\LG-1} $ and set
\BEN\label{wei}
\nu_{ij} =\mu_i - \mu_j\,.
\EEN
Then $\nu_{ij}$, for $i\neq j$, is easily seen to be the weight of ${\cal C}^0_{\mbox{\UB$\GA$}}$-action
on $\HOM(\HH^j, \HH^{i} )$. This gives the following.
\begin{pro} \label{pro-Liei-s}
Let $\GA$ be a simple component in $\CSA$. Then for every $i\neq 0$, one has
{\BM
$$
\GS^i_{\mbox{\UB$\GA$}} = \bigoplus^{\LG-1}_{p=0} \HOM(\HH^p, \HH^{p+i} )\,,
$$
where the direct sum on the right hand side is the weight decomposition of 
$ \GS^i_{\mbox{\UB$\GA$}}$ under the action of 
${\cal C}^0_{\mbox{\UB$\GA$}}$ with the  summand
$ \HOM(\HH^p, \HH^{p+i} )$ being the weight-subsheaf 
corresponding to the weight
$\mbox{\UB$\nu_{p+i,p}$}$ in (\ref{wei}), for $\mbox{\UB$p=0,\ldots, \LG-1$}$.
Furthermore, the summands
$ \HOM(\HH^p, \HH^{p+i} )$  are irreducible
${}^s \GS^0_{\mbox{\UB$\GA$}}$-modules.}
\end{pro}
\begin{pf}
All the assertions are immediate from (\ref{Lie-s2}), the orthogonal decomposition of
$\FT^{\prime}$ in (\ref{ordH-l}) and Proposition \ref{pro-Lie0-s}.
\end{pf}

Substituting the decompositions of Proposition \ref{pro-Liei-s} into (\ref{sLie-grad}) yields
{\BM
\BEN\label{LAG-wd}
\GS_{\mbox{\UB$\GA$}} = \GS^0_{\mbox{\UB$\GA$}} \oplus \left(\bigoplus_{i\neq j} \HOM(\HH^j, \HH^{i} )\right)\,,
\EEN
the decomposition of 
$\LAG$
into the weight-sheaves of ${\cal C}^0_{\mbox{\UB$\GA$}}$}-action.
\begin{rem}\label{quiver-sl}
Observe that the set
\BEN\label{sl-lg-roots}
R_{\LG} = \left\{\left. \nu_{ij} \right| i\neq j \in \{0,\ldots,\LG-1 \} \right\}
\EEN
can also be identified with the set of roots of 
${\bf sl}_{\LG} (\CC)$. Setting
\BEN\label{W-nu}
\mbox{\BM${\cal W}$}_{\nu_{ij}} = \mbox{\BM$\HOM(\HH^j, \HH^{i} )$}
\EEN
the decomposition in (\ref{LAG-wd}) can be rewritten as follows
\BEN\label{LAG-roots}
\mbox{\BM$\GS_{\mbox{\UB$\GA$}} = \GS^0_{\mbox{\UB$\GA$}}$} \oplus \left(\bigoplus_{\nu \in R_{\LG}} \mbox{\BM${\cal W}$}_{\nu} \right)\,.
\EEN

Let us also observe that the appearance of ${\bf sl}_{\LG} (\CC)$ with a distinguished set of roots $R_{\LG}$ (and its polarization) is
determined by the orthogonal decomposition (\ref{ordH-l}) and the triangular decomposition (\ref{d-D}). Indeed, the decomposition 
(\ref{ordH-l}) and the operators $D^{\pm}$ in (\ref{d-D}) can be viewed as the following quiver
\BEN\label{quiv}
\xymatrix{
*=0{\bullet}\ar@/^/[r]^(.0){0} \ar@/_/@{{<}{-}}[r] &*=0 {\bullet}\ar@/^/[r]^(.15){1} \ar@/_/@{{<}{-}}[r] &*=0{\bullet}\ar@{}[r]^(.0){2}&{\cdots}&*=0{\bullet}\ar@/^/[r] \ar@/_/@{{<}{-}}[r]&*=0{\bullet}\ar@{}[r]^(.3){\LG-1}& }
\EEN
The vertices are labeled by integers $\{0,1,\dots,\LG -1\}$ from left to right and represent the summands $\HH^p\,(p=0,\ldots,\LG-1)$
of the decomposition in (\ref{ordH-l}) and the arrows between the neighboring vertices $p$ and $(p+1)$ represent the action of operators
$D^{+}_p$ (the upper arrow) and $D^{-}_{p+1}$ (lower arrow), where $D^{\pm}_p$ is the restriction to $\HH^p$ of the operators 
$D^{\pm}$ in (\ref{d-D}).

Take the ordered set of vertices $\{\{0\},\{1\},\ldots,\{ \LG-1 \} \}$ of the quiver in (\ref{quiv}) and form the vector space
$$
V = \CC\{\{0\},\{1\},\ldots,\{ \LG-1 \} \}.
$$
The vectors $e_p$'s in $C^0$ in (\ref{C0}) can be thought of as endomorphisms of $V$ fixing the $p$-th vertex and annihilating all others.
Thus $C^0$ becomes a distinguished Cartan subalgebra of ${\bf gl}(V)$ while
$$
\goth{h}^0 =\left\{ \left.\sum^{\LG-1}_{p=0} c_p e_p \right| \sum^{\LG-1}_{p=0} c_p =0 \right\}
$$
gives a distinguished Cartan subalgebra of ${\bf sl}(V)$.

From this it follows that that the set
$$ 
R_{\LG} = \left\{ \left.\nu_{ij} \right| i\neq j \in \{0,\ldots,\LG-1 \} \right\}
$$
as in (\ref{sl-lg-roots}) is the set of roots of ${\bf sl}(V)$ with respect to the Cartan subalgebra $\goth{h}^0$. Furthermore, the set
$$
R^{+}_{\LG} = \left\{\left. \nu_{ij} \right| 0\leq i < j\leq \LG-1  \right\} \,\,(resp. \,\,R^{-}_{\LG} = \left\{\left. \nu_{ij} \right| \LG-1 \geq i >j\geq 0 \right\} )
$$
is a subset of positive (resp. negative) roots of $R_{\LG}$, while the roots
$\nu_p = \nu_{p,p+1}$ (resp. $-\nu_p$), for $p=0,1, \ldots,\LG-2$, are positive (resp. negative) simple roots of ${\bf sl}(V)$ with respect to the Cartan subalgebra $\goth{h}^0$. This way the edges of the quiver in (\ref{quiv}) can be identified with preferred generators of the root spaces
$ ({\bf sl}(V))_{\pm\nu_p }\,\,(p=0,1, \ldots,\LG-2)$, while the operators $D^{+}_p$ (resp. $D^{-}_{p+1}$), for
$ p=0,1, \ldots,\LG-2$, become representations of the edges of the quiver in the category of $\OO_{\JABG}$-modules.
\end{rem}
\section{Period maps and Torelli problems}\label{sec-periods}
In this section we take a more geometric point of view on the orthogonal decomposition
\BEN\label{ordH-l1}
 \HT_{-\LG } =\bigoplus^{\LG-1}_{p=0} \HH^p
\EEN
resulting from (\ref{ordH-i}).
Namely, we suggest to view it as a Hodge-like decomposition and view the spaces
$\{\HH^p \ZA \}_{p=0,\ldots, \LG-1 }$, the fibres of the sheaves $\HH^p$'s at
$\ZA \in \JABG$, as periods associated to points of  $\JABG$. This allows us to define
the period map(s) for $\JABG$. Furthermore, the variation of these periods with respect to $[\alpha]$ is related to
the multiplication  in $\HT_{-\LG } $ by local sections of 
$\HT =\HH^0$ (see (\ref{ordH-1}) for this equality). On the other hand, from the study of 
{\BM
$\LAGT$}
in \S\ref{sec-Lie1} we know that the multiplication in
$\HT_{-\LG }$
coincides with the action of the sheaf
{\BM
${\cal C}( \HT)$
of Cartan subalgebras of 
$\LAGT$ (see Remark \ref{act-m}).
Hence the results about 
$\LAGT$}
and its action on
$\HT_{-\LG }$
can be reinterpreted as properties of the aforementioned period map.

Once the period map for $\JABG$ is in place, one can formulate Torelli-type problems.
One of the main results of this section is that these problems have positive solution
precisely over simple components in the sense of Definition \ref{s-c}.

\subsection{Definition of the period map(s) for $\JABG$.}
We begin with the sheaf $\FF^{\prime}$ defined in Corollary \ref{cor-Fpr}
and set
\BEN\label{Fpr-tilde}
\FT^{\prime} =\pi^{\ast} \FF^{\prime}.
\EEN
This together with (\ref{pi-Fpr}) give the equality
\BEN\label{Fpr-til=Hl}
\FT^{\prime} = \HT_{-\LG}\,.
\EEN
Combining this with the orthogonal decomposition (\ref{ordH-l1}) we obtain
\BEN\label{ord-Fpr-til}
\FT^{\prime} = \HT_{-\LG} =\bigoplus^{\LG-1}_{p=0} \HH^p \,.
\EEN

Using the filtration $\FI^{\bullet}$ in (\ref{filtF}), we define
\BEN\label{Fpr-p}
{}^{\prime}{\FI^p} = \FI^p \cap \FT^{\prime}
\EEN
to obtain the following filtration of $\FT^{\prime}$
\BEN\label{filtFpr}
\FT^{\prime} ={}^{\prime}{\FI^0} \supset {}^{\prime}{\FI^1} \supset \ldots \supset {}^{\prime}{\FI^{\LG-1}}
 \supset {}^{\prime}{\FI^{\LG}} =0.
\EEN
From the orthogonal decomposition of  $\FT^{\prime} $ in (\ref{ord-Fpr-til}) and the orthogonal decomposition of
$\FI^p$ in (\ref{ordF-i}) 
one deduces the orthogonal decomposition
\BEN\label{ordPF-p}
\PF^p = \bigoplus^{\LG-1}_{i=p}  \HH^i.
\EEN
Set 
$$
Gr^p (\PF^{\bullet}) = \PF^p / \PF^{p+1}\,,\,\,\,for\,\, p=0,1, \ldots, \LG-1,
$$
to be the associated graded sheaves and observe a natural identification
\BEN\label{Gr=H}
Gr^p (\PF^{\bullet}) \cong \HH^p,
\EEN
for every $p=0,1, \ldots, \LG-1$. This isomorphism together with Remark \ref{rkH} imply
\BEN\label{rkGr}
rk(Gr^p (\PF^{\bullet})) =rk( \HH^p) =h^p_{\GA}\,.
\EEN
Define
\BEN\label{Hv-r}
\overrightarrow{h^{\prime}}_{\GA} =(h^0_{\GA},\ldots, h^{l_{\GA} -1} _{\GA})
\EEN
to be the {\it reduced Hilbert vector of} $\GA$ (compare with the definition of 
$h_{\GA}$ in Lemma \ref{hg=c}).

Consider the scheme
$\mbox{\BM${\FL}$}_{\overrightarrow{h^{\prime}}_{\GA}}$ of relative partial flags of type
$\overrightarrow{h^{\prime}}_{\GA}$ in $\FT^{\prime}$, i.e.
$\mbox{\BM${\FL}$}_{\overrightarrow{h^{\prime}}_{\GA}}$ is the scheme over $\GAB$ with the structure morphism
\BEN\label{fl}
\overrightarrow{Fl}_{\GA} : \FLA \longrightarrow \GAB
\EEN
such that the fibre 
$\FLA ([Z])$ over a closed point $[Z] \in \GAB$ is the variety of partial flags of type
$\overrightarrow{h^{\prime}}_{\GA}$ in the vector space
$\FT^{\prime} ([Z])$, the fibre of $\FT^{\prime}$ at $[Z]$. Recalling the identification
\BEN\label{Fpr-Z=f-Zpr}
\FT^{\prime} ([Z]) = \HO{Z^{\prime}})
\EEN
in Corollary \ref{cor-Fpr}, where $Z^{\prime} =p^{\prime -1}_2([Z])$ is the fibre over $[Z]$ of $p^{\prime }_2$ in 
(\ref{Z-Zpr}), we can describe the set of closed points of 
$\FLA ([Z])$ as follows
\BEN\label{Fla-Z}
\FLA([Z]) =\left\{[F] =\left.[\scriptstyle{\HO{Z^{\prime}}) =F^{0} \supset F^{1} \supset \ldots \supset F^{\LG-1} \supset F^{\LG} =0}] \right| 
\scriptstyle{dim(F^p /F^{p+1}) =h^p_{\GA},\,\,{\rm{for}}\,\, 0\leq p \leq \LG-1 }\right\}.
\EEN
By the universality of 
$\FLA$
we have the morphism
\BEN\label{p}
p_{\GA} : \JABG \longrightarrow \FLA
\EEN
of $\GAB$-schemes which sends closed points 
$\ZA$ of $\JABG$ to the partial flag 
\BEN\label{p-ZA}
[\PF^{\bullet} \ZA]=[\scriptstyle{\HO{Z^{\prime}})={}^{\prime}{\FI^0} \ZA \supset {}^{\prime}{\FI^1} \ZA
\supset \ldots \supset {}^{\prime}{\FI^{\LG-1}} \ZA
 \supset {}^{\prime}{\FI^{\LG}} \ZA=0}]
\EEN
determined by the filtration (\ref{filtFpr}) at $\ZA$.
\begin{defi}\label{period}
The morphism $p_{\GA}$ in (\ref{p}) is called the period map of $\JABG$.
\end{defi}

Set
\BEN\label{rel-tan}
\mbox{\BM${\cal T}$}_{\overrightarrow{Fl}_{\GA}} =\mbox{\BM${\cal T}$}_{\FLA / {\GAB}}
\EEN
to be the relative tangent sheaf of the morphism
$Fl_{\GA}$ in (\ref{fl}). This is a locally free sheaf  on
$\FLA$, since $\overrightarrow{Fl}_{\GA}$ is a smooth morphism.
 Its dual
$\mbox{\BM${\cal T}$}^{\ast}_{\overrightarrow{Fl}_{\GA}}$ is the relative cotangent sheaf of
$\overrightarrow{Fl}_{\GA}$. Invoking our convention in \S\ref{sh-bdl} we denote by 
$\mbox{\BM${T}$}^{\ast}_{\overrightarrow{Fl}_{\GA}}$  (resp.
$\mbox{\BM$ T$}_{\overrightarrow{Fl}_{\GA}}$ ) the relative cotangent (resp. tangent) {\it bundle} over 
$\FLA$. Thus 
$\mbox{\BM${ T}$}^{\ast}_{\overrightarrow{Fl}_{\GA}}$ is the scheme over $\FLA$ with the structure morphism
\BEN\label{rel-cot-bdl}
\sigma: \mbox{\BM${T}$}^{\ast}_{\overrightarrow{Fl}_{\GA}} \longrightarrow \FLA
\EEN
of schemes over $\GAB$. In particular, the fibre of $ \mbox{\BM${T}$}^{\ast}_{\overrightarrow{Fl}_{\GA}}$ over a point  $[Z] \in \GAB$ is
the cotangent bundle
$T^{\ast}_{\FLA([Z])}$ of the variety of partial flags
$\FLA([Z]) $ in (\ref{Fla-Z}).

It is well-known that $T^{\ast}_{\FLA([Z])}$ has the following description
\BEN\label{cot-fl}
T^{\ast}_{\FLA([Z])} =\left\{ \left.([F],x) \in {\scriptstyle{\FLA([Z]) \times End(\HO{Z^{\prime}}))}} \right|
\stackrel{\scriptstyle{[F] =[\HO{Z^{\prime}}) =F^0 \supset  \ldots \supset F^i \supset F^{i+1} \supset \ldots \supset F^{\LG} =0],}} 
{\scriptstyle{x(F^i) \subset F^{i+1}\,,\,\,\forall i\geq 0.} } \right\}.
\EEN

Let
\BEN\label{rel-tan1}
T_{\pi} = T_{{\JABG} / {\GAB}}
\EEN
be the relative tangent bundle\footnote{according to our convention in \S\ref{sh-bdl} the relative tangent {\it sheaf} is denoted by
${\cal{T}}_{\pi}$.}
 of the natural projection
\BEN\label{pi-2}
\pi: \JABG \longrightarrow \GAB\,.
\EEN

We will now define a canonical lifting of the period map $p_{\GA}$ in (\ref{p})
to a morphism
from  $T_{\pi}$ to the relative cotangent bundle $\mbox{\BM${T}$}^{\ast}_{\overrightarrow{Fl}_{\GA}}$.
\begin{pro}\label{pro-p+}
There exists a distinguished morphism
\BEN\label{p+}
p^{+}_{\GA} : T_{\pi} \longrightarrow \mbox{\BM${T}$}^{\ast}_{\overrightarrow{Fl}_{\GA}}
\EEN
for which the diagram
\BEN\label{p-p+}
\xymatrix{
T_{\pi}  \ar[r]^(.45){p^{+}_{\GA}} \ar[d]^{\tau} & {\mbox{\BM${T}$}^{\ast}_{\overrightarrow{Fl}_{\GA}}} \ar[d]^{\sigma} \\
{\JABG} \ar[r]^{p_{\GA}}& {\FLA} }
\EEN
commutes.\footnote{the vertical arrows in the diagram (\ref{p-p+}) are the natural projection.}
\end{pro}
\begin{pf}
First we recall from \RI, Proposition 1.4, that there is a distinguished isomorphism
\BEN\label{M}
M: \HT /{\OO_{\JABG}} \longrightarrow {\cal T}_{\pi}\,.
\EEN
 Also recall the morphism
$D$ in (\ref{m-D}) together with its triangular decomposition (\ref{d-D}). Taking the component
$D^{+}$ gives the morphism
\BEN\label{D+1}
D^{+} : \HT \longrightarrow \mbox{\BM$\GS$}_{\GA} \subset \ENDO(\FT^{\prime})
\EEN
Since $D^{+}$ vanishes\footnote{this vanishing comes from the following two facts: (1) the inclusion
$\OO_{\JABG} \hookrightarrow \HT$ takes the constant section $1_{\JABG}$ of $\OO_{\JABG}$ to the section
$h_0$ of $\HT$ whose value
$h_0 \ZA = 1_Z \in \HT\ZA$ is the constant function of value $1$ on $Z$, for every $\ZA \in \JABG$; (2)
$D^{+} (t)=0$, for any constant function $t$ on $Z$.}
 on $\OO_{\JABG}$, we deduce that 
$D^{+}$ factors through the quotient
$\HT /{\OO_{\JABG}}$. The resulting morphism
\BEN\label{D+2}
D^{+} : \HT /{\OO_{\JABG}} \longrightarrow \mbox{\BM$\GS$}_{\GA} \subset \ENDO(\FT^{\prime})
\EEN 
still will be denoted by $D^{+}$. Composing it with the inverse $M^{-1}$ of $M$ in (\ref{M}) we obtain the morphism
\BEN\label{d+}
d^{+} : {\cal T}_{\pi} \longrightarrow \mbox{\BM$\GS$}_{\GA} \subset \ENDO(\FT^{\prime})\,.
\EEN
Furthermore, the image of $d^{+}$ is contained in the summand
{\BM 
$\GS^1_{\mbox{$\GA$}}$ of the decomposition of
$\LAG$} in (\ref{sLie-grad}), 
i.e. for a tangent vector $v $ in the fibre of
${\cal T}_{\pi}$ at a point $\ZA \in \JABG$, the endomorphism
$$
d^{+}_{\ZA} (v) : \FT^{\prime} ([Z]) =\HO{Z^{\prime}}) \longrightarrow \HO{Z^{\prime}})
$$
has degree $1$ with respect to the grading in (\ref{ordH-l1}), where $d^{+}_{\ZA}$ stands for the restriction of 
$d^{+}$ to the fibre at $\ZA$ and $d^{+}_{\ZA} (v)$ is the value of $d^{+}_{\ZA}$ at $v$. Thus given a closed point
$([Z],[\alpha],v) \in \TPI$, the point
$(p_{\GA} \ZA, d^{+}_{\ZA} (v))$ lies, according to the description in (\ref{cot-fl}), in 
$T^{\ast}_{\FLA ([Z])}$. Hence setting
\BEN\label{p+formula}
p^{+}_{\GA} ([Z],[\alpha],v) =(p_{\GA} \ZA, d^{+}_{\ZA} (v))
\EEN
gives a well-defined map
$$
p^{+}_{\GA} : \TPI \longrightarrow \mbox{\BM${ T}$}^{\ast}_{\overrightarrow{Fl}_{\GA}}\,.
$$
Since the formula in (\ref{p+formula}) depends holomorphically on all parameters in
$([Z],[\alpha],v)$ it follows that $p^{+}_{\GA}$ is a morphism of varieties.

The commutativity of the diagram in (\ref{p-p+}) is part of the definition of
$p^{+}_{\GA}$ in (\ref{p+formula}).
\end{pf}
 
From the identification
in (\ref{Fpr-til=Hl})
it also follows that 
$\FT^{\prime}$ carries the filtration $\HT_{-\bullet}$ in (\ref{filtHT-JG}).
This gives rise to a companion period map which will be denoted
${}^{op}{p_{\GA}}$ and called {\it the opposite of } $p_{\GA}$:
\BEN\label{op-p}
{}^{op}{p_{\GA}} :\JABG \longrightarrow \FLAO\,,
\EEN
where $\FLAO$ is the scheme of relative partial flags in $\FT^{\prime}$ determined by the
{\it opposite} reduced Hilbert vector 
\BEN\label{oHv}
\overleftarrow{h^{\prime}}_{\GA}=(h^{\LG-1}_{\GA},\ldots,h^0_{\GA} ).
\EEN
This means that
$\mbox{\BM${\FL}$}_{\overleftarrow{h^{\prime}}_{\GA}}$ is the scheme over $\GAB$ with the structure morphism
\BEN\label{flo}
\overleftarrow{Fl}_{\GA} : \FLAO \longrightarrow \GAB
\EEN
such that the fibre 
$\FLAO ([Z])$ over a closed point $[Z] \in \GAB$ is the variety of partial flags of type
$\overleftarrow{h^{\prime}}_{\GA}$ in the vector space
$\FT^{\prime} ([Z]) =\HO{Z^{\prime}})$, i.e. 
 the set of closed points of 
$\FLAO ([Z])$ is as follows
\BEN\label{Flao-Z}
\FLAO([Z]) =\left\{[F] =\left.[\scriptstyle{\HO{Z^{\prime}}) =F_{\LG} \supset F_{\LG-1} \supset \ldots \supset F_{1} \supset F_{0} =0}] \right| 
\scriptstyle{dim(F_{p+1} /F_{p}) =h^p_{\GA},\,\,{\rm{for}}\,\, 0\leq p \leq \LG-1 }\right\}.
\EEN

By definition, the map ${}^{op}{p_{\GA}}$ in (\ref{op-p}) is a morphism of schemes over $\GAB$ which sends a closed point $\ZA \in \JABG$ to the partial flag
${}^{op}{p_{\GA}} \ZA$ determined by the filtration $\HT_{-\bullet}$, i.e. we have
\BEN\label{op-p-Z}
{}^{op}{p_{\GA}} \ZA =[\HO{Z^{\prime}}) =\HT_{-\LG} \ZA \supset \ldots \supset \HT_{-1} \ZA \supset \HT_0 \ZA =0]\,.
\EEN

This morphism also admits a distinguished lifting 
\BEN\label{p-}
p^{-}_{\GA} : \TPI \longrightarrow \mbox{\BM${ T}$}^{\ast}_{\overleftarrow{Fl}_{\GA}}
\EEN
defined\footnote{$\mbox{\BM${ T}$}^{\ast}_{\overleftarrow{Fl}_{\GA}}$ stands for the relative cotangent bundle of the structure morphism $\overleftarrow{Fl}_{\GA}$
in (\ref{flo}).}  by the formula
\BEN\label{p-formula}
p^{-}_{\GA}  ([Z],[\alpha],v) =({}^{op}{p_{\GA}} \ZA, d^{-}_{\ZA} (v))
\EEN
for every closed point $ ([Z],[\alpha],v) \in \TPI$. In this formula the morphism
\BEN\label{d-}
 d^{-}: {\cal T}_{\pi} \longrightarrow \mbox{\BM$\GS$}_{\GA} \subset \ENDO(\FT^{\prime})
\EEN
 is defined by composing  the inverse of $M$ in (\ref{M}) with the morphism
\BEN\label{D-}
D^{-} : \HT /{\OO_{\JABG}} \longrightarrow \mbox{\BM$\GS$}_{\GA} \subset \ENDO(\FT^{\prime})
\EEN 
defined in the same manner as $D^{+}$ (see (\ref{D+1}) and (\ref{D+2})) with the only difference
of using the component $D^{-}$ of the triangular decomposition in (\ref{d-D}), instead of $D^{+}$ used in the proof of 
Proposition \ref{pro-p+}.
\begin{rem}\label{val-d}
\begin{enumerate}
\item[1)]
In the sequel the value of $d^{\pm}$ on a tangent vector $v$ at a point $\ZA \in \JABG$ will be denoted by
 $d^{\pm} (v)$ (with the reference to $\ZA$ omitted).
\item[2)]
To calculate 
$d^{\pm} (v)$ 
one takes any lifting of $M^{-1} (v) \in {\HT \ZA}/{\CC \{1_Z \}} $ to a vector 
$\tilde{v}$ in $\HT \ZA$. Then 
\BEN\label{d=D}
d^{\pm} (v) = D^{\pm}(\tilde{v})
\EEN
In the sequel we refer to such $\tilde{v}$ as a lifting of $v$.
\end{enumerate}
\end{rem}

The two period maps, $p_{\GA}$ and ${}^{op}{p_{\GA}} $, are related by taking the orthogonal complement with respect to the quadratic form
$\QT$ defined $\FT$ in (\ref{q}), i.e.
$\PF^p$ in (\ref{ordPF-p}) is the orthogonal complement in $\FT^{\prime} = \HT_{-\LG}$ of
$\HT_{-p} =\bigoplus^{p-1}_{i=0} \HH^i$. This will be expressed as the orthogonality between two morphisms
\BEN\label{morph-ort}
 ( {}^{op}{p_{\GA}} )^{\perp} =  p_{\GA}\,.
\EEN
The liftings $p^{\pm}_{\GA}$ are related by the operation of taking the adjoint (with respect to the quadratic form
$\QT$ in $\FT^{\prime}$)
\BEN\label{adj}
(\cdot)^{\dagger} : \ENDO(\FT^{\prime}) \longrightarrow \ENDO (\FT^{\prime})\,.
\EEN
In particular, for every local section $t$ of $\HT$, the local sections
$D^{\pm} (t)$ of $\ENDO(\FT^{\prime})$ are adjoint to each other (see \RI, \S7, Definition 7.16, for more details), i.e.
$$
(D^{\pm} (t))^{\dagger} = D^{\mp} (t)\,.
$$
This leads to the following identities:
\BEN\label{adj1}
(d^{\pm})^{\dagger} = d^{\mp} \,\,and \,\,(p^{\pm})^{\perp,\dagger} = p^{\mp}\,. 
\EEN
\indent
The definition of the liftings $p^{\pm}_{\GA}$ is based on the triangular decomposition (\ref{d-D}). Hence, these morphisms are
of {\it algebraic} nature. The upshot of the subsequent discussion is to show that they are related in an explicit functorial way to the derivatives of
the period maps $p_{\GA}$ and ${}^{op}{p_{\GA}} $ along the directions of the fibres of the natural projection
$\pi: \JABG \longrightarrow \GAB$.
\subsection{The relative derivatives of $p_{\GA}$ and ${}^{op}{p_{\GA}} $}
Let $d_{\pi} ({}^{op}{p_{\GA}})$ (resp. $d_{\pi} (p_{\GA})$) be the relative differential of 
${}^{op}{p_{\GA}}$ (resp. $p_{\GA}$), i.e. 
$d_{\pi} ({}^{op}{p_{\GA}})$ (resp. $d_{\pi} (p_{\GA})$) is the restriction of the differential
$d ({}^{op}{p_{\GA}})$ (resp. $d (p_{\GA})$) to the relative tangent sheaf 
${\cal T}_{\pi}$ of the natural projection
$$
\pi : \JABG \longrightarrow \GAB\,.
$$
Since the morphisms ${}^{op}{p_{\GA}}$ and $p_{\GA}$ are morphisms of $\GAB$-schemes their relative differentials are morphisms
between the {\it relative tangent bundles} of 
$\JABG$ and $\FLAO$ and $\FLA$ respectively:
\BEN\label{rel-diff-bdl}
d_{\pi} ({}^{op}{p_{\GA}}): \TPI \longrightarrow \mbox{\BM${ T}$}_{\overleftarrow{Fl}_{\GA}},\,\,
 \,\,d_{\pi} (p_{\GA}) : \TPI \longrightarrow \mbox{\BM${ T}$}_{\overrightarrow{Fl}_{\GA}}\,.
\EEN
 Equivalently, on the level of sheaves one has
\begin{eqnarray}\label{rel-diff-sh}
d_{\pi} ({}^{op}{p_{\GA}}) : {\cal T}_{\pi}& \longrightarrow  ({}^{op}{p_{\GA}})^{\ast} \mbox{\BM${\cal T}$}_{\overleftarrow{Fl}_{\GA}}\,, \\ \label{rel-diff-sh1}
d_{\pi} (p_{\GA}):  {\cal T}_{\pi}& \longrightarrow  (p_{\GA})^{\ast} \mbox{\BM${\cal T}$}_{\overrightarrow{Fl}_{\GA}}\,.                                       \label{rel-diff-sh2}
\end{eqnarray}
\begin{pro}\label{pro-Grif-trans}
The relative differentials $d_{\pi} ({}^{op}{p_{\GA}})$ and $d_{\pi} (p_{\GA})$ satisfy Griffiths transversality condition
\begin{eqnarray}\label{Grif-trans}
d_{\pi} ({}^{op}{p_{\GA}}) : {\cal T}_{\pi}& \longrightarrow  
\bigoplus^{\LG-1}_{m=1} \HOM (\HT_{-m} /{ \HT_{-m+1}}, \HT_{-m-1} /{ \HT_{-m} })\,, \\ \nonumber
                                               &                                                                                                            \\ 
d_{\pi} (p_{\GA}):  {\cal T}_{\pi}& \longrightarrow  
\bigoplus^{\LG-1}_{m=1} \HOM (\PF^{m} /{ \PF^{m+1}}, \PF^{m-1} /{ \PF^{m} })\,.       \nonumber
\end{eqnarray}
Furthermore, let
$d_{\pi} ({}^{op}{p_{\GA}})_m$ (resp. $d_{\pi} (p_{\GA})_m$) be the $m$-th component of 
$d_{\pi} ({}^{op}{p_{\GA}})$ (resp. $d_{\pi} (p_{\GA})$). Then for any local section $v$ of $ {\cal T}_{\pi}$ and any local section
$h$ (resp. $\phi$) of $\HT_{-m}$  (resp. $\PF^m$) one has the following:
\begin{enumerate}
\item[a)]
$d_{\pi} ({}^{op}{p_{\GA}})_m (v) (h) \equiv -m  \tilde{v} \cdot h\, (mod\, \HT_{-m})$,
\item[b)]
$d_{\pi} (p_{\GA})_m (v) (\phi) \equiv  m \tilde{v} \cdot \phi \,(mod\, \PF^{m})$,
\end{enumerate}
where $\tilde{v}$ stands for an arbitrary lifting of $M^{-1} (v)$ to a local section of $\HT$
and $\tilde{v} \cdot h $ (resp. $\tilde{v} \cdot \phi$)
 stands for the product in $\FTP$ and it is independent of a lifting chosen after factoring out by $ \HT_{-m}$ (resp. $\PF^{m}$).
\end{pro}
\begin{pf}
It is enough to consider the situation fibrewise. For this fix a closed point $[Z] \in \GAB$ and let
$\JAB_{[Z]}$ be the fibre of 
$\pi: \JABG \longrightarrow \GAB$ over $[Z]$.
Recall that $\JAB_{[Z]}$ is a non-empty Zariski open subset of 
$\PP(\EZ)$. So, given a point $[\alpha] \in \JAB_{[Z]}$, the fibre 
$\TPI \ZA$
of $\TPI$ at $\ZA$ is
the tangent space $T_{\PP(\EZ), [\alpha]}$ of $\PP(\EZ)$ at $[\alpha]$.

Let $v$ be a tangent vector of $\PP(\EZ)$ at a point  $[\alpha] \in \JAB_{[Z]}$  and let
$\alpha(\epsilon)$ be an arc in  $\JAB_{[Z]}$ passing through $[\alpha] $ in the direction of $v$.
Using a canonical identification of the tangent space
$\TPI \ZA=T_{\PP(\EZ), [\alpha]}$ with
$\EZ /{\CC \alpha}$
 we may choose this arc so that it comes from the arc
$$
\tilde{\alpha}(\epsilon) = \alpha +\epsilon \beta
$$
in $\EZ$,where $\epsilon$ is in a small disk around $0$ in $\CC$ and $\beta \in \EZ$ is such that
$$
v\equiv \beta\,(mod\,\CC \alpha).
$$
Using the identifications in (\ref{HT=fr}) and (\ref{a/b1}) we have
\BEN\label{HT-ep}
\HT([Z],\alpha(\epsilon)) =\HT([Z],[\tilde{\alpha}(\epsilon)]) =\frac{\alpha}{\alpha +\epsilon \beta} \HT\ZA =
\frac{1}{1 +\epsilon t} \HT\ZA\,,
\EEN
where $\displaystyle{t =\frac{\beta}{\alpha} \in \HT\ZA}$. Furthermore, from the definition of  the isomorphism $M$ in (\ref{M}) it follows
that $t$ is a lifting of $M^{-1} (v) \in \HT\ZA/{\CC}\{1_Z\}$.

By definition of $\HT_{-m}([Z],\alpha(\epsilon))$ in (\ref{HT-i}) and (\ref{HT-ep}) we have
\BEN\label{HTm-ep}
\HT_{-m}([Z],\alpha(\epsilon)) =  \frac{1}{(1 +\epsilon t)^m} \HT_{-m}\ZA\,. 
\EEN
This implies that for every $h \in \HT_{-m}\ZA $ the expression
\BEN\label{h-ep}
h(\epsilon) =\frac{1}{(1 +\epsilon t)^m} h
\EEN
is a section of $\HT_{-m}$ over the arc $\alpha(\epsilon)$. Hence, by definition, the value
$d_{\pi} (\OPG)(v)$ of the relative differential of $\OPG$ at $\ZA$, along the vector $v$, is given by
the linear maps
$$
\HT_{-m}\ZA \longrightarrow \HO{Z^{\prime}}) /{\HT_{-m}\ZA},\,\,(m=1,\ldots, \LG -1)
$$
which sends $h \in \HT_{-m}\ZA $ to the vector
$$
\left.\frac{d}{d\epsilon} (h(\epsilon)) \right|_{\epsilon=0}\,(mod\,\HT_{-m}\ZA) \equiv - mth\,(mod\,\HT_{-m}\ZA)\,,
$$
where the last equivalence follows from the $\epsilon$-expansion of the right hand side in (\ref{h-ep}).
This proves the formula a) of the proposition and implies that
$d_{\pi} (\OPG)(v)$ restricted to
$\HT_{-m}\ZA$ takes its values in the subspace
$\HT_{-m-1}\ZA /{\HT_{-m}\ZA}$ of \linebreak
$\HO{Z^{\prime}}) /{\HT_{-m}\ZA}$.
Hence the Griffiths transversality condition for $d_{\pi} (\OPG)$.

Turning to the formula b) of the proposition we use the orthogonality relation
$$
\PF^m \ZA =\big{(} \HT_{-m}\ZA \big{)}^{\perp}\,.
$$
Let $\phi \in \PF^m \ZA$. To calculate 
$d_{\pi} (\PG) (\phi)$ we choose a section
$\phi(\epsilon)$ of $\PF^m$ over the arc $\alpha(\epsilon)$ with $\phi(0)=\phi$.
From (\ref{HTm-ep}) it follows
\BEN\label{q-phi-ep}
\QT(\phi(\epsilon), \frac{1}{(1 +\epsilon t)^m} h )=0,
\EEN
for all $h\in \HT_{-m}\ZA$, where $\QT$ stands for the bilinear symmetric form on $\FTP$ defined in (\ref{q}).
Taking the linear term of the $\epsilon$-expansion on the left hand side of (\ref{q-phi-ep}) we obtain
$$
\QT(\phi_1 -mt \phi, h) =0,
$$
for all $h\in \HT_{-m}\ZA$, where $\phi_1 = \left.\frac{d}{d\epsilon} (\phi(\epsilon)) \right|_{\epsilon=0}$. Hence
$$
\phi_1 -mt \phi \equiv 0\,(mod\, \PF^m \ZA)\,.
$$
This implies 
$$
d_{\pi} (\PG) (v) (\phi) \equiv \phi_1 \,(mod\, \PF^m \ZA) \equiv mt \phi \,(mod\, \PF^m \ZA) 
$$
as asserted by the formula b) of the proposition. Since $t \phi \in \PF^{m-1} \ZA$ the above relation also yields Griffiths transversality for
$d_{\pi} (\PG)$.
\end{pf}

In view of the Griffiths transversality condition in (\ref{Grif-trans}) it will be convenient to introduce graded sheaves
associated to $\FT^{\prime}$ relative to two filtrations $\HT_{-\bullet}$ and $\PF^{\bullet}$ defined in
(\ref {filtHT-JG}) and (\ref{filtFpr}) respectively. Thus we set 
\begin{eqnarray}
Gr^{\bullet}_{\HT_{-\bullet}} (\FT^{\prime}) &:=& \bigoplus^{\LG}_{m=1} Gr^m_{\HT_{-\bullet}} (\FT^{\prime})\,,  \\   \label{grH}
Gr^{\bullet}_{\PF^{\bullet}} (\FT^{\prime}) &:=& \bigoplus^{\LG-1}_{m=0} Gr^m_{\PF^{\bullet}} (\FT^{\prime})\,,     \label{grF}
\end{eqnarray}
where the graded pieces are defined as follows
\begin{eqnarray}
Gr^m_{\HT_{-\bullet}} (\FT^{\prime}) & =& \HT_{-m} /{\HT_{-m+1}}\,,  \\    \label{gr-pieceH}
Gr^m_{\PF^{\bullet}} (\FT^{\prime}) &=& \PF^m /{\PF^{m+1}}\,.   \label{gr-pieceF}
\end{eqnarray}

Using the translation functor in the category of graded sheaves we can rewrite the morphisms in (\ref{Grif-trans}) in that category
as follows.
\begin{eqnarray}
d_{\pi} ({}^{op}{p_{\GA}}) : {\cal T}_{\pi}& \longrightarrow  
 \HOM\big {(} Gr^{\bullet}_{\HT_{-\bullet}} (\FT^{\prime}), Gr^{\bullet}_{\HT_{-\bullet}} (\FT^{\prime}) [1] \big{)}\,, \\ \nonumber
                                               &                                                                                                            \\ \label{Grif-trans-gr}
d_{\pi} (p_{\GA}):  {\cal T}_{\pi}& \longrightarrow  
\HOM \big{(} Gr^{\bullet}_{\PF^{\bullet}} (\FT^{\prime}), Gr^{\bullet}_{\PF^{\bullet}} (\FT^{\prime}) [-1]  \big{)}\,,     \nonumber
\end{eqnarray}
where $\HOM$ is taken in the category of graded sheaves.\footnote{for graded modules we always assume that a graded component is zero, if its degree
is not in the range of the grading. Thus, for example, for the graded module $Gr^{\bullet}_{\HT_{-\bullet}} (\FT^{\prime}) [1]$, the component
$Gr^{\LG}_{\HT_{-\bullet}} (\FT^{\prime}) [1] = Gr^{\LG+1}_{\HT_{-\bullet}} (\FT^{\prime}) =0$.}\label{grading-conv}

We aim at relating these morphisms with the 
 morphisms
$d^{\pm}$ defined in (\ref{d+}) and (\ref{d-}) respectively.
Since the latter morphisms are defined algebraically, the virtue of such a relation will be purely algebraic
expressions for the relative differentials of our period maps.

Our first task will be to recast 
$d^{\pm}$ as morphisms of 
${\cal T}_{\pi}$ into the category of graded sheaves as well. The main point here is the orthogonal decomposition
in (\ref{ord-Fpr-til}) which makes $\FT^{\prime}$ itself a graded sheaf. To stress this graded structure of
$\FT^{\prime}$ we will write 
$\FT^{\prime \bullet}$. Thus
$\FT^{\prime \bullet}$ is the sheaf $\FT^{\prime}$ together with additional structure of the grading in (\ref{ord-Fpr-til}).
With this in mind we can rewrite the morphisms $d^{\pm}$ as follows
\BEN\label{d+-gr}
d^{\pm} : {\cal T}_{\pi} \longrightarrow \HOM \big{(} \FT^{\prime \bullet}, \FT^{\prime \bullet} [\pm 1] \big{)}\,,
\EEN
where $\HOM$ is taken in the category of graded sheaves.

The second observation is that the three graded sheaves $\FT^{\prime \bullet}$,
$Gr^{\bullet}_{\HT_{-\bullet}} (\FT^{\prime})$,
$Gr^{\bullet}_{\PF^{\bullet}} (\FT^{\prime})$ involved in our consideration are naturally related.
\begin{lem}\label{two-iso}
There are natural isomorphisms of graded sheaves
\begin{eqnarray*}
\phi:&  \FT^{\prime \bullet} & \longrightarrow Gr^{\bullet}_{\PF^{\bullet}} (\FT^{\prime})\,, \\
\psi:&   \FT^{\prime \bullet} & \longrightarrow Gr^{\bullet}_{\HT_{-\bullet}} (\FT^{\prime}) [1]\,.
\end{eqnarray*}
\end{lem}
\begin{pf}
From the orthogonal decomposition in (\ref{ordPF-p}) one obtains
\BEN\label{PFm-ord}
\PF^m =\PF^{m+1} \oplus \HH^m, \,\,\forall 0\leq m \leq \LG-1.
\EEN
This yields a canonical isomorphism
$$
\phi^m : \HH^m =\big{(} \FT^{\prime \bullet} \big{)}^m \longrightarrow  \PF^m /{\PF^{m+1}} =
 Gr^m_{\PF^{\bullet}} (\FT^{\prime})\,, 
$$
for every $m=0,\ldots, \LG-1$. Putting them together yields the isomorphism
$$
\phi:=\big{(}\oplus^{\LG-1}_{m=0} \phi^m \big{)}: \FT^{\prime \bullet}  \longrightarrow 
Gr^{\bullet}_{\PF^{\bullet}} (\FT^{\prime})
$$
of graded sheaves.

To construct $\psi$ we use the identities
\BEN\label{HTm-ord}
\HT_{-m} =\HT_{-m+1} \oplus \HH^{m-1},\,\,\forall 1\leq m \leq \LG,
\EEN
resulting from Remark \ref{rkH}, (\ref{ordH-i}). This yields a canonical isomorphism
$$
\psi^{m-1} : \HH^{m-1}  \longrightarrow  \HT_{-m} /{\HT_{-m+1}} = Gr^m_{\HT_{-\bullet}} (\FT^{\prime})=
\big{(} Gr^{\bullet}_{\HT_{-\bullet}} (\FT^{\prime}) [1] \big{)}^{m-1}\,,
$$
for  every $m=1,\ldots, \LG$. Putting them together yields the isomorphism
$$
\psi :=\big{(} \oplus^{\LG-1}_{m=0} \psi^m \big{)}: \FT^{\prime \bullet}  \longrightarrow 
Gr^{\bullet}_{\HT_{-\bullet}} (\FT^{\prime}) [1]
$$
of graded sheaves.
\end{pf}

Using the isomorphism $\psi$ (resp. $\phi$) we can relate the sheaves
$\HOM(  \FT^{\prime \bullet}, \FT^{\prime \bullet} [\pm 1])$,
the target of the morphisms $d^{\pm}$, with the targets of the morphisms $d_{\pi} ({}^{op}{p_{\GA}})$ and $d_{\pi} ({p_{\GA}})$, which are
the sheaves
$\HOM(Gr^{\bullet}_{\HT_{-\bullet}} (\FT^{\prime}) ,Gr^{\bullet}_{\HT_{-\bullet}} (\FT^{\prime})[1])$
 and
$\HOM(Gr^{\bullet}_{\PF^{\bullet}} (\FT^{\prime}) , Gr^{\bullet}_{\PF^{\bullet}} (\FT^{\prime}) [-1])$,
respectively.
\begin{lem}\label{two-iso1}
The isomorphism $\psi$ (resp. $\phi$) induces the isomorphism
\begin{eqnarray}
hom(\psi): \HOM( \FT^{\prime \bullet},  \FT^{\prime \bullet} [ 1])&
\longrightarrow
\HOM(Gr^{\bullet}_{\HT_{-\bullet}} (\FT^{\prime}) ,Gr^{\bullet}_{\HT_{-\bullet}} (\FT^{\prime})[1])\,. \\ \label{hom-psi}
\mbox{( resp.}\,\,
hom(\phi) :
 \HOM(  \FT^{\prime \bullet}, \FT^{\prime \bullet} [-1] )&
\longrightarrow
\HOM(Gr^{\bullet}_{\PF^{\bullet}} (\FT^{\prime}) , Gr^{\bullet}_{\PF^{\bullet}} (\FT^{\prime}) [-1])\,.  )    \label{hom-phi}
\end{eqnarray}
\end{lem}
\begin{pf}
To define $hom(\psi)$ take a local section $a=(a^p)$ of 
$\HOM( \FT^{\prime \bullet},  \FT^{\prime \bullet} [ 1])$ and define the local section
$hom(\psi)(a)$
of 
$\HOM(Gr^{\bullet}_{\HT_{-\bullet}} (\FT^{\prime}) ,Gr^{\bullet}_{\HT_{-\bullet}} (\FT^{\prime})[1])$
by requiring the diagram
\BEN\label{a-hom-a}
\xymatrix{
{\HH^{m-1}} \ar[r]^{a^{m-1}} \ar[d]_{\psi^{m-1}} & \HH^m \ar[d]^{\psi^m}\\
Gr^{m}_{\HT_{-\bullet}} (\FT^{\prime}) \ar[r] &Gr^{m+1}_{\HT_{-\bullet}} (\FT^{\prime}) }
\EEN
to commute for all $m\geq 1$, i.e. the bottom horizontal arrow is $(hom(\psi)(a))^m$ given by the formula
$$
(hom(\psi)(a))^m =\psi^m \circ a^{m-1} \circ (\psi^{m-1})^{-1}\,.
$$
This clearly gives an isomorphism. 

The definition of $hom(\phi)$ is completely analogous.
\end{pf}

It turns out that the straightforward relationship
$d_{\pi} ({}^{op}{p_{\GA}}) = hom(\psi) \circ d^{+}$
(resp. $d_{\pi}(p_{\GA}) = hom(\phi) \circ d^{-}$) is incorrect.
The following version of the formulas a) and b) of Proposition \ref{pro-Grif-trans} gives
the correct relations between the relative derivatives of our period maps and their algebraic counterparts - the morphisms
$d^{\pm}$.
\begin{lem}\label{der=d+-}
Let $\ZA$ be a closed point in $\JABG$ and let $v$ be a tangent  vector in the fibre
$({\cal T}_{\pi} )_{\ZA}$ of the relative tangent sheaf ${\cal T}_{\pi}$ at $\ZA$.
Denote by 
$d_{\pi}({}^{op} p_{\GA}) (v)$ (resp. $d_{\pi}( p_{\GA}) (v)$)
the relative differential of ${}^{op} p_{\GA}$ (resp. $p_{\GA}$) at $\ZA$
evaluated on $v$. Let
$d^{\pm} (v)$ be the evaluation of the morphisms $d^{\pm}$ on $v$. Then the following holds.
\begin{enumerate}
\item[a)]
For every $m \in \{1,\ldots, \LG \}$ the diagram
$$
\xymatrix{
{\HH^{m-1} \ZA} \ar[r]^{\scriptstyle{-m d^{+}(v)}} \ar[d]_{\psi^{m-1}} & {\HH^m \ZA} \ar[d]^{\psi^m} \\
{\scriptscriptstyle{Gr^m_{\HT_{-\bullet}} \!(\FT^{\prime}) \ZA }}\ar[r]_{\scriptscriptstyle{
d_{\pi}({}^{\scriptscriptstyle{op}}\! p_{\scriptscriptstyle{\GA}}) (v)}} &
 {\scriptscriptstyle{Gr^{m+1}_{\HT_{-\bullet}} (\FT^{\prime}) \ZA}}  }
$$
commutes.
\item[b)]
For every $m \in \{0,\ldots, \LG -1 \}$ the diagram
$$
\xymatrix{
{\HH^{m} \ZA} \ar[r]^{\scriptstyle{m d^{-}(v)}} \ar[d]_{\phi^{m}} & {\HH^{m-1} \ZA} \ar[d]^{\phi^{m-1}} \\
{\scriptscriptstyle{Gr^m_{\PF^{\bullet}} (\FT^{\prime})\ZA}} \ar[r]^{\scriptscriptstyle{d_{\pi}( p_{\GA}) (v)}}   &
{\scriptscriptstyle{ Gr^{m-1}_{\PF^{\bullet}} (\FT^{\prime})\ZA}} }
$$
commutes.
\end{enumerate}
(The vertical arrows in the above diagrams are the components of the isomorphisms
$\phi$ and $\psi$ in Lemma \ref{two-iso}; 
all graded components are assumed to be zero, whenever their degree is not in the range of a grading, see the footnote on page \pageref{grading-conv}).
\end{lem}
\begin{pf}
Let $t$ be an element of $\HT\ZA$ lifting $M^{-1} (v) \in \HT\ZA /{\CC}\{1_Z \}$, where
$M$ is the isomorphism in (\ref{M}) and $1_Z$ is the constant function of value $1$ on $Z$. Then the formula a) in Proposition \ref{pro-Grif-trans} reads as follows
$$
d_{\pi} (\OPG)(v)(h) \equiv -mth\,(mod\, \HT_{-m} \ZA)\,,
$$
for any $h \in \HH^{m-1} \ZA \subset \HT_{-m} \ZA$, where the inclusion comes from (\ref{HTm-ord}).

Using the triangular decomposition (\ref{d-Dt}) of the operator $D(t)$ of  multiplication by $t$, we write $th$ in the above formula as follows
\BEN\label{th-tri}
th =D^{-}(t)h + D^{0}(t)h + D^{+}(t)h.
\EEN
Furthermore, the terms
$D^{-}(t)h$ and $ D^{0}(t)h$ are in $\HT_{-m} \ZA$, for $h\in \HT_{-m} \ZA$. Combining this together with the definition of the isomorphism
$\psi$ in the proof of Lemma \ref{two-iso}, we obtain the following equality
$$
 d_{\pi} (\OPG)(v)(\psi^{m-1}(h)) = -m \psi^m (D^{+}(t)h) = -m \psi^m (d^{+}(v)h)\,,
$$
where the last equality follows from Remark \ref{val-d}, (\ref{d=D}).
Hence the commutativity of the first diagram.

For the second diagram of the lemma, we use  the formula b) in Proposition \ref{pro-Grif-trans} to obtain
$$
d_{\pi} (\PG)(v)(h) \equiv mth \,(mod \PF^m \ZA)\,,
$$
for any $h \in \HH^m \ZA \subset \PF^m \ZA$, where the inclusion comes from (\ref{PFm-ord}). Using the decomposition for $th$ in (\ref{th-tri}) once
again and observing that
$D^{+}(t)h$ and $ D^{0}(t)h$ are in $\PF^{m} \ZA$, for all $h\in \HH^{m} \ZA$, we deduce
$$
d_{\pi} (\PG)(v)(h) \equiv mD^{-}(t)h \,(mod\, \PF^m \ZA)\,.
$$
By  the definition of the isomorphism $\phi$ in Lemma \ref{two-iso} this yields 
$$
d_{\pi} (\PG)(v)(\phi^m(h)) = m\phi^{m-1}(D^{-}(t)h) = m\phi^{m-1}(d^{-}(v)h)\,,
$$
for all $h \in \HH^m \ZA$, and where the last equality is again Remark \ref{val-d}, (\ref{d=D}). Hence the commutativity of the second diagram of the lemma.
\end{pf}

 From the diagrams in Lemma \ref{der=d+-} we see that 
$d_{\pi}({}^{op} p_{\GA})$ (resp. $d_{\pi}( p_{\GA})$) ,  up to the canonical  identification provided by
$\psi$ (resp. $\phi$), coincides not with 
$d^{+}$ (resp. $d^{-}$) as could be naively expected,  but with its {\it scaled} version, 
where the scaling keeps track of the degree of the grading of sheaves. The following general result gives a functorial
procedure for such scalings in the category of graded sheaves.
\begin{lem}\label{lem-sc}
Let ${\cal M}^{\bullet} =\bigoplus_{p\in {\bf Z}} {\cal M}^p$
be a ${\bf Z}$-graded sheaf of modules and let
$c_{{\cal M}^{\bullet}} \in Hom({\cal M}^{\bullet}, {\cal M}^{\bullet})$
be the ``counting" endomorphism, i.e. the components of $c_{{\cal M}^{\bullet}}$
are given by the formula
\BEN\label{morph-c}
c^m_{{\cal M}^{\bullet}} = \frac{m(m+1)}{2} id_{{\cal M}^m}\,.
\EEN
Given a morphism of sheaves
$$
a: {\cal A} \longrightarrow \HOM({\cal M}^{\bullet}, {\cal M}^{\bullet} [n])\,,
$$
for some integer $n$  and with $\HOM$ functor taken in the category of graded sheaves, one sets
\BEN\label{sc}
c(a):= -ad(c) \circ a\,.
\EEN
Then the resulting morphism
$$
c(a): {\cal A} \longrightarrow \HOM({\cal M}^{\bullet}, {\cal M}^{\bullet} [n])
$$
has the components given by the following formula
\BEN\label{sc1}
c(a)^m = -(nm + \frac{n(n+1)}{2}) a^m, \,\,\forall m \in {\bf Z}\,.
\EEN
\end{lem}
\begin{pf}
This is an elementary straightforward calculation.
\end{pf}

Applying the above formalism to $d^{\pm}$ we obtain the morphisms
\BEN\label{d-sc}
c(d^{\pm}): {\cal T}_{\pi} \longrightarrow \HOM(\FT^{\prime \bullet}, \FT^{\prime \bullet} [\pm 1] )
\EEN
subject to the scaling of the components of $d^{\pm}$ in the diagrams of Lemma \ref{der=d+-}.
Thus we obtain the following ``algebraic" formulas for the relative differentials of the period maps
$p_{\GA}$ and ${}^{op}p_{\GA}$. 
\begin{pro}\label{pro-der=c(d)}
The relative differentials
$d_{\pi}({}^{op} p_{\GA})$ and $d_{\pi}( p_{\GA})$) are related to the morphisms
$d^{\pm}$ as follows:
\begin{eqnarray}\label{der=c(d)}
d_{\pi}({}^{op} p_{\GA})&=& hom(\psi) \circ c(d^{+})\,,\\  \nonumber
& &\\       
d_{\pi}( p_{\GA})&=& hom(\phi) \circ c(d^{-})\,,   \nonumber
\end{eqnarray}
where $hom(\psi)$ and $hom(\phi)$ are the isomorphisms from Lemma \ref{two-iso1}.
Equivalently, the following diagrams commute
\begin{enumerate}
\item[a)]
$$
\xymatrix{
   & {\scriptstyle{\HOM(\FT^{\prime \bullet}, \FT^{\prime \bullet} [ 1])}} \ar[dd]^{hom(\psi)} \\
{\cal T}_{\pi} \ar[ru]^{c(d^{+})} \ar[rd]^{\scriptstyle{d_{\pi}({}^{op} p_{\GA})}}& \\
  &{\scriptstyle{{\HOM(Gr^{\bullet}_{\HT_{-\bullet}} (\FT^{\prime}) ,Gr^{\bullet}_{\HT_{-\bullet}} (\FT^{\prime})[1])}}}  }
$$
\item[b)]
$$
\xymatrix{
   & {\scriptstyle{{\HOM(\FT^{\prime \bullet}, \FT^{\prime \bullet} [ -1])}}} \ar[dd]^{hom(\phi)} \\
{\cal T}_{\pi} \ar[ru]^{c(d^{-})} \ar[rd]^{\scriptstyle{d_{\pi}( p_{\GA})}}& \\
  &{\scriptstyle{{\HOM(Gr^{\bullet}_{\PF^{\bullet}} (\FT^{\prime}) ,Gr^{\bullet}_{\PF^{\bullet}} (\FT^{\prime})[-1])}}}  }
$$
\end{enumerate}
\end{pro}
\begin{rem}\label{ci=}
The identities in (\ref{der=c(d)}) establish the equality of relative differentials of
$p_{\GA}$ and ${}^{op} p_{\GA}$ with the morphisms $d^{\pm}$, up to the canonical identifications determined by
canonical isomorphisms $\phi$ and $\psi$ in Lemma \ref{two-iso}.  We agree on these canonical identifications (ci) and write
the identities in (\ref{der=c(d)}) as follows:
\begin{eqnarray}\label{der=c(d)1}
d_{\pi}({}^{op} p_{\GA})&\stackrel{ci}{=}&  c(d^{+})\,,\\  \nonumber
& &\\       
d_{\pi}( p_{\GA})&\stackrel{ci}{=}&  c(d^{-})\,.   \nonumber
\end{eqnarray}
\end{rem}

Next we return to the liftings
$$
p^{+}_{\GA} : T_{\pi} \longrightarrow \mbox{\BM${T}$}^{\ast}_{\overrightarrow{Fl}_{\GA}} \,\,\,and \,\,\,
p^{-}_{\GA} : T_{\pi} \longrightarrow \mbox{\BM${T}$}^{\ast}_{\overleftarrow{Fl}_{\GA}}
$$
of $p_{\GA}$ and ${}^{op} p_{\GA}$ introduced in (\ref{p+}) and (\ref{p-}) respectively.
Since they are determined in terms of the morphisms $d^{\pm}$ we should be able to relate them to
the relative differentials as well. Indeed, first observe that we can define the morphisms
\BEN\label{c-p+-}
c(p^{+}_{\GA}) : \TPI \longrightarrow  \mbox{\BM${T}$}^{\ast}_{\overrightarrow{Fl}_{\GA}} \,\,\,and \,\,\,
c(p^{-}_{\GA}) : T_{\pi} \longrightarrow \mbox{\BM${T}$}^{\ast}_{\overleftarrow{Fl}_{\GA}}
\EEN
simply by replacing $d^{\pm}$ used in formulas (\ref{p+formula}), (\ref{p-formula}) by
$c(d^{\pm})$. Following this by the operation of taking adjoint
$(\cdot)^{\dagger}$ (see (\ref{adj})) and using the notation introduced in Remark \ref{ci=}
we deduce the following {\it vector bundle} version of the identities in (\ref{der=c(d)1}).
\begin{cor}\label{der=c(p+-)}
\begin{eqnarray*}
d_{\pi}({}^{op} p_{\GA})&\stackrel{ci}{=}&  - c(p^{-})^{\dagger}\,,\\ 
& &\\       
d_{\pi}( p_{\GA})&\stackrel{ci}{=}& - c(p^{+})^{\dagger}\,.   
\end{eqnarray*}
\end{cor}
\begin{pf}
The proofs of two identities are very similar, so will do only the first one.

Let $([Z],[\alpha],v)$ be a point in $\TPI$. By definition the value of 
$c(p^{-})$ at $([Z],[\alpha],v)$ is as follows
\BEN\label{c(p-)Zav}
c(p^{-})([Z],[\alpha],v) = ( \OPG \ZA, c(d^{-})(v))\,.
\EEN
Applying the operation of taking the adjoint yields
\BEN\label{adj-c(p-)Zav}
c(p^{-})^{\dagger}([Z],[\alpha],v) = ( \OPG \ZA, c(d^{-})^{\dagger}(v))\,.
\EEN
So it remains to calculate $c(d^{-})^{\dagger}(v)$. From (\ref{sc1}) it follows
\BEN\label{adj-c(p-)Zav1}
c(d^{-})(v) =\sum^{\LG-1}_{p=1} pd^{-}_p (v)\,,
\EEN
where $d^{-}_p (v)$ is the $p$-th component of $d^{-} (v)$, i.e. $d^{-}_p (v)$ is the restriction of $d^{-} (v)$
to the summand $\HH^p \ZA$. Applying $(\cdot)^{\dagger}$ yields
$$
c(d^{-})^{\dagger}(v) =\sum^{\LG-1}_{p=1} p(d^{-}_p)^{\dagger} (v) =\sum^{\LG-1}_{p=1} pd^{+}_{p-1} (v) = 
\sum^{\LG-2}_{p=0} (p+1)d^{+}_{p} (v) =-c(d^{+})(v)\stackrel{ci}{=} -d_{\pi} (\OPG) (v)\,,
$$
where the last equality follows from the first identity in (\ref{der=c(d)1}).
\end{pf}

\subsection{Torelli problems for maps $ p_{\GA}$ and ${}^{op} p_{\GA}$}\label{sec-Torelli}

Once the maps $\PG$ and $\OPG$ are in place, it is natural to consider Torelli-type properties for them.
\begin{defi}\label{Torelli-pr}
\begin{enumerate}
\item[1)]
We say that Torelli property holds for a component $\GA$ of $\CSA$, if the morphism
$\PG$ or, equivalently, $\OPG$ is an embedding.
\item[2)]
We say that Infinitesimal Torelli property holds for $\GA$ if the differential of 
$\PG$ or, equivalently, the differential of $\OPG$ is injective at every point of $\JABG$.
\end{enumerate}
\end{defi}

The work done in the previous section allows to reformulate the Infinitesimal Torelli property as purely algebraic property of injectivity
of the morphisms $d^{\pm}$ defined in (\ref{d+}) and (\ref{d-}) respectively.

\begin{pro}\label{InfTor-equiv}
Let $\GA$ be a component in $\CSA$. Then
the following properties are equivalent:
\begin{enumerate}
\item[1)]
the  Infinitesimal Torelli property holds for $\GA$,
\item[2)]
the morphism
$$
d^{+} : {\cal T}_{\pi} \longrightarrow \HOM(\FTB, \FTB[1])
$$
is injective,
\item[3)]
the morphism
$$
d^{-}: {\cal T}_{\pi} \longrightarrow \HOM(\FTB, \FTB[-1])
$$
is injective.
\end{enumerate}
\end{pro}
\begin{pf}
By definition $\PG$ (resp. $\OPG$) (see (\ref{p}) and (\ref{op-p})) is a morphism of $\GAB$-schemes. This implies that
the differential of $\PG$ (resp. $\OPG$) is injective if and only if the relative differential 
$d_{\pi} (\PG)$
 (resp. $d_{\pi} (\OPG)$) is injective. From Proposition \ref{pro-der=c(d)} it follows that the injectivity of
$d_{\pi} (\PG)$ (resp. $d_{\pi} (\OPG)$) is equivalent to the injectivity of $d^{-}$ (resp. $d^{+}$).
This gives the equivalence between 1) and 2) (resp. 1) and 3)). The equivalence between 2) and 3) is assured by 
the fact that $d^{+}$ and $d^{-}$ are adjoint to each other (see (\ref{adj1})).
\end{pf}

The proposition above establishes a link between the Infinitesimal Torelli property for an admissible component $\GA$ in $\CSA$ and
the sheaf of Lie algebras 
{\BM
$\LAGT$}.
 Our next result shows that the kernel of $d^{\pm}$ or, equivalently, the failure of the Infinitesimal Torelli property is controlled by
the center
{\BM
 $\CG$ of $\LAGT$.
More precisely, recall that the sheaf of Cartan subalgebras
${\cal C}(\HT)$ of $\LAGT$
}
 (see  Proposition \ref{Cartan}) is naturally identified with
$\FT^{\prime}$ (Claim \ref{Cartan=f}). Furthermore, the center, via this identification, is isomorphic to the subsheaf
$\pi^{\ast} \FF^{\prime}_c$ of $\HT$ constructed in Proposition \ref{sh-c}.
On the other hand $\HT$ is related to the relative tangent sheaf ${\cal T}_{\pi}$ by the isomorphism
$$
M:  \HT / {\OO_{\JABG}} \longrightarrow {\cal T}_{\pi} 
$$
which was recalled in (\ref{M}).
\begin{pro}\label{ker-der=center}
One has equalities
$$
ker(d^{+}) = ker(d^{-}) =ker(d_{\pi} (\PG)) = ker (d_{\pi} (\OPG))
$$
 and isomorphisms
\BEN\label{iso-ker-center}
\mbox{\BM$\CG$} /{\OO_{\JABG}} \cong \pi^{\ast} \FT^{\prime}_c /{\OO_{\JABG}} \cong ker(d_{\pi} (\PG))\,.
\EEN
\end{pro}
\begin{pf}
The first assertion is a restatement of Proposition \ref{InfTor-equiv}. For the second assertion observe that the first isomorphism
is given by the morphism $D$ (see Proposition \ref{sh-c}). To see the second isomorphism in (\ref{iso-ker-center})
consider a local section
$v$ of $ker(d^{\pm})$. Choose a lifting $\tilde{v}$ of $M^{-1} (v)$ to a local section of $\HT$. Then we have
$$
D^{\pm}(\tilde{v}) =0
$$
in the triangular decomposition (\ref{d-Dt}) and hence $D(\tilde{v})= D^0 (\tilde{v})$. This implies
\BEN\label{com-D0-v}
[ D^0 (\tilde{v}), D(t)]=[ D (\tilde{v}), D(t)]= 0,
\EEN
for every local section $t$ of $\HT$, and where the last equality above is the commutativity of the multiplication in
$\FT^{\prime}$.

Substituting $D(t) = D^{-} (t) + D^{0} (t) + D^{+} (t)$ into (\ref{com-D0-v}) and decomposing according to the grading of
{\BM
$\LAGT$}
in (\ref{Lie-grad}) yield
$$
[ D (\tilde{v}), D^{\pm}(t)] = [ D (\tilde{v}), D^{0}(t)] = 0\,,
$$
for every local section $t$ of $\HT$. Since 
{\BM
$\LAGT$}
is generated (as a Lie algebra) by $D^{\pm}(t),D^0 (t)$,
we deduce that $ D (\tilde{v})$ is a local section of the center
{\BM
$\CG$}. This gives an inclusion
\BEN\label{inc1}
M^{-1} (ker(d^{\pm})) \subset \pi^{\ast}(\FT^{\prime}_c )/{\OO_{\JABG}}\,.
\EEN
On the other hand, from Proposition \ref{pro-c-basis}, 2) and 4), it follows that local sections of the center
{\BM
$\CG$}
preserve the orthogonal decomposition $\FT^{\prime}$ in (\ref{ord-Fpr-til}). Hence
$D^{\pm} (t) =0$, for every local section  $t$ of $\pi^{\ast}(\FT^{\prime}_c )$. This gives an inclusion
$$
\pi^{\ast}(\FT^{\prime}_c ) /{\OO_{\JABG}} \subset M^{-1} (ker(d^{\pm}))\,.
$$
Combining this with (\ref{inc1}) yields the equality
$$
\pi^{\ast}(\FT^{\prime}_c ) /{\OO_{\JABG}} = M^{-1} (ker(d^{\pm}))\,.
$$
This completes the proof of the second assertion.
\end{pf}
\begin{cor}\label{InfTorelli=s}
The Infinitesimal Torelli property holds for a component $\GA$ in $\CSA$ if and only if
$\GA$ is simple (Definition \ref{s-c}).
\end{cor}
\begin{pf}
Follows immediately from
Definition \ref{s-c} and Proposition \ref{ker-der=center}.
\end{pf}

One of the features of our period map(s) is that the Torelli property turns out to be equivalent to the Infinitesimal Torelli property.
\begin{thm}\label{Torelli=InfTor}
Let $\GA$ be a component in $\CSA$.
Then the following statements are equivalent:
\begin{enumerate}
\item[1)]
the Torelli property holds for $\GA$,
\item[2)]
the Infinitesimal Torelli property holds for $\GA$,
\item[3)]
$\GA$ is simple.
\end{enumerate}
\end{thm}
\begin{pf}
By Corollary \ref{InfTorelli=s} the statements 2) and 3) are equivalent. 
We consider the equivalence 1) and 2).
It is also clear that 1) implies 2). So we turn to the implication in the other direction.

We will work with the morphism $\OPG$. Since it is a $\GAB$-morphism, the Torelli property holds if and only if
 $\OPG$ is an embedding on each fibre of $\pi: \JABG \longrightarrow \GAB$. 

 Fix
$[Z] \in \GAB$ and consider
$\JAB_Z =\pi^{-1} ([Z])$, the fibre of $\JABG$ over $[Z]$. This gives the morphism
$$
\OPG(Z) : \JAB_Z  \longrightarrow  \FLAO ([Z])
$$
which is the restriction of $\OPG$ to $\JAB_Z $ and where $\FLAO ([Z])$ is the fibre $\FLAO$ over $[Z]$.

Assume 2) holds. Then $\OPG(Z)$ is an immersion. So to see that it is an embedding it is enough to show that it is injective on closed points of 
$\JAB_Z $. We argue by contradiction.

Let $[\alpha]$ and $ [\beta]$ be two distinct points of
$\JAB_Z $ and assume $\OPG(Z,[\alpha]) =\OPG(Z,[\beta])$.
By definition of $\OPG$ in (\ref{op-p-Z}) the two filtrations
$\HT_{-\bullet} \ZA$ and $\HT_{-\bullet} ([Z],[\beta])$ coincide. This implies in particular
that
$ \HT ([Z],[\beta]) = \HT \ZA $.
Combining this with (\ref{a/b1}) we obtain
\BEN\label{a=b}
\frac{\alpha}{\beta} \HT \ZA  =  \HT \ZA\,.  
\EEN
This implies that $t = \frac{\alpha}{\beta}$ belongs to $\HT\ZA$ (since $1_Z \in \HT\ZA$) and the multiplication by $t$
preserves $\HT \ZA  $. But then the multiplication by $t$ preserves the whole filtration
$\HT_{-\bullet} \ZA$, as follows from the definitions in (\ref{HT-i}). Hence $D^{+} (t) =0$. By assumption  $t$ is not a constant function
on $Z$, so its projection $\overline{t}$ in 
$ \HT\ZA / {\CC\{1_Z\}}$ is non-zero. This gives the non-zero tangent vector
$v =M(\overline{t})$ in $\TPI$ at $\ZA$ such that
$d^{+} (v) =D^{+}(t) =0$, where the first equality is the definition of $d^{+} (v) $ (see Remark \ref{val-d}).
By Proposition \ref{InfTor-equiv} this is equivalent to the failure of the Infinitesimal Torelli property. But this contradicts the assumption
that 2) holds.
\end{pf}

Once we know that the Torelli property (Definition \ref{Torelli-pr}) is equivalent to the injectivity of
morphisms
\BEN\label{d+-1}
d^{\pm} : {\cal T}_{\pi} \longrightarrow \HOM (\FTB, \FTB[\pm1]) = \bigoplus^{\LG-1}_{p=0} \HOM(\HH^p, \HH^{p \pm 1})\,,
\EEN
we can ask the same question for its components
\BEN\label{d+-p}
d^{\pm}_p : {\cal T}_{\pi} \longrightarrow \HOM (\HH^p, \HH^{p\pm 1})\,,
\EEN
where $p=0, \ldots, \LG-1$, with the understanding that $d^{-}_0 =d^{+}_{\LG-1} =0$.
Observe that the operation of taking adjoint in (\ref{adj}) interchanges $d^{+}_p$ with
$d^{-}_{p+1}$, for all $p\in \{0,\ldots,\LG-2\}$. So it is enough to consider the the morphisms of the same sign, say the $d^{+}_p$'s.

From the definition of the filtration
$\HH_{-\bullet}$ it follows easily that the conditions of failure of 
$d^{+}_p$ to be injective for various values of $p$ are not independent. Namely, one has the inclusion
\BEN\label{Tp-Tp+1}
ker(d^{+}_p) \subset ker(d^{+}_{p+1})\,,
\EEN
for all $p \in \{0,\ldots,\LG-2\}$ (see \S\ref{sec-strat}, Lemma \ref{lem-T(p)}, for a proof).
This suggests the following terminology.
\begin{defi}\label{def-TI}
For a component $\GA$ in $\CSA$ set
$$
\tau_{\GA} = \min \left\{ \left. p\in \{0,\ldots,\LG-2 \} \right|  ker(d^{+}_p) \neq 0 \right\}
$$
and call it Torelli index of $ \GA$.
\end{defi}
\begin{rem}\label{TP-TI}
In view of the inclusions in (\ref{Tp-Tp+1}) the two notions - Torelli index and Torelli property - can become different only if the
length $\LG$ of the filtration $\HT_{-\bullet}$ is $\geq 3$.
\end{rem}

Using the terminology of Definition \ref{def-TI} and the results obtained so far one deduces the following numerical criterion.
\begin{cor}\label{ti=0}
A component $\GA$ is simple if and only if
the Torelli index $\tau_{\GA} >0$.
\end{cor}

The non-vanishing of $ker(d^{+}_{\tau_{\GA} })$, for the values of $\tau_{\GA}$ in the range\footnote{we assume $\LG \geq 3$, see Remark \ref{TP-TI}.}
$[1,\LG-2]$, leads to the geometric properties of configurations similar to the ones
given in Corollary \ref{Z-c-dec} (see the proof of Proposition \ref{gp-stTor}) and hence to a certain hierarchy among simple components. We will not pursue this discussion here
except of distinguishing components which are on the top of this hierarchy.
\begin{defi}\label{st-T}
\begin{enumerate}
\item[1)]
A component $\GA$ in $\CSA$ is said to have strong Torelli property if
its Torelli index $\tau_{\GA} = \LG-2$, i.e. the morphisms $d^{+}_p$ are injective for all
$p\in \{0, \ldots, \LG-2 \}$.
\item[2)]
A configuration $Z$ on $X$ is said to have strong Torelli property if there is a component
$\GA $ containing $[Z]$ and having strong Torelli property.
\end{enumerate}
\end{defi}

It turns out that the configurations in general position with respect to the adjoint linear system
$\left| K_X +L\right|$
satisfy strong Torelli property. Namely, the following holds.
\begin{pro}\label{gp-stTor}
Let $Z$ be a configuration of $d$ points on $X$ 
 such that the index of $L$-speciality
$\delta(L,Z) =r+1 \geq 2$ and $d \geq r+2$. Assume $Z$ to be in general position with respect to the adjoint linear system
$\left| K_X + L \right|$.
Then $Z$ satisfies strong Torelli property.
\end{pro}
\begin{pf}
Let $\GA$ be an admissible component in $\CS$ containing $[Z]$ and let $\JAB_{Z}$ be the fibre of $\JABG$ over $[Z]$. Fix 
$[\alpha] \in \JAB_{Z}$ and consider the morphism
$$
\kappa \ZA : Z \longrightarrow \PP(\HT\ZA^{\ast})
$$
as in (\ref{kappa}). From the proof of Corollary 7.13, \RI, it follows that
$\kappa \ZA$ is an embedding and its image, which we continue to denote by $Z$, is the set of $d$ distinct points in
general position in $\PP(\HT\ZA^{\ast})$. Thus $\HT_{-\LG} = \HO Z)$ and we consider its orthogonal decomposition
\BEN\label{ord-ZA}
\HO Z) = \bigoplus^{\LG-1}_{p=0} \HH^p \ZA\,.
\EEN
This comes along with the linear map
$$
d^{+} \ZA : \TPI \ZA \longrightarrow \bigoplus^{\LG-2}_{p=0} Hom (\HH^p \ZA , \HH^{p+1} \ZA)
$$
which is the morphism $d^{+}$ in (\ref{d+-1}) at $\ZA$. Our objective is to show that the $p$-th component
$d^{+}_p \ZA $ of $d^{+} \ZA $
$$
d^{+}_p \ZA : \TPI \ZA \longrightarrow  Hom (\HH^p \ZA , \HH^{p+1} \ZA)
$$
is injective, for all $p\in \{0,\ldots,\LG-2 \}$.

From Theorem \ref{Torelli=InfTor}, Corollary \ref{s-gp} and the inclusions (\ref{Tp-Tp+1}) it follows that
the component $d^{+}_0 \ZA$ is injective. So we may assume that $\LG \geq 3$ and let $p_0$ be the smallest index for which
$d^{+}_{p_0} \ZA$ fails to be injective. Thus from what is said above $p_0 \in [1,\LG-2]$ and we set
\BEN\label{Tp0}
T^{(p_0)} \ZA = ker( d^{+}_{p_0} \ZA )\,.
\EEN
Using the isomorphism $M$ in (\ref{M}) and Remark \ref{val-d}, we obtain the subspace 
$\tilde{T}^{(p_0)} \ZA$ of $\HT\ZA$, composed of elements $t$ such that
$
D^{+}_{p_0} (t) =0
$,
where $D^{+}_{p} (t)$ is the restriction to $\HH^p \ZA$ of the operator $D^{+} (t)$ in (\ref{d-Dt}).
This implies that the multiplication by $t$, for all $t \in \tilde{T}^{(p_0)} \ZA$, preserves the subspace $\HT_{-p_0 -1} \ZA$ of the filtration 
$\HT_{-\bullet}$ in (\ref{filtHT-JG}) at $\ZA$ as well as
the summands $\HH^p \ZA$ in (\ref{ord-ZA}), for all $p\geq p_0 +1$. From this point on our considerations are similar to the ones in the study
of the center of 
{\BM
$\LAGT$} in \S\ref{Center}. Namely, we consider the weight decomposition of $\HO Z)$ under the action (by multiplication) of 
$\tilde{T}^{(p_0)} \ZA$:
\BEN\label{wd-act-Tp0}
\HO Z) = \bigoplus_{\lambda \in (\tilde{T}^{(p_0)} \ZA )^{\ast}} V^{(p_0)}_{\lambda}\,,
\EEN
where $V^{(p_0)}_{\lambda}$ is the weight space of $\tilde{T}^{(p_0)} \ZA $-action corresponding to a weight $\lambda$.
As in Proposition \ref{pro-wd} one shows:
\begin{enumerate}
\item[1)]
$V^{(p_0)}_{\lambda}$ is an ideal of $\HO Z)$, for every weight $\lambda$ occurring in (\ref{wd-act-Tp0}),
\item[2)]
$ V^{(p_0)}_{\lambda} \cdot V^{(p_0)}_{\mu} = 0$, for all $\lambda \neq \mu$.
\end{enumerate}
Furthermore, since the action of $\tilde{T}^{(p_0)} \ZA$ preserves $\HT_{-p_0 -1} \ZA$ and the summands $\HH^p \ZA$, for $p\geq p_0 +1$,
each weight space $V^{(p_0)}_{\lambda}$ admits the following orthogonal decomposition
\BEN\label{wsp-ord}
V^{(p_0)}_{\lambda} = V^{\leq p_0}_{\lambda} \oplus \left(\bigoplus_{p > p_0} V^{p}_{\lambda} \right)\,,
\EEN
where $V^{\leq p_0}_{\lambda} = V^{(p_0)}_{\lambda} \bigcap \HT_{-p_0 -1} \ZA $ and 
$V^{p}_{\lambda} =V^{(p_0)}_{\lambda} \bigcap \HH^p \ZA$, for $p>p_0$.

We now turn to the geometric interpretation of the weight decomposition in (\ref{wd-act-Tp0}). Set $A^{\lambda}$ to be the subscheme
of $Z$ corresponding to the ideal 
$V^{(p_0)}_{\lambda}$ and let $A_{\lambda}$ be the subscheme of $Z$ complementary to $A^{\lambda}$. This gives the decomposition
of $Z$ 
\BEN\label{Z-decAlam}
Z =\bigcup_{\lambda} A_{\lambda}
\EEN
into disjoint union of subconfigurations $A_{\lambda}$'s of $Z$, where the union is taken over the weights $\lambda$ occurring in
(\ref{wd-act-Tp0}). As in the proof of Lemma \ref{lem-Vlam} one shows that the functions in 
$\tilde{T}^{(p_0)} \ZA$
are constant on every $A_{\lambda}$. In particular, each $A_{\lambda}$ determines the codimension one subspace
$H_{\lambda} \subset \tilde{T}^{(p_0)} \ZA$ of functions vanishing on $A_{\lambda}$. Recalling that $Z$ is identified via the map
$\kappa\ZA$ with the subset of $d$ points in $\PP(\HT\ZA^{\ast})$, we deduce that each subconfiguration $A_{\lambda}$ is contained in a hyperplane in $\PP(\HT\ZA^{\ast})$. Furthermore, since $Z$ is in general position in $\PP(\HT\ZA^{\ast})=\PP^r$ it follows 
\BEN\label{degAlam}
deg (A_{\lambda} )\leq r,\,\,\forall \lambda\,.
\EEN
We now claim that this leads to a contradiction. Indeed,
consider the last summand \linebreak $\HH^{\LG-1} \ZA$ of the decomposition in 
(\ref{ord-ZA}). We can find a weight $\lambda$ such that
$$
V^{\LG-1}_{\lambda} = V^{(p_0)}_{\lambda} \bigcap (\HH^{\LG-1} \ZA ) \neq 0\,.
$$
 By definition functions in $V^{\LG-1}_{\lambda} $ vanish on 
$A^{\lambda}$ and hence they have their support in $A_{\lambda}$. Let $x$ be a non-zero element of $ V^{\LG-1}_{\lambda}$ and let $a$ be a point in $A_{\lambda}$ such that
$x(a) \neq 0$. From the inequality in (\ref{degAlam}) it follows that there is $t_a \in \HT\ZA$ such that $t_a (a) \neq 0$ and
$t_a (b) =0$, for all $b\in A_{\lambda} \setminus {a}$. This implies that the support of
$x \cdot t_a $ is $\{a\}$. Thus the delta-function $\delta_a$ is in $\HH^{\LG-2} \ZA \oplus \HH^{\LG-1} \ZA $. But the following claim
implies that this is possible
only if $\LG =2$, which is contrary to the assumption that $\LG \geq 3$.
\begin{cl}\label{delta-ord}
Let $Z$ be as above and let
 $\delta_a$ be the delta-function on $Z$ supported at $a \in Z$.
Let
\BEN\label{delta-ord1}
\delta_a =\sum^{\LG-1}_{p=0} \delta^{(p)}_a
\EEN
be the decomposition of $\delta_a$ according to the direct sum in (\ref{ord-ZA}). Then $ \delta^{(0)}_a \neq 0$.
\end{cl}
{\it Proof of {\bf Claim \ref{delta-ord}.}}\,
Let $\delta^{(m_0)}_a$ be the first (from the left) non-zero component in the decomposition (\ref{delta-ord1}).
To show that it lies in $\HT\ZA$ it is enough to check that $D^{-} (t) (\delta^{(m_0)}_a ) =0$, for all
$t\in \HT\ZA$, - this is a characterization of elements in $\HH^0 \ZA$ proved in \RI, Remark 7.8. For this we 
multiply $\delta_a$ by $t\in \HT\ZA$ and decompose it according to the direct sum in (\ref{ord-ZA})
\BEN\label{t-del-a}
t\cdot \delta_a = D^{-} (t) (\delta^{(m_0)}_a ) + \delta^{\prime}_a\,,
\EEN
where $\delta^{\prime}_a$ is the component of $ t\cdot \delta_a$ contained in 
$\PF^{m_0} \ZA = \bigoplus_{p\geq m_0} \HH^p \ZA$. But $t\cdot \delta_a =t(a)\delta_a $ is a scalar multiple of $\delta_a$,
for all $t\in \HT\ZA$. Hence $t\cdot \delta_a \in \PF^{m_0} \ZA$. This and the decomposition in (\ref{t-del-a}) imply
$$
D^{-} (t) (\delta^{(m_0)}_a ) =0\,,
$$
for all $t\in \HT\ZA$. 
This completes the proof of Claim \ref{delta-ord} as well as the proof of the proposition.
%{\it Proof of surjectivity in (\ref{surj}).}
%Since $A_{\lambda}$ contains at most $r$ points, for every $a\in A_{\lambda}$ we can find $t_a \in \HT\ZA$ such that its restriction to
%$A_{\lambda}$ has support at $a$ only. This implies that 
%\BEN\label{t-a-f}
%t_a \cdot f = t_a (a) f(a) \delta_a ,\,\,\forall f\in V^{(p_0)}_{\lambda}
%\EEN
%Hence $\delta_a $ belongs to the image of the map in (\ref{surj}) as long as there is $f\in V^{\leq p_0}_{\lambda}$ (see (\ref{wsp-ord}) for %%notation))
%which does not vanish at $a$. If this is the case for all $\lambda$ and all $a \in A_{\lambda}$, then, in view of the disjoint sum in
%(\ref{Z-decAlam}), we are done.
%\\
%\indent
%Now assume that there is $\lambda$ and a point $a \in A_{\lambda}$ such that all $f$ in $ V^{\leq p_0}_{\lambda}$ vanish at $a$.
%This implies that there is $x \in V^{p}_{\lambda}$, for some $p > p_0$, which does not vanish at $a$ (otherwise all $f$ in 
%$ V^{(p_0)}_{\lambda}$ vanish at $a$ thus forcing it to be in $A^{\lambda}$, the complement of $A_{\lambda}$ which is impossible). %%        Substituting
%$x$ for $f$ in (\ref{t-a-f}) implies that the delta-function 
%$\delta_a$ lies in $\HH^{p-1} \ZA \oplus \HH^{p} \ZA \oplus \HH^{p+1}$. From Claim \ref{delta-ord} it follows that
%$p\leq 1$. This together with the inequality $p_0 < p$ yield $p_0 =0$, which contradicts the assumption that
%$p_0 \geq 1$. This completes the proof of surjectivity in (\ref{surj}) as well as the proof of the proposition.
\end{pf} 

\section{${\bf sl_2}$-structures on $\FTP$}\label{sec-sl2}
The morphism $d^{+}$ (resp. $d^{-}$) considered in \S\ref{sec-periods}, (\ref{d+}) (resp. (\ref{d-})), attaches intrinsically the nilpotent
endomorphism $d ^{+} (v)$ (resp. $d ^{-} (v)$) to every tangent vector $v$ in $\TPI$. We have seen their importance with respect to
the period maps defined for $\JABG$, for every component $\GA$ in $\CSA$. In this section we explore more subtle
representation theoretic aspects of this assignment by completing $d^{+} (v)$ to an ${\bf sl_2}$-triple. This is made possible by
the well-known Jacobson-Morozov theorem.\footnote{for this and other standard facts about such triples we refer to
\cite{[Kos]}.} Once such a triple is chosen, we look at its representation on the fibre of the sheaf
$\FTP$ at the point of $\JABG$ underlying the tangent vector $v$. This yields further, finer, decomposition
of the orthogonal decomposition of $\FTP$ in  (\ref{ord-Fpr-til}). The resulting structure is somewhat reminiscent of the linear algebra data arising in the theory of Mixed Hodge structure.
\\
\\
\indent
We fix a component $\GA$ in $\CSA$ and assume it to be simple.\footnote{see Definition \ref{s-c}; from the results in
\S\ref{sec-ses}, Theorem \ref{ty-dec-cont}, it follows that this assumption is inessential.} As before the morphism
$\pi: \JABG \longrightarrow \GAB$ stands for the natural projection and $\TPI$ denotes its relative tangent bundle
(recall our convention in \S\ref{sh-bdl} of distinguishing locally free sheaves and the corresponding vector bundle).
We begin by considering the situation at a closed point of $\TPI$ and then give a sheaf version of our construction.

\subsection{Constructions on a fibre of $\FTP$}\label{const-fib}
Fix a point $\ZA \in \JABG$ and let $v$ be a tangent vector in $\TPI$ lying over $\ZA$.
Evaluating the morphism $d^{+}$ (resp. $d^{-}$) in (\ref{d+}) (resp. (\ref{d-})) at the point
$([Z],[\alpha],v) \in \TPI$, we obtain the endomorphism $d^{+} (v)$ (resp. $d^{-} (v)$) of
$\FTP ([Z])$, the fibre\footnote{recall, by (\ref{Fpr-tilde}), the fibre $\FTP ([Z])$ does not depend on $[\alpha]$.}
 of $\FTP$ at $\ZA$. By definition of 
{\BM
$\LAG$}
the endomorphisms $d^{\pm} (v)$ are nilpotent elements of
{\BM
$\LAG$}$\ZA$,
the fibre of 
{\BM
$\LAG$}
at $\ZA$. 
By Jacobson-Morozov  theorem, 
$d^{+} (v)$ (resp. $d^{-} (v)$) can be completed to an ${\bf sl_2}$-triple
$\{d^{+} (v), h, y \}$ (resp. $ \{y^{\prime}, h^{\prime}, d^{-} (v) \}$),
where $h$ (resp. $h^{\prime}$) is a semismple element of
{\BM
$\LAG$}$\ZA$
subject to the standard relations
\begin{eqnarray*}
[h,d^{+} (v)] = 2 d^{+} (v)& [h,y]=-2y&[d^{+} (v),y]=h\,,\\
\mbox{(resp.}\,\, [h^{\prime}, d^{-} (v)] = -2 d^{-} (v)& [h^{\prime} ,y^{\prime}]=2y^{\prime} &[y^{\prime},d^{-} (v)]=h^{\prime}\,\,)\,.
 \end{eqnarray*}

It is well-known  that semisimple elements coming along with
$d^{+} (v)$ (resp. $d^{-} (v)$) in an ${\bf sl_2}$-triple form a homogeneous space modeled on the nilpotent Lie algebra
\BEN\label{g+-}
{\bf g}^{+} (v) = ker (ad(d^{+} (v))) \cap im(ad(d^{+} (v)))
\,\,
(\mbox{resp.}\,\, {\bf g}^{-} (v) = ker (ad(d^{-} (v))) \cap im(ad(d^{-} (v)))\,)\,,
\EEN
i.e. two choices for a semisimple element in an ${\bf sl_2}$-triple for
$d^{\pm} (v)$
differ by an element in 
${\bf g}^{\pm} (v)$.
In fact it is known that if $h$ (resp. $h^{\prime}$) is a semisimple element which goes along with 
$d^{+} (v)$ (resp. $d^{-} (v)$) in an ${\bf sl_2}$-triple, then any other semisimple element
$\tilde{h}$ (resp. $\tilde{h^{\prime}}$) for $d^{+} (v)$ (resp. $d^{-} (v)$) can be taken to be of the form
$$
\tilde{h} =\exp(ad(w))h\,\,
(\mbox{resp.}\,\,
\tilde{h^{\prime}} =\exp(ad(w^{\prime}))h^{\prime})\,,
$$
for some $w \in {\bf g}^{+} (v)$ (resp. $w^{\prime} \in {\bf g}^{-} (v)$).
Thus the set of semisimple elements for
$d^{+} (v)$ (resp. $d^{-} (v)$) in an ${\bf sl_2}$-triple is a principle homogeneous space for the action of the unipotent
group
\begin{eqnarray*}
 G^{+} (v) =&\{ exp(ad(w)) \mid w \in {\bf g}^{+} (v) \} \\
(\mbox{resp.}\,\, G^{-} (v) =&\{ exp(ad(w^{\prime})) \mid w^{\prime} \in {\bf g}^{-} (v) \}\,)
\end{eqnarray*}
(see \cite{[Kos]}, Theorem 3.6). 

Next we bring in the grading of
{\BM
$\LAG$}
in (\ref{sLie-grad}). This gives the grading 
\BEN\label{gr-gZ}
\mbox{\BM$\LAG $}\ZA = \bigoplus^{\LG-1}_{i=-(\LG-1)} \mbox{\BM$\LAG^i $}\ZA
\EEN
on the fibre 
$\mbox{\BM$\LAG $}\ZA$ of 
{\BM
$\LAG$}
at $\ZA$. The fact that 
$d^{+} (v)$ (resp. $d^{-} (v)$) is of degree $1$ (resp. $(-1)$) with respect to this grading allows us to choose
$h$ (resp. $h^{\prime}$) to be in the summand 
$ \mbox{\BM$\LAG^0 $}\ZA$ and $y$ (resp. $y^{\prime}$) in 
$ \mbox{\BM$\LAG^{-1} $}\ZA$ (resp. 
$ \mbox{\BM$\LAG^{1} $}\ZA$).
We will always assume such a choice of $h$ and $y$ (resp. $h^{\prime}$ and  $y^{\prime}$).

In this graded version of  ${\bf sl_2}$-triples associated to 
$d^{+} (v)$ (resp. $d^{-} (v)$) the corresponding semisimple elements form the homogeneous subset 
${\bf h^0}(d^{+} (v))$ (resp. ${\bf h^0}(d^{-} (v))$) of \label{set-h0}
$\mbox{\BM$\LAG^0 $}\ZA$
modeled on the nilpotent Lie subalgebra
\BEN\label{g0+-}
{\bf g^0_{+}} ([Z],[\alpha],v) =\mbox{\BM$\LAG^0 $}\ZA \bigcap {\bf g}^{+} (v)\,\,(resp.\,\, 
{\bf g^0_{-}} ([Z],[\alpha],v) =\mbox{\BM$\LAG^0 $}\ZA \bigcap {\bf g}^{-} (v) )\,,
\EEN
 where ${\bf g}^{\pm} (v)$ are nilpotent Lie subalgebras defined in (\ref{g+-}).
Furthermore, the unipotent groups
\BEN\label{G0+-}
 G^0_{\pm} ([Z],[\alpha],v) =\{ \exp(ad(w)) \mid w \in {\bf g^0}_{\pm} ([Z],[\alpha],v)\}
\EEN 
act simply transitively on
${\bf h^0}(d^{\pm} (v))$
respectively (this can be seen by adapting the argument in the proof of Theorem 3.6, p.987, \cite{[Kos]}, to the graded situation at hand).
Thus
${\bf h^0}(d^{+} (v))$ (resp. ${\bf h^0}(d^{-} (v))$
is a principal homogeneous space for the unipotent group
$ G^0_{+} ([Z],[\alpha],v)$ (resp. $ G^0_{-} ([Z],[\alpha],v)$).

We will now fix an ${\bf sl_2}$-triple
$\{ d^{+} (v), h, y \}$ with $h \in  {\bf h^0}(d^{+} (v))$ and $y \in  \mbox{\BM$\LAG^{-1} $}\ZA$,
and consider its action on 
$\FTP([Z])$.
This gives the weight decomposition
\BEN\label{sl-wd}
\FTP([Z]) =\bigoplus_{n \in {\bf Z}} W(n)
\EEN
under the action of $h$, i.e. $W(n)$ is the eigen space of $h$ corresponding to the eigen value $n$ of $h$.
Since
$(d^{+} (v))^{\LG} =0$ it follows that the weights occurring in
(\ref{sl-wd}) are in the set
\BEN\label{spectr-h}
\{-(\LG-1), \ldots, (\LG-1) \}\,.
\EEN
Following the convention in Hodge theory, we shift the grading of weights to the right by $(\LG-1)$ to obtain
\BEN\label{sl-wd-shift}
\FTP([Z]) =\bigoplus^{2(\LG-1)}_{n =0} V(n)\,,
\EEN
where $V(n) =W(n-\LG+1)$, for $n=0, \ldots, 2(\LG-1)$. Abusing the language it will be called weight decomposition
of $\FTP([Z]) $ as well.

A choice of $h$ in $\mbox{\BM$\LAG^0 $}\ZA$ implies that
$h$ preserves the orthogonal decomposition in (\ref{ord-Fpr-til}). Thus each summand
$\HH^p \ZA$ admits the weight decomposition
\BEN\label{Hp-wd}
\HH^p \ZA  = \bigoplus^{2(\LG-1)}_{n =0} \HH^p \ZA (n)\,,
\EEN
where $\HH^p \ZA (n) = \HH^p \ZA \bigcap V(n)$. We will adopt the notation of Hodge theory by writing the 
double grading above as follows
\BEN\label{double}
\HH^{p,n-p} ([Z],[\alpha],v,h) :=\HH^p \ZA (n)\,. 
\EEN
The following result gives a more precise version of the decomposition in (\ref{Hp-wd}) in this bigraded form.
\begin{lem}\label{Hp-dgr}
The weight decomposition in (\ref{Hp-wd}) has the following form 
$$
\HH^p \ZA  = \bigoplus^{p +\LG-1}_{n =p} \HH^{p,n-p} ([Z],[\alpha],v,h)\,. 
$$
\end{lem}
\begin{pf}
We need to establish the possible range of weights of $h$ on $\HH^p \ZA$.
Let $m$ be a positive weight of $h$ on $\HH^p \ZA$.
From the properties of ${\bf sl_2}$-representations it follows that the vectors in 
$\HH^p \ZA$ of weight $m$ come from some vectors of weight $-m$ upon
applying $(d^{+} (v))^m$. Such vectors are situated in the summand 
$\HH^{p-m} \ZA$ of the orthogonal decomposition 
in (\ref{ord-Fpr-til}). Since $(p-m)$ is non-negative we obtain the upper bound
\BEN\label{up}
m\leq p,
\EEN
for the weights of $h$ on $\HH^p \ZA$. By our shifting convention
$m=n-\LG+1$, for some $n\in \{0, \ldots, 2(\LG-1) \}$. Combining this with
(\ref{up}) yields the asserted upper bound
$$
n \leq p + \LG -1.
$$

Let $-m\, (m>0)$ be a negative weight of $h$ on $\HH^p \ZA$. Then we can push vectors
in $\HH^p \ZA$ having this weight to the weight space of weight $m$ by applying
$(d^{+} (v))^m$. This will take the vectors from
$\HH^p \ZA$ to $\HH^{p+m} \ZA$. Since $p+m \leq \LG-1$ we obtain
$$
-m =n-\LG+1 \geq p-\LG+1,
$$
where the equality is the shift convention. This implies $n\geq p$ as asserted in the lemma.
\end{pf}
The above considerations can be turned around by saying that every weight space
$V(n)$ in (\ref{sl-wd-shift}) admits the orthogonal decomposition induced by (\ref{ord-Fpr-til})
$$
V(n) = \bigoplus^{\LG-1}_{i=0} \HH^{i,n-i} ([Z],[\alpha],v,h)\,.
$$
This follows from the fact that $h$ is in $\mbox{\BM$\LAG^0$}\ZA$. From Lemma \ref{Hp-dgr} this decomposition has the following form
\BEN\label{Vn-dgr}
V(n) = \bigoplus_{i+j=n \atop 0 \leq i,j \leq \LG-1} \HH^{i,j} ([Z],[\alpha],v,h)\,.
\EEN
Putting together this double grading with the weight decomposition in (\ref{sl-wd-shift}) gives a bigrading
on $\FTP([Z])$. This fact and  various properties of the weight and double gradation of $\FTP([Z])$
are summarized below.
\begin{pro}\label{pro-dgr}
\begin{enumerate}
\item[1)]
The weight and orthogonal decompositions of $\FTP([Z])$ in (\ref{sl-wd-shift}) and (\ref{ord-Fpr-til}) respectively
define a bigrading  on  $\FTP([Z])$:
\BEN\label{bi-gr-Z}
\FTP([Z]) = \bigoplus _{(p,q)} \HH^{p,q}([Z],[\alpha],v,h)\,,
\EEN
where the direct sum is taken over the points $(p,q) \in {\bf Z^2}$ lying in the square with vertices
$(0,0), \,( \LG-1, 0), \,(\LG-1, \LG-1),\, (0,\LG-1)$.  (Once $([Z],[\alpha],v)$ is fixed and $h$ is chosen the reference to these 
parameters will be omitted and we simply write
$\HH^{p,q}$ instead of $\HH^{p,q}([Z],[\alpha],v,h)$.)

\item[2)]
The endomorphism $d^{+}(v)$ (resp. $h$ and $y$) is of type $(1,1)$ (resp. $(0,0)$ and $(-1,-1)$)
with respect to the bigrading in (\ref{bi-gr-Z}), i.e.
$$
d^{+}(v) (\HH^{p,q}) \subset \HH^{p+1,q+1}\,,
$$
for every $(p,q)$ occurring in (\ref{bi-gr-Z}) (with the understanding that a graded piece
$\HH^{m,n}$ equals zero, unless $(m,n)$ lies in the square described in 1) of the proposition).
\item[3)]
For every integer $k\in \{0,\ldots, \LG-1\}$ the endomorphism
$d^{+}(v) $ induces the isomorphism
$$
(d^{+}(v) )^k : V(\LG-1-k) \longrightarrow V(\LG-1+k) 
$$
which is of type $(k,k)$ with respect to the bigrading in (\ref{Vn-dgr}). In particular, one has isomorphisms
$$
(d^{+}(v) )^k : \HH^{p,q} \longrightarrow \HH^{p+k,q+k}\,,
$$
for every $(p,q)$ such that $p+q =\LG -1-k$ and $p,q \geq 0$.
\end{enumerate}
\end{pro}
\begin{pf}
The first assertion follows from (\ref{Vn-dgr}) and the range for the values of $n$ in the weight decomposition
(\ref{sl-wd-shift}). For the second assertion observe that 
$d^{+}(v)$ raises a weight by $2$ and has degree $1$ with respect to the orthogonal decomposition in 
(\ref{ord-Fpr-til}), i.e.
$$
d^{+}(v) (V(n)) \subset V(n+2)\,\,\mbox{and}\,\, d^{+}(v)(\HH^p \ZA) \subset \HH^{p+1} \ZA\,.
$$
This implies 
$$
d^{+}(v) (\HH^{p,q}) =d^{+}(v) (\HH^p \ZA \bigcap V(p+q)) \subset \HH^{p+1} \ZA \bigcap V(p+q+2) = \HH^{p+1, q+1}\,.
$$
The third assertion follows from the properties of 
${\bf sl_2}$-representations and part 2) of the proposition.
\end{pf}

The weight spaces $V(n)$ in (\ref{sl-wd-shift}) depend not only on 
$d^{+}(v)$ but also on the choice of $h$ in ${\bf h^0} (d^{+}(v))$, the space defined on
p.\pageref{set-h0}. However, the subspaces
\BEN\label{W-n}
W^n = \bigoplus _{ m\geq n} V(m)
\EEN
are independent of $h$ (this follows from the fact that the group
$G^0_{+} (v)$ in (\ref{G0+-})  preserves these subspaces and acts transitively on ${\bf h^0} (d^{+}(v))$).
This way one obtains a filtration of 
$\FTP([Z])$:
\BEN\label{w-filt-v}
\FTP([Z]) =W^0 \supset W^1 \supset \ldots \supset W^{2(\LG-1)} \supset W^{2\LG-1} =0
\EEN
intrinsically associated to $d^{+}(v)$. This filtration will be called
the {\it weight filtration of $d^{+}(v)$ on $\FTP([Z])$} and denoted by
$W^{\bullet} (\FTP([Z]), v)$. Its properties, summarized below, are easily derived from
Proposition \ref{pro-dgr}.
\begin{pro}\label{pro-w-filt-v}
The weight filtration 
$W^{\bullet} (\FTP([Z]), v)$ in (\ref{w-filt-v}) has the following properties.
\begin{enumerate}
\item[1)]
 The homomorphism $d^{+}(v)$ takes $W^k$ to $W^{k+2}$, i.e. one has
$$
d^{+}(v) (W^k) \subset W^{k+2}\,,
$$
for every $k\geq 0$.
\item[2)]
Set
$$
Gr^k_{W^{\bullet}}  (\FTP([Z]))= W^k / W^{k+1}\,.
$$
 The orthogonal decomposition
(\ref{ord-Fpr-til}) induces a grading on  $Gr^k_{W^{\bullet}} (\FTP([Z]))$
$$
Gr^k_{W^{\bullet}}  (\FTP([Z])) = \bigoplus^{\LG-1}_{i=0} Gr^{i,k-i}_{W^{\bullet}}  (\FTP([Z]))\,, 
$$
where the summands are defined as follows
$$
Gr^{i,k-i}_{W^{\bullet}}  (\FTP([Z])) = \HH^i \ZA \bigcap W^k /{\HH^i \ZA \bigcap W^{k+1}} =
\big{(} \bigoplus_{m\geq k} \HH^{i,m-i} \big{)} /\big{(} \bigoplus_{m\geq k+1} \HH^{i,m-i} \big{)}\,.
$$
\item[3)]
The associated graded vector space
$$
Gr_{W^{\bullet}} (\FTP([Z]))  = \bigoplus^{2(\LG-1)}_{k=0} Gr^k_{W^{\bullet}}  (\FTP([Z])) =\bigoplus^{2(\LG-1)}_{k=0} W^k / W^{k+1}
$$
of $\FTP([Z])$ with respect to the weight filtration $W^{\bullet}$ in (\ref{w-filt-v}), together with the grading in 2), acquires the bigrading
\BEN\label{Gr-dgr}
Gr_{W^{\bullet}} (\FTP([Z]))  = \bigoplus_{(p,q)}  Gr^{p,q}_{W^{\bullet}}  (\FTP([Z]))\,, 
\EEN
where the direct sum is taken over the points $(p,q) \in {\bf Z^2}$ lying in the square described in Proposition \ref{pro-dgr}, 1).
\item[4)]
The endomorphism $d^{+}(v)$ induces the homomorphism
$$
gr(d^{+}(v)) : Gr_{W^{\bullet}} (\FTP([Z])) \longrightarrow Gr_{W^{\bullet}} (\FTP([Z]))
$$
which has type (1,1) with respect to the double grading in 3), i.e. one has
$$
gr(d^{+}(v)) (Gr^{p,q}_{W^{\bullet}}  (\FTP([Z])) ) \subset  Gr^{p+1,q+1}_{W^{\bullet}}  (\FTP([Z]))\,,
$$
for every $(p,q)$ occurring in the decomposition in (\ref{Gr-dgr}).
\item[5)]
For every integer $k \in \{0,\ldots, \LG-1 \}$ the endomorphism $gr(d^{+}(v))$ induces the isomorphism
$$
(gr(d^{+}(v)))^k : Gr^{\LG-1-k}_{W^{\bullet}}  (\FTP([Z])) \longrightarrow Gr^{\LG-1+k}_{W^{\bullet}}  (\FTP([Z]))
$$
which is of type $(k,k)$ with respect to the bigrading in 3). In particular, one has isomorphisms
$$
(gr(d^{+}(v)))^k : Gr^{p,q}_{W^{\bullet}}  (\FTP([Z])) \longrightarrow Gr^{p+k,q+k}_{W^{\bullet}}  (\FTP([Z]))\,,
$$
for every $(p,q)$ such that 
 $p+q =\LG -1-k$ and $p,q \geq 0$.
\end{enumerate}
\end{pro}

\subsection{Sheaf version of the constructions in \S\ref{const-fib}}\label{const-sh} 
We will be working on the relative tangent bundle $\TPI$ of the natural projection
$$
\pi: \JABG \longrightarrow \GAB\,.
$$
Thus $\TPI$ is a scheme over $\JABG$ with the structure projection
\BEN\label{tau}
\tau: \TPI \longrightarrow \JABG\,,
\EEN
whose fibre $\TPI \ZA$ over a point $\ZA \in \JABG$ is the vector space
of vertical\footnote{this terminology is explained in the footnote on page \pageref{footnote:vv}.} tangent vectors of $\JABG$ at $\ZA$.

We begin by building the scheme parametrizing the ${\bf sl_2}$-triples discussed in \S\ref{const-fib}.
For this consider the pullback 
{\BM
$\mbox{\UB$\tau^{\ast}$}\LAG$
of the sheaf of semisimple algebras\footnote{since the component $\GA$ is assumed to be simple, the sheaf 
{\BM
$\LAG$}, by Corollary \ref{cor-s},
is actually a sheaf of simple Lie algebras.
But this will not matter in the constructions below.} $\LAG$}.
The morphisms $d^{\pm}$ defined in (\ref{d+}) and (\ref{d-}) give rise to two distinguished sections of
{\BM
$\mbox{\UB$\tau^{\ast}$}\LAG$}. We will denote them by
$d^{\pm}_{\diamond}$. Thus the value of $d^{\pm}_{\diamond}$ at a closed point
$([Z],[\alpha],v) \in \TPI$ are the nilpotent elements $d^{\pm}(v)$ in 
{\BM
$\LAG$}$\ZA$ considered in \S\ref{const-fib}.

The sections $d^{\pm}_{\diamond}$ define the morphisms of sheaves
\BEN\label{ad-d+-}
ad(d^{\pm}_{\diamond}): \tau^{\ast} \mbox{\BM$\LAG$} \longrightarrow \tau^{\ast} \mbox{\BM$\LAG$}\,,
\EEN
whose value at $([Z],[\alpha],v) \in \TPI$ is
\BEN\label{ad-d+-v}
ad(d^{\pm}(v)): \mbox{\BM$\LAG$}\ZA \longrightarrow \mbox{\BM$\LAG$}\ZA\,.
\EEN

From now on we consider the section $d^{+}_{\diamond}$ only. The morphism
$ad(d^{+}_{\diamond})$ gives rise to two subsheaves 
$ker(ad(d^{+}_{\diamond}))$ and $im(ad(d^{+}_{\diamond}))$ of
{\BM
$\mbox{\UB$\tau^{\ast}$}\LAG$}.
Set
\BEN\label{sh-g+-}
\mbox{\BM$\GS$}( d^{+}) = ker(ad(d^{+}_{\diamond})) \bigcap im(ad(d^{+}_{\diamond}))\,.
\EEN
Taking account of the grading of 
{\BM
$\LAG$}
in (\ref{sLie-grad}) and the fact that $d^{+}_{\diamond}$ is a section of
{\BM
$\mbox{\UB$\tau^{\ast}$}\LAG^1$}
  we obtain a grading on 
$\mbox{\BM$\GS$}( d^{+})$:
\BEN\label{sh-g+-gr}
\mbox{\BM$\GS$}( d^{+}) = \bigoplus^{\LG-1}_{i=-(\LG -1)} \mbox{\BM$\GS^i$}( d^{+})\,, 
\EEN
where $ \mbox{\BM$\GS^i$}( d^{+}) =\mbox{\BM$\GS$}( d^{+}) \bigcap  \tau^{\ast}\mbox{\BM$\LAG^i$}$.
We will be interested in degree zero part
\BEN\label{sh-g0+-}
\mbox{\BM$\GS^0$}( d^{+}) = \mbox{\BM$\GS$}( d^{+}) \bigcap  \tau^{\ast}\mbox{\BM$\LAG^0$}\,.
\EEN
This is a sheaf of nilpotent Lie algebras whose fibre at a point
$([Z],[\alpha],v) \in \TPI$ is the subalgebra
${\bf g^0_{+}} ([Z],[\alpha],v)$ defined in (\ref{g0+-}).

Taking 
$Exp =\exp \circ ad$ of $\mbox{\BM$\GS^0$}( d^{+}) $ we obtain the subsheaf of
${\cal A}{\it ut} (\mbox{\BM$\GS$}_{\GA} )$ which will be denoted by
$\mbox{\BM$G^0$}(d^{+}) $. This is a sheaf of unipotent groups whose fibre 
at $([Z],[\alpha],v) \in \TPI $ is the group
$G^0_{+} ([Z],[\alpha],v)$ defined in (\ref{G0+-}).
We can now define a scheme over 
$\TPI$ which parametrizes ${\bf sl_2}$-triples having the values of 
$d^{+}$ as nil-positive elements. This scheme is denoted by
$\mbox{\BM$h^0$}(d^{+})$ and it is defined by the following incidence relation.
\begin{equation}\label{sc-sl2}
\mbox{\BM$h^0$}(d^{+}) =\left\{ \left.{\scriptstyle{([Z],[\alpha],v,h) \in \TPI \times_{\JABG} \mbox{\BM$\LAG^0$}}} \right|
\stackrel{\mbox{$\scriptstyle{\{ d^{+}(v), h, y \}}$ \scriptsize{is an} ${\scriptstyle{\bf sl_2}}$\scriptsize{-triple,}
\scriptsize{where $\scriptstyle{d^{+}(v)}$ is nil-positive,}}}
{\mbox{$\scriptstyle{h \in \mbox{\BM$\LAG^0$}\ZA}$ \scriptsize{is semisimple, and $\scriptstyle{y \in \mbox{\BM$\LAG^{-1}$}\ZA}$ is nil-negative elements}}}
\right\} \,, 
\end{equation}
where $\mbox{\BM$\LAG^0$}$ is considered here as a scheme over $\JABG$ and 
$\TPI \times_{\JABG} \mbox{\BM$\LAG^0$}$ is the fibre product  of 
$\JABG$-schemes $\TPI$ and $\mbox{\BM$\LAG^0$}$.

The projection
\BEN\label{pr-1}
 \eta^{+}_1: \mbox{\BM$h^0$}(d^{+}) \longrightarrow \TPI
\EEN
makes 
$\mbox{\BM$h^0$}(d^{+}) $ a $\mbox{\BM$G^0$}(d^{+}) $-principal homogeneous fibration over
$\TPI$, whose fibre over a closed point
$([Z],[\alpha],v)$ is the homogeneous space
${\bf h^0}(d^{+}(v))$ introduced in \S\ref{const-fib}, p.\pageref{set-h0}.

Set
\BEN\label{tau+}
\tau^{+} =\tau \circ {\bf \eta^{+}_1} : \mbox{\BM$h^0$}(d^{+})
 \stackrel{{\bf \eta^{+}_1}}{\longrightarrow} \TPI 
\stackrel{\tau}{\longrightarrow} \JABG
\EEN
and consider the pullback 
$(\tau^{+})^{\ast}\mbox{\BM$\LAG$}$. From the definition of 
the scheme 
$\mbox{\BM$h^0$}(d^{+})$ in (\ref{sc-sl2}) it follows that 
$(\tau^{+})^{\ast}\mbox{\BM$\LAG$}$ comes equipped with three distinguished sections
which will be denoted by
$s^{+}, s^0, s^{-}$.
Their values at a closed point
$([Z],[\alpha],v,h) \in \mbox{\BM$h^0$}(d^{+})$ are respectively:

$
s^{+} ([Z],[\alpha],v,h) = d^{+} (v)$,  
$s^0 ( [Z],[\alpha],v,h) =h$ and 
$s^{-} ([Z],[\alpha],v,h)$ is uniquely determined by the requirement that
\BEN\label{sl-univ}
\{ d^{+} (v), h, s^{-} ([Z],[\alpha],v,h) \}
\EEN
is an ${\bf sl_2}$-triple with 
$ d^{+} (v)$ and $h$ being respectively its nil-positive and semisimple elements.
Thus $\mbox{\BM$h^0$}(d^{+})$ carries a ``universal'' ${\bf sl_2}$-triple spanned
by the sections $s^{+}, s^0, s^{-}$.
\begin{rem}\label{s+=pb}
By construction, the section $s^{+}$ is the pullback by $\eta^{+}_1$ of the section
$d^{+}_{\diamond}$, i.e. we have
$$
s^{+} = (\eta^{+}_1)^{\ast} d^{+}_{\diamond}\,. 
$$
\end{rem}

To define a sheaf version of the weight decomposition in (\ref{sl-wd}) consider the pullback
$(\tau^{+})^{\ast}\FTP$. By definition the sheaf
{\BM
$\LAG$}
comes together with its faithful representation on $\FTP$, i.e. it is defined as a subsheaf\footnote{from Corollary \ref{cor-s} we know that 
$\mbox{\BM$\LAG$} = {\bf sl}(\FTP)$.} of
$\ENDO(\FTP)$. In particular, we think of sections 
$s^{\pm},s^0$ as endomorphisms of the sheaf $(\tau^{+})^{\ast}\FTP$.
Furthermore, the endomorphism
\BEN\label{s0}
s^0: (\tau^{+})^{\ast}\FTP \longrightarrow  (\tau^{+})^{\ast}\FTP
\EEN
is semisimple with integer eigen values and its eigen sheaves provide the sheaf version of the weight decomposition in (\ref{sl-wd}).
More precisely, define the weight sheaf ${\cal W}(n)$ corresponding to a weight $n \in {\bf Z}$ by the formula
\BEN\label{w-sh}
{\cal W}(n) = ker (s^0 - n {\bf id}_{ (\tau^{+})^{\ast}\FTP} )\,.
\EEN
This gives a decomposition of 
$(\tau^{+})^{\ast}\FTP$  into the direct sum of the subsheaves
${\cal W}(n) \,(n\in  {\bf Z})$. From the
orthogonal decomposition (\ref{ord-Fpr-til}) and the condition that the values of 
$s^{+}$ are endomorphisms of degree $1$ with respect to this decomposition, it follows that $(s^{+})^{\LG} =0$.
This implies that
the weights of $s^0$ belong to the set $\{-(\LG-1), \ldots, \LG-1 \}$. This yields the following
sheaf version of (\ref{sl-wd}):
\BEN\label{sl-wd-sh}
(\tau^{+})^{\ast}\FTP = \bigoplus^{\LG-1}_{n= -(\LG-1)} {\cal W}(n)\,.
\EEN
Setting
${\cal V}(n):= {\cal W}(n-\LG+1)$, for $n=0,1, \ldots, 2(\LG-1)$, gives the decomposition
\BEN\label{sh-shifted}
(\tau^{+})^{\ast}\FTP = \bigoplus^{2(\LG-1)}_{n= 0} {\cal V}(n)\,.
\EEN
This decomposition will be called the (shifted) weight decomposition of 
$(\tau^{+})^{\ast}\FTP $.

The sheaf 
$(\tau^{+})^{\ast}\FTP $
continues to have the orthogonal decomposition
\BEN\label{ord-tau}
(\tau^{+})^{\ast}\FTP =\bigoplus^{\LG-1}_{p= 0} (\tau^{+})^{\ast}\HH^p
\EEN
and the sections $s^{\pm},\,s^0$ have respectively degrees $\pm1,\,0$, relative to this decomposition.

Define the sheaves $\HH^{p,q}$ as follows:
\BEN\label{Hpq}
\HH^{p,q} = (\tau^{+})^{\ast}\HH^p \bigcap {\cal V}(p+q)\,.
\EEN
Then the sheaf version of Lemma \ref{Hp-dgr} takes the form
\BEN\label{sh-Hp-dgr}
(\tau^{+})^{\ast}\HH^p =\bigoplus^{p+\LG-1}_{n= p} \HH^{p,n-p}\,,
\EEN
for every $p=0,\ldots, \LG-1$. This implies the sheaf version of the bigraded decomposition in (\ref{Vn-dgr}):
\BEN\label{sh-Vn-gr}
{\cal V}(n) = \bigoplus_{i+j=n \atop 0\leq i,j \leq \LG-1} \HH^{i,j}\,.
\EEN
 Finally, the sheaf analogue of Proposition \ref{pro-dgr}
is as follows.
\begin{pro}\label{pro-dgr-sh}
\begin{enumerate}
\item[1)]
The weight and the orthogonal decompositions of $(\tau^{+})^{\ast}\FTP$ in (\ref{sh-shifted}) and (\ref{ord-tau}) respectively
define a bigrading  on  $(\tau^{+})^{\ast}\FTP$:
\BEN\label{bi-gr-tau}
(\tau^{+})^{\ast}\FTP = \bigoplus _{(p,q)} \HH^{p,q}\,,
\EEN
where the sheaves $\HH^{p,q}$ are as in (\ref{Hpq}) and where the 
direct sum is taken over the points $(p,q) \in {\bf Z^2}$ lying in the square with vertices
$(0,0), \,( \LG-1, 0), \,(\LG-1, \LG-1),\, (0,\LG-1)$. 
\item[2)]
The sections $s^{\pm}$ (resp. $s^0$ ) viewed as endomorphisms of 
$(\tau^{+})^{\ast}\FTP$
are of type $(\pm 1, \pm 1)$ (resp. $(0,0)$ )
with respect to the bigrading in (\ref{bi-gr-tau}), i.e.
$$
s^{\pm} (\HH^{p,q}) \subset \HH^{p\pm 1,q \pm 1}\,\,(\mbox{resp.} \,\,
s^{0} (\HH^{p,q}) \subset \HH^{p,q})\,,
$$
for every $(p,q)$ occurring in (\ref{bi-gr-tau}), with the understanding that a graded piece
$\HH^{m,n}$ equals zero, unless $(m,n)$ lies in the square described in 1) of the proposition.
\item[3)]
For every integer $k\in \{0,\ldots, \LG-1\}$ the endomorphism
$s^{+} $ induces the isomorphism
$$
(s^{+} )^k : {\cal V}(\LG-1-k) \longrightarrow {\cal V}(\LG-1+k) 
$$
which is of type $(k,k)$ with respect to the bigrading in (\ref{bi-gr-tau}). In particular, one has isomorphisms
$$
(s^{+} )^k : \HH^{p,q} \longrightarrow \HH^{p+k,q+k}\,,
$$
for every $(p,q)$ such that $p+q =\LG -1-k$ and $p,q \geq 0$.
\end{enumerate}
\end{pro}
Next we turn to the sheaf analogue of the weight filtration (\ref{w-filt-v}).
For this we set
\BEN\label{w-filt-sh0}
{\cal V}^n = \bigoplus_{m\geq n} {\cal V}(m)\,.
\EEN
We know that ${\cal V}^n$ is constant along the fibres of $\eta^{+}_1$ in (\ref{pr-1}) (this follows from the definition
of $\mbox{\BM$h^0$}(d^{+})$ and the fact that the spaces in (\ref{W-n}) are independent of a semisimple element
which is used to define them).
This implies that there is a sheaf
${\cal W}^n$ on $\TPI$ such that
\BEN\label{Vn=pb}
{\cal V}^n = (\eta^{+ }_1)^{\ast} {\cal W}^n\,.
\EEN
This way we obtain a filtration of $\tau^{\ast} \FTP$:
\BEN\label{w-filt-sh}
\tau^{\ast} \FTP ={\cal W}^0 \supset {\cal W}^1 \supset \ldots \supset {\cal W}^{2(\LG-1)} \supset {\cal W}^ {2\LG-1} =0
\EEN
which will be called the weight filtration of 
$\tau^{\ast} \FTP$ associated to
$d^{+}$ or, simply, the weight filtration of $d^{+}$. 

This filtration and the orthogonal decomposition 
\BEN\label{ord-tau1}
\tau^{\ast} \FTP = \bigoplus^{\LG-1}_{p=0} \tau^{\ast} \HH^p\,,
\EEN
obtained as the pullback by $\tau$ of (\ref{ord-Fpr-til}), give the following
  sheaf version of Proposition \ref{pro-w-filt-v}.
\begin{pro-defi}\label{pro-w-filt-sh}
The sheaf $\tau^{\ast} \mbox{\BM$\LAG$}$ comes along with a distinguished section
$d^{+}_{\diamond}$. This section, viewed as an endomorphism of the sheaf
$\tau^{\ast} \FTP$, gives rise to the following structures. 
\begin{enumerate}
\item[1)]
$d^{+}_{\diamond}$ defines the weight filtration 
${\cal W}^{\bullet}$ as in (\ref{w-filt-sh}) and  $d^{+}_{\diamond}$ acts on it by shifting the index of the filtration by $2$, i.e. one has
$$
d^{+}_{\diamond} ({\cal W}^k) \subset {\cal W}^{k+2}\,,
$$
for every $k\geq 0$.
\item[2)]
Set
$$
Gr^k_{{\cal W}^{\bullet}}  (\tau^{\ast} \FTP)= {\cal W}^k /{\cal W}^{k+1}\,.
$$
 The orthogonal decomposition in (\ref{ord-tau})
 induces a grading on  $Gr^k_{{\cal W}^{\bullet}} (\tau^{\ast} \FTP)$
$$
Gr^k_{{\cal W}^{\bullet}}  ( \tau^{\ast} \FTP) = \bigoplus_i Gr^{i,k-i}_{{\cal W}^{\bullet}}  ( \tau^{\ast} \FTP)\,, 
$$
where the summands are defined as follows
$$
Gr^{i,k-i}_{{\cal W}^{\bullet}}  (\tau^{\ast} \FTP) = \tau^{\ast} (\HH^i ) \bigcap {\cal W}^k /{\tau^{\ast} (\HH^i)  \bigcap {\cal W}^{k+1}} \,.
$$ 
\item[3)]
The associated graded sheaf
$$
Gr_{{\cal W}^{\bullet}} (\tau^{\ast} \FTP)  = \bigoplus^{2(\LG-1)}_{k=0} Gr^k_{{\cal W}^{\bullet}}  (\tau^{\ast} \FTP) =
\bigoplus^{2(\LG-1)}_{k=0} {\cal W}^k /{ \cal W}^{k+1}
$$
of $\tau^{\ast} \FTP$, defined with respect to the weight filtration ${\cal W}^{\bullet}$,  together with the grading
 in 2), acquires the bigrading
\BEN\label{Gr-dgr1}
Gr_{{\cal W}^{\bullet}} (\tau^{\ast} \FTP)  = \bigoplus_{(p,q)}  Gr^{p,q}_{{\cal W}^{\bullet}}  (\tau^{\ast} \FTP)\,, 
\EEN
where the direct sum is taken over the points $(p,q) \in {\bf Z^2}$ lying in the square described in Proposition \ref{pro-dgr-sh}, 1).
\item[4)]
The endomorphism $d^{+}_{\diamond}$ induces the endomorphism
$$
gr(d^{+}_{\diamond}) : Gr_{{\cal W}^{\bullet}} (\tau^{\ast} \FTP) \longrightarrow Gr_{{\cal W}^{\bullet}} (\tau^{\ast} \FTP)
$$
which has type (1,1) with respect to the double grading in 3), i.e. one has
$$
gr(d^{+}_{\diamond}) \big{(} Gr^{p,q}_{{\cal W}^{\bullet}}  (\tau^{\ast} \FTP ) \big{)} \subset  Gr^{p+1,q+1}_{{\cal W}^{\bullet}}  (\tau^{\ast} \FTP)\,,
$$
for every $(p,q)$ occurring in the decomposition in (\ref{Gr-dgr}).
\item[5)]
For every integer $k \in \{0,\ldots, \LG-1 \}$ the endomorphism $gr(d^{+}_{\diamond})$ induces the isomorphism
$$
(gr(d^{+}_{\diamond}))^k : Gr^{\LG-1-k}_{{\cal W}^{\bullet}}  (\FTP) \longrightarrow Gr^{\LG-1+k}_{{\cal W}^{\bullet}}  (\FTP([Z]))
$$
which is of type $(k,k)$ with respect to the bigrading in 3). In particular, one has isomorphisms
$$
(gr(d^{+}_{\diamond}))^k : Gr^{p,q}_{{\cal W}^{\bullet}}  (\tau^{\ast} \FTP) \longrightarrow Gr^{p+k,q+k}_{{\cal W}^{\bullet}}  (\tau^{\ast} \FTP)\,,
$$
for every $(p,q)$ such that 
 $p+q =\LG -1-k$ and $p,q \geq 0$.
\end{enumerate}
The data of $\tau^{\ast} \FTP$ together with $d^{+}_{\diamond}$ and its weight filtration
${\cal W}^{\bullet}$, the orthogonal decomposition in (\ref{ord-tau1}), the associated sheaf
$Gr_{{\cal W}^{\bullet}} (\tau^{\ast} \FTP) $ in 3) with its bigrading in (\ref{Gr-dgr}) and the endomorphism
$gr(d^{+}_{\diamond})$ will be called {\rm ${\bf sl_2}$-structure of $\FTP$ associated to $d^{+}$} or
{\rm the natural positive ${\bf sl_2}$-structure of $\FTP$}. 
\end{pro-defi}
\begin{rem}
The terminology `natural' and `positive' refers implicitly to a possibility of varying our construction. This is obvious for `positive', since replacing $d^{+}$ by $d^{-}$, one arrives to the {\rm{negative}} ${\bf sl_2}$-structure on $\FTP$ (see \S\ref{sec-inv} for details).

The adjective `natural' comes from the fact that we can modify $d^{+}_{\diamond}$ by introducing parameters in the definition of $d^{+}$ (resp.
$d^{-}$). Namely, one considers the morphism
\BEN\label{d+-deform}
d^{+}_{\bf z} : {\cal T}_{\pi} \longrightarrow \ENDO(\FT)\,,
\EEN
where ${\bf z} =(z_0,z_1, \ldots, z_{\LG-2} ) \in \CC^{\LG-1}$ and $d^{+}_{\bf z} = \sum^{\LG-2}_{p=0} z_p d^{+}_p$. Thus one has `moduli' of
positive (resp. negative) ${\bf sl_2}$-structures on $\FTP$. 

Besides this simple variation of $d^{\pm}_{\diamond}$ one can also vary it by ``turning" on the action of the neutral summand
$\mbox{\BM$\GS^0_{\GA}$}$ in the decomposition of $\mbox{\BM$\GS_{\GA}$}$ in (\ref{sLie-grad}). More precisely, for a vertical tangent vector
$v$ at a point $\ZA \in \JABG$, instead of taking
$d^{\pm} (v)$ one can take
$ d^{\pm} (x_v (v))$, where $x_v$ is an element of $\mbox{\BM$\GS^0_{\GA}$}\ZA$ naturally associated to $v$ (e.g. there is a canonical way 
to lift $v$ to an element $\tilde{v} \in \HT\ZA$, so one could take $x_v = D^{0} (\tilde{v})$). Once the neutral element $x_v$ is chosen,
all edges of the graph in (\ref{graph}) are colored by operators in $\mbox{\BM$\GS_{\GA}$}\ZA$ - the vertical edges are colored by $x_v$ and
the right (resp. left) handed arrows in (\ref{graph}) are colored by $d^{+} (v)$ (resp. $d^{-} (v)$). This allows us to associate operators
with every path of the graph in (\ref{graph}). This is what we called path-operators in \RI, \S6. In particular, if a path is a loop, then the corresponding operator is of degree $0$ and we can modify our tangent vector $v$ by the action of such loop-operators. Thus ${\bf sl_2}$-structures of $\FTP$ have not only continuous parameters but also the discrete ones indexed by loops of the trivalent graph in (\ref{graph}).
We will address these variational aspects of ${\bf sl_2}$-structures elsewhere.
\end{rem}  

\section{${\bf sl_2}$-structures on {\BM$\LAG$}}\label{sec-sl2bis}

The construction of ${\bf sl_2}$-structure in the previous section could have been done for any locally free sheaf of graded modules equipped
with a compatible graded action of the sheaf of graded Lie algebras
{\BM
$\LAG$, where the grading on
$\LAG$}
is as described in \S\ref{sec-gr}, (\ref{sLie-grad}). We have chosen to do it for the sheaf $\FTP$, rather than to give a functorial treatment,
since this sheaf arises naturally, it comes with the natural {\it graded} action of 
{\BM
$\LAG$}
and this action is closely tied with the geometry of $X$. However, there is another natural choice - the sheaf of graded Lie algebras
{\BM
$\LAG$}
itself, equipped with the adjoint action. This object is `principal' in the category  of graded sheaves of modules with {\it graded}
{\BM
$\LAG$}-action in a sense that all other objects in this category as well as their ${\bf sl_2}$-structures are representations of this one.
Thus it is important to have a good description of ${\bf sl_2}$-structures on 
{\BM
$\LAG$}. This is the subject of the present section.

We begin by recall the grading
{\BM
\BEN\label{sLie-grad1}
\LAG = \bigoplus^{\LG - 1}_{i=-(\LG-1)} \GS^i_{\mbox{\UB$\GA$}}
\EEN}
 introduced in \S\ref{sec-gr}, (\ref{sLie-grad}).
\begin{rem}\label{tw-ord}
This grading can be viewed as {\rm{twisted}} orthogonal decomposition of 
{\BM
$\LAG$}
in the following sense.

The sheaf 
{\BM
$\LAG$
can be equipped with the quadratic form determined by the fibrewise Killing form.  We will call it the Killing form of
$\LAG$. With respect to this form
the summands 
$\GS^i_{\mbox{\UB$\GA$}}$ and 
$\GS^j_{\mbox{\UB$\GA$}}$ are orthogonal, unless $j=-i$. In this, latter case, the Killing form induces a non-degenerate pairing
\BEN\label{Kf}
\GS^i_{\mbox{\UB$\GA$}} \times \GS^{-i}_{\mbox{\UB$\GA$}} \longrightarrow  \mbox{\UB$\OO_{\JABG}$}\,.
\EEN
Next recall that
$\LAG$ carries an involution
\BEN\label{inv}
(\cdot)^{\dagger}: \LAG \longrightarrow \LAG
\EEN 
defined by taking the adjoint $x^{\dagger}$ of a local section $x$ of $\LAG$.

This involution interchanges the summands in (\ref{sLie-grad1}) of opposite sign, i.e. one has
\BEN\label{gp-inv}
 (\GS^i_{\mbox{\UB$\GA$}})^{\dagger}= \GS^{-i}_{\mbox{\UB$\GA$}}\,,
\EEN}
for every $i \in \{-(\LG-1), \ldots, (\LG-1) \}$. Thus the decomposition (\ref{sLie-grad1}) is orthogonal, up to the twist by the involution
$(\cdot)^{\dagger}$.
\end{rem}

Our construction of a natural (positive) ${\bf sl_2}$-structure on 
{\BM
$\LAG$ follows the same steps as the considerations in \S\ref{const-sh}. 

We begin by taking the pullback 
$\mbox{\UB$\tau^{\ast}$}\LAG$ of the sheaf $\LAG$ }
via the natural projection
$\tau$ in (\ref{tau}). As it was pointed out in \S\ref{const-sh}, the sheaf 
{\BM
$\mbox{\UB$\tau^{\ast}$}\LAG$}
 comes together with the tautological section $d^{+}_{\diamond}$. 
Then we go further by lifting 
{\BM
$\LAG$}
 from $\JABG$ to ${\bf h^0}(d^{+})$, the scheme defined in (\ref{sc-sl2}).

The sheaf
{\BM
$\mbox{\UB$(\tau^{+})^{\ast}$}\LAG$
 (see (\ref{tau+}), for the definition of $\tau^{+}$) continues to have the twisted orthogonal decomposition
\BEN\label{ord-tau-g}
\mbox{\UB$(\tau^{+})^{\ast}$}\LAG = \bigoplus^{\LG - 1}_{i=-(\LG-1)} \mbox{\UB$(\tau^{+})^{\ast}$} \GS^i_{\mbox{\UB$\GA$}}\,.
\EEN}
But in addition, it is equipped with
the `universal' ${\bf sl_2}$-triple 
$\{ s^{+}, s^{0}, s^{-} \}$, where the  sections 
$s^{\pm},\, s^0$ are defined in (\ref{sl-univ}).
Observe that these sections are the sections of the summands 
{\BM
$\mbox{\UB$(\tau^{+})^{\ast}$}\LAG^{\pm1}$ and $\mbox{\UB$(\tau^{+})^{\ast}$}\LAG^0$}
respectively.
We consider the decomposition of 
{\BM
$\mbox{\UB$(\tau^{+})^{\ast}$}\LAG$}
under the adjoint action of this ${\bf sl_2}$-triple.

The operator $ad(s^0)$ acts semisimply and gives the decomposition
{\BM
\BEN\label{w-g}
\mbox{\UB$(\tau^{+})^{\ast}$}\LAG = \bigoplus^{2(\LG-1)}_{n=-2(\LG-1)} \LAG (n)
\EEN
into the weight subsheaves
$\LAG (n)$}
of $ad(s^0)$. The range for the weights in (\ref{w-g}) comes from the fact 
\BEN\label{i-n}
(ad(s^{+}))^{2\LG-1} =0
\EEN
which in turn is the result of $s^{+}$ being a section of 
{\BM
$\mbox{\UB$(\tau^{+})^{\ast}$}\LAG^1$ and the grading in (\ref{ord-tau-g}).

The weight decomposition (\ref{w-g}) together with the twisted orthogonal decomposition in 
(\ref{ord-tau-g}) gives rise to a bigrading on 
$\mbox{\UB$(\tau^{+})^{\ast}$}\LAG$
\BEN\label{tau+bigr}
\mbox{\UB$(\tau^{+})^{\ast}$}\LAG = \bigoplus_{(p,q)} \LAG^{p,q}\,,
\EEN
where
$\LAG^{p,q} = \mbox{\UB$(\tau^{+})^{\ast}$}\LAG^p \bigcap \LAG(p+q)$.}
The result below establishes the range for $(p,q)$ of this bigrading.
\begin{lem}\label{range-p,q}
\begin{enumerate}
\item[1)]
For every $p$ in the twisted orthogonal decomposition in (\ref{ord-tau-g}) the following holds
{\BM
$$
\mbox{\UB$(\tau^{+})^{\ast}$}\LAG^p = \bigoplus^{\LG-1}_{q=-(\LG-1)} \LAG^{p,q}\,.
$$
\item[2)]
The direct sum
$$
\mbox{\UB$(\tau^{+})^{\ast}$}\LAG = \bigoplus_{(p,q)} \LAG^{p,q}
$$}
is taken over the points 
$(p,q) \in {\bf Z^2}$
 lying in the square with vertices
$(-(\LG-1),-(\LG-1)),\,\,((\LG-1),-(\LG-1)),\,\,((\LG-1), (\LG-1)), (-(\LG-1), (\LG-1))$.
\item[3)]
For every $n\in \{-2(\LG-1), \ldots, 2(\LG-1) \}$, the weight-subsheaf
{\BM
$\LAG(n)$ in (\ref{w-g}) decomposes as follows:
$$
\LAG(n) = \bigoplus_{p+q=n \atop 0\leq \mid p \mid,\,\mid q \mid \leq \LG-1} \LAG^{p,q}\,.
$$}
\end{enumerate}
\end{lem}
\begin{pf}
From (\ref{sh-Hp-dgr}) we know that the weights of $s^0$ on
$(\tau^{+})^{\ast} \HH^m$ are in the set
\BEN\label{sp-Hm}
\{ m-(\LG-1),\ldots, m \}\,.
\EEN
A local section of 
{\BM
$\mbox{\UB$(\tau^{+})^{\ast}$}\LAG^p$}
takes 
$(\tau^{+})^{\ast} \HH^m$ 
to
$(\tau^{+})^{\ast} \HH^{m+p}$.
Thus the weights are shifted from the set in (\ref{sp-Hm}) to the set
$\{ m+p-(\LG-1),\ldots, m+p \}$. This implies that the possible weights of 
$ad(s^0)$ on $\mbox{\UB$(\tau^{+})^{\ast}$}\LAG^p$ are in the set
\BEN\label{sp-gp}
\{p-(\LG-1), \ldots, p+(\LG-1) \}\,.
\EEN
Thus the range for the index $q$ in
{\BM
$\LAG^{p,q} = \mbox{\UB$(\tau^{+})^{\ast}$}\LAG^p \bigcap \LAG(p+q)$}
is in the set
$\{-(\LG-1), \ldots, (\LG-1) \}$ as asserted.

Combining the first assertion with the range of the grading in (\ref{ord-tau-g})
we deduce the second statement. The third assertion follows from the definition of the bigraded pieces
{\BM
$\LAG^{p,q}$}
and 2) of the lemma.
\end{pf}

Next proposition summarizes some basic properties of the double grading of 
{\BM
$
\mbox{\UB$(\tau^{+})^{\ast}$}\LAG
$
in (\ref{tau+bigr}).
\begin{pro}\label{Lie-bigr}
\begin{enumerate}
\item[1)]
For every $(p,q)$ the summand
$\LAG^{p,q}$ in (\ref{tau+bigr}) acts on
$\mbox{\UB$(\tau^{+})^{\ast}\FTP$}$ by endomorphisms of type 
$(p,q)$ with respect to the bigrading of 
$\mbox{\UB$(\tau^{+})^{\ast}\FTP$}$ in (\ref{bi-gr-tau}), i.e. one has
$$
\LAG^{p,q} \otimes \HH^{s,t} \longrightarrow \HH^{s+p,t+q}\,,
$$
for every $(s,t)$ in (\ref{bi-gr-tau}).
\item[2)]
The bigrading in (\ref{tau+bigr}) turns $(\tau^{+})^{\ast} \LAG$ into a sheaf of bigraded Lie algebras, i.e.
the Lie algebra bracket in $(\tau^{+})^{\ast} \LAG$ satisfies the following relations
$$
[ \LAG^{p,q}, \LAG^{p^{\prime},q^{\prime}}] \subset  \LAG^{p+p^{\prime},q+q^{\prime}}\,,
$$
for all $(p,q)$ and $(p^{\prime},q^{\prime})$ in (\ref{tau+bigr}).
\item[3)]
Summands 
$ \LAG^{p,q}$ and $ \LAG^{p^{\prime},q^{\prime}}$ are orthogonal with respect to the Killing form on 
$(\tau^{+})^{\ast} \LAG$, unless $(p^{\prime},q^{\prime}) =- (p,q)$. In this case the pairing
$$
\LAG^{p,q} \times \LAG^{-p,-q} \longrightarrow \mbox{\UB$\OO_{{\bf h^0}(d^{+})}$}\,,
$$
induced by the Killing form, is non-degenerate and defines an isomorphism
$$
\LAG^{-p,-q} \cong (\LAG^{p,q} )^{\ast}\,.
$$
\end{enumerate}
\end{pro}
\begin{pf}
The properties 1) and 2) are immediate from the definitions of the double gradings in question. The property 3)
follows from the (twisted) orthogonality of the decompositions in (\ref{ord-tau-g}) and (\ref{w-g}).
\end{pf}
}
From the definition of the bigrading in (\ref{tau+bigr}) it also follows that
$s^{+}$ (resp. $s^{0}$ and $s^{-}$) is a section of 
{\BM
$\LAG^{1,1}$
(resp. $\LAG^{0,0}$ and $\LAG^{-1,-1}$).
}
This fact and other basic properties of $ad(s^{+})$ following from the theory of ${\bf sl_2}$-representations are recorded below.
\begin{pro}\label{pro-ads+}
\begin{enumerate}
\item[1)]
The adjoint action of $s^{+}$ (resp. $s^{0}$ and $s^{-}$) induces the morphism
$$
ad(s^{+})\,\, (resp. \,\,ad(s^{0}),\,ad(s^{-}))\, : (\tau^{+})^{\ast} \mbox{\BM$\GS$}_{\GA} \longrightarrow
  (\tau^{+})^{\ast} \mbox{\BM$\GS$}_{\GA}
$$
which has type $(1,1)$ (resp. $(0,0)$ and $(-1,-1)$) with respect to the bigrading in (\ref{tau+bigr}).
\item[2)]
For every integer $p \in [0, 2(\LG-1)]$, one has an isomorphism
$$
(ad(s^{+}))^p : \mbox{\BM$\GS$}_{\GA} (-p) \longrightarrow \mbox{\BM$\GS$}_{\GA} (p)
$$
which is of type $(p,p)$ with respect to the bigrading in Lemma \ref{range-p,q}, 3). In particular, it induces isomorphisms
on the bigraded components
{\BM
$$
\mbox{\UB$(ad(s^{+}))^p$} : \LAG^{s,-s-p}  \longrightarrow \LAG^{s+p,-s}\,, 
$$}
for every integer $s \in [-(\LG-1), (\LG-1)-p]$.
\end{enumerate}
\end{pro}

Next we define the weight filtration associated to the weight decomposition in (\ref{w-g}) by setting
{\BM
\BEN\label{w-filt-eta}
{\hat{\GS}}^{[n]}_{\mbox{\UB$\GA$}} = \bigoplus_{m\geq n} \LAG (m)\,.
\EEN}
These subsheaves are sheaves of Lie subalgebras of 
{\BM
$\mbox{\UB$(\tau^{+})^{\ast}$} \LAG$}. Furthermore,
as we observed in \S\ref{const-sh} these sheaves are constant along the fibres of the projection $\eta^{+}_1$
in (\ref{pr-1}) and hence they descend to $\TPI$ to give {\it the weight filtration of
{\BM
$\mbox{\UB$\tau^{\ast}$}\LAG$} associated to} $d^{+}_{\diamond}$:
{\BM
\BEN\label{w-filt-g}
\mbox{\UB$\tau^{\ast}$}\LAG =\GS^{[-2(\LG-1)]}_{\mbox{\UB$\GA$}} \supset \ldots \supset \GS^{[n]}_{\mbox{\UB$\GA$}}\supset 
\GS^{[n+1]}_{\mbox{\UB$\GA$}}
\supset \ldots
 \supset \GS^{[2(\LG-1)]}_{\mbox{\UB$\GA$}} \supset \GS^{[2\LG-1]}_{\mbox{\UB$\GA$}} =0\,,
\EEN
where the subsheaves $\GS^{[n]}_{\mbox{\UB$\GA$}}$}'s above are such that their pullback by $\eta^{+}_1$ gives the subsheaves
{\BM
${\hat{\GS}}^{[n]}_{\mbox{\UB$\GA$}}$ in (\ref{w-filt-eta}):
\BEN\label{w-filt-eta=pb}
\mbox{\UB$(\eta^{+}_1)^{\ast}$} \GS^{[n]}_{\mbox{$\GA$}} = {\hat{\GS}}^{[n]}_{\mbox{\UB$\GA$}} = \bigoplus_{m\geq n} \LAG(m)\,.
\EEN

The following properties of the filtration
$\GS^{[\bullet]}_{\mbox{\UB$\GA$}}$ }
in (\ref{w-filt-g}) are immediate from its definition and properties of the weight decomposition
in (\ref{w-g}).
{\BM
\begin{lem}\label{w-filt-g1}
\begin{enumerate}
\item[1)]
Each sheaf $\GS^{[n]}_{\mbox{\UB$\GA$}}$ in the filtration (\ref{w-filt-g}) is a subsheaf of Lie subalgebras of
$\mbox{\UB$\tau^{\ast}$}\LAG$.
\item[2)]
The Lie bracket in $\mbox{\UB$\tau^{\ast}$}\LAG$ satisfies the following relations
$$
[\GS^{[m]}_{\mbox{\UB$\GA$}}, \GS^{[n]}_{\mbox{\UB$\GA$}}] \subset \GS^{[m+n]}_{\mbox{\UB$\GA$}}\,,
$$
for all $m,n$.
\end{enumerate}
\end{lem}
}
Define the associated graded sheaf
{\BM
\BEN\label{g-gr}
Gr(\mbox{\UB$\tau^{\ast}$}\LAG) = \bigoplus^{2(\LG-1)}_{n=- 2(\LG-1)} Gr^n(\mbox{\UB$\tau^{\ast}$}\LAG)\,,
\EEN
where $ Gr^n(\mbox{\UB$\tau^{\ast}$}\LAG) = \GS^{[n]}_{\mbox{$\GA$}} / {\GS^{[n+1]}_{\mbox{$\GA$}}}$.
\begin{pro}\label{pro-g-gr}
\begin{enumerate}
\item[1)]
  $Gr(\mbox{\UB$\tau^{\ast}$}\LAG)$ is a sheaf of graded Lie algebras.
\item[2)]
The Killing form on 
$\mbox{\UB$\tau^{\ast}$}\LAG$ induces a non-degenerate quadratic form on 
 $Gr(\mbox{\UB$\tau^{\ast}$}\LAG)$ which coincides with the Killing form on
$Gr(\mbox{\UB$\tau^{\ast}$}\LAG)$. In particular, the direct sum in (\ref{g-gr}) is a twisted orthogonal decomposition
in the sense of Remark \ref{tw-ord}.
\end{enumerate}
\end{pro}
}
\begin{pf}
{\BM
Lemma \ref{w-filt-g1}, 2), implies that the Lie bracket on 
$\mbox{\UB$\tau^{\ast}$}\LAG$ descends to
$Gr(\mbox{\UB$\tau^{\ast}$}\LAG)$.

The part 2) of the proposition follows readily from 1). Indeed,
let $x$ be a local section of 
$Gr^n(\mbox{\UB$\tau^{\ast}$}\LAG)$, for some integer $n\in [-2(\LG-1),2(\LG-1)]$,
and let $\tilde{x}$ be an arbitrary lifting of $x$ to a local section of $\GS^{[n]}_{\mbox{\UB$\GA$}}$.
Then $ad( \tilde{x})$ induces an action on the graded sheaf
$Gr(\mbox{\UB$\tau^{\ast}$}\LAG)$ which, by 1) of the proposition, coincides with $ad(x)$.
This implies that the Killing form on
$\mbox{\UB$\tau^{\ast}$}\LAG$ descends to the Killing form on
$Gr(\mbox{\UB$\tau^{\ast}$}\LAG)$.

To see the twisted orthogonality of the grading in (\ref{g-gr}) we take local sections 
$x$ and $y$ of 
$Gr^m(\mbox{\UB$\tau^{\ast}$}\LAG)$ and
$Gr^n(\mbox{\UB$\tau^{\ast}$}\LAG)$
respectively. Then $ad(x)\circ ad(y)$ shifts the degree of the grading in (\ref{g-gr}) by $(m+n)$. In particular,
$ad(x)\circ ad(y)$ is nilpotent, unless $m+n =0$. This and the definition of the Killing form imply
$$
\langle x,y \rangle=Tr(ad(x)\circ ad(y)) =0\,,
$$
for all local sections $x$ and $y$ as above, unless $m+n =0$. 
In this latter case, for $x \neq 0$ choose a lifting $\tilde{x}$ of $x$ to a local section of
$\GS^{[m]}_{\mbox{\UB$\GA$}}$. Then we know that $\tilde{x}$ is not a local section of $\GS^{[m+1]}_{\mbox{\UB$\GA$}}$. 
Then from Proposition \ref{Lie-bigr}, 3) and 4), it follows
that there exists a local section
$\tilde{y}$ of $\GS^{[-m]}_{\mbox{\UB$\GA$}}$ such that 
$$
\langle\tilde{x}, \tilde{y} \rangle \neq 0\,,
$$
where $\langle \cdot, \cdot \rangle$ stands for the pairing given by the Killing form in
$\mbox{\UB$\tau^{\ast}$}\LAG$. Let $y$ be the image of $\tilde{y}$ under the natural projection
$$
 \GS^{[-m]}_{\mbox{\UB$\GA$}} \longrightarrow  \GS^{[-m]}_{\mbox{\UB$\GA$}} /{\GS^{[-m+1]}_{\mbox{\UB$\GA$}}} =
 Gr^{-m}(\mbox{\UB$\tau^{\ast}$}\LAG)\,.
$$
Then from the first part of 2) we deduce
$$
\langle x, y \rangle =\langle\tilde{x}, \tilde{y} \rangle \neq 0
$$
which proves that the pairing
$$
Gr^{m}(\mbox{\UB$\tau^{\ast}$}\LAG) \times Gr^{-m}(\mbox{\UB$\tau^{\ast}$}\LAG) \longrightarrow \mbox{\UB$\OO_{\TPI}$}
$$
induced by the Killing form in 
$\mbox{\UB$\tau^{\ast}$}\LAG$}
is non-degenerate.
\end{pf}

{\BM
Next we take into account the twisted orthogonal decomposition of 
$\mbox{\UB$\tau^{\ast}$}\LAG$ 
\BEN\label{sLie-grad2}
\mbox{\UB$\tau^{\ast}$}  \LAG = \bigoplus^{\LG - 1}_{i=-(\LG-1)} \mbox{\UB$\tau^{\ast}$} \GS^i_{\mbox{\UB$\GA$}}
\EEN
obtained by applying $\mbox{\UB$\tau^{\ast}$}$ to the decomposition in (\ref{sLie-grad1}). This together with
the filtration 
$\GS^{[\bullet]}_{\mbox{\UB$\GA$}}$ in (\ref{w-filt-g}) define
a bigrading on 
$Gr(\mbox{\UB$\tau^{\ast}$}\LAG)$.
Namely, set
\BEN\label{gr-g-bigr}
Gr^{p,q}(\mbox{\UB$\tau^{\ast}$}\LAG) = \mbox{\UB$\tau^{\ast}$}\LAG^p \bigcap \GS^{[p+q]}_{\mbox{\UB$\GA$}} /
                       \mbox{\UB$\tau^{\ast}$}\LAG^p \bigcap \GS^{[p+q+1]}_{\mbox{\UB$\GA$}} \,.
\EEN
From Lemma \ref{range-p,q} we obtain that
$Gr(\mbox{\UB$\tau^{\ast}$}\LAG)$ admits the following bigrading
\BEN\label{gr-g-bigr1}
Gr(\mbox{\UB$\tau^{\ast}$}\LAG) = \bigoplus_{p,q} Gr^{p,q}(\mbox{\UB$\tau^{\ast}$}\LAG)\,,
\EEN}
where the points $(p,q)\in {\bf Z^2} $ lie in the square described in Lemma \ref{range-p,q}, 2).

Recall that the sheaf
{\BM
$\mbox{\UB$\tau^{\ast}$}\LAG$}
is equipped with the section
 $d^{+}_{\diamond}$ and  
its pullback via $\eta^{+}_1$ in (\ref{pr-1}) gives the section $s^{+}$ of 
{\BM
$\mbox{\UB$(\tau^{+})^{\ast}$}\LAG$} (see Remark \ref{s+=pb}). 
This implies that all the properties of
 $s^{+}$ described in Proposition \ref{pro-ads+} have their analogues for $d^{+}_{\diamond}$.
We summarize this together with the preceding discussion below.
\begin{pro-defi}\label{pro-gr-g-bigr}
The sheaf 
{\BM
$\mbox{\UB$\tau^{\ast}$}\LAG$}
comes along with a distinguished section $d^{+}_{\diamond}$ which gives rise to the following structures.
\begin{enumerate}
\item[1)]
$d^{+}_{\diamond}$ defines the weight filtration
{\BM
$\GS^{[\bullet]}_{\mbox{$\GA$}}$
of
$\mbox{\UB$\tau^{\ast}$}\LAG$ as in (\ref{w-filt-g}) and which has the properties described in Lemma \ref{w-filt-g1}.
Furthermore, $\mbox{\UB$d^{+}_{\diamond}$}$ is a section of $\GS^{[2]}_{\mbox{\UB$\GA$}}$}
 and its adjoint action  
  shifts the index of the filtration by $2$, i.e. one has
{\BM
$$
\mbox{\UB$ad ( d^{+}_{\diamond})$}(\GS^{[n]}_{\mbox{\UB$\GA$}}) \subset \GS^{[n+2]}_{\mbox{\UB$\GA$}}\,.
$$
\item[2)]
The filtration $\GS^{[\bullet]}_{\mbox{\UB$\GA$}}$ together with the twisted orthogonal decomposition
of $\mbox{\UB$\tau^{\ast}$}\LAG$} in (\ref{sLie-grad2}) gives rise to the bigrading
{\BM
$$
Gr(\mbox{\UB$\tau^{\ast}$}\LAG) = \bigoplus_{p,q} Gr^{p,q}(\mbox{\UB$\tau^{\ast}$}\LAG)
$$
of the associated graded sheaf
$Gr(\mbox{\UB$\tau^{\ast}$}\LAG) $, where 
$Gr^{p,q}(\mbox{\UB$\tau^{\ast}$}\LAG)$ are defined in (\ref{gr-g-bigr}) and the direct sum is taken over
the points 
$(p,q)\in {\bf Z^2} $ lying in the square described in Lemma \ref{range-p,q}, 2). This bigrading turns 
$Gr(\mbox{\UB$\tau^{\ast}$}\LAG)$ into a sheaf of bigraded Lie algebras with respect to the Lie bracket induced
from
$\mbox{\UB$\tau^{\ast}$}\LAG$}.
\item[3)]
Let 
$\overline{d^{+}_{\diamond}}$ be the image of $d^{+}_{\diamond}$ under the natural projection
{\BM
 $$
\GS^{[2]}_{\mbox{\UB$\GA$}} \longrightarrow {\GS^{[2]}_{\mbox{\UB$\GA$}}}/{\GS^{[3]}_{\mbox{\UB$\GA$}}}\,.
$$
}
Then $\overline{d^{+}_{\diamond}}$ is a section of
{\BM 
$Gr^{1,1}(\mbox{\UB$\tau^{\ast}$}\LAG)$}
 and its adjoint action has the following properties.
\begin{enumerate}
\item[(i)]
{\BM
$\mbox{\UB$ad(\overline{d^{+}_{\diamond}})$}(Gr^{p,q}(\mbox{\UB$\tau^{\ast}$}\LAG)) \subset Gr^{p+1,q+1}(\mbox{\UB$\tau^{\ast}$}\LAG)$\,,
for every $(p,q)$ occurring in the double graded direct sum in 2).
\item[(ii)]
For every integer $k\in [0,2(\LG-1)]$, one has the isomorphism
$$
\mbox{\UB$(ad(\overline{d^{+}_{\diamond}}))^k$} : Gr^{-k}(\mbox{\UB$\tau^{\ast}$}\LAG) \longrightarrow Gr^k(\mbox{\UB$\tau^{\ast}$}\LAG)
$$
which has type $(k,k)$ with respect to the bigrading in 2). In particular, it induces isomorphisms
$$
\mbox{\UB$(ad(\overline{d^{+}_{\diamond}}))^k$} : Gr^{p,q}(\mbox{\UB$\tau^{\ast}$}\LAG) \longrightarrow
 Gr^{p+k,q+k} (\mbox{\UB$\tau^{\ast}$}\LAG)\,,
$$
}
for every $(p,q)$ such that $p+q=-k$ and $p \in \{-(\LG-1),\ldots, (\LG-1) -k \}$.
\end{enumerate}
\end{enumerate}

The data of
{\BM 
$\mbox{\UB$\tau^{\ast}$}\LAG$ together with the section $\mbox{\UB$d^{+}_{\diamond}$}$ and its weight filtration
$\GS^{[\bullet]}_{\mbox{$\GA$}}$, the twisted orthogonal decomposition in (\ref{sLie-grad2}), the associated graded sheaf
$Gr(\mbox{\UB$\tau^{\ast}$}\LAG)$
 in (\ref{g-gr}) with its bigrading in 2) and the section
$\mbox{\UB$\overline{d^{+}_{\diamond}}$}$
will be called ${\bf sl_2}$-structure of 
$\LAG$ associated to $\mbox{\UB$d^{+}$}$.
}
\end{pro-defi}

The essential point in defining the ${\bf sl_2}$-structure of 
{\BM
$\LAG$}
 was to view the morphism
$$   
d^{+}: {\cal T}_{\pi} \longrightarrow \mbox{\BM$\GS$}_{\GA}
$$
defined in (\ref{d+}) as a distinguished section of 
$\tau^{\ast} \mbox{\BM$\GS$}_{\GA}$, where 
$$
\tau: \TPI \longrightarrow \JABG
$$
is the natural projection. But of course the morphism contains more information.
Indeed, by going to a point $([Z],[\alpha], v)$ of $\TPI$ we have taken just a single value
$d^{+}(v)$ of $d^{+}$ at $\ZA$. Though clearly intrinsically attached to
$([Z],[\alpha], v)$, it ignores completely all other values of $d^{+}$ at $\ZA$.
We will now remedy this by taking the pullback by $\tau$ of the whole morphism $d^{+}$.
In other words we consider the morphism
\BEN\label{tau-d+}
\tau^{\ast}(d^{+}) :  \tau^{\ast} {\cal T}_{\pi} \longrightarrow \tau^{\ast} \mbox{\BM$\GS$}_{\GA}\,.
\EEN
\begin{rem} \label{geom-taud+}
Besides the above, formal, reason there is also a geometric reason to consider the morphism
$\tau^{\ast}(d^{+})$. This is based on the geometric point of view on the morphism
$d^{+}$. Namely, consider 
$ \mbox{\BM$\GS$}_{\GA}$ as a vector bundle over $\JABG$ and view $d^{+}$
as a map of bundles over $\JABG$. This gives the following commutative diagram
of morphisms
\BEN\label{d+-bundle}
\xymatrix{
{\TPI} \ar[rr]^{d^{+}} \ar[dr]_{\tau} & &{ \mbox{\BM$\GS$}_{\GA}} \ar[dl]^{\gamma} \\
           &{\JABG} }
\EEN
where $\tau$ and $\gamma$ are the natural projections.

Consider now the relative (with respect to the projection $\tau$) differential 
$d_{\tau} (d^{+})$ of the morphism
$d^{+}$. This gives the morphism
\BEN\label{rel-diff-tau}
d_{\tau} (d^{+}) : {\cal T}_{\tau} \longrightarrow (d^{+})^{\ast} {\cal T}_{\gamma}\,,
\EEN
where ${\cal T}_{\tau}$ (resp. ${\cal T}_{\gamma}$) is the relative tangent sheaf of $\tau$ (resp. $\gamma$).
Since $\TPI$ and $ \mbox{\BM$\GS$}_{\GA}$ are vector bundles over $\JABG$ and $\tau$ (resp. $\gamma$)
is its natural projection we have
\BEN\label{rel-tan=pb}
{\cal T}_{\tau} = \tau^{\ast}({\cal T}_{\pi}) \,\,(\mbox{resp.} \,\, {\cal T}_{\gamma} =\gamma^{\ast}( \mbox{\BM$\GS$}_{\GA}))
\EEN
and the morphism in (\ref{rel-diff-tau}) becomes
\BEN\label{rel-diff-tau1}
d_{\tau} (d^{+}) : \tau^{\ast}({\cal T}_{\pi}) \longrightarrow (d^{+})^{\ast}(\gamma^{\ast}( \mbox{\BM$\GS$}_{\GA})) =
                        \tau^{\ast}(\mbox{\BM$\GS$}_{\GA})\,,
\EEN
where the last equality comes from the commutativity of the diagram in (\ref{d+-bundle}). Finally, since 
$d^{+}$ is linear along the fibres of the projections, this is nothing but the morphism
$ \tau^{\ast} (d^{+})$ in (\ref{tau-d+}). Thus we have
\BEN\label{rel-diff-tau=pb}
d_{\tau} (d^{+}) = \tau^{\ast} (d^{+})\,.
\EEN
In view of this equality we will denote this morphism by $d^{+}_{\tau}$.
\end{rem}

We have seen that the sheaf $\tau^{\ast}(\mbox{\BM$\GS$}_{\GA})$ together with the section $d^{+}_{\diamond}$ carries a rich
structure described in Proposition-Definition \ref{pro-gr-g-bigr}. A part of this structure is the weight filtration 
in (\ref{w-filt-g}).
Using the morphism
$d^{+}_{\tau}$ this filtration can be transplanted to $\tau^{\ast}({\cal T}_{\pi})$.  This becomes especially pertinent
in case of $\GA$ being simple, since in this case the Infinitesimal Torelli property holds (see \S\ref{sec-Torelli},
Corollary \ref{InfTorelli=s}) and it makes the morphism $d^{+}_{\tau}$  injective ( see Proposition \ref{InfTor-equiv}).
Thus, in a nut-shell, we have that
 the relative tangent space 
$\TPI \ZA$ of $\JABG$ at a point 
$\ZA \in \JABG$ acquires a filtration intrinsically associated with every tangent vector
$v \in  \TPI \ZA$ or, even better, with every point of the projective space $\PP(\TPI \ZA)$. 
In what follows we establish some basic properties of these filtrations.

First we recall that the morphism $d^{+}_{\tau}$ takes values in the summand
{\BM
$\mbox{\UB$\tau^{\ast}$}\GS^1_{\mbox{\UB$\GA$}}$ of the twisted orthogonal decomposition in (\ref{sLie-grad2}).
The weight filtration 
$\GS^{[\bullet]}_{\mbox{\UB$\GA$}}$ in (\ref{w-filt-g}) induces the weight filtration on 
$\mbox{\UB$\tau^{\ast}$}\GS^1_{\mbox{\UB$\GA$}}$. Due to Lemma \ref{range-p,q}, 1), this gives the following
\BEN\label{w-filt-g-d1} 
\mbox{\UB$\tau^{\ast}$}\GS^1_{\mbox{\UB$\GA$}} = \GS^{1,[2-\LG]}_{\mbox{\UB$\GA$}} \supset \GS^{1,[3-\LG]}_{\mbox{\UB$\GA$}} \supset
    \ldots \supset \GS^{1,[\LG]}_{\mbox{\UB$\GA$}} \supset \GS^{1,[\LG +1]}_{\mbox{\UB$\GA$}} =0\,,
\EEN
where 
$ \GS^{1,[n]}_{\mbox{\UB$\GA$}} = \mbox{\UB$\tau^{\ast}$}\GS^1_{\mbox{\UB$\GA$}} \bigcap  \GS^{[n]}_{\mbox{\UB$\GA$}}$.
} 
Transferring this filtration to
$\tau^{\ast}({\cal T}_{\pi})$ via the morphism 
$d^{+}_{\tau}$ gives the filtration
${\cal W}^{\bullet}(\tau^{\ast}({\cal T}_{\pi}), d^{+})$ on $\tau^{\ast}({\cal T}_{\pi})$, where one sets
{\BM
\BEN\label{Wk-TPI}
\mbox{\UB${\cal W}^{n}(\tau^{\ast}({\cal T}_{\pi}), d^{+}) = (d^{+}_{\tau} )^{-1}$} (\GS^{1,[n]}_{\mbox{\UB$\GA$}})\,.
\EEN}
\begin{defi}\label{w-filt-TPI}
The filtration ${\cal W}^{\bullet}(\tau^{\ast}({\cal T}_{\pi}), d^{+})$ will be called the weight filtration of 
$\tau^{\ast}({\cal T}_{\pi})$ associated to $d^{+}$ or simply, the $d^{+}$-filtration of 
$\tau^{\ast}({\cal T}_{\pi})$.
\end{defi}
\begin{pro}
The $d^{+}$-filtration of 
$\tau^{\ast}({\cal T}_{\pi})$ has the following form:
\BEN\label{pro-w-filt-TPI}
\tau^{\ast}({\cal T}_{\pi}) ={\cal W}^0(\tau^{\ast}({\cal T}_{\pi}), d^{+}) \supset \ldots 
\supset {\cal W}^{\LG}(\tau^{\ast}({\cal T}_{\pi}), d^{+}) \supset {\cal W}^{\LG+1}(\tau^{\ast}({\cal T}_{\pi}), d^{+})=0\,.
\EEN
\end{pro}
\begin{pf}
The assertion is equivalent to saying that the image of 
$d^{+}_{\tau}$ is contained in
{\BM 
$\GS^{1,[0]}_{\mbox{\UB$\GA$}}$
}
of the filtration in (\ref{w-filt-g-d1}). To see this observe that the sheaf
$\tau^{\ast}({\cal T}_{\pi})$ comes with the tautological section which will be denoted by
$i_0$, i.e. for a point $([Z],[\alpha],v) \in \TPI$ the value  $i_0 ([Z],[\alpha],v) =v$.
Consider its image 
$d^{+}_{\tau} (i_0)$. This is now a section of 
{\BM
$\mbox{\UB$\tau^{\ast}$}\GS^1_{\mbox{\UB$\GA$}}$.}
 This is of course the distinguished section $d^{+}_{\diamond}$ of
{\BM
$\mbox{\UB$\tau^{\ast}$}\LAG$}
 defining the ${\bf sl_2}$-structure  
in Proposition-Definition \ref{pro-gr-g-bigr},
so from now on we use the notation $d^{+}_{\diamond}$ instead of $d^{+}_{\tau} (i_0)$. 

We know that $d^{+}_{\diamond}$ is a section
of 
{\BM
$ \GS^{1,[2]}_{\mbox{\UB$\GA$}}$}
 of the weight filtration (\ref{w-filt-g-d1}). Furthermore, the image of
$d^{+}_{\tau}$ is a subsheaf of commutative Lie subalgebras\footnote{recall: (1) by Remark \ref{val-d} the values of  $d^{+}_{\tau}$ are the same as the values of the morphism
$D^{+}$ in the triangular decomposition in (\ref{d-D}), and (2) the image of $D^{+} (\HT)$, by 
\RI, Lemma 7.6, is a subsheaf of abelian Lie subalgebra of 
$\mbox{\BM$\GS$}_{\GA}$.} of 
$\tau^{\ast} \mbox{\BM$\LAG$}$. In particular, the image of
$d^{+}_{\tau}$ commutes with the section
$d^{+}_{\diamond}$. Thus we have an inclusion
\BEN\label{im-ker}
im(d^{+}_{\tau} ) \subset ker (ad(d^{+}_{\diamond}))\,.
\EEN
From the properties of ${\bf sl_2}$-representations and the definition of the filtration
{\BM
$\GS^{[\bullet]}_{\mbox{\UB$\GA$}}$ (see (\ref{w-filt-eta}), (\ref{w-filt-eta=pb})) it follows that
$\mbox{\UB$ker (ad(d^{+}_{\diamond})) \subset$} \GS^{[0]}_{\mbox{\UB$\GA$}}$
}.
\end{pf}

\section{Involution on {\BM$\LAG$}}\label{sec-inv}
In the  two preceding sections we considered
${\bf sl_2}$-triples of 
{\BM
$\LAG$}
arising from nilpotent elements $d^{+}(v)$, for $v$ being tangent vectors in the relative tangent sheaf
${\cal T}_{\pi}$ and $d^{+}$ is the morphism defined in (\ref{d+}). However, it was explicitly indicated in the beginning of
\S\ref{sec-sl2} that we also have a choice of working with elements $d^{-} (v)$, the values of the morphism
$d^{-}$ in (\ref{d-}). Such a choice leads to 
${\bf sl_2}$-triples
$\{e,h,d^{-}(v)\}$, where
$d^{-}(v)$ is nil-negative, $h$ is semisimple and $e$ is nil-positive elements of a triple. In addition, the grading of 
{\BM
$\LAG$
in (\ref{sLie-grad}) allows to take $\mbox{\UB$h$}$ (resp. $\mbox{\UB$e$}$) in 
$\LAG^0$ (resp. $\LAG^1$).} With these choices in mind we go through the constructions analogous to the one's in
\S\ref{const-sh}. There is nothing conceptually new:  the objects have the same meaning as in
\S\ref{const-sh} with only the notational difference - we switch the sign `$+$' to `$-$'. Thus the scheme parametrizing
all ${\bf sl_2}$-triples with values of $d^{-}$ as nil-negative elements is denoted by
${\bf h^0}( d^{-})$. It is defined similar to 
${\bf h^0}( d^{+})$  in (\ref{sc-sl2}), i.e. this is the incidence correspondence in 
$\TPI \times _{\JABG} \mbox{\BM$\GS^0$}_{\GA}$, parametrizing the quadruples
$([Z],[\alpha],v,h)$ such that the triple
$\{e,h,d^{-}(v)\}$ is an  ${\bf sl_2}$-triple with
$d^{-}(v)$ nil-negative, $h \in  \mbox{\BM$\GS^0$}_{\GA} \ZA$ semisimple 
and $e \in \mbox{\BM$\GS^1$}_{\GA} \ZA$ nil-positive elements respectively.
The projection of 
${\bf h^0}( d^{-})$ on the first factor of 
$\TPI \times _{\JABG} \mbox{\BM$\GS^0$}_{\GA}$
will be denoted by $\eta^{-}_1$. Composing it with the projection $\tau$ gives the morphism
\BEN\label{tau-}
\tau^{-} =\tau \circ {\bf \eta^{-}_1} : \mbox{\BM$h^0$}(d^{-})
 \stackrel{{\bf \eta^{-}_1}}{\longrightarrow} \TPI 
\stackrel{\tau}{\longrightarrow} \JABG\,.
\EEN
As in \S\ref{const-sh} we take 
$(\tau^{-} )^{\ast} \mbox{\BM$(\LAG)$}$ and denote by 
$(s^{\prime})^{\pm},\,(s^{\prime})^0$ its three distinguished sections whose values at a closed point 
$([Z],[\alpha],v,h)$ of ${\bf h^0}( d^{-})$ are given by the following identities:
$$
(s^{\prime})^{-} ([Z],[\alpha],v,h) =d^{-} (v),\,\, (s^{\prime})^0 ([Z],[\alpha],v,h)=h
$$
and $ (s^{\prime})^{+} ([Z],[\alpha],v,h)$ is uniquely determined by $h$, $d^{-}(v)$ and the requirement
that
$$
\{(s^{\prime})^{+} ([Z],[\alpha],v,h), h, d^{-}\}
$$ 
is an ${\bf sl_2}$-triple with
$h$ semisimple, $d^{-} (v)$ (resp. $(s^{\prime})^{+} ([Z],[\alpha],v,h)$) nil-negative (resp. positive) elements respectively.

The next step is to consider the defining representation of 
this universal (negative) triple on
$(\tau^{-})^{\ast} \FTP$. This gives the weight decomposition of 
$(\tau^{-})^{\ast} \FTP$ under the action of the semisimple element
$(s^{\prime})^0$:
\BEN\label{wd-d-}
(\tau^{-})^{\ast}\FTP = \bigoplus^{\LG-1}_{n= -(\LG-1)} {\cal W^{\prime}}(n)\,.
\EEN
Applying the conventional shift in grading we obtain
\BEN\label{sh-shifted1}
(\tau^{-})^{\ast}\FTP = \bigoplus^{2(\LG-1)}_{n= 0} {\cal V^{\prime}}(n)\,,
\EEN
where
${\cal V^{\prime}}(n):= {\cal W^{\prime}}(n-\LG+1)$, for $n=0,1, \ldots, 2(\LG-1)$.

The main goal of this section is to relate this weight decomposition with the one in \S\ref{const-sh}, (\ref{sh-shifted}).
The key point here is the involution 
$(\cdot)^{\dagger}$ of taking the adjoint which we have already encountered in (\ref{inv}). In particular, 
the involution switches from $d^{+}$ to $d^{-}$ and vice verse, i.e.
$$
(d^{\pm})^{\dagger} = d^{\mp}\,,
$$
something we have already seen in (\ref{adj1}). This implies that 
$(\cdot)^{\dagger}$ induces an isomorphism
\BEN\label{inv1}
i: {\bf h^0}(d^{+}) \longrightarrow {\bf h^0}(d^{-})
\EEN
given by the formula
\BEN\label{i-formula}
i([Z],[\alpha],v,h) = ([Z],[\alpha],v,h^{\dagger})\,.
\EEN
In particular, $i$ commutes with the projections $\eta^{\pm}_1$ making the following diagram
\BEN\label{diag-i}
\xymatrix{
 {\bf h^0}(d^{+}) \ar[rr]^{i} \ar[dr]_{\eta^{+}_1}& & {\bf h^0}(d^{-}) \ar[dl]^{\eta^{-}_1} \\
& \TPI& }
\EEN
commutative

Applying the involution to the universal ${\bf sl_2}$-triple
$\{ s^{+}, s^{0}, s^{-} \}$ in \S\ref{const-sh} yields another ${\bf sl_2}$-triple
\BEN\label{sl2-adj}
\{ (s^{-})^{\dagger}, (s^{0})^{\dagger}, (s^{+})^{\dagger} \}
\EEN
 of sections of $(\tau^{+})^{\ast} \mbox{\BM$\LAG$}$.
This triple is related to the universal negative triple
$\{( s^{\prime})^{+}, ( s^{\prime})^{0}, ( s^{\prime})^{-} \}$
of sections of 
$(\tau^{-})^{\ast} \mbox{\BM$\LAG$}$ as follows.
\begin{lem}\label{two-sl2}
The ${\bf sl_2}$-triple 
$\{ (s^{-})^{\dagger}, (s^{0})^{\dagger}, (s^{+})^{\dagger} \}$
is equal to the pullback by $i$ in (\ref{inv1}) of the 
${\bf sl_2}$-triple 
$\{( s^{\prime})^{+}, ( s^{\prime})^{0}, ( s^{\prime})^{-} \}$. More precisely, one has
\BEN\label{rel-two-sl2}
i^{\ast}(( s^{\prime})^{\pm}) =(s^{\mp})^{\dagger} \,\,and\,\, i^{\ast}(( s^{\prime})^0) =(s^0)^{\dagger}\,.
\EEN
\end{lem}
\begin{pf}
Let 
$([Z],[\alpha],v,h)$ be a closed point of 
${\bf h^0}(d^{+})$. Then by definition
\BEN\label{lhs}
(s^{0})^{\dagger} ([Z],[\alpha],v,h) = h^{\dagger}\,.
\EEN
On the other hand, evaluating 
$i^{\ast}(( s^{\prime})^0)$ at $([Z],[\alpha],v,h)$ gives
$$
i^{\ast}(( s^{\prime})^0)([Z],[\alpha],v,h) =( s^{\prime})^0(i([Z],[\alpha],v,h))=( s^{\prime})^0 ([Z],[\alpha],v,h^{\dagger}) =h^{\dagger}\,.
$$
Combining this with (\ref{lhs}) yields the second equality in (\ref{rel-two-sl2}).

Turning to the first equality in (\ref{rel-two-sl2}) we have
$$
i^{\ast}(( s^{\prime})^{-})([Z],[\alpha],v,h) =( s^{\prime})^{-} ([Z],[\alpha],v,h^{\dagger}) =d^{-}(v) =(d^{+}(v))^{\dagger} =
(s^{+})^{\dagger}([Z],[\alpha],v,h)\,.
$$
Hence the equality $i^{\ast}(( s^{\prime})^{-}) =(s^{+})^{\dagger}$.
The uniqueness of $(s^{-})^{\dagger}$ (resp. $(s^{\prime})^{+}$) guaranties the remaining equality
$i^{\ast}(( s^{\prime})^{+}) = (s^{-})^{\dagger}$.
\end{pf}

From the diagram in (\ref{diag-i}) and the definition of $\tau^{+}$ (resp. $\tau^{-}$) in (\ref{tau+}) (resp. (\ref{tau-})) one obtains
\BEN\label{FP+-}
(\tau^{+})^{\ast} \FTP = i^{\ast} ((\tau^{-})^{\ast} \FTP)\,.
\EEN
Combining this with the decompositions in (\ref{sh-shifted}) and (\ref{sh-shifted1}) yields two gradings on 
$(\tau^{+})^{\ast} \FTP$
\BEN\label{two-gr}
(\tau^{+})^{\ast} \FTP = \bigoplus^{2(\LG-1)}_{n=0} {\cal V}(n) =\bigoplus^{2(\LG-1)}_{n=0} i^{\ast} ({\cal V^{\prime}}(n) )\,.
\EEN
They are related by the orthogonality conditions below.
\begin{pro}\label{pro-two-gr}
For every $n \in \{0,1,\ldots, 2(\LG-1) \}$ the following holds.
\begin{enumerate}
\item[1)]
$( {\cal V}(n))^{\perp} =\bigoplus_{m \neq n} i^{\ast} ({\cal V^{\prime}}(m) )$.
\item[2)]
The quadratic form $\QT$ in $\FTP$ induces non-degenerate  bilinear pairing
$$
 {\cal V}(n) \times i^{\ast}{\cal V^{\prime}}(n) \longrightarrow \OO_{{\bf h^0}(d^{+})}
$$
yielding an isomorphism
$$
i^{\ast} {\cal V^{\prime}}(n) =({\cal V}(n))^{\ast}\,.
$$
\end{enumerate}
\end{pro}
\begin{pf}
Let $([Z],[\alpha],v,h)$ be a closed point of ${\bf h^0}(d^{+})$. Then the fibre of
$(\tau^{+})^{\ast} \FTP$ at this point is $\FTP \ZA$, the fibre of $\FTP$ at $\ZA \in \JABG$.
From the equality in (\ref{two-gr}) it follows that this vector space admits two weight decompositions
$$
\FTP\ZA =\bigoplus^{2(\LG-1)}_{n=0} {\cal V}(n) ([Z],[\alpha],v,h)=
\bigoplus^{2(\LG-1)}_{n=0} {\cal V^{\prime}}(n) ([Z],[\alpha],v,h^{\dagger}) 
$$
corresponding to the action of $h$ and $h^{\dagger}$ respectively. Our task will be to compare them.

Fix a summand
${\cal V}(n) ([Z],[\alpha],v,h)$. According to the shifting convention this is the weight space of $h$ corresponding to the weight
$n^{\prime} =n-\LG+1$. For every $x \in {\cal V}(n) ([Z],[\alpha],v,h)$ and 
$y \in  {\cal V^{\prime}}(m) ([Z],[\alpha],v,h^{\dagger}) $ we have
$$
n^{\prime} \QT(x,y) =\QT(h(x),y)=\QT(x, h^{\dagger}(y)) = m^{\prime} \QT(x,y)\,,
$$
where $m^{\prime} =m -\LG+1$. This implies that
$\QT(x,y) =0$, for all $m^{\prime} \neq n^{\prime}$. This gives an inclusion
\BEN\label{Vpr-inc-V}
\bigoplus_{m\neq n} {\cal V^{\prime}}(m) ([Z],[\alpha],v,h^{\dagger})  \subset ({\cal V}(n) ([Z],[\alpha],v,h))^{\perp}\,.
\EEN
On the other hand the bilinear pairing 
$$
\FTP\ZA \times \FTP\ZA \longrightarrow \CC
$$
induced by $\QT$ is non-degenerate. This together with the inclusion in (\ref{Vpr-inc-V}) imply that the pairing
\BEN\label{Vpr-V-pair}
{\cal V}(n) ([Z],[\alpha],v,h) \times {\cal V^{\prime}}(n) ([Z],[\alpha],v,h^{\dagger}) \longrightarrow \CC
\EEN
induces an injection
\BEN\label{Vpr-inc-V1}
{\cal V}(n) ([Z],[\alpha],v,h) \subset ({\cal V^{\prime}}(n))^{\ast} ([Z],[\alpha],v,h^{\dagger})\,.
\EEN
Interchanging the roles of 
${\cal V}(n)$ and ${\cal V^{\prime}}(n)$, we obtain an inclusion 
$$
{\cal V^{\prime}}(n) ([Z],[\alpha],v,h^{\dagger}) \subset ({\cal V}(n))^{\ast} ([Z],[\alpha],v,h)\,.
$$
Combining this inclusion with the one in (\ref{Vpr-inc-V1})
 implies that
that both inclusions are  equalities and the pairing in (\ref{Vpr-V-pair}) is perfect. Hence part 2) of the proposition is proved.

To complete the first part of the proposition we need to show that the inclusion in (\ref{Vpr-inc-V}) is an equality.
But this follows from the equality
$$
dim\left({\cal V}(n) ([Z],[\alpha],v,h) \right)= dim \left(({\cal V^{\prime}}(n))^{\ast} ([Z],[\alpha],v,h^{\dagger})\right)
$$
and the fact that the pairing in (\ref{Vpr-V-pair}) is perfect.
\end{pf}

As in the case of $d^{+}$ we go on to define the bigrading on 
$(\tau^{-})^{\ast} \FTP$. This is done by setting
\BEN\label{bigr-pr}
{}^{\prime} \HH^{p,q} ={\cal V^{\prime}}(p+q) \bigcap (\tau^{-})^{\ast} \HH^q\,.
\EEN
This yields
\BEN\label{Hq-pr-bigr}
(\tau^{-})^{\ast} \HH^q = \bigoplus^{\LG-1}_{p=0} {}^{\prime} \HH^{p,q}
\EEN
and hence a bigrading of $(\tau^{-})^{\ast} \FTP$
\BEN\label{tau-bigr-pr}
(\tau^{-})^{\ast} \FTP = \bigoplus_{(p,q)} {}^{\prime} \HH^{p,q}\,,
\EEN
where the sum is taken over the points $(p,q) \in {\bf Z^2}$ lying in the same square as in Proposition \ref{pro-dgr-sh}.

Observe that the section $(s^{\prime})^{-}$ is the pullback by $\eta^{-}_1$ of $d^{-}_{\diamond}$, the section $d^{-}_{\diamond}$ of 
$\tau^{\ast} \mbox{\BM$\LAG$}$ canonically attached to the morphism $d^{-}$.
All the properties of Proposition \ref{pro-dgr-sh} hold for the bigrading in (\ref{tau-bigr-pr}) with the only difference
that the role of 
$s^{+}$ there is taken over by
$(s^{\prime})^{-}$. This is an endomorphism of type
$(-1,-1)$ with respect to this bigrading, i.e. $(s^{\prime})^{-}$ induces the morphisms
\BEN\label{s-}
(s^{\prime})^{-}:  {}^{\prime} \HH^{p,q} \longrightarrow  {}^{\prime} \HH^{p-1,q-1}\,,
\EEN
for every $(p,q)$ occurring in (\ref{tau-bigr-pr}). Thus 3) of Proposition \ref{pro-dgr-sh}
reads as follows:

for every integer $k\in \{0,\ldots, \LG-1\}$, the endomorphism
$(s^{\prime})^{-} $ induces the isomorphism
\BEN\label{s-k}
((s^{\prime})^{-}  )^k : {\cal V^{\prime}}(\LG-1+k) \longrightarrow {\cal V^{\prime}}(\LG-1-k) 
\EEN
which is of type $(-k,-k)$ with respect to the bigrading in (\ref{tau-bigr-pr}). In particular, one has isomorphisms
\BEN\label{s-k-pq}
((s^{\prime})^{-}  )^k : {}^{\prime}\HH^{p,q} \longrightarrow  {}^{\prime}\HH^{p-k,q-k}\,,
\EEN
for every $(p,q)$ such that $p+q =\LG -1+k$ and $p,q \geq 0$.

The bigradings (\ref{bi-gr-tau}) and (\ref{tau-bigr-pr}) are also related via the pullback by the isomorphism $i$ in (\ref{inv1}).
The following can be viewed as a sort of Hodge-Riemann bilinear relations between the two bigradings.
\begin{pro}\label{Riemann}
The sheaves
$\HH^{p,q}$ in (\ref{bi-gr-tau}) and $i^{\ast} ({}^{\prime}\HH^{p^{\prime},q^{\prime}})$
are orthogonal with respect to the symmetric bilinear pairing on $(\tau^{+})^{\ast} \FTP$  induced by the quadratic form
$\QT$ on $\FTP$, unless $p=q^{\prime}$ and $q=p^{\prime}$. In this latter case the pairing
$$
\HH^{p,q} \times i^{\ast} ({}^{\prime}\HH^{q,p}) \longrightarrow \OO_{\bf h^0 (d^{-})}
$$
is perfect.
\end{pro}
\begin{pf}
By definition of the bigraded sheaves in (\ref{bigr-pr}) one has
$$
i^{\ast} ({}^{\prime}\HH^{p^{\prime},q^{\prime}}) =
 i^{\ast} \left({\cal V^{\prime}}(p^{\prime} +q^{\prime}) \bigcap (\tau^{-})^{\ast}\HH^{q^{\prime}} \right)=
 i^{\ast} \left({\cal V^{\prime}}(p^{\prime} +q^{\prime}) \right) \bigcap i^{\ast} \left((\tau^{-})^{\ast} \HH^{q^{\prime}} \right)=
 i^{\ast} \left({\cal V^{\prime}}(p^{\prime} +q^{\prime}) \right) \bigcap (\tau^{+})^{\ast} \HH^{q^{\prime}}\,.
$$
On the other hand by (\ref{Hpq})
$$
\HH^{p,q} = {\cal V}(p+q) \bigcap (\tau^{+})^{\ast} \HH^p\,.
$$
Using the orthogonality of the summands
$(\tau^{+})^{\ast} \HH^{q^{\prime}}$ and $(\tau^{+})^{\ast} \HH^p$, for all $q^{\prime} \neq p$ and Proposition \ref{pro-two-gr}, we obtain
$$
i^{\ast} ({}^{\prime}\HH^{p^{\prime},q^{\prime}}) \perp \HH^{p,q}
$$
unless $q^{\prime} =p$ and $p^{\prime} +q^{\prime} =p+q$. Hence the first assertion.

The second one follows from the first and non-degeneracy of the pairing
$$
(\tau^{+})^{\ast} (\FTP) \times (\tau^{+})^{\ast} (\FTP)  \longrightarrow \OO_{{\bf h^0} (d^{+})}
$$
induced by $\QT$.
\end{pf}

Next we turn to the weight filtration associated to $d^{-}$.
For every $n \in \{ 0, \ldots, 2(\LG-1) \}$ set
\BEN\label{w-filt-pr}
{\cal V^{\prime}}_n = \bigoplus_{m\leq n} {\cal V^{\prime}} (m)\,.
\EEN
These sheaves are constant along the fibres of the  projection
$$
\eta^{-}_1 :{\bf h^0} (d^{-}) \longrightarrow \TPI
$$
for the same reason as for the sheaves ${\cal V}^n$ in (\ref{w-filt-sh0}) and the projection $\eta^{+}_1$.
Thus there are sheaves ${\cal W^{\prime}}_n$ on $\TPI$ such that
\BEN\label{Wpr}
(\eta^{-}_1)^{\ast} {\cal W^{\prime}}_n = {\cal V^{\prime}}_n
\EEN
and the increasing filtration 
\BEN\label{w-filt-sh1}
0={\cal W^{\prime}}_{-1} \subset {\cal W^{\prime}}_0 \subset \ldots {\cal W^{\prime}}_{2(\LG-1)}= \tau^{\ast} \FTP
\EEN
 which will be called the weight filtration of $\tau^{\ast} \FTP$ associated to $d^{-}$.

The properties of this filtration are similar to the ones of ${\cal W}^{\bullet}$ in Proposition-Definition \ref{pro-w-filt-sh}
with the difference that $d^{+}_{\diamond}$ there should be replaced by $d^{-}_{\diamond}$ and the arrows of the corresponding morphisms should be reversed.
Thus one obtains the following.
\begin{pro-defi}\label{pro-w-filt-sh1}
The sheaf $\tau^{\ast} \mbox{\BM$\LAG$}$ comes along with a distinguished section
$d^{-}_{\diamond}$. This section viewed as an endomorphism of the sheaf
$\tau^{\ast} \FTP$ gives rise to the following structures. 
\begin{enumerate}
\item[1)]
$d^{-}_{\diamond}$ defines the weight filtration 
${\cal W^{\prime}}_{\bullet}$ as in (\ref{w-filt-sh1}) and  $d^{-}_{\diamond}$ acts on it by shifting the index of the filtration by $(-2)$, i.e. one has
$$
d^{-}_{\diamond} ({\cal W^{\prime}}_k) \subset {\cal W^{\prime}}_{k-2}\,,
$$
for every $k\geq 0$.
\item[2)]
Set
$$
Gr^k_{{\cal W^{\prime}}_{\bullet}}  (\tau^{\ast} \FTP)= {\cal W^{\prime}}_k /{\cal W^{\prime}}_{k-1}\,.
$$
 The orthogonal decomposition in (\ref{ord-tau1})
 induces a grading on  $Gr^k_{{\cal W^{\prime}}_{\bullet}} (\tau^{\ast} \FTP)$
$$
Gr^k_{{\cal W^{\prime}}_{\bullet}}  ( \tau^{\ast} \FTP) = \bigoplus Gr^{k-p,p}_{{\cal W^{\prime}}_{\bullet}}  ( \tau^{\ast} \FTP)\,, 
$$
where the summands are defined as follows
$$
Gr^{k-p,p}_{{\cal W^{\prime}}_{\bullet}}  (\tau^{\ast} \FTP) =  
\big{(}{\cal W^{\prime}}_k \bigcap \tau^{\ast} (\HH^p ) \big{)} / \big{(}{\cal W^{\prime}}_{k-1} \bigcap \tau^{\ast} (\HH^p ) \big{)}\,.
$$ 
\item[3)]
The associated graded sheaf
$$
Gr_{{\cal W^{\prime}}_{\bullet}} (\tau^{\ast} \FTP)  = 
\bigoplus^{2(\LG-1)}_{k=0} Gr^k_{{\cal W^{\prime}}_{\bullet}}  (\tau^{\ast} \FTP) =
\bigoplus^{2(\LG-1)}_{k=0} {\cal W^{\prime}}_k /{ \cal W^{\prime}}_{k-1}
$$
of $\tau^{\ast} \FTP$ defined with respect to the weight filtration ${\cal W^{\prime}}_{\bullet}$  together with the grading
 in 2) acquires the bigrading
\BEN\label{Gr-dgr2}
Gr_{{\cal W}^{\prime}_{\bullet}} (\tau^{\ast} \FTP)  = \bigoplus_{(p,q)}  Gr^{p,q}_{{\cal W}^{\prime}_{\bullet}}  (\tau^{\ast} \FTP) \,,
\EEN
where the direct sum is taken over the points $(p,q) \in {\bf Z^2}$, lying in the square described in Proposition \ref{pro-dgr-sh}, 1).
\item[4)]
The endomorphism $d^{-}_{\diamond}$ induces the endomorphism
$$
gr(d^{-}_{\diamond}) : Gr_{{\cal W}^{\prime}_{\bullet}} (\tau^{\ast} \FTP) \longrightarrow Gr_{{\cal W}^{\prime}_{\bullet}} (\tau^{\ast} \FTP)
$$
which has type (-1,-1) with respect to the double grading in 3), i.e. one has
$$
gr(d^{-}_{\diamond}) (Gr^{p,q}_{{\cal W}^{\prime}_{\bullet}}  (\tau^{\ast} \FTP ) \subset 
 Gr^{p-1,q-1}_{{\cal W}^{\prime}_{\bullet}}  (\tau^{\ast} \FTP)\,,
$$
for every $(p,q)$ occurring in the decomposition in (\ref{Gr-dgr2}).
\item[5)]
For every integer $k \in \{0,\ldots, \LG-1 \}$, the endomorphism $gr(d^{-}_{\diamond})$ induces the isomorphism
$$
(gr(d^{-}_{\diamond}))^k : Gr^{\LG-1+k}_{{\cal W}^{\prime}_{\bullet}}  (\FTP) \longrightarrow 
Gr^{\LG-1-k}_{{\cal W}^{\prime}_{\bullet}}  (\FTP)
$$
which is of type $(-k,-k)$ with respect to the bigrading in 3). In particular, one has isomorphisms
$$
(gr(d^{-}_{\diamond}))^k : Gr^{p,q}_{{\cal W}^{\prime}_{\bullet}}  (\tau^{\ast} \FTP) \longrightarrow
 Gr^{p-k,q-k}_{{\cal W}^{\prime}_{\bullet}}  (\tau^{\ast} \FTP)\,,
$$
for every $(p,q)$ such that 
 $p+q =\LG -1+k$ and $0\leq p,q \leq \LG-1$.
\end{enumerate}

The data of $\tau^{\ast} \FTP$ together  with $d^{-}_{\diamond}$ and its weight filtration
${\cal W}^{\prime}_{\bullet}$, the orthogonal decomposition in (\ref{ord-tau1}), the associated sheaf
$Gr_{{\cal W}^{\prime}_{\bullet}} (\tau^{\ast} \FTP) $ in 3) with its bigrading in (\ref{Gr-dgr2}) and the endomorphism
$gr(d^{-}_{\diamond})$ will be called {\rm ${\bf sl_2}$-structure of $\FTP$ associated to $d^{-}$} or
{\rm the natural negative ${\bf sl_2}$-structure of $\FTP$}. 
\end{pro-defi}

We now have two weight filtrations 
${\cal W}^{\bullet}$ and 
${\cal W}^{\prime}_{\bullet}$
of $\tau^{\ast} \FTP$ associated to $d^{+}$ and $d^{-}$ respectively. As in the case of  the corresponding weight decompositions in 
Proposition \ref{pro-two-gr} they are related by orthogonality.
\begin{pro}\label{pro-two-filt}
For every $n \in \{0,1,\ldots, 2(\LG-1) \}$ the following holds.
\begin{enumerate}
\item[1)]
${\cal W}^{\prime}_{n} =({\cal W}^{n+1})^{\perp}$\,.
\item[2)]
The bilinear pairing on 
$\tau^{\ast} \FTP$, induced by the quadratic form $\QT$ on $\FTP$, defines a perfect pairing
$$
{\cal W}^{n+1} \times \big{(} \tau^{\ast} \FTP / {{\cal W}^{\prime}_{n}} \big{)} \longrightarrow \OO_{\TPI}\,.
$$
\end{enumerate}
\end{pro}
\begin{pf}
By the defining property of the filtration ${\cal W}^{\bullet}$ in (\ref{Vn=pb}) we have
$$
(\eta^{+}_1)^{\ast} {\cal W}^{n+1} = {\cal V}^{n+1} =\bigoplus_{m\geq n+1} {\cal V}(m)\,,
$$
where the last equality is the definition in (\ref{w-filt-sh0}). Applying Proposition \ref{pro-two-gr}, 1), to
${\cal V}^{n+1}$ yields
$$
((\eta^{+}_1)^{\ast} {\cal W}^{n+1})^{\perp} = ( {\cal V}^{n+1})^{\perp} = \bigoplus_{s\leq n} i^{\ast} {\cal V}^{\prime} (s) =
 i^{\ast} ({\cal V}^{\prime}_n) = i^{\ast} ((\eta^{-}_1)^{\ast} {\cal W}^{\prime}_n )= (\eta^{+}_1)^{\ast} {\cal W}^{\prime}_n \,,
$$
where the third and the fourth equalities come from (\ref{w-filt-pr}) and (\ref{Wpr}) respectively, while the last one is the commutativity
of the diagram in (\ref{diag-i}). Thus we obtain the first assertion.

Turning to the second part observe that 1) implies that the pairing
$$
  {\cal W}^{n+1} \times \tau^{\ast} \FTP \longrightarrow \OO_{\TPI}
$$
defined by $\QT$ descends to 
${\cal W}^{n+1} \times \big{(} \tau^{\ast} \FTP / {{\cal W}^{\prime}_{n}} \big{)}$
and the induced pairing
$$
{\cal W}^{n+1} \times \big{(} \tau^{\ast} \FTP / {{\cal W}^{\prime}_{n}} \big{)} \longrightarrow \OO_{\TPI}
$$
gives an injection
$$
{\cal W}^{n+1} \hookrightarrow \big{(} \tau^{\ast} \FTP / {{\cal W}^{\prime}_{n}} \big{)}^{\ast}\,.
$$
From Proposition \ref{pro-two-gr}, 2), the ranks of both sheaves above are equal. Hence the above inclusion is an isomorphism.
\end{pf}

We have now two {\it a priori} distinct 
${\bf sl_2}$-structures on
$\tau^{\ast} \FTP$: the one is positive, associated to $d^{+}$ in Proposition-Definition \ref{pro-w-filt-sh}, and the other one is negative, 
associated to $d^{-}$ in Proposition-Definition \ref{pro-w-filt-sh1}. They give rise to two bigraded sheaves of modules
\begin{eqnarray}
Gr_{{\cal W}^{\bullet}} (\tau^{\ast} \FTP) & = &\bigoplus_{(p,q)}  Gr^{p,q}_{{\cal W}^{\bullet}}  (\tau^{\ast} \FTP)\,, \\   \label{bigr1}
Gr_{{\cal W}^{\prime}_{\bullet}} (\tau^{\ast} \FTP) & =& 
\bigoplus_{(p,q)}  Gr^{p,q}_{{\cal W}^{\prime}_{\bullet}}  (\tau^{\ast} \FTP)\,.    \label{bigr2}
\end{eqnarray}
These two bigradings are related as follows.
\begin{cor}\label{cor-two-bigr}
\begin{enumerate}
\item[1)]
For every $n \in \{0,1, \ldots, 2(\LG-1) \}$, there is a natural perfect pairing
\BEN\label{gr-pair}
{\cal W}^n /{\cal W}^{n+1} \times {\cal W}^{\prime}_n /{\cal W}^{\prime}_{n-1} \longrightarrow \OO_{\TPI}
\EEN
inducing the perfect pairing of graded sheaves
\BEN\label{gr-pair1}
 Gr_{{\cal W}^{\bullet}} (\tau^{\ast} \FTP) \times Gr_{{\cal W}^{\prime}_{\bullet}} (\tau^{\ast} \FTP)  \longrightarrow \OO_{\TPI}\,.
\EEN
This gives rise to the duality isomorphism
\BEN\label{gr-dual}
Gr_{{\cal W}^{\bullet}} (\tau^{\ast} \FTP) \cong \big{(} Gr_{{\cal W}^{\prime}_{\bullet}} (\tau^{\ast} \FTP) \big{)}^{\ast}\,.
\EEN
\item[2)]
The pairing in (\ref{gr-pair1}) induces pairings of the bigraded sheaves
$$
Gr^{p,q}_{{\cal W}^{\bullet}}  (\tau^{\ast} \FTP) \times  Gr^{p^{\prime},q^{\prime}}_{{\cal W}^{\prime}_{\bullet}}  (\tau^{\ast} \FTP)
\longrightarrow \OO_{\TPI}
$$
which are identically zero, unless $p^{\prime} =q$ and $q^{\prime} =p$ and in that case
the pairing
$$
Gr^{p,q}_{{\cal W}^{\bullet}}  (\tau^{\ast} \FTP) \times  Gr^{q,p}_{{\cal W}^{\prime}_{\bullet}}  (\tau^{\ast} \FTP)
\longrightarrow \OO_{\TPI}
$$
is perfect. In particular, the duality isomorphism in (\ref{gr-dual}) is a ``braided" isomorphism of bigraded sheaves, i.e.
the isomorphism in (\ref{gr-dual}) maps 
$Gr^{p,q}_{{\cal W}^{\bullet}}  (\tau^{\ast} \FTP) $ isomorphically onto
$\big{(} Gr^{q,p}_{{\cal W}^{\prime}_{\bullet}}  (\tau^{\ast} \FTP) \big{)}^{\ast}$.
\end{enumerate}
\end{cor}
\begin{pf}
Applying Proposition \ref{pro-two-filt} to $n$ and $n+1$ one deduces the perfect pairing asserted in (\ref{gr-pair}).
Putting these pairings together for all $n$ yields the perfect pairing in (\ref{gr-pair1}).

The part 2) follows from the first part and the definitions of graded sheaves 
$Gr^{p,q}_{{\cal W}^{\bullet}}  (\tau^{\ast} \FTP) $ (resp.
$Gr^{p^{\prime},q^{\prime}}_{{\cal W}^{\prime}_{\bullet}}  (\tau^{\ast} \FTP)$)
in Proposition-Definition \ref{pro-w-filt-sh}, 3)
(resp. Proposition-Definition \ref{pro-w-filt-sh1}, 3)).
\end{pf}

\section{Stratification of $\TPI$}\label{sec-strat}
In this section we use the morphism $d^{+}$ defined  in (\ref{d+}) for geometric purposes. 
Throughout this section we fix an admissible component $\GA$ in $\CS$ and assume that it is not quasi-abelian 
(see Definition \ref{qa-c}). As before $\TPI$ denotes the relative tangent bundle of $\JABG$ over $\GAB$. 
\subsection{Linear stratification of $\TPI$}
This stratification arises, when we consider the morphism $d^{+}$ written as in (\ref{d+-1}). For every 
$p\in \{0,\ldots, \LG-2 \}$ set
\BEN\label{d+p}
d^{+}_p : {\cal T}_{\pi} \longrightarrow \HOM (\HH^p , \HH^{p+1})
\EEN
to be the $p$-th component of $d^{+}$ and define
\BEN\label{T(p)}
{\cal T}^{(p)}_{\pi}  = ker (d^{+}_p )\,.
\EEN
The following result gives some basic properties of these subsheaves.
\begin{lem}\label{lem-T(p)}
For every $p\in \{0,\ldots, \LG-2 \}$ the following holds
\begin{enumerate}
\item[1)]
${\cal T}^{(p)}_{\pi}  \subset {\cal T}^{(p+1)}_{\pi}$\,,
\item[2)]
 ${{\cal T}^{(\LG-2)}_{\pi}} \neq \TPI$\,.
\end{enumerate}
\end{lem}
\begin{pf}
Let $v$ be a local section of ${\cal T}^{(p)}_{\pi} $.
From Remark \ref{val-d} it follows that a local section
$\tilde{v}$ of $\HT$ lifting $M^{-1} (v)$, where the isomorphism $M$ is as in (\ref{M}), is subject to the following property
\BEN\label{D+v=0}
D^{+}_p ( \tilde{v}) =0.
\EEN
This implies that the multiplication by $\tilde{v}$ preserves the subsheaf
$\HT_{-p-1}$ of the filtration 
$\HT_{-\bullet}$ in (\ref{filtHT-JG}). From the definition of this filtration in (\ref{HT-i}) one deduces 
that the multiplication by $\tilde{v}$ preserves $\HT_{-q}$, for all $q\geq p+ 1$. This in turn is equivalent to
$$
D^{+}_q ( \tilde{v}) =0,\,\,\forall q \geq p\,.
$$
Using Remark \ref{val-d} once again, we obtain
$$
d^{+}_q  (v) = D^{+}_q ( \tilde{v}) =0,\,\,\forall q \geq p\,.
$$
Hence $v$ is a local section of ${\cal T}^{(q)}_{\pi}$, for all 
$q \geq p$, yielding the inclusion 
$$
{\cal T}^{(p)}_{\pi}  \subset {\cal T}^{(p+1)}_{\pi} \,.
$$
This proves the first assertion of the lemma. To see the second one we assume the equality
${\cal T}^{(\LG-2)}_{\pi} = {\cal T}_{\pi} $. Then the argument of the proof of the first part implies that the multiplication by any local section of $\HT$ preserves 
$\HT_{-{\LG +1}}$. But this means that the filtration in (\ref{filtHT-JG})
stabilizes at $ \HT_{-{\LG +1}}$. This contradicts the fact that the length of the filtration  is $\LG$.
\end{pf}

 We switch now to the geometric point of view of bundles over $\JABG$:
 the subsheaves ${\cal T}^{(p)}_{\pi}  $ in the above lemma   correspond to proper subschemes of $\TPI$. We denote them by 
$T^{(p)}_{\pi}$. Then the inclusions in Lemma \ref{lem-T(p)}, 1), give the stratification
\BEN\label{linear-strat}
  T^{(0)}_{\pi} \subset T^{(1)}_{\pi} \subset \cdots \subset 
T^{(\LG-2)}_{\pi} \subset T^{(\LG-1)}_{\pi} = \TPI
\EEN
which we refer to as {\it linear stratification} of  $\TPI$. 
\begin{rem}
An admissible component $\GA$ is simple (Definition \ref{s-c}) if and only if  the stratum $T^{(0)}_{\pi}$ is the zero-section of $\TPI$. 
More generally, if $\tau_{\GA}$ is the Torelli index of $\GA$, then 
$$
   T^{(0)}_{\pi} = \cdots =   T^{(\tau_{\GA} -1)}_{\pi} 
$$
is the zero-section of $\TPI$.

In the extreme case of $\GA$, being subject to the strong Torelli property
(see Definition \ref{st-T}), the linear filtration is reduced to
$$T^{(0)}_{\pi}  =   T^{(\LG -2)}_{\pi} \subset \TPI\,, $$
where $T^{(0)}_{\pi}$ is the zero-section of $\TPI$.
\end{rem}
\subsection{${\bf sl_2}$-stratification of $\TPI$ }
In this subsection we use ${\bf sl_2}$-structures discussed in 
\S\ref{sec-sl2} to produce a ``non-linear" stratification of $\TPI$.

From the assumption that $\GA$ is not quasi-abelian it follows that  
the sheaf 
{\BM
$\LAG$} is non-zero. We know that its pullback $\tau^{\ast}\mbox{\BM$\LAG$}$ comes along with
two distinguished sections $d^{\pm}_{\diamond}$. From \S\ref{sec-inv} they are exchanged by the involution of taking the adjoint
and their ${\bf sl_2}$-structures are related via orthogonality. 
So it will be enough to consider one of them, say $d^{+}_{\diamond}$. Considering the natural (defining) representation of
{\BM
$\LAG$}
on $\FTP$,
we view $d^{+}_{\diamond}$ as an endomorphism of $\tau^{\ast} \FTP$
\BEN\label{d+0-endo}
d^{+}_{\diamond} : \tau^{\ast} \FTP \longrightarrow \tau^{\ast} \FTP
\EEN
and consider its kernel
\BEN\label{ker+}
{\cal K}^{+} := ker (d^{+}_{\diamond})\,.
\EEN
From the point of view of ${\bf sl_2}$-representations this subsheaf parametrizes the highest weight vectors. 
Hence the importance of this subsheaf. 

Since $d^{+}_{\diamond}$ is of degree $1$ with respect to the orthogonal decomposition of $\tau^{\ast} \FTP$ in (\ref{ord-tau1})
 its kernel inherits this decomposition, i.e. we have
\BEN\label{ord-K+}
{\cal K}^{+} =\bigoplus^{\LG-1}_{p=0} ({\cal K}^{+})^p\,,
\EEN
where $({\cal K}^{+})^p =\tau^{\ast} \HH^p \bigcap {\cal K}^{+}$.
Furthermore, the weight filtration 
${\cal W}^{\bullet}$ defined by $d^{+}_{\diamond}$ (see Proposition-Definition \ref{pro-w-filt-sh}, 1)) restricts to each
summand
$({\cal K}^{+})^p$ and hence produces the weight filtration 
\BEN\label{w-filt-Kp}
({\cal K}^{+})^p =({\cal K}^{+})^{p,[0]} \supset ({\cal K}^{+})^{p,[1]} \supset \ldots \supset ({\cal K}^{+})^{p,[2(\LG-1)]}
 \supset  ({\cal K}^{+})^{p,[2\LG-1]} =0\,,
\EEN
where 
$({\cal K}^{+})^{p,[m]} = ({\cal K}^{+})^p \bigcap {\cal W}^m$.
The properties of ${\bf sl_2}$-representations allow us to be more precise about the form of these filtrations.
\begin{lem}\label{lem-w-filt-Kp}
The weight filtration $({\cal K}^{+})^{p,[\bullet]}$ has the following form
\BEN\label{w-filt-Kp1}
({\cal K}^{+})^p =({\cal K}^{+})^{p,[\LG-1]} \supset ({\cal K}^{+})^{p,[\LG]} \supset \ldots \supset ({\cal K}^{+})^{p,[p +\LG-1]}
 \supset  ({\cal K}^{+})^{p,[p+\LG]} =0\,.
\EEN
\end{lem}
\begin{pf}
Take the pullback 
$(\eta^{+}_1)^{\ast}({\cal K}^{+})^p$ and examine the weights of
$s^{0}$ on it.\footnote{we use the notation of \S\ref{const-sh}, Proposition \ref{pro-dgr-sh}.} Since 
$({\cal K}^{+})^p$  is generated by highest weight vectors of the representation of the
${\bf sl_2}$-triple 
$\{s^{+},s^0,s^{-} \}$, these weights must be positive. On the other hand from
Lemma \ref{Hp-dgr} it follows that the weights of $s^0$ on the summand
$(\tau^{+})^{\ast} \HH^p$ in (\ref{ord-tau}) are at most $p$. Keeping in mind the shift convention we deduce
$$
({\cal K}^{+})^{p}  \subset {\cal W}^{\LG-1} \,\,and\,\, ({\cal K}^{+})^{p,[m]} =0,\,\,\forall m\geq p+\LG\,.
$$
Hence the assertion of the lemma.
\end{pf}
The orthogonal decomposition (\ref{ord-K+}) together with the weight filtration in (\ref{w-filt-Kp1}) give rise to the associated
bigraded sheaf of modules
\BEN\label{K+-bigr}
Gr_{{\cal W}^{\bullet}} ({\cal K}^{+}) =\bigoplus^{\LG-1}_{p=0} 
\left( \bigoplus^{p+\LG-1}_{q=\LG-1} Gr^{p,q}_{{\cal W}^{\bullet}} ({\cal K}^{+})\right)\,,
\EEN
where
$Gr^{p,q}_{{\cal W}^{\bullet}} ({\cal K}^{+}) = ({\cal K}^{+})^{p,[q]} / ({\cal K}^{+})^{p,[q+1]}$.
Thus , for every $(p,q)$ occurring in (\ref{K+-bigr}), the bigraded sheaf
$Gr^{p,q}_{{\cal W}^{\bullet}} ({\cal K}^{+})$
parametrizes the highest weight vectors of weight $(q-\LG+1)$ (shifting convention) which are located in the $p$-th summand of the orthogonal decomposition
in (\ref{ord-tau1}). Equivalently,  $Gr^{p,q}_{{\cal W}^{\bullet}} ({\cal K}^{+})$ parametrizes 
irreducible ${\bf sl_2}$-representations of weight $(q-\LG+1)$ which are contained in the range
$[p-q+\LG-1, p]$ of the grading in (\ref{ord-tau1}).
Hence the bigraded sheaf encapsulates the decomposition of 
$(\tau^{+})^{\ast} \FTP$ into the direct sum of isotypical irreducible ${\bf sl_2}$-submodules,
with  the additional feature of keeping track of where these submodules are located with respect to the grading 
of $(\tau^{+})^{\ast} \FTP$ in (\ref{ord-tau}). 

With the representation theoretic meaning of 
the bigraded sheaf $Gr_{{\cal W}^{\bullet}} ({\cal K}^{+}) $ clarified, we turn to its geometric aspect.
Namely, we use the decomposition in (\ref{K+-bigr}) to stratify 
$\TPI$ according to the ranks of its bigraded summands. To do this we need some notation.

For a point $([Z], [\alpha],v) \in \TPI$, we denote as before $d^{+}(v)$  to be the value of 
$d^{+}_{\diamond}$ at  $([Z], [\alpha],v)$ and set
\BEN\label{K+p-v}
{\cal K}^{+} (v) =ker(d^{+}(v) ) =\bigoplus^{\LG-1}_{p=0} ({\cal K}^{+})^p (v)
\EEN
to be the fibre of ${\cal K}^{+}$ as well as of its orthogonal summands $({\cal K}^{+})^p $'s at $([Z], [\alpha],v) $.
Analogously, ${\cal W}^{\bullet} (v)$ will denote the fibre of 
${\cal W}^{\bullet}$ at $([Z], [\alpha],v)$. This is the weight filtration of 
$\FTP([Z])$ associated to the nilpotent endomorphism $d^{+}(v)$. Set 
\BEN\label{Grpq-v}
Gr^{p,q}_{{\cal W}^{\bullet} (v)} ({\cal K}^{+}(v)) = ({\cal K}^{+})^{p,[q]} (v) / ({\cal K}^{+})^{p,[q+1]} (v)
\EEN
to be the associated bigraded vector spaces, where
$({\cal K}^{+})^{p,[q]} (v)  = ({\cal K}^{+})^p (v) \bigcap {\cal W}^{q} (v)$.
\begin{rem}\label{bigr=multm}
The representation theoretic meaning of these spaces is quite transparent:
it parametrizes irreducible ${\bf sl_2}$-modules of weight $(q-\LG+1)$
which are contained in the submodule
$$
\bigoplus^p_{k=p-q+\LG-1} \HH^k \ZA \subset \FTP([Z]) =\bigoplus^{\LG-1}_{k=0} \HH \ZA\,.
$$
The ${\bf sl_2}$-action on $ \FTP([Z])$ is the one given by an 
${\bf sl_2}$-triple in ${\bf h^0}(d^{+}(v))$, the fibre of 
${\bf h^0}(d^{+})$ in (\ref{sc-sl2}) lying over $([Z], [\alpha],v) $.
If $\{d^{+}(v), h, y \}$ is such an ${\bf sl_2}$-triple, then the above can be stated as follows.

Define the subspace
$$
M^p (v,h) = \bigoplus_{k\geq0} y^k \left( ({\cal K}^{+})^{p} (v)\right)
$$
of $ \FTP([Z])$ and observe that it is invariant under the action of $\{d^{+}(v), h, y \}$. Thus it is an ${\bf sl_2}$-submodule of
$\FTP([Z])$ contained in $\HT_{-p-1} \ZA$. Also set ${\bf S^n}=Sym^n(\CC^2)$ to be the standard irreducible 
${\bf sl_2}$-module of weight $n$, then we have a (non-canonical) isomorphism
\BEN\label{bigr=mult-formula}
Gr^{p,q}_{{\cal W}^{\bullet}(v)} ({\cal K}^{+}(v)) \cong Hom_{\bf sl_2} ({\bf S^{q-\LG+1}},M^p (v,h))\,,
\EEN
where $Hom_{\bf sl_2}$ is the $Hom$-functor in the category of ${\bf sl_2}$-modules.
\end{rem}

Using the notation in (\ref{Grpq-v}) we set
\BEN\label{mu-qp}
\mu_{qp} (v) =dim \left(Gr^{p,\LG-1 +q}_{{\cal W}^{\bullet}(v)} ({\cal K}^{+}(v))\right), \,\,0\leq q \leq p \leq \LG-1\,.
\EEN
From (\ref{bigr=mult-formula}) it follows that these integers are the multiplicities
of irreducible ${\bf sl_2}$-submodules of $\FTP([Z])$  having the weight $q$ and contained
in the range $[p-q,p]$ of the grading
$$
\FTP([Z]) =\bigoplus^{\LG-1}_{m=0} \HH^m \ZA\,.
$$

We arrange the integers $\mu_{qp} (v)$ in $\LG \times \LG$ upper triangular matrix
\BEN\label{matr-mu}
M(v^{+}) = (M(v^{+})_{qp})
\EEN
with zeros below the main diagonal, i.e.
$M(v^{+})_{qp} =\mu_{qp},\,\,\forall 0\leq q \leq p \leq \LG-1$ and 
$M(v^{+})_{qp} = 0,\,\,\forall 0\leq p <q \leq \LG-1$.
\begin{defi}\label{mult-mat-def}
The matrix $M(v^{+})$ in (\ref{matr-mu}) will be called 
the multiplicity matrix of $d^{+}$ at $([Z], [\alpha],v) $ or,
simply, $d^{+}$-multiplicity matrix of $([Z], [\alpha],v) $. Its entries in the upper triangle will be called
the $d^{+}$-multiplicities of $([Z], [\alpha],v) $.
(If no ambiguity is likely we omit the reference to $\ZA$ in the above terminology and simply speak about
 $d^{+}$-multiplicity matrix and  $d^{+}$-multiplicities of $v$).
\end{defi}
\begin{pro}\label{pro-strat}
Let $\GA$ be an admissible component in $\CS$ and assume it is not quasi-abelian.
Let $\TPI$ be the relative tangent bundle of the natural projection
$$
\pi: \JABG \longrightarrow \GAB\,.
$$
Then $\TPI$ admits a partition into a finite union of locally closed subsets on each of which the multiplicity
matrix of
$d^{+}$ in (\ref{matr-mu}) is constant.
\end{pro}
\begin{pf}
The finiteness follow from the fact that the set of matrices in (\ref{matr-mu}) is finite. Indeed, all the matrix entries
$M(v^{+})_{qp}$ are non-negative and
\BEN\label{col-M}
M^p (v^{+}) :=\sum^p_{q=0} M(v^{+})_{qp} =\sum^p_{q=0} \mu_{qp} = dim(({\cal K}^{+})^p (v)) \leq  rk(\HH^p) =h^p_{\GA}\,.
\EEN

Next we turn to defining the strata of the asserted partition. For this consider the stratification of 
$\TPI$ according to the rank of 
${\cal K}^{+}$:
\BEN\label{str-T}
\TPI=T^1_{\pi} \supset T^2_{\pi} \supset \ldots \supset T^{d^{\prime}_{\GA} -1}_{\pi} \supset
 T^{d^{\prime}_{\GA} }_{\pi} \supset T^{d^{\prime}_{\GA} +1}_{\pi} =\emptyset\,,
\EEN
where $d^{\prime}_{\GA} =rk (\FTP)$ and where each stratum $ T^s_{\pi} $, set theoretically, is defined as follows
\BEN\label{str-Ts}
T^s_{\pi} =\left\{\left. ([Z], [\alpha],v) \in \TPI \right| dim({\cal K}^{+} (v)) = d^{\prime}_{\GA} -rk(d^{+} (v)) \geq s \right\}\,. 
\EEN
Thus (\ref{str-T}) is the stratification of $\TPI$ by the degeneracy loci of the morphism
$d^{+}_{\diamond}$ in (\ref{d+0-endo}) and each stratum is a closed subscheme of $\TPI$ which we consider with its reduced structure.

For every locally closed part
$$
 {\stackrel{\circ}{T^s}}_{\pi} = T^s_{\pi} \setminus T^{s+1}_{\pi}
$$
which is non-empty, consider the set 
$Irr({\stackrel{\circ}{T^s}}_{\pi})$ of its irreducible components.
For every component $\Sigma^s \in Irr({\stackrel{\circ}{T^s}}_{\pi})$ the restriction
${\cal K}^{+} \otimes \OO_{\Sigma^s}$ is a locally free sheaf of rank $s$. The components in $Irr({\stackrel{\circ}{T^s}}_{\pi})$
will be called admissible strata of rank $s$.

Choose an admissible stratum $\Sigma^s$ of rank $s$. In view of the orthogonal decomposition
of ${\cal K}^{+}$ in (\ref{ord-K+}) the summands 
$({\cal K}^{+})^p \otimes \OO_{\Sigma^s}$ are locally free, whenever they are non-zero.
Setting 
\BEN\label{kp-sigma}
k^p (\Sigma^s) = rk (({\cal K}^{+})^p \otimes \OO_{\Sigma^s})\,,
\EEN
we obtain the vector
\BEN\label{c-sigma}
c(\Sigma^s) =(k^0 (\Sigma^s), k^1 (\Sigma^s), \ldots, k^{\LG-1} (\Sigma^s))\,.
\EEN
This is a composition of $s$, i.e.
\BEN\label{c-sigma1}
\left| c(\Sigma^s) \right| = \sum^{\LG-1}_{p=0}  k^p (\Sigma^s) =s
\EEN
Thus we obtain a map 
\BEN\label{c-map}
c:Irr({\stackrel{\circ}{T^s}}_{\pi}) \longrightarrow C_{\LG} (s)\,,
\EEN
where $ C_{\LG} (s)$ is the set of compositions of $s$ having $\LG$ parts (the parts here are allowed to be zero).

Let $K=(k^0,k^1,\ldots, k^{\LG-1})$ be a composition lying in the image of the map $c$ in (\ref{c-map}) and let 
$\Sigma^s$ be an admissible stratum of rank $s$ with
$c(\Sigma^s) =K$. We can now stratify such a $\Sigma^s$ according to the ranks of
the subsheaves
$({\cal K}^{+})^{p,[\LG-1 +q]} \otimes  \OO_{\Sigma^s}$ of the 
filtration in (\ref{w-filt-Kp1}) restricted to $\Sigma^s$, for $q=0,\ldots,p$, i. e. 
 $\Sigma^s$ can be partitioned into the disjoint union of locally closed strata on each of which the 
ranks of the sheaves $({\cal K}^{+})^{p,[\LG-1 +q]} \otimes  \OO_{\Sigma^s}$ are constant for all 
$0\leq q \leq p \leq \LG-1$. This in turn implies that on these strata the ranks of the quotient sheaves
$$
Gr^{p,\LG-1+q}_{{\cal W}^{\bullet}} ({\cal K}^{+}) =({\cal K}^{+})^{p,[\LG-1 +q]} / ({\cal K}^{+})^{p,[\LG +q]}
$$
are constant for all $0\leq q \leq p \leq \LG-1$. These ranks arranged in an upper triangular matrix yield the 
$d^{+}$-multiplicity matrices  which are constant on each of these strata.
\end{pf} 
\begin{rem}\label{rem-strat}
\begin{enumerate}
\item[1)]
For every point
$ ([Z], [\alpha],v) \in \TPI$ 
the endomorphism
$d^{+} (v)$ annihilates the summand
$\HH^{\LG-1}$ and the subspace $\tilde{\bf c}$ of $\HT \ZA =\HH^0 \ZA$
described in Proposition \ref{pro-c-basis}. This implies that
$$
T^s_{\pi} = \TPI\,,
$$
 for all $s \leq h^{\LG-1}_{\GA} + rk(\mbox{\BM$\CG$})$, where
$ h^{\LG-1}_{\GA} = rk(\HH^{\LG-1})$.
\item[2)]
A component $\GA$ is simple (Definition \ref{s-c}) if and only if the stratum 
$T ^{d^{\prime}_{\GA}}_{\pi}$ is the zero section of $\TPI$ (this is seen by combining
Corollary \ref{InfTorelli=s} and Proposition \ref{InfTor-equiv}).
\item[3)]
The strata in Proposition \ref{pro-strat} are indexed by $d^{+}$-multiplicity matrices
(Definition \ref{mult-mat-def}). Given such a matrix, say $A=(A_{qp})$, denote by
$\Sigma(A)$ a stratum of $\TPI$ on which the $d^{+}$-multiplicity matrix is constant and equals the matrix $A$.
From this matrix we can read off the following:
\begin{enumerate}
\item[a)]
the ranks of the bigraded summands of the bigraded sheaf 
$Gr_{{\cal W}^{\bullet}} ({\cal K}^{+}) \otimes \OO_{\Sigma(A)}$ in (\ref{K+-bigr}), i.e.
$$
rk (Gr^{p,\LG-1 +q}_{{\cal W}^{\bullet}} ({\cal K}^{+}) \otimes \OO_{\Sigma(A)}) = A_{qp}\,,
$$
for all $0\leq q \leq p \leq \LG-1$.
\item[b)]
The rank of the stratum $\Sigma(A)$ (see the proof of Proposition \ref{pro-strat}) which is given by 
$\left| A\right| =\sum_{q,p} A_{qp}$, the sum of all entries of $A$, i.e. one has an inclusion
$$
\Sigma(A) \subset  {\stackrel{\circ}{T}}^{\mid A\mid }_{\pi} = T^{\mid A\mid }_{\pi} \setminus T^{\mid A\mid +1}_{\pi}\,.
$$
\item[c)]
The composition $c(\Sigma(A) )$ (see (\ref{c-sigma})) or, equivalently, the ranks of the sheaves
$({\cal K}^{+})^p \otimes \OO_{\Sigma(A)}$, for $p=0,1,\ldots, \LG-1$:
$$
k^p (\Sigma(A)) = rk (({\cal K}^{+})^p \otimes \OO_{\Sigma(A)}) = \sum_{q} A_{qp}\,.
$$
\end{enumerate}
The $d^{+}$-multiplicity matrix $A$ of a stratum $\Sigma(A)$ can be viewed as a bigraded version of the Hilbert vector
$h_{\GA}$ (see Lemma \ref{hg=c}) of the admissible component underlying $\Sigma(A)$. 
The passage from the $d^{+}$-multiplicity matrix $A$ to the
composition
$$
c(\Sigma(A))=(k^0(\Sigma(A)),\ldots, k^{\LG-1}(\Sigma(A))
$$
 gives the obvious inequalities
$$
k^p(\Sigma(A)) \leq h^p_{\GA},\,\,\forall 0\leq p \leq \LG-1\,.
$$
Furthermore, $k^{\LG-1}(\Sigma(A))=  h^{ \LG-1}_{\GA}$ and the other inequalities are strict unless 
$\Sigma(A)$ is contained in one of the ``linear" strata $\TPI^{(p)}$ of $\TPI$ described in (\ref{linear-strat}).
\end{enumerate}
\end{rem}

 Let ${\bf M}_{\LG} ({\bf Z_{+}})$ be the semi-ring of $\LG \times \LG$-matrices
with non-negative integer coefficients and let 
${\bf B}^{+}_{\LG} ({\bf Z_{+}})$ the subset of ${\bf M}_{\LG} ({\bf Z_{+}})$ consisting of 
upper triangular matrices. The result of Proposition \ref{pro-strat} suggests the following terminology. 
\begin{defi}\label{matr-adm}
A matrix $A$ in 
${\bf B}^{+}_{\LG} ({\bf Z_{+}})$
is called admissible if there is a stratum in Proposition \ref{pro-strat}, whose
$d^{+}$-multiplicity matrix is $A$. The set of all admissible matrices in
${\bf B}^{+}_{\LG} ({\bf Z_{+}})$ will be denoted by
${\bf B}^{+adm}_{\LG} ({\bf Z_{+}})$.
\end{defi}
\begin{rem}
Replacing $d^{+}_{\diamond}$ by $d^{-}_{\diamond}$ gives us the kernel
\BEN\label{K-}
{\cal K}^{-} = ker (d^{-}_{\diamond}) = \bigoplus^{\LG-1}_{p=0} ({\cal K}^{-})^p\,.
\EEN
Using the filtration ${\cal W}^{\prime}_{\bullet}$ in (\ref{w-filt-sh1}) and keeping in mind that the weights of the negative ${\bf sl_2}$-structure are distributed in the same way as for the positive one (see Proposition \ref{Riemann}) , we obtain the following induced filtration on each summand $ ({\cal K}^{-})^p$
\BEN \label{K-p-filt}
0=({\cal K}^{-})^{[p-1],p} \subset ({\cal K}^{-})^{[p],p} \subset ({\cal K}^{-})^{[p+1],p} \subset \cdots \subset 
({\cal K}^{-})^{[\LG-1],p} =({\cal K}^{-})^p\,,
\EEN
where $({\cal K}^{-})^{[m],p} =  ({\cal K}^{-})^p \bigcap {\cal W}^{\prime}_{m}$. This leads to the associated bigraded sheaf of modules
\BEN\label{K-bigr}
Gr_{{\cal W}^{\prime}_{\bullet}} ({\cal K}^{-}) = \bigoplus^{\LG-1}_{p=0} \left( \bigoplus^{\LG-1}_{q=p} 
Gr^{q,p}_{{\cal W}^{\prime}_{\bullet}}  ({\cal K}^{-}) \right)\,.
\EEN
Considering all the above at a point $([Z],[\alpha],v)$ of $\TPI$, gives us bigraded modules
$$
Gr_{{\cal W}^{\prime}_{\bullet} (v)} ({\cal K}^{-} (v)) = \bigoplus^{\LG-1}_{p=0} \left( \bigoplus^{\LG-1}_{q=p} 
Gr^{q,p}_{{\cal W}^{\prime}_{\bullet} (v)}  ({\cal K}^{-} (v)) \right)\,,
$$
for every vertical tangent vector $v$ in $\TPI \ZA$, the fibre of $\TPI$ over $\ZA \in \JABG$.

Similar to $d^{+} (v)$ we set
\BEN\label{mult-v-}
\mu^{\prime}_{qp} (v) =dim (Gr^{q,p}_{{\cal W}^{\prime}_{\bullet} (v)} ) ({\cal K}^{-} (v)))
\EEN
and call these numbers $d^{-}$-multiplicities of $v$. As in (\ref{matr-mu}) we arrange them in a matrix, but this time the {\rm lower}
triangular one,
\BEN\label{matr-mu-}
M(v^{-}) = ( M(v^{-})_{qp})
\EEN
with $M(v^{-})_{qp} =\mu^{\prime}_{qp} (v)$, for all $0\leq p \leq q \leq \LG-1$, and zeros above the main diagonal. The matrix
$M(v^{-})$ will be called the multiplicity matrix of $d^{-}$ at $([Z],[\alpha],v)$ or, simply, 
$d^{-}$-multiplicity matrix of $v$.  

Using bilinear pairings in Proposition \ref{Riemann} together with properties of ${\bf sl_2}$-representations it is not difficult to see
that $d^{-}$- and $d^{+}$-multiplicities of $v$ are related as follows
\BEN\label{mult+-}
\mu^{\prime}_{qp} (v) = \mu_{(\LG-1 -q)(\LG-1 -q +p)}\,,
\EEN
for all $0\leq p \leq q \leq \LG-1$.
\end{rem}

We will now give an alternative description of the stratification in Proposition \ref{pro-strat}.
It is based on the well-known correspondence between nilpotent elements of 
${\bf sl_n (C)}$ and partitions of $n$. Recall that given a nilpotent element $x$ in 
 ${\bf sl_n (C)}$ one assigns to it the partition of $n$ by taking the Jordan form of $x$ and arranging the sizes of its Jordan blocks in the 
decreasing order. In our situation to every point $([Z],[\alpha],v)$ in $\TPI$ we assign the nilpotent
element $d^{+}(v)$ of the Lie algebra
$\mbox{\BM$\LAG$}\ZA$. Using the defining representation of $\mbox{\BM$\LAG$}\ZA$
on $\FTP([Z])$ we view
$d^{+}(v)$ as a nilpotent element of ${\bf sl}(\FTP([Z]))$ and we let $\lambda(v)$ to be the partition of
$d^{\prime}_{\GA} =rk(\FTP)$ associated to it. Below we show that this partition is completely determined by the
$d^{+}$-multiplicity matrix $M(v^{+})$ in (\ref{matr-mu}). To do this we introduce some additional notation.

Given the matrix $M(v^{+})$ we assign  partitions to its column vectors by setting
\BEN\label{part-p}
\lambda^{(p)} (v) =(1^{\mu_{0p} (v)} 2^{\mu_{1p} (v)} \ldots (p+1)^{\mu_{pp} (v)}), \,\,for \,\,p=0, \ldots, \LG-1\,,
\EEN
where the notation
$m^k$ means that the part $m$ occurs in a given partition $k$ times.\footnote{one calls $k$ the multiplicity of the part $m$ in a partition;
for this and other standard notation and facts about partitions our reference is \cite{[Mac]}.}
\begin{lem}\label{part-union}
The partition $\lambda(v)$ corresponding to $d^{+}(v)$ is the union of partitions $\lambda^{(p)} (v)$'s in 
(\ref{part-p})
\BEN\label{part-union1}
\lambda(v) = \bigcup^{\LG-1}_{p=0} \lambda^{(p)} (v)
\EEN
(recall: the union $\lambda \cup \mu$ of two partitions is defined by the partition whose parts are those of $\lambda$
and $\mu$, arranged in the decreasing order).
\end{lem}
\begin{pf}
Observe that for every $p=0,1,\ldots, \LG-1$, the partition $\lambda^{(p)} (v)$ records the dimensions of the summands of the
 graded module
\BEN\label{gr-K+p}
Gr_{{\cal W}^{\bullet}(v)} (({\cal K}^{+})^p (v) ) =
 \bigoplus^{p}_{q=0} Gr^{p,\LG-1 + q}_{{\cal W}^{\bullet} (v)} ({\cal K}^{+} (v))\,.
\EEN
  We claim that this module contributes to the Jordan form of $d^{+}(v)$ the Jordan blocks of sizes prescribed by the  
partition $\lambda^{(p)} (v)$.
Indeed, from the properties of ${\bf sl_2}$-representations it follows that the summand
$ Gr^{p,\LG-1 + q}_{{\cal W}^{\bullet} (v)} ({\cal K}^{+} (v))$
parametrizes irreducible ${\bf sl_2}$-submodules of $\FTP([Z])$ of weight $q$ whose highest weight vectors are contained in $\HH^p \ZA$.
Each such submodule contributes a Jordan block of size $(q+1)$ in the Jordan form of
$d^{+}(v)$. Hence we have
$\mu_{qp} (v) = dim  Gr^{p,\LG-1 + q}_{{\cal W}^{\bullet} (v)} ({\cal K}^{+} (v))$
Jordan blocks of size $(q+1)$. Varying $q$ from $0$ to $p$ gives Jordan blocks of the Jordan form of
$d^{+}(v)$ prescribed by the partition $\lambda^{(p)} (v)$. Varying $p$ from $0$ to $\LG-1$
gives all Jordan blocks in the Jordan form of $d^{+}(v)$.
\end{pf}

The following lemma records some properties of the partitions
 $\lambda^{(p)} (v)$'s.
\begin{lem}\label{part-lam-v}
The partitions $\lambda^{(p)} (v)$, for $p=0,\ldots, \LG-1$, are subject to the following properties.
\begin{enumerate}
\item[1)]
$\lambda^{(p)} (v) =(1^{\mu_{0p} (v)} 2^{\mu_{1p} (v)} \ldots (p+1)^{\mu_{pp} (v)})$.
Its weight $\left|\lambda^{(p)} (v) \right|$ (the sum of all its parts) and its length
$l(\lambda^{(p)} (v))$ (the number of parts) are as follows
$$
\left|\lambda^{(p)} (v) \right| =\sum^p_{q=0} (q+1) \mu_{qp} (v), \,\,
l(\lambda^{(p)} (v)) = \sum^p_{q=0}  \mu_{qp} (v) = dim (({\cal K}^{+})^p (v) )\,.
$$
\item[2)]
The components $h^p_{\GA}$, for $p\leq \LG-1$, of the Hilbert vector  $h_{\GA}$ (see Lemma \ref{hg=c} for definition) of $\GA$  are as follows
$$
h^p_{\GA}  = \sum_{t\geq p} \sum^t_{s=t-p} \mu_{st} (v) = \sum_{t\geq p} (\lambda^{(p)} (v))^{\prime}_{t-p+1}\,,
$$
where 
$\lambda^{\prime}$ denotes the partition conjugate of a partition $\lambda$.
\item[3)]
$\sum^{\LG-1}_{p=0} \left|\lambda^{(p)} (v) \right| =d^{\prime}_{\GA}$
where $d^{\prime}_{\GA}$ is the rank of $\FTP$.
\end{enumerate}
\end{lem}
\begin{pf}
The part 1) is the definition of $\lambda^{(p)} (v) $ and all the formulas follow directly from it.
The part 3) is a consequence of Lemma \ref{part-union}. To see part 2) we use the diagrammatic representation of partitions
following the convention of \cite{[Mac]}.

Consider the vertical levels numbered from left to right by integers
$\{0,\ldots,\LG-1\}$. This should be thought as a visualization of the orthogonal decomposition of 
$\FTP$ in (\ref{ord-Fpr-til}). Place the diagram $\lambda^{(t)} (v) $ so that its first column is at the vertical level labeled by $t$.
Reflect the diagram with respect to this level (the rows of the reflected diagram now go from the right, starting at the level $t$, to the left).
The value $h^p_{\GA}$, which is the rank of $\HH^p$, has now the following pictorial description:
it is the total number of boxes of the reflected diagrams which one finds on the level $p$. It is clear that such boxes are contributed by the 
rows of the partitions  $\lambda^{(t)} (v) $'s, for $t\geq p$, which arrive to the level $p$. Hence, for $t\geq p$, the contribution of 
$\lambda^{(t)} (v) $ is given by parts which are $\geq (t-p+1)$. This gives the first equality in 2).

The second equality in 2) is an obvious transposition: the boxes of rows of $\lambda^{(t)} (v)\,\,(t\geq p) $ counted above form
the $(t-p+1)$-st column of $\lambda^{(t)} (v) $ or, equivalently, the $(t-p+1)$-st row of the transposed partition
$(\lambda^{(t)} (v))^{\prime}$.
\end{pf}

The partition $\lambda(v)$ corresponding to $d^{+} (v)$ reflects the presence of grading given by the orthogonal decomposition
of $\FTP$ in (\ref{ord-Fpr-til}). So it will be useful to formalize some of the properties of 
$\lambda(v)$ by a notion of {\it graded} partition.
\begin{defi}\label{gr-part}
Let $n$ be a positive integer and let
$\overrightarrow{h}=(h_0,\ldots, h_{l-1})$ be a composition of $n$ whose parts $h_p$ are positive for all $p\in \{0,\ldots, l-1\}$.
A partition $\lambda$ of $n$ is said to be $\overrightarrow{h}$-graded, if a sequence of partitions
$\{\lambda^{(p)} \}_{p=0,\ldots, l-1}$ is given such that
\begin{enumerate}
\item[a)]
$\lambda = \bigcup^{l-1}_{p=0} \lambda^{(p)}$\,,
\item[b)]
the diagram of each partition $\lambda^{(p)}$ is contained in the diagram of the partition
$((p+1)^{h_p})$, which is an $h_p \times (p+1)$ rectangle.
\end{enumerate}
The set of all $\overrightarrow{h}$-graded partitions of $n$ will be denoted by $P_n (\overrightarrow{h})$.
\end{defi}
\begin{rem}\label{part-lam=gr}
The partition $\lambda(v)$ associated to $d^{+} (v)$ comes together with the sequence of partitions 
$\{\lambda^{(p)} (v) \}_{p=0,\ldots, \LG-1}$ defined in Lemma \ref{part-lam-v}, 1).
Taking the integer $d^{\prime}_{\GA}$ and its composition
$\overrightarrow{h^{\prime}}_{\GA} = (h^0_{\GA}, \ldots, h^{\LG-1}_{\GA})$, and using 
Lemmas \ref{part-union}, \ref{part-lam-v}, we obtain the conditions a) and b) of Definition \ref{gr-part}.
Thus the partitions 
$\lambda(v)$ belong to $P_{d^{\prime}_{\GA}} (\overrightarrow{h^{\prime}}_{\GA})$ and the assignment of 
$\lambda(v)$ to $([Z],[\alpha],v) \in \TPI$ defines a map
\BEN\label{map-gr-part}
\TPI \longrightarrow P_{d^{\prime}_{\GA}} (\overrightarrow{h^{\prime}}_{\GA})\,.
\EEN
 The image of this map will be denoted by
$P^a_{d^{\prime}_{\GA}} (\overrightarrow{h^{\prime}}_{\GA})$ and its elements will be called {\rm{admissible}}
$\overrightarrow{h^{\prime}}_{\GA}$-graded partitions of $d^{\prime}_{\GA}$.
\end{rem}

From Lemma \ref{part-union} it follows that we have a bijection between the set of admissible partitions 
$P^a_{d^{\prime}_{\GA}} (\overrightarrow{h^{\prime}}_{\GA})$ and the set of admissible matrices
${\bf B}^{+adm}_{\LG} ({\bf Z_{+}})$ defined in 
Definition \ref{matr-adm}. In particular, given an admissible graded partition
$\lambda \in P^a_{d^{\prime}_{\GA}} (\overrightarrow{h^{\prime}}_{\GA})$, we denote by
$M_{\lambda}$ the corresponding admissible matrix. In view of this identification
Proposition \ref{pro-strat} can be reformulated as follows.
\begin{pro}\label{pro-strat1}
Let $\GA$ and $\TPI$ be as in Proposition \ref{pro-strat}. 
Then there is a stratification of $\TPI$ given by the finite union of locally closed sets 
$T^{\lambda}_{\pi}$ indexed by the admissible partitions $\lambda \in P^a_{d^{\prime}_{\GA}} (\overrightarrow{h^{\prime}}_{\GA})$.
Each stratum
$T^{\lambda}_{\pi}$ parametrizes points 
$([Z],[\alpha],v)$ of $\TPI$ for which the Jordan form of $d^{+} (v)$ is constant and prescribed by the partition 
$\lambda$. This stratum is the same as the one in Proposition \ref{pro-strat} and it is labeled by the $d^{+}$-multiplicity matrix
$M_{\lambda}$.
\end{pro}
\begin{cor}\label{g-part-mat}
There is a unique stratum $ T^{\lambda}_{\pi}$ in Proposition \ref{pro-strat1}, which is a non-empty 
Zariski open subset of  $\TPI$. The corresponding partition in 
$P^a_{d^{\prime}_{\GA}} (\overrightarrow{h^{\prime}}_{\GA})$ will be denoted by
$\lambda_{\GA}$ and the corresponding $d^{+}$-multiplicity matrix $M_{\lambda_{\GA}}$ will be denoted by
$M_{\GA}$.
\end{cor}
\begin{pf}
This follows immediately from the irreducibility of $\TPI$.
\end{pf}

\section{Configurations and theirs equations}\label{sec-equations}

In previous sections we have seen how the nilpotent elements of 
{\BM
$\LAG$},
given by the values of morphisms $d^{\pm}$ (see (\ref{d+}) and (\ref{d-})), give rise to very rich algebraic and geometric structures
on $\JA$ (to be more precise on the relative tangent sheaf ${\cal T}_{\pi}$). These are ${\bf sl_2}$-structures in 
\S\S\ref{sec-sl2},\ref{sec-sl2bis} and  the stratification in \S\ref{sec-strat} respectively.
However, it is not clear yet that these nilpotent elements are useful for elucidating the properties of configurations on $X$,
as it is the case with central elements of
{\BM
$\CG$, the center of $\LAGT$} (see, for example, Corollary \ref{Z-c-dec}, Theorem \ref{th-Lie-dec}).
This section addresses this question and it can be viewed as a concrete application of the theory developed so far.
 Namely, we show how to use ${\bf sl_2}$-subalgebras associated to the nilpotent elements, given by the values of morphisms
$d^{\pm}$, to write down equations defining configurations arising from geometric considerations.

The way to produce these equations is somewhat evocative of the classical method of Petri, 
 which gives explicit equations of hypersurfaces of degree 2 (quadrics) and 3 (cubics) through
a canonical curve.\footnote{see e.g., Mumford's survey, \cite{[Mu]}, 
and the references therein for more details.} We also give explicit equations for hypersurfaces (of all degrees) passing through a given
configuration.  The equations, in general, might be quite complicated and not very illuminating. What is essential and different 
%to be able to say, paraphrasing Mumford, 
%``I have seen every geometrically  interesting configuration of points on an algebraic surface once" (compare to Mumford's
%``I have seen every curve once" in  \cite{[Mu]}, p.17).  
 in our approach
is that the main ingredient in getting those equations is representation theoretic. Namely, we exploit the decomposition
of the space of functions on a configuration into the irreducible ${\bf sl_2}$-submodules under the action of 
${\bf sl_2}$-subalgebras, associated to the nilpotent elements of 
{\BM
$\LAG$}
defined by the values of $d^{\pm}$.  This point of view on obtaining equations of projective embeddings, to our knowledge, is new
and seems to be quite fruitful for gaining insight into projective properties of configurations of points on surfaces, geometry of curves on surfaces
and surfaces themselves (see \S\S\S\ref{sec-m00}, \ref{sec-K3}, 9.6). 

\subsection{Geometric set-up}\label{set-up}
In this subsection we recall a geometric context of our constructions.
Let $\GA$ be a component in $\CSA$ and let $\ZA$ be a point in $\JABG$ (see \S\ref{breve} for notation). 
This gives a short exact sequence of sheaves on $X$
\BEN\label{exs-again}
\xymatrix@1{
0 \ar[r] & {\OO_X} \ar[r]& {\SE}_{[\alpha]} \ar[r] & {\ID_{Z}} (L) \ar[r]& 0 }
\EEN
corresponding to the extension class 
$\alpha \in \EZ= Ext^1 (\ID_{Z} (L), \OO_X )$, where $\ID_Z$ is the sheaf of ideals of $Z$ on $X$.
The sheaf 
${\SE}_{[\alpha]}$ sitting in the middle of (\ref{exs-again}) is locally free of rank $2$ with Chern classes
\BEN\label{ch-cl}
c_1 ({\SE}_{[\alpha]}) =L\,\,and \,\,c_2 ({\SE}_{[\alpha]}) =d\,.
\EEN

From (\ref{exs-again}) it follows that
${\SE}_{[\alpha]}$ comes with a distinguished global section which we call $e$. This is the image of 
$1 \in \HO X)$ under the monomorphism in (\ref{exs-again}). The epimorphism in that sequence can be now identified with
the exterior product with the section $e$, i.e. the sequence (\ref{exs-again}) can be alternatively viewed as the Koszul sequence for the pair
$({\SE}_{[\alpha]},e)$
\BEN\label{exs-Koszul}
\xymatrix@1{
0 \ar[r] & {\OO_X} \ar[r]^e& {\SE}_{[\alpha]} \ar[r]^{\wedge e} & {\ID_{Z}} (L) \ar[r]& 0\,. }
\EEN

We will now give a description of the space $\EZ$ in terms of geometry related to 
$(\SE_{[\alpha]},e)$. For this we assume that it has another global section, say $e^{\prime}$, such that the subscheme
\BEN\label{C-curve}
C=(e \wedge e^{\prime} =0)
\EEN
is a smooth irreducible curve.\footnote{this holds if $\SE_{[\alpha]}$ is, for example, generated by its global sections.} This situation is given by the 
following exact sequence of sheaves 
\BEN\label{C-sheaf}
\xymatrix@1{
0 \ar[r] & {\OO_X} \oplus {\OO_X} \ar[r]^(0.55){(e, e^{\prime})}& {\SE}_{[\alpha]} \ar[r] & {\OO_C (\left.L\right|_{C} - Z)} \ar[r]& 0\,, }
\EEN
where $ Z$ and $\left.L\right|_{C}$, the restriction of $L$ to $C$, are viewed as divisors on $C$.
\begin{lem}\label{ext=ls}
Let $(\SE_{[\alpha]},e),\,Z,$ and $C$ be as above and assume $X$ to be a regular surface, i.e.
the irregularity $q(X) =h^1 (\OO_X) =0$. Then one has the following. 
\begin{enumerate}
\item[1)]
The restriction of sections $e,\, e^{\prime}$ to $C$ give rise to two sections $s,\, s^{\prime}$ of the line bundle $\OO_C (Z)$ on $C$.
The subspace $P(s,\, s^{\prime})$ of $H^0(\OO_C (Z))$ spanned by these sections generates $\OO_C (Z)$, i.e. the linear pencil
$\left| P(s,\, s^{\prime})\right|$ on $C$ is base point free.
\item[2)]
There is a natural identification
\BEN\label{ext=ls1}
\EZ \cong  H^0(\OO_C (Z)) / {\CC s}\,.
\EEN
\end{enumerate}
\end{lem}

\noindent
{\it Proof.} Taking the restriction of the sequence in (\ref{C-sheaf}) to $C$ gives the sequence
\BEN\label{E-on-C}
\xymatrix@1{
0 \ar[r] &  {\OO_C (Z)} \ar[r]& {\SE}_{[\alpha]} \otimes \OO_C \ar[r] & {\OO_C (\left.L\right|_{C} - Z)} \ar[r]& 0 }
\EEN
from which it follows that the monomorphism in (\ref{C-sheaf}) factors through $ {\OO_C (Z)}$. Furthermore, the resulting
morphism
$$
{\OO_X} \oplus {\OO_X} \longrightarrow \OO_C (Z)
$$
is surjective. Denoting the image of $(1,0) \in H^0 (\OO_X) \oplus H^0 (\OO_X)$ (resp. $ (0,1)$) by $s$ (resp. $ s^{\prime}$),
we obtain the first part of the lemma.

The second part can be seen as follows. Put together the Koszul sequence in (\ref{exs-Koszul}) and the sequence (\ref{C-sheaf}) to obtain
the following commutative diagram
\BEN\label{Koszul-C-sheaf}
\xymatrix{
0 \ar[r] & {\OO_X} \ar[r]^e \ar[d]& {\SE}_{[\alpha]} \ar[r]^{\wedge e} \ar@{{}{=}{}}[d]& {\ID_{Z}} (L) \ar[r] \ar[d]& 0\\
0 \ar[r] & {\OO_X} \oplus {\OO_X} \ar[d] \ar[r]^(0.55){(e, e^{\prime})}& {\SE}_{[\alpha]} \ar[r] & {\OO_C (\left.L\right|_{C} - Z) }\ar[r]& 0\\
  &{\OO_X}& & &  }
\EEN
This yields the following short exact sequence
\BEN\label{Z-C}
\xymatrix@1{
0\ar[r] &{\OO_X} \ar[r]^(0.45){w}&  {\ID_{Z}} (L) \ar[r]&  {\OO_C (\left.L\right|_{C} - Z)} \ar[r]& 0\,, }
\EEN
where $w =e \wedge e^{\prime}$ is viewed as a section of $\OO_X (L)$ defining the curve $C$.

Tensoring with the canonical bundle 
$\OO_X (K_X)$ of $X$ we obtain
\BEN\label{Z-C-K}
\xymatrix@1{
0\ar[r] &{\OO_X (K_X)} \ar[r]^(0.45){w}&  {\ID_{Z}} (L+K_X) \ar[r]&  {\OO_C (L+K_X)\otimes \OO_C ( - Z)} \ar[r]& 0\,. }
\EEN

By the adjunction formula $\OO_C (L+K_X) = \Omega_C$ is the canonical bundle of $C$, so the last sheaf in 
(\ref{Z-C-K}) is $\Omega_C \otimes \OO_C ( - Z)$. With this in mind we consider the long exact sequence of cohomology groups of
(\ref{Z-C-K}) to obtain
$$
\xymatrix@1{
0 \ar[r]&  H^1(\ID_{Z} (L+K_X) )\ar[r]&  H^1(\Omega_C \otimes \OO_C ( - Z)) \ar[r]& H^2(\OO_X (K_X)) \ar[r] &0\,, }
$$
where the injectivity on the left is the assumption of $X$ being regular and the surjectivity on the right comes from our basic vanishing
conditions on $\OO_X (L)$ in (\ref{vc}). Dualizing this last sequence and using Serre duality give the following
\BEN\label{coh-Z-C}
\xymatrix{
0 & {\EZ} \ar[l]& H^0(\OO_C (Z)) \ar[l]& H^0( \OO_X) \ar[l]& 0 \ar[l]\,, }
\EEN
where the injection on the right is given by the multiplication by the section $s$ of $\OO_C (Z)$ defined in the part 1) of the lemma.
This yields the asserted identification
 $$\EZ \cong  H^0(\OO_C (Z)) / {\CC s}\,. \qquad \qquad \Box$$
\\
\\
\indent
From the identification (\ref{ext=ls1}) it follows that a point $\ZA$ in $\JABG$ with the sheaf
$\SE_{[\alpha]}$ subject to the hypothesis of Lemma \ref{ext=ls} is equivalent to the following geometric set-up:
\begin{center}
\begin{enumerate}
\item[1)]
 a smooth irreducible curve $C$ in the linear system $\left| L \right|$ together with a line bundle $\OO_C (D)$
\item[2)]
 A flag $[s] \in \left| P(s) \right|  \subset \left| \OO_C (D) \right|$ in the linear system
of $\OO_C (D)$ subject to the following properties:
\begin{enumerate}
\item[(i)]
 $s$ is a global section of $\OO_C (D)$ such that $(s=0) =Z$
\item[(ii)]
 $\left| P(s) \right| $ is a
base point free linear pencil  
  such that
the line $P(s) /{\CC s}$  in $H^0 (\OO_C (D)) /{\CC s}$ corresponds to the line $\CC\alpha$ in 
$\EZ$ under the isomorphism in Lemma \ref{ext=ls}.
\end{enumerate}
\end{enumerate}
\end{center}
Equivalently, these data defines the morphism
\BEN\label{kappa-C}
\kappa_C : C \longrightarrow \PP(H^0(\OO_C (D))^{\ast})
\EEN
together with a pencil $\left| P(s) \right|$ of hyperplanes with the one, $H_{s}$, corresponding to $s$. In particular, the
morphism 
$\kappa \ZA$ in (\ref{kappa}) is nothing but the restriction $\left.\kappa_C \right|_{Z}$ of $\kappa_C$ to $Z$
\BEN\label{kappa-C-Z}
\left.\kappa_C \right|_{Z} : Z \longrightarrow  H_s = \PP(\big{(}H^0(\OO_C (D))/{\CC s} \big{)}^{\ast})
\EEN
together with the explicit identifications
\BEN\label{two-ident}
\HT\ZA =\left\{ \left.\frac{\gamma}{\alpha} \right| \gamma \in \EZ \right\} =\left\{ \left.\frac{t}{s^{\prime}} \,(mod\, I_Z) \right| t\in H^0(\OO_C (D)) \right\}= 
H^0(\OO_C (D))/{\CC s}\,,
\EEN
where $s^{\prime}$ is a lifting to $H^0(\OO_C (D))$ of the generator of the line
$P(s) /{\CC s}$ corresponding to $\alpha$ under the isomorphism in Lemma \ref{ext=ls}, 2), 
 and $I_Z$ stands for the ideal of rational functions on $C$ vanishing on $Z$.
The first equality in (\ref{two-ident}) comes from (\ref{HT=fr}), while the second comes from the isomorphism in 
Lemma \ref{ext=ls}, 2). 

The identifications in (\ref{two-ident}) allow us to view the symmetric algebra
$S^{\bullet}( \HT\ZA)$ as the ring of (non-homogeneous) polynomials on $(H^0(\OO_C (D))/{\CC s})^{\ast}$,
while the filtration $\HT_{-\bullet} \ZA$ contains all the information about the quotient of this polynomial ring by the ideal of
polynomials vanishing on $Z^{\prime}$, the image of $Z$ under the morphism $\left.\kappa_C \right|_{Z}$ in 
(\ref{kappa-C-Z}). Hence the relevance of our filtration 
$\HT_{-\bullet}$ in (\ref{filtHT-JG}) to the morphism
$\kappa_C$ in (\ref{kappa-C}) and to the geometry of curves on $X$. 
\begin{rem}
One can go from the geometric set-up above to the bundle (of rank 2) point of view by taking a two dimensional subspace $P$
of $H^0 (\OO_C (D))$ such that the corresponding linear subsystem $\left| P \right|$ of $\left| \OO_C (D) \right|$ is base point free.
This is done by looking at the morphism of sheaves on $X$
\BEN\label{ls-bdl}
P \otimes \OO_X \longrightarrow \OO_C (D)\,.
\EEN
The kernel of this morphism is a sheaf of rank $2$ on $X$ and the surjectivity in (\ref{ls-bdl}) guaranties that the kernel, call it ${\cal S}$, is
locally free. Thus we have a short exact sequence 
\BEN\label{ls-bdl1}
\xymatrix{
0 \ar[r]& {\cal S} \ar[r] & P \otimes \OO_X \ar[r]& \OO_C (D) \ar[r]& 0\,,}
\EEN
whose dual is the sequence analogous to the one in (\ref{C-sheaf}). In particular, the dual 
${\cal S}^{\ast}$ of ${\cal S}$ is a locally free sheaf on $X$ having rank $2$ and the Chern invariants
$L$ and $d =deg(D)$. Furthermore, it comes with a distinguished two-dimensional subspace $P^{\ast}$ of 
sections. 

There seem to be an incongruity here, since taking a pair $({\cal S}^{\ast}, P^{\ast})$ and going back to the geometric set-up,
the space $P^{\ast}$ is identified with a subspace of 
$H^0 (\OO_C (D))$ which should be the one we started with. However, this space is $P$. This seeming contradiction is settled by 
the fact that $dim(P)=2$ and $P^{\ast}$ can be identified with $P$ upon choosing an isomorphism
$\bigwedge^2 P \cong \CC$.

The construction of  the vector bundle ${\cal S}$ in (\ref{ls-bdl}) is a two-dimensional analogue of Lazarsfeld's construction in
\cite{[Laz]} and in Tyurin's work, \cite{[Ty]}. In this two dimensional form it also was used by Donagi and Morrison in \cite{[D-Mo]}.
\end{rem} 
\subsection{${\bf sl_2}$-basis of $\HO{Z^{\prime}})$}\label{sec-sl2-basis}
Let $\GA$ be a component of $\CSA$ and assume it to be simple (Definition \ref{s-c}). From 
\S\ref{sec-Lie1}, Theorem \ref{ty-dec-cont}, it follows that this assumption is inessential and we make it to simplify the discussion only.

Let $\ZA$ be a point of $\JABG$ and consider the filtration
$\HT_{-{\bullet}}$ in (\ref{filtHT-JG}) at $\ZA$. We know that in the resulting filtration
$\HT_{-{\bullet}} \ZA$ of $\HO Z)$ the subspace
$\HT_{-\LG} \ZA$ is a subring and, by Corollary \ref{cor-Fpr}, it depends on $[Z]$ only. This subspace
is isomorphic to the space of functions on $Z^{\prime}$, the image of $Z$ under the morphism
\BEN\label{kappa1}
\kappa\ZA: Z \longrightarrow  Z^{\prime} \subset \PP({\HT\ZA}^{\ast})
\EEN
 in Remark \ref{kap}, (\ref{kappa}),
 i.e. one has
\BEN\label{Zpr-bis}
  \HO{Z^{\prime}}) \cong \HT_{-\LG} \ZA \subset \HO Z)\,,
\EEN
where the isomorphism is given by the pullback $\kappa^{\ast}\ZA$.
 
We aim at writing down equations defining the image $Z^{\prime}$ in the projective space \linebreak$\PP({\HT\ZA}^{\ast})$.
 The main ingredient of our approach is representation theoretic. It consists of using
 our Lie algebraic considerations from \S\ref{sec-sl2} to construct a particular basis for 
$\HO{Z^{\prime}})$. Namely, let $v$ be a {\it non-zero} vertical tangent vector\footnote{`vertical' as usual refers to being tangent along the fibres of the natural projection
$\pi:\JABG \longrightarrow \GAB$.} of $\JABG$ at $\ZA$
 and consider the endomorphism
\BEN\label{d+v-Z}
d^{+}(v): \HO{Z^{\prime}})\cong \bigoplus^{\LG-1}_{p=0} \HH^p \ZA \longrightarrow \bigoplus^{\LG-1}_{p=0} \HH^p \ZA \cong \HO{Z^{\prime}})
\EEN
which is the value of the morphism $d^{+}$ in (\ref{d+}) at the point $([Z],[\alpha],v)$ of the relative tangent bundle $\TPI$ of $\pi$.

We know that $d^{+}(v)$ is nilpotent. Let $\lambda (v)$ be the partition of $d^{\prime}_{\GA} = deg(Z^{\prime})$ corresponding to $d^{+}(v)$.
From Lemma \ref{part-union} we know that it can be written as follows
$$
\lambda (v) = \bigcup^{\LG-1}_{p=0} \lambda^{(p)} (v)\,,
$$
where each $\lambda^{(p)} (v) =(1^{\mu_{0p} (v)} 2^{\mu_{1p} (v)} \ldots (p+1)^{\mu_{pp} (v)})$ is as in Lemma \ref{part-lam-v}, 1).
Since $v$ is fixed in this discussion we will often omit it in the notation above.

From the properties of ${\bf sl_2}$-representations it follows that for every $q\leq p$ with $\mu_{qp} \neq 0$,
we can choose elements 
$y^{(1)}_{qp}, \ldots, y^{(\mu_{qp})}_{qp}$ in the summand $\HH^{p-q} \ZA$ of the orthogonal decomposition of $ \HO{Z^{\prime}})$
in (\ref{d+v-Z}) such that the family of vectors
\BEN\label{set-Bqp}
B_{qp} (v) =\left\{ \left.(d^{+}(v))^m y^{(s)}_{qp} \right| 1\leq s \leq \mu_{qp},\,\, 0\leq m \leq q \right\}
\EEN
is linearly independent in $ \HO{Z^{\prime}})$. Taking the union
\BEN\label{union-set-Bqp}
B(v) = \bigcup_{(q,p)} B_{qp} (v)
\EEN
over all $(q,p)$ subject to $0\leq q \leq p \leq \LG-1$, with the convention that $B_{qp} =\emptyset$, whenever $\mu_{qp} =0$,
we obtain a basis of 
$\HO{Z^{\prime}})$.

In the next step of our construction we modify our basis by replacing the operator $d^{+} (v)$ by multiplication by an appropriate  element
of $\HT\ZA$. More precisely, recall the isomorphism $M$ in (\ref{M}) of the relative tangent sheaf 
${\cal T}_{\pi}$ and $\HT /{\OO_{\JABG}}$. We take $M^{-1} (v) \in \HT\ZA /{\CC\{1_Z\}}$ and let $\tilde{v}$ be its lifting to
$\HT\ZA$. By Remark \ref{val-d} 
$$
d^{+} (v) =D^{+} (\tilde{v})\,,
$$
where $D^{+} (\tilde{v})$ is the positive component of the multiplication operator $D(\tilde{v})$ as in (\ref{d-Dt}). Replacing the operator
$d^{+} (v)$ by the multiplication by $\tilde{v}$ we obtain the families of vectors
\BEN\label{set-Bqp-til}
{\tilde{B}}_{qp}(\tilde{v}) = \left\{ \left.(\tilde{v})^m y^{(s)}_{qp} \right| 1\leq s \leq \mu_{qp},\,\, 0\leq m \leq q \right\}
\EEN
which are still linearly independent in $\HO{Z^{\prime}})$. Taking the union 
\BEN\label{union-Bqp-til}
\tilde{B} (\tilde{v}) = \bigcup_{(q,p)} {\tilde{B}}_{qp} (\tilde{v})
\EEN
gives us a basis of $\HO{Z^{\prime}})$. Furthermore, by construction the basis is adapted to the filtration 
$\HT_{-\bullet} \ZA$ of  $\HO{Z^{\prime}})$ in a sense that the set
\BEN\label{set-Bp-til}
{\tilde{B}}_p (\tilde{v}) = \left\{ \left.(\tilde{v})^m y^{(s)}_{qp^{\prime}} \in B_{qp^{\prime}}  (\tilde{v}) 
 \right| m+p^{\prime} -q \leq p,\,\, \mu_{qp^{\prime}} \neq 0,\,\,1\leq s \leq \mu_{qp^{\prime}} \right\}
\EEN
is a basis of the subspace 
$\HT_{-p-1} \ZA$, for every $p=0,\ldots,\LG-1$. In particular, the set
\BEN\label{B0}
{\tilde{B}}_0 (\tilde{v}) =\left\{ \left.y^{(s)}_{pp} \right| 1\leq s \leq \mu_{pp},\,\, 0\leq p \leq \LG-1, \,\, \mu_{pp} \neq 0 \right\}
\EEN
is a basis of $\HT_{-1} \ZA = \HT \ZA$ (this last equality comes from the definition of the filtration 
$\HT_{-\bullet}$ in  (\ref{HT-i})). 
\begin{rem}\label{rem-B0}
To simplify the notation the elements $y^{(s)}_{pp}$ will be denoted by $y^{(s)}_p$. Thus in the sequel
the basis in (\ref{B0}) will be given in the following form
\BEN\label{B0-1}
{\tilde{B}}_0 (\tilde{v}) =\left\{ \left. y^{(s)}_{p} \right| 1\leq s \leq \mu_{pp},\,\, 0\leq p \leq  \LG-1, \,\, \mu_{pp} \neq 0 \right\}\,.
\EEN
\end{rem}

We summarize the above discussion in the following statement.
\begin{pro-defi}\label{sl2-basis}
Let $\GA$ be a simple component in $\CSA$ and let $\ZA$ be a point in $\JABG$
together with the filtration $\HT_{-\bullet} \ZA$ and the orthogonal decomposition
$$
\HT_{-\LG} ([Z]) = \bigoplus^{\LG-1}_{p=0} \HH^p \ZA\,.
$$ 
Let $t$ be an element of $\HT \ZA$ which is non-constant, viewed as a function on $Z$,
and let $D^{+}(t)$ be the positive component in the decomposition of the operator of multiplication by $t$
(see (\ref{d-Dt})). Let $\lambda(t) =\bigcup^{\LG-1}_{p=0} \lambda^{(p)}(t)$ be the partition associated to $D^{+} (t)$, where
the partitions
$\lambda^{(p)}(t)$ are defined  (see Lemma \ref{part-lam-v}) as follows\footnote{if no ambiguity is likely, the parameter $t$ will be omitted from the above notation.}
$$
\lambda^{(p)}(t) = (1^{\mu_{0p} (t)} 2^{\mu_{1p} (t)} \ldots (p+1)^{\mu_{pp} (t)}),\,\,\mbox{ for every}\,\,
p \in \{0,\ldots,\LG-1 \}\,.
$$
Then for every $(q,p)$, with $0\leq q \leq p\leq \LG-1$ and $\mu_{qp} \neq 0$, there exist elements
$y^{(1)}_{qp},\ldots, y^{(\mu_{qp})}_{qp} \in \HH^{p-q} \ZA$ subject to the following properties:
\begin{enumerate}
\item[1)]
the set of elements
$$
{\tilde{B}}_{qp}(t) = \left\{ \left.t^m y^{(s)}_{qp} \right| 1\leq s \leq \mu_{qp},\,\, 0\leq m \leq q \right\}
$$
is linearly independent in $\HT_{-\LG} ([Z])$ and is contained in $\HT_{-p-1} \ZA$;
 \item[2)]
the union 
$$
\tilde{B} (t) =\bigcup_{(q,p)} {\tilde{B}}_{qp}(t)
$$
is a basis of $\HT_{-\LG} ([Z])$;
\item[3)]
for every $p \in \{0,\ldots,\LG-1 \}$, the set
$$
{\tilde{B}}_p (t) = \left\{ \left.t^m y^{(s)}_{qp^{\prime}} \in B_{qp^{\prime}}  (t) 
 \right| m+p^{\prime} -q \leq p,\,\, \mu_{qp^{\prime}} \neq 0,\,\,1\leq s \leq  \mu_{qp^{\prime}} \right\}
 $$
is a basis of the subspace 
$\HT_{-p-1} \ZA$. In particular, the set
$$
{\tilde{B}}_0 (t) =\left\{ \left.y^{(s)}_{p} \right| 1\leq s \leq \mu_{pp},\,\, 0\leq p \leq  \LG-1, \,\, \mu_{pp} \neq 0 \right\}
$$
is a basis of $\HT \ZA$.
\end{enumerate}

The basis ${\tilde{B}} (t)$ in 2) will be called an ${\bf sl_2}$-basis of $\HT_{-\LG} ([Z])$ associated to $D^{+} (t)$. In view of the identification
of $\HT_{-\LG} ([Z])$ with $\HO{Z^{\prime}})$ in (\ref{Zpr-bis}) this basis will be also called an ${\bf sl_2}$-basis of 
$\HO{Z^{\prime}})$ associated to $D^{+} (t)$.
\end{pro-defi}

\subsection{Equations defining $Z^{\prime}$}\label{sec-conf-eqns}
We now return to the morphism $\kappa \ZA$ in (\ref{kappa1}) and show how to use the basis 
${\tilde{B}} (t)$ in Proposition-Definition \ref{sl2-basis} for writing down equations defining the image
$Z^{\prime}$ of $\kappa \ZA$. The main idea is very simple. The basis ${\tilde{B}}_0 (t)$ can be used to construct
the monomial basis for the symmetric algebra
$S^{\bullet}( \HT\ZA )$. Restricting a monomial to $Z^{\prime}$ gives an element of
$\HO{Z^{\prime}})$. If it is non-zero we can express it uniquely with respect to the basis 
${\tilde{B}} (t)$. Such expressions lead to non-homogeneous  equations defining
$Z^{\prime}$ in $\HT\ZA^{\ast}$. Furthermore, using the identification  
of $\HT\ZA$ with $\EZ$ in  (\ref{HT=Ext}) one obtains homogeneous equations defining
$Z^{\prime}$ in the projective space $\PP((\EZ)^{\ast})$.

To realize the strategy outlined above we fix the set of indeterminates
\BEN\label{set-indet}
{\bf Y} =\left\{ \left.Y_{sp} \right|  0\leq p \leq \LG-1,\,\,\mu_{pp} \neq 0,\,\,1\leq s \leq \mu_{pp} \right\}\,.
\EEN
Observe that the indexing corresponds to the basis elements in ${\tilde{B}}_0 (t)$, the basis of
$\HT\ZA$ in Proposition-Definition \ref{sl2-basis}, 3). So one should think of
${\bf Y}$ as a basis for the space of linear functionals on $(\HT\ZA)^{\ast}$. Thus when we evaluate them on 
$Z^{\prime}$ we obtain a basis of $\HT\ZA$. We agree on the following matching
\BEN\label{Y-y}
 \left. Y_{sp}\right|_{Z^{\prime}} =y^{(s)}_p\,,
\EEN
for all $ Y_{sp} \in {\bf Y}$.

The element $t$ lies in $\HT\ZA$ and hence can be expressed uniquely in terms of the basis ${\tilde{B}}_0 (t)$
$$
t=\sum_{sp} c_{sp} y^{(s)}_p\,.
$$
This implies that $t$ is the restriction to $Z^{\prime}$ of the linear function
\BEN\label{T-t}
T=\sum_{sp} c_{sp} Y_{sp}\,.
\EEN

 The monomials in the indeterminates $Y_{sp}$ give a basis for the algebra of polynomial functions on
$\HT\ZA^{\ast}$.
Let ${\bf m}=(m_{sp}) $ be a multi-degree (the indexing is the same as for indeterminates) and let 
\BEN\label{mon-m}
Y^{\bf m} = \prod_{s,p} Y^{m_{sp}}_{sp}
\EEN
be the corresponding monomial. Substituting for $Y_{sp}$ the elements $y^{(s)}_p$'s of the basis  ${\tilde{B}}_0 (t)$ in (\ref{B0-1}),
  we obtain the element
\BEN\label{y-m}
y^{\bf m} = \prod_{s,p} (y^{(s)}_p)^{m_{sp}}
\EEN
in $\HO{Z^{\prime}})$. If we set 
$$
\left| {\bf m} \right| =\sum_{s,p} m_{sp}
$$
to be the total degree of $Y^{\bf m}$ in (\ref{mon-m}), then we have
$y^{\bf m} \in \HT_{-\left| {\bf m} \right|} \ZA$. If this element is zero, then of course 
$Y^{\bf m} =0$ is already an equation for $Z^{\prime}$ in 
$\HT\ZA^{\ast}$. Otherwise, $y^{\bf m}$ can be expressed uniquely in terms of the basis
${\tilde{B}}_{\mid {\bf m} \mid -1} (t)$ (see Proposition-Definition \ref{sl2-basis})
\BEN\label{y-m1}
y^{\bf m} = \sum a^{ks}_{ji} ({\bf m}) t^k y^{(s)}_{ji}\,,
\EEN
where the sum is taken over $(j,i,k,s)$ such that $k+i-j \leq \left| {\bf m} \right| -1$ and $1\leq s \leq \mu_{ji}$.

For every element $y^{(s)}_{qp}$ of the basis ${\tilde{B}}(t)$
choose its lifting  $P^{(s)}_{qp}$ in the symmetric algebra 
$S^{\bullet}( \HT\ZA)$, i.e. 
$P^{(s)}_{qp}$ is a polynomial function in the indeterminates of the set
${\bf Y}$ which,  when restricted to $Z^{\prime}$, is equal to $y^{(s)}_{qp}$. Since
by definition $y^{(s)}_{qp}$ lies in $\HH^{p-q} \ZA \subset \HT_{-(p-q)-1}$ (for this inclusion see (\ref{ordH-i})),
we can choose
$P^{(s)}_{qp}$ in $S^{p-q+1}( \HT\ZA)$. With such a choice made once and for all and the lifting $T$ of $t$
 in (\ref{T-t}), we obtain 
\BEN\label{eqns}
F({\bf m}) = Y^{\bf m} - \sum a^{ks}_{qp} ({\bf m}) T^k P^{(s)}_{qp}
\EEN
a polynomial vanishing on $Z^{\prime}$.

At this stage the polynomials $F({\bf m})$ are non-homogeneous. To have homogeneous polynomials
vanishing on $Z^{\prime}$, recall that the elements of 
$\HT \ZA$ can be identified as fractions of elements in $\EZ$ (see (\ref{HT=fr})). So to homogenize our polynomials we
set
\BEN\label{set-indet1}
{\bf T} =\left\{ \left.T_{sp} \right|  0\leq p \leq \LG-1,\,\,\mu_{pp} \neq 0,\,\,1\leq s \leq \mu_{pp} \right\}
\EEN
to be the basis of linear functionals on $(\EZ)^{\ast}$ corresponding to the basis ${\bf Y}$ in (\ref{set-indet})
under the isomorphism in (\ref{HT=Ext}).
Since $\alpha \in \EZ$, it also can be viewed as a linear functional on $(\EZ)^{\ast}$ and as such it will be denoted by $T_{\alpha}$.
This and (\ref{HT=fr}) yield the following fractional form for the elements in ${\bf Y}$ 
\BEN\label{Y-T}
Y_{sp} =\frac{T_{sp}}{T_{\alpha}}\,.
\EEN
Substituting into (\ref{eqns}) and multiplying by $T^{\mid {\bf m} \mid}_{\alpha}$ yields
\BEN\label{hom-eqns}
H({\bf m}) =T^{\left| {\bf m} \right|}_{\alpha} F({\bf m}) = T^{\bf m} -
 \sum a^{ks}_{qp} ({\bf m}) {\tilde{T}}^k  Q^{(s)}_{qp}  T^{\left| {\bf m} \right| -k -p+q-1}_{\alpha}\,,
\EEN 
where 
$T^{\bf m} =\prod_{(s,p)} (T_{sp})^{m_{sp}}$ is the monomial of multi-degree ${\bf m}$ in the set of indeterminates 
${\bf T}$,  $\tilde{T} $ is the linear form corresponding to $T$, i.e. $T=\frac{ \tilde{T}}{ T_{\alpha}}$,
and $Q^{(s)}_{qp} =T^{p-q+1}_{\alpha}P^{(s)}_{qp}$ the homogenized form of 
$P^{(s)}_{qp}$. These $H({\bf m})$'s are  now homogeneous forms (in the set of
indeterminates ${\bf T} $) of degree $\left| {\bf m} \right|$ vanishing on $Z^{\prime}$. 

Of course in such a generality this is of limited use. However, the above considerations
 give an algebro-geometric interpretation of the multiplicities $\mu_{qp}$ in the definition of the partitions
$\lambda^{(p)} (t)$ or, equivalently, of the corresponding multiplicity matrix (see Definition \ref{mult-mat-def}).
\begin{pro}\label{eqns-mult}
Let $\ZA$ and $t$ be as in Proposition-Definition \ref{sl2-basis} and let $\mu_{qp} (t)$ be the multiplicities associated
to the nilpotent endomorphism $D^{+} (t)$ (see Proposition-Definition \ref{sl2-basis} for notation).

Let $\mu_{qp} (t) \neq 0$ and let $Q^{(s)}_{qp}\ ,(s=1,\ldots, \mu_{qp} (t))$  be the homogeneous polynomials of degree
$(p-q+1)$ appearing in the equations (\ref{hom-eqns}). Then there exist homogeneous forms 
$A^{(s)}_{qp}$ subject to the following properties.
\begin{enumerate}
\item[1)]
The degree  $\nu^{(s)}_{qp}= deg( A^{(s)}_{qp})$ is at most $( p+1)$.
\item[2)]
Let $\tilde{T}$ and $T_{\alpha}$ be as in (\ref{hom-eqns}) then the homogeneous form
\BEN\label{eqns-mult1}
G^{(s)}_{qp} =( \tilde{T})^{q+1} Q^{(s)}_{qp} -T^{p+2-\nu^{(s)}_{qp}} _{\alpha} A^{(s)}_{qp}
\EEN
of degree $(p+2)$ vanishes on $Z^{\prime}$, for every $s \in \{1,\ldots, \mu_{qp} (t) \}$.
Furthermore, these forms are linearly independent in
$S^{p+2} \EZ$. In particular,
$dim (I_{p+2} (Z^{\prime})) \geq \mu_{qp}$, for all $0\leq q \leq p$, where
$I_n (Z^{\prime})$ stands for the subspace of homogeneous polynomials of degree $n$ on
$(\EZ)^{\ast}$ vanishing on $Z^{\prime}$.
\end{enumerate}
\end{pro}
\begin{pf}
By definition of the elements $y^{(s)}_{qp}$ in the set $B_{qp}$ in (\ref{set-Bqp}) we have
\BEN\label{eq1}
(D^{+} (t))^{q+1} (y^{(s)}_{qp}) =0,\,\,\forall 1\leq s \leq \mu_{qp}\,.
\EEN
Replacing $D^{+} (t)$ by the multiplication by $t$, yields the following 
\BEN\label{eq2}
 t^{q+1} y^{(s)}_{qp} \in \HT_{-p-1} \ZA,\,\,\forall 1\leq s \leq \mu_{qp}\,.
\EEN
Using the notation in (\ref{eqns}) we obtain polynomials in the set of indeterminates ${\bf Y}$ (see (\ref{set-indet}))
\BEN\label{eq3}
T^{q+1} P^{(s)}_{qp} - B^{(s)}_{qp},\,\,for\,\, s\in \{1,\ldots, \mu_{qp} \}
\EEN
which vanish on $Z^{\prime}$ (here $T$ and $P^{(s)}_{qp}$ are as in (\ref{eqns})) and where
$B^{(s)}_{qp}$ is a polynomial whose degree is at most $(p+1)$. Homogenizing, as it was done
 in (\ref{hom-eqns}), gives the forms $G^{(s)}_{qp}$, as asserted in (\ref{eqns-mult1}), where $A^{(s)}_{qp}$, the homogenization of the polynomials
$B^{(s)}_{qp}$, are subject to the asserted properties.

To see the linear  independence of the forms $G^{(s)}_{qp}\,(s=1,\ldots, \mu_{qp} (t))$ in (\ref{eqns-mult1})
consider the relation in $S^{p+2} \EZ$
$$
\sum^{\mu_{qp} (t)}_{s=0} c_s G^{(s)}_{qp} = 0\,,
$$
for some constants $c_s \in \CC$.
Substituting the expressions of $G^{(s)}_{qp} $ from (\ref{eqns-mult1}) yields
\BEN\label{eq4}
( \tilde{T})^{q+1} \sum^{\mu_{qp} (t)}_{s=0} c_s Q^{(s)}_{qp} - 
\sum^{\mu_{qp} (t)}_{s=0} c_s T^{p+2-\nu^{(s)}_{qp}} _{\alpha} A^{(s)}_{qp} =0\,.
\EEN
From this it follows that
\BEN\label{eq5}
\sum^{\mu_{qp} (t)}_{s=0} c_s Q^{(s)}_{qp} =T^m_{\alpha} F\,,
\EEN
for some integer $m\geq 1$ and some homogeneous polynomial $F$ of degree $(p-q+1 -m)$. 
Dehomogenizing, i.e. dividing by $T^{p-q+1}_{\alpha}$, and evaluating on $Z^{\prime}$ yields the following
\BEN\label{eq6}
\sum^{\mu_{qp} (t)}_{s=0} c_s y^{(s)}_{qp} \in \HT_{-(p-q+1-m)}\ZA \subset  \HT_{-(p-q)}\ZA\,.
\EEN
But by construction the elements $y^{(s)}_{qp} \,(s=1,\ldots, \mu_{qp} (t))$ are linearly independent in 
$\HH^{p-q} \ZA$ which is orthogonal to $\HT_{-(p-q)}\ZA$ in 
$\HT_{-(p-q)-1}\ZA$ (see (\ref{ordH-i})). This implies that the constants $c_s$ must be all equal to $0$.
\end{pf}
\begin{rem}\label{eqns-geo}
 Assume that $\ZA$ satisfies the geometric set-up described in \S\ref{set-up}. This means that $Z$ lies on a smooth curve
$C$ in the linear system $\left| L \right|$ and there is a base point free linear pencil
$\left| P(\sigma,\sigma^{\prime}) \right|$ in  $\left| L \right|$ generated by two global sections $\sigma,\sigma^{\prime}$ of
$\OO_C (Z)$ subject to the following conditions
\begin{enumerate}
\item[1)]
$Z =(\sigma=0)$,
\item[2)]
the quotient space $P(\sigma,\sigma^{\prime}) /{\CC\sigma} \subset H^0 (\OO_C (Z)) /{\CC\sigma}$ is identified with the line
$\CC\alpha \subset \EZ$ under the isomorphism (\ref{ext=ls1}) in Lemma \ref{ext=ls}, and we assume that the extension
class $\alpha$ corresponds to the coset $\sigma^{\prime}\,mod(\CC\sigma)$ under this identification.
\end{enumerate}
In particular, in this set-up we have the following description of $\HT\ZA$
\BEN\label{EZ=ls}
\HT \ZA \cong \left\{ \left.\frac{x}{\sigma^{\prime}} \,(mod \,I_Z) \right| x\in H^0 (\OO_C (Z)) \right\}\,,
\EEN
where $I_Z$ stands for the ideal of rational functions on $C$ vanishing on $Z$. Thus the morphism $\kappa \ZA$ in (\ref{kappa1}) is the restriction
of the morphism $\kappa_C$ in (\ref{kappa-C}) to $Z$. Equivalently, the configuration
$Z^{\prime} =\kappa_C (Z)$ is the hyperplane section of the curve $C^{\prime}$, the image of $\kappa_C$, cut out by the hyperplane
$H_{\sigma}$ of $\PP(H^0 (\OO_C (Z))^{\ast})$, corresponding to the section $\sigma$.
 Thus the equations in (\ref{hom-eqns}) are equations defining a
hyperplane section of the curve $C^{\prime}$. More precisely, the set ${\bf T}$ in (\ref{set-indet1}) is a set of homogeneous coordinates
in 
$H_{\sigma} =\PP((H^0 (\OO_C (Z))/ {\CC\sigma})^{\ast})$. In particular, the linear form
$T_{\alpha}$ in (\ref{hom-eqns}) can be identified with the linear form $\overline{\sigma^{\prime}}$,
the restriction of $\sigma^{\prime}$ to the subspace $(H^0 (\OO_C (Z))/ {\CC\sigma})^{\ast}$ of
$H^0 (\OO_C (Z))^{\ast}$ (here we view $\sigma^{\prime}$ as a linear function on 
$H^0 (\OO_C (Z))^{\ast}$). With this in mind the homogeneous forms in (\ref{hom-eqns}) take the following form
\BEN\label{hysec-eqns}
H({\bf m}) = T^{\bf m} -
 \sum a^{ks}_{qp} ({\bf m}) {\tilde{T}}^k  Q^{(s)}_{qp}  {\overline{\sigma^{\prime}}}^{\left| {\bf m} \right| -k -p+q-1}\,.
\EEN 
\end{rem}

\subsection{$\mu_{00}$ and multi-secant planes} \label{sec-m00}
In Proposition \ref{eqns-mult} we have given an interpretation of the multiplicities
$\mu_{qp} (t)$ in terms of equations defining $Z^{\prime}$ (all the notation in that proposition are preserved and used freely here).
In this subsection we give an interpretation of the first multiplicity $\mu_{00} (t) $ in terms of geometric properties of 
the configuration $Z^{\prime}$ in the projective space
$\PP(\HT\ZA^{\ast})$.

We take a non-constant function $t$ in $\HT\ZA$ and consider the endomorphism
\BEN\label{D+t}
D^{+} (t) : \HO{Z^{\prime}}) \longrightarrow \HO{Z^{\prime}})\,.
\EEN
Its kernel $K^{+}(t)$ is the same as the fibre of the sheaf ${\cal K}^{+}$ in (\ref{ker+}) at the point
$([Z],[\alpha], M(\overline{t}))$ of $\TPI$, where $\overline{t}$ is the projection of $t$ to 
$\HT\ZA /{\CC}$ and
$M(\overline{t})$ is the image of $\overline{t}$ under the isomorphism $M$ in (\ref{M}). In particular,
$K^{+}(t)$ admits the orthogonal decomposition
$$
K^{+}(t) =\bigoplus^{\LG-1}_{p=0} (K^{+}(t) )^p\,,
$$
where $(K^{+}(t) )^p$ is the fibre of the sheaf $({\cal K}^{+})^p$ (see (\ref{ord-K+}) for notation) at $([Z],[\alpha], M(\overline{t}))$.
From Lemma \ref{lem-w-filt-Kp} it follows that the multiplicity $\mu_{00} (t)$ is equal to the  dimension of the summand
$(K^{+}(t) )^0$.
\begin{lem}\label{t=D0t}
The multiplication by $t$ restricted to $(K^{+}(t) )^0$ coincides with the operator $D^0 (t)$
in the triangular decomposition (\ref{d-Dt}) and hence gives rise to the linear map
$$
D(t): (K^{+}(t) )^0 \longrightarrow \HT\ZA\,.
$$
\end{lem}
\begin{pf}
The operator $D^{-} (t)$ vanishes on the summand $\HH^0 \ZA$. Hence $D(t)$ restricted to $(K^{+}(t) )^0$ equals $D^0 (t)$.
\end{pf}

The assertion of the above lemma written out explicitly gives the following relation
\BEN\label{t=D0t1}
tx =D^0 (t)(x),\,\,\forall x\in  (K^{+}(t) )^0\,. 
\EEN

We will now apply this relation to a particular choice of $t$ to derive some geometric consequences.
First we observe that $\mu_{00} (t)$ is an upper continuous function of $t\in \HT\ZA$ and it takes integer values in
the interval\footnote{recall: $\HH^0 \ZA =\HT\ZA$ and $dim \HT\ZA =r+1$ is the index of $L$-speciality of $Z$. Hence the upper bound of the interval
of $[1, r+1]$.
The lower bound comes from the fact that $D^{+} (t)$ annihilates the subspace $\CC\{1_Z \} \subset \HT\ZA$ of constant functions.}
 $[1, r+1]$. Also observe that the upper bound is achieved precisely when $t \in \CC1$. This follows from the assumption of $\GA$ being simple.
Set
\BEN\label{m00-Z}
\mu_{00} \ZA =min\left\{\left.\mu_{00} (t) \right|  t\in \HT\ZA \right\}\,.
\EEN
From what is said above this constant is the {\it generic} value of $d^{+}$-multiplicity $\mu_{ 00} (t)$ as $t$ varies through
$\HT\ZA$, i.e. this value is taken on a non-empty Zariski open subset of $\HT\ZA$. 
The following result gives a geometric meaning of this number.
\begin{lem}\label{m00-geom}
 Let $Z^{\prime}$ be the image of the morphism $\kappa\ZA$ as in (\ref{kappa1}). Then $Z^{\prime}$ admits a decomposition
$$
Z^{\prime} = Z^{\prime}_1 +Z^{\prime}_2
$$
subject to the following properties.
\begin{enumerate}
\item[1)]
$Z^{\prime}_1$ and $Z^{\prime}_2$ are disjoint.
\item[2)]
$Z^{\prime}_1$ spans a hyperplane in $\PP( \HT\ZA^{\ast})$.
\item[3)]
$Z^{\prime}_2$  spans the projective subspace  of $\PP( \HT\ZA^{\ast})$ whose dimension
is at most $r+1-\mu_{00} \ZA$. 
\end{enumerate}
\end{lem}
\begin{pf}
Choose a subset $Z^{\prime}_0$ of $Z^{\prime}$ consisting of $r$ points spanning a hyperplane in $\PP( \HT\ZA^{\ast})$ and let
$t$ be a linear function on  $\HT\ZA^{\ast}$ defining the span of $Z^{\prime}_0$. Restricting $t$ to $Z^{\prime}$ gives a non-constant
function on $Z^{\prime}$ which we continue to denote by $t$. By construction it belongs to $\HT\ZA$ and it vanishes on the subset
$Z^{\prime}_0$. 
 We now apply the relation
(\ref{t=D0t1}) to it. This implies that 
$D^0 (t) (x)$ vanishes on $Z^{\prime}_0$  as well. Hence it must be a scalar multiple of $t$ and we obtain
\BEN\label{scalar}
   D^0 (t) (x) =\xi(x) t, \,\,\forall x\in (K^{+} (t))^0\,,
\EEN
where $\xi(x)$ is a scalar. From this it follows that $\xi$ is a linear function on $ (K^{+} (t))^0$. Hence
$$
 D^0 (t) (x) =0,\,\,\forall x\in ker(\xi)\,.
$$
Substituting this in (\ref{t=D0t1}) we obtain
\BEN\label{tx=0}
tx=0,\,\,\forall x\in ker(\xi)\,.
\EEN
This relation yields the asserted decomposition. Indeed, define $Z^{\prime}_1$ to be the subset of $Z^{\prime}$, where
$t$ vanishes. Observe that $Z^{\prime}_1$ is a proper subset of $Z^{\prime}$ containing $Z^{\prime}_0$. In particular,
the span of $Z^{\prime}_1$ is a hyperplane in $\PP( \HT\ZA^{\ast})$.

Define $Z^{\prime}_2$ to be the complement of $Z^{\prime}_1$ in $Z^{\prime}$. From (\ref{tx=0}) it follows that
$x\in \HT\ZA$ vanishes on $Z^{\prime}_2$ if and only if $x\in ker(\xi)$. Hence the span of 
$Z^{\prime}_2$ in $\PP( \HT\ZA^{\ast})$ is a projective subspace of dimension $(r+1-\mu_{00} (t))$.
\end{pf}
\begin{cor}\label{cor-m00}
Let $t$ be as in the proof of Lemma \ref{m00-geom} and let
$$
Z^{\prime} =Z^{\prime}_1 +Z^{\prime}_2
$$
be the decomposition of $Z^{\prime}$ corresponding to $t$.
Then one obtains the decomposition
$$
Z=Z_1 + Z_2
$$
of $Z$, where $Z_i =(\kappa\ZA)^{-1} (Z^{\prime}_i)\,(i=1,2)$. Furthermore, the indexes of 
$L$-speciality of these subconfigurations are as follows
$$
\delta(L,Z_1) =\mu_{00} (t)-1 \,\,and\,\,\delta(L,Z_2) =1\,.
$$
\end{cor}
\begin{pf}
The argument is the same as in the proof of Corollary \ref{Z-c-dec}. Namely, set $\JJ_{Z_i}\,(i=1,2)$ to be the ideal sheaves on $Z$ 
of subconfigurations $Z_i\,(i=1,2)$. Then we have the direct sum decomposition
$$
\HO Z) =H^0 (\JJ_{Z_1}) \oplus H^0 (\JJ_{Z_2}) 
$$
and the isomorphisms 
$$
\HO{Z_1}) \cong H^0 (\JJ_{Z_2}),\,\,\,\,\, \HO{Z_2}) \cong H^0 (\JJ_{Z_1})\,. 
$$
All this follows from a diagram similar to the one in (\ref{Zlam-Zlamc}).
Continuing to argue as in the proof of Corollary \ref{Z-c-dec}, (ii), we have that the index of $L$-speciality of $Z_1$ (resp. $Z_2$)
is computed as the dimension of the space 
$H^0 (\JJ_{Z_2}) \bigcap \HT\ZA$ (resp. $H^0 (\JJ_{Z_1}) \bigcap \HT\ZA$). From the proof of 
Lemma \ref{m00-geom} it follows
$$
 H^0 (\JJ_{Z_2}) \bigcap \HT\ZA =ker(\xi) \,\,and\,\,H^0 (\JJ_{Z_1}) \bigcap \HT\ZA= \CC \{t\}\,.
$$
Hence the assertion of the corollary.
\end{pf}
Applying the above result to the geometric set-up in \S\ref{set-up} we obtain the following.
\begin{cor}\label{m00-setup}
Let
$C$ be a smooth curve in the linear system $\left|L \right|$ and let $Z$ be a configuration of degree $d$ on $C$ such that
the line bundle $\OO_C (Z)$ is special and base point free.  Then we can choose a base point free linear pencil
$\left| P(\sigma, \sigma^{\prime})\right| $ in the linear system $\left| H^0(\OO_C (Z)) \right|$ which defines a point 
$\ZA$ of the Jacobian $\JA$ as described in Remark \ref{eqns-geo} and identifies the space
$\HT\ZA$ with the space of fractions of the form
$$
f=\left.{\frac{x}{\sigma^{\prime}}} \right|_Z,\,\, for \,\, x\in H^0(\OO_C (Z))
$$
(see (\ref{EZ=ls})).
Furthermore, there exists a non-zero section $x \in H^0(\OO_C (Z))$
such that the element
$t=\displaystyle{\left.{\frac{x}{ \sigma^{\prime}}}\right|_Z} \in \HT\ZA$
defines the decomposition 
$$
Z=Z_1 +Z_2
$$
as in Corollary \ref{cor-m00}, where $Z_i \,(i=1,2)$ are special divisors on $C$
such that
$h^0 (\OO_C (Z_1)) =\mu_{00} (t)$ and $h^0 (\OO_C (Z_2)) = 2$.
\end{cor}
\begin{pf}
Only the last assertion needs to be proved. The speciality of $Z_i \,(i=1,2)$ follows from the speciality of $Z$.
To see the  formulas for $h^0 (\OO_C (Z_i))\,(i=1,2)$, use the exact sequence analogous to the one in 
(\ref{Z-C-K}) with $Z_{1,2}$ in place of $Z$ there. The associated cohomology sequence gives
\BEN\label{m00-setup1}
\xymatrix@1{
0\ar[r]& H^1 (\ID_{Z_i} (L+K_X) \ar[r]& H^1 ( \Omega_C \otimes \OO_C (-Z_i)) \ar[r]& H^2 (\OO_X (K_X)) \ar[r]& 0\,, }
\EEN
for $i=1,2$, where we used the assumption of the geometric set-up of $X$ being a regular surface (see Lemma \ref{ext=ls}).
This exact sequence together with Serre duality and  the definition of index of $L$-speciality (see (\ref{is})) gives the formula
$$
h^0 (\OO_C (Z_i)) = \delta (L,Z_i) +1\,,
$$
for $i=1,2$. Substituting the values of $\delta (L,Z_i) $ from Corollary \ref{cor-m00} yields the asserted formulas.
\end{pf}

The next result gives the value of generic multiplicity $\mu_{00} \ZA$ (see (\ref{m00-Z}))
 in the case $Z$ is in general position with respect to the adjoint linear
system
$\left| L+K_X \right|$.
\begin{cor}\label{m00-gp}
Let $\ZA$ be a point of $\JA$, where
 $Z$ is a configuration of $d$ points in general position with respect to the adjoint linear
system
$\left| L+K_X \right|$.  Assume that the index of $L$-speciality
$\delta(L,Z) =r+1 \geq 2$ and $deg Z \geq 2r+1$. Then
$\mu_{00} \ZA =1$.
\end{cor}
\begin{pf}
From \RI, Corollary 7.13, it follows that the map
$\kappa \ZA$ in (\ref{kappa1}) is an embedding and its image
$Z^{\prime}$ is a set of $d$ distinct points in general position in the projective space 
$\PP(\HT\ZA^{\ast})$.

Assume $\mu_{00} \ZA \geq 2$. Applying Lemma \ref{m00-geom} yields a decomposition
\BEN\label{m00-geom-dec}
Z^{\prime} =Z^{\prime}_1 + Z^{\prime}_2\,,
\EEN
where $Z^{\prime}_1$ and $Z^{\prime}_2$ span respectively a hyperplane and a subspace of dimension
$\leq (r+1- \mu_{00} \ZA) \leq r-1$ in $\PP(\HT\ZA^{\ast})$, where the second inequality is the consequence of the assumption 
$\mu_{00} \ZA \geq 2$. The fact that $Z^{\prime}$ is in general position in $\PP(\HT\ZA^{\ast})$ implies 
that $deg Z^{\prime}_i \leq r$, for $i=1,2$. This together with decomposition in (\ref{m00-geom-dec})
imply
$$
d=deg Z^{\prime} =deg Z^{\prime}_1 + deg Z^{\prime}_2 \leq 2r
$$
contrary to the hypothesis that $d\geq 2r+1$.
\end{pf}

\subsection{Complete intersections on $K3$-surfaces}\label{sec-K3}
In this subsection we  apply our theory to configurations which are complete intersections on a $K3$-surface.
In particular, we give a complete set of very simple explicit quadratic equations defining such configurations - 
the quadrics in question are of rank $\leq 4$. This in turn leads to recovering
quadrics through canonical curve which have a much more simpler form then the ones obtained by Petri's method (see \cite{[Mu]}).

Let $X$ be a $K3$-surface and let $\OO_X (L)$ be a very ample line bundle on $X$. Consider a configuration $Z$ on $X$ which is a 
complete intersection of two smooth curves $C_1$ and $C_2$ in the linear system $\left| L \right|$. Let $\gamma_i \,(i=1,2)$ be sections of 
$H^0 (\OO_X (L)$ defining the curves $C_i \,(i=1,2)$, i.e. $C_i =(\gamma_i =0)$, for $i=1,2$. 
 It has been shown in \RI, \S5.2, that the space of extensions $\EZ$ is identified as follows
\BEN\label{ext-ci}
\EZ = H^0 (\OO_X (L)) / {\CC\{\gamma_1, \gamma_2\}}\,.
\EEN
It has been also shown that the orthogonal decomposition of $\HO Z)$ at a point $\ZA$ for a general choice of 
$\alpha$ has the following form
\BEN\label{ord-ci}
\HO Z) = \HH^0 \ZA \oplus  \HH^1 \ZA \oplus  \HH^2 \ZA\,,
\EEN
where $dim  \HH^0 \ZA =\frac{L^2 }{2},\,\, dim  \HH^1 \ZA =\frac{L^2 }{2} -1\,\,dim\HH^2 \ZA =1$.

Let $Z$ be a general complete intersection as above. Then it is well-known that it is in general position  in 
$\PP((\EZ)^{\ast})$.
\begin{lem}\label{part-ci}
Let $([Z],[\alpha],v)$ be a  closed point of the relative tangent sheaf of $\JA$ with respect to the projection $\pi$ in (\ref{pi1}),
with $Z$ being a general complete intersection as above. Then, for a sufficiently general vector $v$, the partition 
$\lambda (v)$ associated to $d^{+} (v)$ (see \S\ref{sec-strat} for details and notation) 
is as follows
$$
\lambda (v) =( 3 \,2^{\frac{d}{2} -2}\, 1)\,,
$$
where $d=L^2 =degZ$. Furthermore, the partitions $\lambda^{(p)}$'s $(p=0,1,2)$ in the decomposition (\ref{part-union1}) are as follows
$$
\lambda^{(0)} =(1), \,\,\lambda^{(1)}=(2^{\frac{d}{2} -2}),\,\,\lambda^{(2)}= (3)\,.
$$
\end{lem}
\begin{pf}
The weight of the orthogonal decomposition, i.e. the number of summands in it, is $3$. So
$(d^{+} (v))^3 =0$ and the parts of $\lambda (v)$ are at most $3$. Since $dim\HH^2 \ZA =1$ it follows that the multiplicity of $3$ 
can be at most $1$ and, for a general $v$, it must be $1$, because $(d^{+} (v))^2 \neq 0$, for a general $v$.
 
To see the multiplicities of $1$ and $2$ we consider the restriction of $d^{+} (v)$ to the summand 
$\HH^0 \ZA$
\BEN\label{d+v0}
d^{+} (v) : \HH^0 \ZA \longrightarrow \HH^1 \ZA\,.
\EEN
From Corollary \ref{m00-gp} it follows that  the kernel 
$(K^{+}(v))^0$ of this homomorphism is precisely the subspace of constants $\CC\{1_Z \}$.
Hence multiplicity of $1$ in the partition $\lambda (v)$ is $1$. Since the dimension of
$\HH^0 \ZA$ equals $\frac{d}{2}$, it follows that the part $2$ occurs in 
$\lambda (v)$ with multiplicity $(\frac{d}{2} -2)$.

The last assertion follows immediately from the definition of the partitions $\lambda^{(p)}$'s in (\ref{part-p}) and the first part of the proof.
\end{pf}

We aim at writing down quadratic equations defining $Z$ in the projective space
$\PP((\EZ)^{\ast}) =\PP((H^0(\OO_X (L) / {\CC\{\gamma_1, \gamma_2\} })^{\ast})$.
Our guide is the general strategy outlined in \S\ref{sec-conf-eqns}. So we begin by writing quadratic
equations determining the image of $\kappa\ZA$ in $\PP(\HT\ZA^{\ast})$ and then use the explicit identifications in
Lemma \ref{ext=ls} to pass to the equations in \linebreak $\PP((H^0(\OO_X (L) / {\CC\{\gamma_1, \gamma_2\} })^{\ast})$.
 We also know from the same lemma that  $\kappa\ZA$, in the case at hand, is an embedding, so we will not distinguish $Z$ and its image
under $\kappa\ZA$.

Set $g=\frac{d}{2} +1$. This is the genus of smooth curves in the linear system
$\left| L \right|$. Fix $g-1$ distinct points $z_1, \ldots, z_{g-1}$ in $Z$ and view them as linear functionals on
$\HT\ZA$. The assumption that $Z$ is in general position is equivalent to 
$\{z_1, \ldots, z_{g-1} \}$ being a basis of
$\HT\ZA^{\ast}$. Let $\{x_1, \ldots, x_{g-1} \}$ be the basis of 
$\HT\ZA$ dual to $\{z_1, \ldots, z_{g-1} \}$, i.e.
\BEN\label{x-z}
x_i (z_j) =  \delta_{ij},\,\,\forall i,j\,.
\EEN

Let $\HH\ZA$ be the subspace of $\HT\ZA$ orthogonal to the constant $1_Z\in \HT\ZA$.
We examine the operators $D^{+} (x_i)$ restricted to $\HH\ZA$.
\begin{lem}\label{D+H} 
$
D^{+}(x_i) : \HH\ZA \longrightarrow \HH^1 \ZA
$
is an isomorphism, for all $i$.
\end{lem} 
\begin{pf}
From the dimensions of the summands in 
(\ref{ord-ci}) it follows that
the spaces $\HH\ZA$ and $\HH^1 \ZA$ have the same dimension $(g-2)$. So it is enough to show the injectivity of
$D^{+}(x_i)$. This in turn is deduced from the proof of Corollary \ref{m00-gp} and the fact that $Z$ is in general position in
$\PP(\HT\ZA^{\ast})$.
\end{pf}
\begin{lem}\label{eq-quad}
For every $i\neq j$ and for every $k$, there exists an element $h_{ijk} \in \HT\ZA$ such that
$q_{ijk} =x_i x_j - x_k h_{ijk}$, viewed as a quadratic polynomial on $\HT\ZA^{\ast}$), vanishes
on $Z$.
\end{lem}
\begin{pf}
From Lemma \ref{D+H} it follows that for all triples $i,j,k$ there exists a unique element
$h^{\prime}_{ijk} \in \HH\ZA$ such that
$$
D^{+}(x_i) (x_j) = D^{+}(x_k)(h^{\prime}_{ijk})\,.
$$
Replacing the operator $D^{+} (x_i)$ (resp. $D^{+} (x_k)$) by multiplication by $x_i$ (resp. $x_k$) we obtain
\BEN\label{eq-quad1}
\overline{x_i x_j - x_k  h^{\prime}_{ijk}} =\overline{m}_{ijk}\,,
\EEN
for some element $m_{ijk} \in \HT\ZA$, where $\overline {a}$, for $a \in S^{\bullet} (\HT\ZA)$, stands for the restriction of
$a$, viewed as a polynomial function on $\HT\ZA^{\ast}$, to $Z$.

Let $i \neq j$. Then the left hand side in (\ref{eq-quad1}) vanishes on the set
$\{z_s \mid s \neq k\}$. Hence $m_{ijk}$ vanishes on this set as well.
Since this set spans the hyperplane in 
$\HT\ZA^{\ast}$ corresponding to the linear functional $x_k$ we obtain
$$
m_{ijk} =c_{ijk} x_k\,,
$$
for some $c_{ijk} \in \CC$. Substituting this in (\ref{eq-quad1}) gives
$$
q_{ijk} = x_i x_j -x_k (h^{\prime}_{ijk} + c_{ijk})
$$
which vanishes on $Z$. Setting $h_{ijk} =h^{\prime}_{ijk} + c_{ijk}$ yields the assertion.
\end{pf}

Fix $k$, say $k=1$, and consider the set of quadratic polynomials
\BEN\label{set-eqns-quad}
Q \ZA =\left\{ \left. q_{ij} =x_i x_j - x_1 h_{ij1} \right| 2\leq i <j \leq g-1 \,\,\right\}\,,
\EEN
where $q_{ij} =q_{ij1}$ are as in Lemma \ref{eq-quad}. This set is non-empty if $g\geq 4$, and then, it gives us $\binom{g-2}{2}$ linear independent polynomials in
$S^2 (\HT\ZA)$. Thus we proved the following.
\begin{lem}\label{Q-k}
Let $g\geq 4$. Then the image $\kappa\ZA: Z\hookrightarrow \PP(\HT\ZA^{\ast})$
lies on 
${g-2 \choose 2}$ linearly independent quadrics given by the set $Q \ZA$ in 
(\ref{set-eqns-quad}).
\end{lem}

For the rest of this discussion we assume $g\geq 4$ and use the considerations of \S\ref{set-up} to identify the map
$\kappa\ZA$ with a hyperplane section of the canonical embedding of one of the smooth curves in 
$\left| L \right|$ passing through $Z$.

Recall that $Z$ is a complete intersection of two smooth curves $C_1$ and $C_2$ in $\mid L \mid$.
Fix $C_1 =(\gamma_1 =0)$, where $\gamma_i\,(i=1,2)$ are sections corresponding to the divisors $C_i\,(i=1,2)$. 
Then by definition $Z$ lies on $C_1$ and the line bundle
$\OO_{C_1} (Z) =\OO_{C_1} (C_2) =\OO_{C_1} (L) = \Omega_{C_1}$ is the canonical
line bundle of $C_1$.

Let $\overline{\gamma_2}$ be the restriction of $\gamma_2$ to $C_1$. Then it is a section of 
$\Omega_{C_1}$ defining $Z$. Applying the identification (\ref{ext=ls1}) to this situation yields
\BEN\label{ext=ls2}
\EZ \cong H^0 (\Omega_{C_1}) /{\CC\overline{\gamma_2}}\,.
\EEN
Then a choice of an extension class $\alpha$ is determined by a choice of another section,
say $\omega_0 \in H^0 (\Omega_{C_1})$, such that the linear pencil generated by 
$\omega_0 $ and $\overline{\gamma_2}$ is base point free and the coset
$\omega_0\,(mod\,\CC\overline{\gamma_2})$ in
$H^0 (\Omega_{C_1}) /{\CC\overline{\gamma_2}}$
goes over to $\alpha$ under the isomorphism in (\ref{ext=ls2}). With this in mind the identification in 
(\ref{two-ident}) yields
\BEN\label{HT=om-fr}
\HT\ZA =\left\{\left. \left.\frac{\omega}{\omega_0} \right|_{Z} \,\, \right| \omega \in H^0 (\Omega_{C_1})\right\}\,,
\EEN
where $\left.\frac{\omega}{\omega_0} \right|_{Z}$ stands for the restriction of the rational function
$\frac{\omega}{\omega_0}$ on $C_1$ to $Z$. In particular, the basis 
$\{x_1, \ldots, x_{g-1} \} $ has the form
$$
x_i =\left.\frac{\omega_i}{\omega_0} \right|_{Z}\,,
$$
where $\omega_1, \ldots, \omega_{g-1}$ are linearly independent sections  of $\Omega_{C_1}$.
Similarly, $h_{ijk}$ in Lemma \ref{eq-quad}, viewed as a function on $Z$, will be of the form
$$
h_{ijk} =\left.\frac{\omega_{ijk}}{\omega_0} \right|_{Z}\,,
$$
for some sections $\omega_{ijk}$ of $\Omega_{C_1}$. Substituting all this in the expressions of
$q_{ij}$ in (\ref{set-eqns-quad}) we obtain the set of quadratic polynomials
\BEN\label{I2Z}
Q_2 (Z) =\left\{\left. \omega_i \omega_j - \omega_1 \omega_{ij1}  \in S^2  H^0 (\Omega_{C_1}) \right| 2\leq i <j \leq g-1,\, \right\}
\EEN
vanishing on $Z$. We now claim that this set gives a basis of the space of quadratic forms vanishing on $Z$.
Indeed, let 
$$
\kappa_{C_1} : C_1 \longrightarrow \PP( H^0 (\Omega_{C_1})^{\ast}) =\PP^{g-1}
$$
be the canonical map of $C_1$. In our case it is an embedding, so we identify $C_1$ and its image under
$\kappa_{C_1}$. Then the  configuration $Z$ is the hyperplane section of $C_1$ obtained by intersecting
$C_1$ with the hyperplane
$H_Z$ corresponding to the section $\overline{\gamma_2} \in H^0 (\Omega_{C_1})$ defining $Z$.
Thus we have
$$
Z \subset H_Z = \PP( (H^0 (\Omega_{C_1}) /{\CC\overline{\gamma_2}})^{\ast})\,.
$$

Set $\ID_Z$ to be the sheaf of ideals of $Z$ in $H_Z$ and denote by
$\overline{\omega_i}$ the restriction of $\omega_i$, viewed as a section of $\OO_{\PP^{g-1}} (1)$,
to the hyperplane $H_Z$.
\begin{pro}\label{quadrics}
The set of quadrics 
$$
\overline{Q}_2 (Z) =
\left\{ \left.\overline{\omega}_i \overline{\omega}_j - \overline{\omega}_1 \overline{ \omega}_{ij1}  \in H^0 (\OO_{H_Z} (2))
 \right| 2\leq i <j \leq g-1\, \right\}
$$
forms a basis of $H^0 (\ID_Z (2))$, the space of quadrics in $H_Z$ vanishing on $Z$.
\end{pro}
\begin{pf}
From Lemma \ref{Q-k} it follows that the set $\overline{Q}_2 (Z)$ consists of
${g-2 \choose 2}$ linearly independent quadrics. The assertion now follows from the fact 
that the  dimension  of $H^0 (\ID_Z (2))$ is equal to
${g-2 \choose 2}$. Indeed, let $\ID_{C_1}$ be the sheaf of ideals of $C_1$ in $\PP^{g-1}$.
Then we have the following exact sequence relating the ideal sheaves of $C_1$ and its hyperplane section $Z$
$$
\xymatrix@1{
0\ar[r]& {\ID_{C_1} (-1)} \ar[r]^{\overline{\gamma}_2}& \ID_{C_1} \ar[r]& \ID_Z  \ar[r]&0 }\,,
$$
where the monomorphism is the multiplication by $\overline{\gamma}_2$, viewed here as a section of 
$\OO_{\PP^{g-1}} (1)$, and the epimorphism is the restriction to the hyperplane $H_Z$.
Tensoring with $\OO_{\PP^{g-1}} (2)$ and taking the associated sequence of cohomology groups gives
the isomorphism
\BEN\label{IC1-IZ}
 H^0 (\ID_{C_1} (2)) \cong  H^0 (\ID_{Z} (2))\,.
\EEN
Now a classical result of Max Noether (see  e.g. \cite{[G-H]}, p.253) yields the count 
$$
h^0 (\ID_{C_1} (2)) = {g-2 \choose 2}\,.
$$
\end{pf}
\begin{rem}\label{rk4}
Observe that the quadrics in $\overline{Q}_2 (Z)$ are all of rank $\leq 4$.
Thus we recover a hyperplane section version of Mark Green's theorem on the generation of the ideal of a canonical curve
by quadrics of rank 4, \cite{[Gr]}.
\end{rem}

Next we lift quadrics from the hyperplane $H_Z$ to $\PP^{g-1}$ and obtain quadrics passing through the curve $C_1$.
More precisely, we go back to the set of sections
$\omega_1, \ldots, \omega_{g-1}$ of $\Omega_{C_1}$. Adding to them the section $\overline{\gamma}_2$
gives a basis for $H^0 (\Omega_{C_1})$.
\begin{cor}\label{quad-C}
Let $\{ \overline{\gamma}_2, \omega_1, \ldots, \omega_{g-1}\}$ be a basis of $H^0 (\Omega_{C_1})$ as above
and let $\omega_{ij} =\omega_{ij1}$, where $\omega_{ij1}$ are as in (\ref{I2Z}). Then there is a unique choice of sections
$\omega^{\prime}_{ij}$ in  $H^0 (\Omega_{C_1})$ such that the quadratic polynomials
$$
Q_2 (C_1) =\left\{\left.\tilde{q}_{ij} = \omega_i \omega_j - \omega_1 \omega_{ij} + \overline{\gamma}_2 \omega^{\prime}_{ij} 
\in S^2  H^0 (\Omega_{C_1}) \right| 2\leq i <j \leq g-1\,\right\}
$$
form a basis of $H^0 (\ID_{C_1} (2))$, the space of quadratic forms vanishing on $C_1$.
\end{cor}
\begin{pf}
Consider $\overline{\gamma}_2$ as a section of $\OO_{\PP^{g-1}} (1)$ defining the hyperplane $H_Z$. This gives the following
exact sequence of sheaves on $\PP^{g-1}$  
$$
\xymatrix@1{
0 \ar[r] & {\OO_{\PP^{g-1}} (-1)} \ar[r]^{\overline{\gamma}_2} & {\OO_{\PP^{g-1}} } \ar[r] &{\OO_{H_Z}} \ar[r]&0\,.   }
$$
Tensoring it with $\OO_{\PP^{g-1}} (2)$ and passing to the cohomology sequence yields 
\BEN\label{rest-quadr}
\xymatrix@1{
0 \ar[r] & H^0(\OO_{\PP^{g-1}} (1)) \ar[r]^{\overline{\gamma}_2} & H^0(\OO_{\PP^{g-1}} (2))  \ar[r] & H^0(\OO_{H_Z} (2))  \ar[r]&0\,.  }
\EEN
The isomorphism in (\ref{IC1-IZ}) implies that for every $2\leq i<j \leq g-1$ there exists unique polynomial
$\tilde{q}_{ij}$ in $H^0(\ID_{C1} (2))$ whose restriction to the hyperplane $H_Z$ gives the polynomial
$(\overline{\omega}_i \overline{\omega}_j - \overline{\omega}_1 \overline{ \omega}_{ij})$.
On the other hand the polynomial
$ \omega_i \omega_j - \omega_1 \omega_{ij}$ has the same restriction to $H_Z$ as 
$\tilde{q}_{ij}$. From the exact sequence (\ref{rest-quadr}) it follows that they differ by a multiple of $\overline{\gamma}_2$ and
this multiple is unique.
\end{pf}

We can lift our equations further, to obtain quadratic equations through the surface itself. Namely, choose a lifting 
$\{\gamma_{i+2}\}_{i=1,\ldots, g-1}$ to
$H^0(\OO_X (L))$ of the elements 
$\{\omega_1,\ldots,\omega_{g-1} \}$ in 
$H^0(\Omega_{C_1}) =H^0 (\OO_{C_1} (L)) =H^0(\OO_X (L)) /{\CC\gamma_1}$. Completing it by $\gamma_1$ and
$\gamma_2$ we obtain a basis
$\{\gamma_1, \gamma_2, \gamma_3, \ldots, \gamma_{g+1} \}$
of $H^0(\OO_X (L))$.
\begin{cor}\label{Q-K3}
Let $\{\gamma_1, \gamma_2, \gamma_3, \ldots, \gamma_{g+1} \}$ be a basis
of $H^0(\OO_X (L))$ as above and let 
$\gamma_{ij}$ and $\gamma^{\prime}_{ij}$ be liftings to $H^0(\OO_X (L))$ of sections
$\omega_{ij}$ and $\omega^{\prime}_{ij}$  in Corollary \ref{quad-C}. Then there is a unique choice of sections
$\gamma^{\prime\prime}_{ij}$ in  $H^0(\OO_X (L))$ such that the homogeneous quadratic polynomials
$$
Q_2 (X) =\left\{\left. 
\gamma_{i+2} \gamma_{j+2} - \gamma_3 \gamma_{ij} + \gamma_2 \gamma^{\prime}_{ij} +\gamma_1 \gamma^{\prime\prime}_{ij}
\in S^2  H^0 (\OO_X (L)) \right| 2\leq i <j \leq g-1\, \right\}
$$
form a basis of  the space of quadratic forms vanishing on $X$.
\end{cor}
\begin{pf}
The argument is analogous to the one in the proof of Corollary \ref{quad-C}, i.e. we relate the embedding of $X$
in the projective space 
$\PP^g=\PP(H^0(\OO_X (L))^{\ast})$ to its hyperplane section $C_1$, determined by the section $\gamma_1$.
This gives the following exact sequence relating the ideal sheaf 
$\ID_X$ of $X$ in $\PP^g$ with the ideal sheaf $\ID_{C_1}$ of $C_1$ in the hyperplane
$H_{\gamma_1} = \PP((H^0(\OO_X (L))/{\CC\gamma_1} )^{\ast})$
$$
\xymatrix@1{
0\ar[r]& {\ID_{X} (-1)} \ar[r]^{\gamma_1}& \ID_{X} \ar[r]& \ID_{C_1}  \ar[r]&0\,. }
$$
Tensoring with $\OO_{\PP^g} (2)$ and taking the resulting cohomology sequence yields an isomorphism
$$
H^0(\ID_{X} (2)) \cong H^0(\ID_{C_1} (2))\,.
$$
This implies that the polynomials $\tilde{q}_{ij}$ in Corollary \ref{quad-C} are the restrictions of the quadrics
$Q_{ij} (1\leq i<j\leq g-1)$ in $S^2  H^0 (\OO_X (L))$, forming a basis of
$H^0(\ID_{X} (2))$. On the other hand, from the expressions of $\tilde{q}_{ij}$ in Corollary \ref{quad-C}, we see that they are the restrictions
the quadratic polynomials
\BEN\label{Q-X}
 \gamma_{i+2} \gamma_{j+2} - \gamma_3 \gamma_{ij} + \gamma_2 \gamma^{\prime}_{ij}\,.
\EEN 
From the exact sequence
$$
\xymatrix@1{
0 \ar[r] & H^0(\OO_{\PP^{g}} (1)) \ar[r]^{{\gamma}_1} & H^0(\OO_{\PP^{g}} (2))  \ar[r] & 
H^0(\OO_{H_{\gamma_1}} (2))  \ar[r]&0  }
$$
relating quadrics in $\PP^{g}$ with quadrics in the hyperplane $H_{\gamma_1}$ it follows that the difference between
$Q_{ij}$ and the polynomial in (\ref{Q-X}) is a multiple of $\gamma_1$. Hence the assertion of the corollary.
\end{pf}
\subsection{Adjoint linear system $\left| L+K_X \right|$ and geometry of $Z$}
It was mentioned that the orthogonal decomposition in (\ref{d+v-Z})
should contain information not only about geometry of the morphism $\kappa\ZA$ in (\ref{kappa1}) but also about geometry of 
configurations with respect to the adjoint linear system  $\left| L+K_X \right|$ (see the discussion in the end of \S\ref{sec-ord}). In this section we take up the considerations
of an ${\bf sl_2}$-triple
associated to the operator $d^{-} (v)$, for $([Z],[\alpha],v)$ in the relative tangent bundle $\TPI$ of $\JABG$, for an admissible component
$\GA$ in $\CS$. It turns out that geometry behind the action of such an ${\bf sl_2}$-triple concerns the image of $Z$ with respect to 
the adjoint linear system  $\left| L+K_X \right|$. Throughout this discussion we assume
\BEN\label{gp-assumption}
Z \,\,\mbox{is in general position with respect to}\,\, \left| L+K_X \right| \,\,\mbox{and}\,\, d > r+2\,.
\EEN
These assumptions imply that the linear system $\left| L+K_X \right|$ defines an embedding
\BEN\label{sysadj-emb}
Z \hookrightarrow \PP((H^0 (L+K_X) / H^0 (\ID_Z (L+K_X) ))^{\ast}) :=\PP^{d-r-2}_Z\,,
\EEN
where $\ID_Z$ is the sheaf of ideals of $Z$ on $X$. So we will not make a distinction between $Z$ and its
image in $\PP^{d-r-2}_Z$. Thus $Z$ will be viewed here, simultaneously, as a configuration of $d$ points on $X$ and 
a configuration of $d$ points in general position in the projective space $\PP^{d-r-2}_Z$.
 
Recall the filtration
\BEN\label{filtFd-Z}
\FT(L)\ZA =\HO Z (L+K_X))=\FI_0 \ZA \supset \FI_1 \ZA \supset \cdots \supset \FI_{\LG} \ZA =0\,,
\EEN
discussed in \S\ref{sec-ord}, (\ref{filtFd}). By construction 
$$
\FI_1 \ZA  = H^0 (L+K_X) / H^0 (\ID_Z (L+K_X) )\,,
$$
 so the filtration in (\ref{filtFd-Z}) gives a filtration of 
$H^0 (L+K_X) / H^0 (\ID_Z (L+K_X) )$
\BEN\label{filtFd-Z1}
H^0 (L+K_X) / H^0 (\ID_Z (L+K_X) ) = \FI_1 \ZA \supset \cdots \supset \FI_{\LG} \ZA =0\,.
\EEN

Arguing as in the case of $d^{+} (v)$ in \S\ref{sec-sl2-basis}, we choose a basis of 
$H^0 (L+K_X) / H^0 (\ID_Z (L+K_X) )$ adapted to the action of an ${\bf sl_2}$-triple associated to $d^{-} (v)$.
More precisely, let
\BEN\label{part-v}
\lambda (v) =(\lambda_1 (v) , \lambda_2 (v) ,  \ldots , \lambda_s (v) )
\EEN
be the partition associated to $d^{-} (v)$ viewed as a nilpotent endomorphism of $\HO Z)$.
We use the usual diagrammatic representation of $\lambda (v)$ as an array of boxes aligned in horizontal rows, from top to bottom, with
$\lambda_i (v)$ boxes in the $i$-th row (see \cite{[Mac]} for details). However, since the operator $d^{-} (v)$ moves the index of the
grading
\BEN\label{gradHOZ}
\HO Z) = \bigoplus^{\LG-1}_{p=0} \HH^p \ZA
\EEN
from right to left, it will be convenient for our purposes to think of boxes in the diagram of $\lambda (v)$ running from
{\it right to left} as well. Thus we fix an ${\bf sl_2}$-triple $\goth{s}$ associated to $d^{-} (v)\,(v\neq 0)$ as described in \S\ref{sec-inv}
and define a basis of $\HO Z)$ adapted to this ${\bf sl_2}$-triple as follows.

We fill in the  boxes in the first column of $\lambda (v)$ with highest weight vectors $f_1, \ldots,f_s$ of the action of  $\goth{s}$ on 
$\HO Z)$, where the vector $f_i$ is placed in the box of the $i$-th row. Then the remaining boxes of $\lambda (v)$ are filled with vectors
of the form $(d^{-} (v))^m (f_i)$, with $m$ running from $1$ to $\lambda_i (v) -1$ in the rows with $\lambda_i (v) \geq 2$.
Furthermore, we can choose vectors $f_i$'s to be homogeneous with respect to the grading in (\ref{gradHOZ}), i.e. each $f_i$ 
belongs to a particular summand, say 
$\HH^{p_i} \ZA$, of the decomposition in (\ref{gradHOZ}). Hence, setting $deg(f_i)=p_i$, the vectors $(d^{-} (v))^m (f_i)$ acquire grading
\BEN\label{grad-basis}
deg((d^{-} (v))^m (f_i) ) =deg (f_i) - m =p_i -m\,.
\EEN
This degree will be attached to the corresponding box of $\lambda (v)$. Thus the degree of boxes in the $i$-th row runs
(from right to left) from $deg(f_i)$ to $deg(f_i) - \lambda_i (v) +1$.

Next we replace the operator $d^{-} (v)$ by the multiplication by an element $\tilde{v}$ in $\HT\ZA$ lifting $v$ (see Remark \ref{val-d}, 2)).
This gives us a basis of $\HO Z)$ composed of vectors 
${\tilde{v}}^m f_i \,(i=1,\ldots,s;\,m=0,\ldots,\lambda_i (v) -1)$, with vector ${\tilde{v}}^m f_i$ sitting in the $i$-th row and 
$(m+1)$-st column (counting from the right).

We now use the identification of $\HO Z)$ with $\HO Z (L+K_X)$ provided by the morphism in (\ref{mor-RT}) (see the discussion following
(\ref{mor-RT})). For every $i \in \{1,\ldots,s\}$, denote by $\phi_i$ the vector in $\HO Z (L+K_X)$ corresponding to $f_i$ 
under this identification. Then the elements
\BEN\label{basis-adjunc}
 {\tilde{v}}^m \phi_i \,(i=1,\ldots,s;\,m=0,\ldots,\lambda_i (v) -1) 
\EEN
form a distinguished basis of $\HO Z (L+K_X))$. Furthermore, this basis is adapted to the filtration $\FI_{\bullet} \ZA$ in (\ref{filtFd-Z}) 
in a sense that ${\tilde{v}}^m \phi_i $ lies in the subspace $\FI_{p_i -m} \ZA$, for every $i$ and $m$ in (\ref{basis-adjunc}). This follows
from the fact that under the identification of $\HO Z)$ with $\HO Z (L+K_X)$ the filtration
$\FI^{\bullet} \ZA$, the fibre at $\ZA$ of the filtration in (\ref{filtF}), corresponds to the filtration  $\FI_{\bullet} \ZA$ 
(\RI, Lemma 2.1).

Since we are interested only in a basis of $\FI_1 \ZA =H^0 (L+K_X) / H^0 (\ID_Z (L+K_X) )$, the basis in (\ref{basis-adjunc}) 
has to be modified by suppressing the elements sitting in the boxes of $\lambda (v)$ having degree $0$.
This means that we need to leave out all those elements in (\ref{basis-adjunc}) which project to a basis of
$\FI_{0} \ZA  / \FI_{1} \ZA \cong \HH^0 \ZA$.
 
On the level of partitions the desired modification is achieved by erasing  the boxes\footnote{the number of boxes erased is $r+1$, the dimension of
$\HH^0 \ZA$.} of  $\lambda (v)$ of degree $0$.
Denote the resulting partition of $d-r-1$ by $\hat{\lambda}(v)$:
\BEN\label{part-tranc}
\hat{\lambda}(v) =(\hat{\lambda}_1 (v) , \hat{\lambda}_2 (v) , \ldots, \hat{\lambda}_{s^{\prime}} (v))
\EEN
and call it truncation of $\lambda (v)$. From the definition it follows that the parts of $\hat{\lambda}(v)$ and $\lambda (v)$ are related as follows
$$
\lambda_i (v) - \hat{\lambda}_i (v) = 0\,\,or\,\,1,\,\,\, \mbox{for every $1\leq i \leq s^{\prime}$}\,,
$$
where $s^{\prime}$ is the length of $\hat{\lambda}(v)$. That length is related to the length $s$ of $\lambda (v)$ by the formula
\BEN\label{length-part}
 s^{\prime} = s- \mu^{\prime}_{(\LG-1)0} (v) = s-\mu_{00} (v)\,,
\EEN
where the first equality is the definition of the $d^{-}$-multiplicities\footnote{the $d^{-}$-multiplicity $\mu^{\prime}_{(\LG-1)0} (v)$ counts
precisely the number of rows of $\lambda(v)$ consisting of a single box having degree $0$.} in (\ref{mult-v-}), while the second comes from (\ref{mult+-}). 

From the considerations above one deduces the following basis for $H^0 (L+K_X) / H^0 (\ID_Z (L+K_X) )$.
\begin{lem}\label{lem-basis-adjunc}
Set $\phi_{im} ={\tilde{v}}^m \phi_i $, for $i=1,\ldots,s^{\prime}$ and $m=0,\ldots, \hat{\lambda}_i (v) -1$, where
${\tilde{v}}^m \phi_i $'s are as in (\ref{basis-adjunc}). Then the  set of elements
\BEN\label{B(v)-adjunc}
B(v) =\left\{ \left.\phi_{im} \in H^0 (L+K_X) / H^0 (\ID_Z (L+K_X) ) \right| i=1,\ldots,s^{\prime};\,\,m=0,\ldots, \hat{\lambda}_i (v) -1 \right\}
\EEN
is a basis of $H^0 (L+K_X) / H^0 (\ID_Z (L+K_X) )$ such that for every $i,m$ in (\ref{B(v)-adjunc}),
$\phi_{im}$ lies in the subspace $\FI_{p_i -m} \ZA$, where $p_i$ is as in (\ref{grad-basis}).
\end{lem}

We will also view elements $\phi_{im}$ in (\ref{B(v)-adjunc}) as homogeneous coordinates of the projective space
$\PP^{d-r-2}_Z$ (see (\ref{sysadj-emb}) for notation). As such they will be denoted by $X_{im}$. Thus thinking of $Z$
as a configuration of points in the projective space $\PP^{d-r-2}_Z$ we have 
\BEN\label{X-phi}
\left.{X_{im}}\right|_Z = \phi_{im},\,\,\forall i,m\,.
\EEN
Our objective now is to use homogeneous coordinates $X_{im}$ to write down equations of subvarieties in 
$\PP^{d-r-2}_Z$ passing through $Z$. The following result illustrates what we have in mind.
\begin{lem}\label{scroll-Z}
Let
$$
B^{\ast}(v) =\left\{\left. V_{im} \in \left(H^0 (L+K_X) / H^0 (\ID_Z (L+K_X) )\right)^{\ast} \right| 
i=1,\ldots,s^{\prime};\,\,m=0,\ldots, \hat{\lambda}_i (v) -1 \right\}
$$
be the basis of $\left(H^0 (L+K_X) / H^0 (\ID_Z (L+K_X) )\right)^{\ast}$ dual to the basis
$\{X_{im} \}_{{i=1,\ldots,s^{\prime}} \atop {m=0,\ldots, \hat{\lambda}_i (v) -1}} $.
\\
\indent
For every part $\hat{\lambda}_i (v) \geq 2$ of the partition $\hat{\lambda}(v)$ in (\ref{part-tranc})
the configuration $Z \subset \PP^{d-r-2}_Z $ is contained  in the cone over a rational normal curve $C_i$ lying in the projective space
\BEN\label{span-Ci}
P_i = Span\left\{\left.V_{ij} \right| j=0,\ldots,\hat{\lambda}_i (v) -1 \right\}
\EEN
with the vertex of the cone being the complementary projective subspace $\Pi_i$ in $\PP^{d-r-2}_Z $ cut out by the hyperplanes
$$
X_{ij} =0, \,\, for \,\,j=0,\ldots,\hat{\lambda}_i (v) -1. 
$$
\end{lem}
\begin{pf}
If $\hat{\lambda}_i (v) =2$, then $P_i = \PP^1$ and $\Pi_i$ is a complementary projective subspace of codimension $2$.
Then the cone in the assertion is the whole projective space $\PP^{d-r-2}_Z $. So we assume $\hat{\lambda}_i (v) \geq 3$
and write the $2\times (\hat{\lambda}_i (v) -1)$-matrix
\BEN\label{matr-cone-Z}
M_i (v) =\left( \begin{array}{ccc}
           X_{i(\hat{\lambda}_i (v) -2)} &\cdots&X_{i0} \\
           X_{i(\hat{\lambda}_i (v) -1)} &\cdots&X_{i1}
               \end{array} \right)\,.
\EEN
The $2\times 2$-minors of this matrix are homogeneous quadratic polynomials in 
$X_{i0}, \ldots,  X_{i(\hat{\lambda}_i (v) -1)}$. Setting them to be equal to zero gives 
${\hat{\lambda}_i (v) -1 \choose 2}$ quadrics in
$\PP^{d-r-2}_Z $ which are all singular along the subspace $\Pi_i$, while their restrictions to the subspace
$P_i$ in (\ref{span-Ci}) cut out a rational normal curve in $P_i$. Hence the $2\times 2$-minors of $M_i (v)$ in 
(\ref{matr-cone-Z}) cut out the cone asserted in the lemma. 

It remains to check that the cone contains $Z$ or, equivalently, that the $2\times 2$-minors of $M_i (v)$
vanish on $Z$. This is insured by the relations in (\ref{X-phi}) and the definition of $\phi_{im}$ in 
Lemma \ref{lem-basis-adjunc}. Indeed, let $z\in Z$. Then its homogeneous coordinates in 
$\PP^{d-r-2}_Z $  are given by
\BEN\label{X(z)}
X_{im} (z) =\phi_{im} (z) =(\tilde{v} (z))^m  \phi_i (z)\,,
\EEN
where the first equality is (\ref{X-phi}) and the second comes from the definition of $\phi_{im}$ in Lemma \ref{lem-basis-adjunc}.
From (\ref{X(z)}) it follows
$$
(X_{im} X_{i(m^{\prime}+1)} -  X_{i m^{\prime} } X_{i(m+1)} ) (z) =
(\tilde{v} (z) )^{m +m^{\prime}+1} \phi^2_i (z) - (\tilde{v} (z) )^{m +m^{\prime}+1} \phi^2_i (z) =0\,,
$$
for every $m\neq m^{\prime}$.
\end{pf}

To complete our considerations we introduce the following geometric realization of a partition.
Let $\mu = (\mu_1 ,\mu_2 , \ldots , \mu_t )$ be a partition with $\mu_t \geq 1$.
To such a partition we associate the vector bundle over $\PP^1$
\BEN\label{part-bundle}
\GS_{\mu} = \bigoplus^t_{i=1} \OO_{\PP^1} (\mu_i -1)\,.
\EEN
Denote by $\PP(\mu)$ the projectivization of the dual $\GS^{\ast}_{\mu}$ of $\GS_{\mu}$. It comes with the natural projection
$$
\pi_{\mu} : \PP(\mu) \longrightarrow \PP^1\,.
$$
On $\PP(\mu)$ we choose $\OO_{\PP(\mu)} (1)$ so that the direct image
$$
 \pi_{\mu \ast} (\OO_{\PP(\mu)} (1) ) =\GS_{\mu}\,.
$$
Observe that $\GS_{\mu}$ is generated by its global sections, so $\OO_{\PP(\mu)} (1)$ defines a morphism 
\BEN\label{morph-part-bund}
\psi_{\mu} :  \PP(\mu) \longrightarrow \PP^{\left| \mu \right| -1}\,,
\EEN
where $\left| \mu \right| = \sum^t_{i=1} \mu_i$ is the weight of the partition $\mu$. Denote by $Y_{\mu}$ the image of 
$\psi_{\mu}$ and call it {\it $\mu$-scroll}.
\begin{rem}\label{mu-scroll}
\begin{enumerate}
\item[1)]
If $\mu_t \geq 2$, then $\OO_{\PP(\mu)} (1)$ is very ample and $Y_{\mu}$ is a rational normal scroll of dimension $t$
in $\PP^{\left| \mu \right| -1}$.
\item[2)]
If $\mu$ contains $1$ with multiplicity $m_1 \geq 1$, then one has the following possibilities
$$
\mu =\left \{ \begin{array}{cc}
(1^t),&if\,\,m_1 =t, \\
(\mu_1 ,\ldots ,\mu_{t-m_1},\underbrace{1,\ldots,1}_{m_1 -times} ), &if\,\,1\leq m_1 <t.
     \end{array} \right.
$$
In the first case $ Y_{\mu} =\PP^{t-1}$ and in the second
$Y_{\mu}$ is the cone over a rational normal scroll $Y_{\mu^{\prime}}$ with the vertex
$Sing( Y_{\mu}) =\PP^{m_1 -1}$ and where 
$$
\mu^{\prime} =(\mu_1 , \ldots , \mu_{t-m_1} )\,.
$$
\end{enumerate}
\end{rem}

If $\mu =(\left| \mu \right| ), \,\,\left| \mu \right|  \geq 2$, then $Y_{\mu}$ is a rational normal curve in 
$\PP^{\left| \mu \right|  -1}$. So a $\mu$-scroll is a natural generalization of a rational normal curve.

A well-known result of the old Italian school of algebraic geometry states that $d$ points in general position in 
$\PP^{d-3}$ lie on a rational normal curve (see e.g. \cite{[G-H]}). The preceding considerations lead to the following generalization
of this classical result.
\begin{thm}\label{italian-gen}
Let $Z$ be a configuration of $d$ points on $X$ with $d >r+2$ and $r \geq 1$ ,where $r +1$
is the index of $L$-speciality of $Z$ (see (\ref{is})). Assume $Z$ to be in general position with respect to the adjoint linear system
$\left| L +K_X \right|$.

Let $\GA$ be an admissible component of $\CS$ containing $[Z]$ and let
$([Z],[\alpha],v)$, with $v\neq 0$, be a point of the relative tangent bundle $\TPI$ of 
$\JABG$ over $\GAB$. Let 
$$
\lambda (v) =( \lambda_1 (v), \ldots , \lambda_s (v))
$$
be the partition of $d$ associated to $d^{-} (v)$ and let
$$
\hat{\lambda}(v) =(\hat{\lambda}_1 (v) , \hat{\lambda}_2 (v) ,\ldots \hat{\lambda}_{s^{\prime}} (v))  
$$
be the truncation of $\lambda (v)$ defined in (\ref{part-tranc}). Then the image of $Z$ with respect to 
$\left| L +K_X \right|$ lies on a $\hat{\lambda}(v) $-scroll in the projective space
$\PP^{d-r-2}_Z$ (see (\ref{sysadj-emb}) for notation).
\end{thm}
\begin{pf}
If $\hat{\lambda}(v) =(1^{s^{\prime}})$ with $s^{\prime} = d-r-1$, then the $\hat{\lambda}(v)$-scroll
$Y_{\hat{\lambda}(v)} =\PP^{d-r-2}_Z$ and the assertion of the theorem is trivial.

Assume $\hat{\lambda}(v) \neq (1^{s^{\prime}})$ and write 
\BEN\label{part-trunc1}
\hat{\lambda}(v) = (\hat{\lambda}_1 (v) , \ldots , \hat{\lambda}_{s^{\prime}-m_1} (v),\underbrace{1,\ldots,1}_{m_1 -times} )\,,
\EEN
where $m_1$ is the multiplicity of $1$ in $\hat{\lambda}(v)$.

Applying Lemma \ref{scroll-Z} to every $\hat{\lambda}_i (v)$, for $i=1,\ldots,s^{\prime}-m_1$, we obtain
$Y_{\hat{\lambda}(v)}$ as the cone over a rational normal scroll
$Y_{\hat{\lambda}^{\prime} (v)}$ 
 with the vertex of the cone
$Sing (Y_{\hat{\lambda}(v)}) = \PP^{m_1 -1}$ and where
$$
\hat{\lambda}^{\prime} (v) =(\hat{\lambda}_1 (v) , \ldots , \hat{\lambda}_{s^{\prime}-m_1} (v) )\,.
$$
The rational normal scroll
$Y_{\hat{\lambda}^{\prime} (v)}$ is contained in the projective subspace $P$ of $\PP^{d-r-2}_Z$ spanned by the set of points
$$
\left\{ \left.V_{ij} \in B^{\ast} (v) \right| i =1,\ldots, s^{\prime}-m_1,\,j= 0,\ldots,\hat{\lambda}_i (v) -1 \right\}\,,
$$
where $B^{\ast} (v)$ is as defined in Lemma \ref{scroll-Z}, and the vertex of the cone $Sing (Y_{\hat{\lambda}(v)})$ is the complementary
subspace cut out by the linear equations
$$
X_{ij} = 0, \,\,i =1,\ldots, s^{\prime}-m_1,\,j= 0,\ldots,\hat{\lambda}_i (v) -1\,.
$$ 
\end{pf}
\begin{rem}
\begin{enumerate}
\item[1)]
The case $\hat{\lambda}(v) =(1^{s^{\prime}})$ in the proof of Theorem \ref{italian-gen} holds if and only if
$\LG=2$, $h^1_{\GA} = r$ and $d=2r +1$.
\item[2)]
Let $\LG \geq 3$. Then the multiplicity $m_1 (\hat{\lambda}(v))$ of $1$ in $\hat{\lambda}(v)$ is given by the following formula
\BEN\label{m1-mult}
m_1 (\hat{\lambda}(v)) =\mu^{\prime}_{(\LG-1)1} (v) + \mu^{\prime}_{(\LG-2)0} (v)\,,
\EEN
where $\mu^{\prime}_{qp}$ are the $d^{-}$-multiplicities of $v$ defined in (\ref{mult-v-}).
\item[3)]
Let $Z$ be a configuration of $d$ points with $d\geq 4$ and let
$[Z]$ be in $\stackrel{\circ}{\GA^1_d}$, the first non-trivial (with respect to our constructions) stratum of the stratification in (\ref{strH}).
 Assume, in addition, that $Z$ is in general position with respect to the adjoint linear system
 $\left| L+K_X \right|$.
Then for every point $([Z],[\alpha],v)$ of the relative tangent sheaf of $\JAB$ over $\XD$ with $v\neq 0$, the partition
$\lambda (v)$ has the following form
$$
 \lambda (v) =(d-1,1)\,.
$$
Its truncation
$\hat{\lambda}(v) =(d-2)$. Then Theorem \ref{italian-gen} implies that the image of $Z$ under the linear system $\left| L+K_X \right|$
lies on a rational normal curve in $\PP^{d-3}_Z$ (see the notation in (\ref{sysadj-emb})), thus recovering the classical result.
Furthermore, the rational normal curve acquires an additional meaning - it can be recovered as
 the closure of the image of the period map $p_{\GA}$ 
(see \S\ref{sec-periods}, (\ref{p}), for notation) over 
$[Z] \in \GAB$, where $\GA$ is an admissible component in $C^1 (L,d)$ containing $[Z]$.
\end{enumerate}
\end{rem}

\section{Representation theoretic constructions}\label{sec-RT-const}
The preceding sections show that the Lie theoretic aspects of the Jacobian $\JA$ provide new methods and insights in the study of 
geometry of surfaces. Starting from this section we change the logic of our investigations - we make use of the sheaves of Lie algebras
{\BM
$\LAG$}, for admissible components $\GA \in \CS$,
to construct various objects (sheaves, complexes of sheaves, constructible functions), either on $\JA$ or on the Hilbert scheme
$\XD$, which can serve as new invariants for vector bundles on $X$ as well as for $X$ itself. Our basic tool for this will be
the morphisms $d^{\pm}$ encountered in \S4, (\ref{d+}), (\ref{d-}). These morphisms relate our Jacobian
to such fundamental objects in the Geometric representation theory as nilpotent orbits, Springer resolution and Springer fibres.

\subsection{Basic set-up}
Let $\GA$ be an admissible component in $\CS$ and assume it to be simple (Definition \ref{s-c}). Thus we tacitly assume that the set of simple components in $\CS$ is nonempty. From \S\ref{sec-ses}, Theorem \ref{ty-dec-cont}, it follows that this is the essential case to consider.

 By Corollary \ref{cor-s}
 the sheaf of Lie algebras
{\BM
$\LAG$ attached to $\JABG$ (see \S\ref{Lie} for notation) has the following description
\BEN\label{LAG=sl}
\LAG = \mbox{\UB$\pi^{\ast} {\bf sl} (\FF^{\prime})$}\,,
\EEN}
where $\pi$ is the natural projection
\BEN\label{proj-pi}
\pi : \JABG \longrightarrow \GAB
\EEN\
and ${\bf sl}( \FF^{\prime})$ stands for the sheaf of germs of traceless endomorphisms of $\FF^{\prime}$ (see Corollary \ref{cor-Fpr} for the 
definition of $\FF^{\prime}$).

Let ${\cal T}_{\pi}$ be the relative tangent sheaf of the morphism $\pi$ and let
$\TPI$ be the corresponding vector bundle, i.e. it is a fibre space over  $\JABG$ with the natural projection
\BEN\label{tau1}
\tau : \TPI \longrightarrow \JABG
\EEN
whose fibre over a point $\ZA \in \JABG$ is the space 
$\TPI \ZA ={\cal T}_{\pi} \ZA$ of vertical 
 tangent vectors of $\JABG$ at $\ZA$. 
 
In \S\ref{sec-periods}, (\ref{d+}) and (\ref{d-}), we defined morphisms of sheaves
$$
d^{\pm} : {\cal T}_{\pi} \longrightarrow \mbox{\BM$\LAG$}\,.
$$
 In this section  
{\BM
$\LAG$}
will be often viewed as a vector bundle over $\JABG$. Then $d^{\pm}$ can be viewed as morphisms of  $\JABG$-schemes
\BEN\label{d+-morph}
\xymatrix{
{\TPI} \ar[rr]^{d^{\pm}} \ar[dr]_{\tau} & &{ \mbox{\BM$\GS$}_{\GA}} \ar[dl]^{\gamma} \\
           &{\JABG} }
\EEN
 as it has been already done in Remark \ref{geom-taud+}. The fact that 
{\BM$\LAG$}
is the pullback of ${\bf sl} ( \FF^{\prime})$ (see (\ref{LAG=sl})) allows further to associate to
$d^{\pm}$ morphisms of schemes over $\GAB$. More precisely, set
\BEN\label{LAGpr}
\mbox{\BM$\GS$}^{\prime}_{\GA} ={\bf sl} (\FF^{\prime})
\EEN
and view it 
as a bundle over $\GAB$ with the natural projection
\BEN\label{LAGpr-proj}
\gamma^{\prime} : \mbox{\BM$\GS$}^{\prime}_{\GA} \longrightarrow \GAB\,.
\EEN
 We can now identify 
{\BM
$\LAG$}
as the fibre-product
\BEN\label{LAG-fb-prod}
\mbox{\BM$\LAG =\GS$}^{\prime}_{\GA} \times_{\GAB} \JABG\,.
\EEN
Composing (\ref{d+-morph}) with the projection $\pi$ in (\ref{proj-pi}) yields the
 commutative diagram
\BEN\label{prd+-}
\xymatrix{
{\TPI} \ar[rr]^{{}^{\prime}d^{\pm}} \ar[dr]_{\tilde{\pi}} & &{ \mbox{\BM$\GS$}^{\prime}_{\GA}} \ar[dl]^{\gamma^{\prime}} \\
           &{\GAB} }
\EEN
where the notation is as follows
\begin{enumerate}
\item[a)]
$\tilde{\pi}: \TPI \longrightarrow \GAB$ is the composition $\tilde{\pi} =\pi \circ \tau$,
\item[b)]
${}^{\prime}d^{\pm} : \TPI \longrightarrow \mbox{\BM$\GS$}^{\prime}_{\GA}$
is the composition of $d^{\pm}$ in (\ref{d+-morph}) with the projection
$$
\mbox{\BM$\LAG =\LAG^{\prime}$} \times_{\GAB} \JABG \longrightarrow \mbox{\BM$\LAG^{\prime}$}\,.
$$
\end{enumerate}

Denote by
\BEN\label{cone-n}
\mbox{\BM${\cal N}$}^{\prime}_{\GA} =\mbox{\BM${\cal N} (\GS^{\prime}_{\mbox{\UB$\GA$}} )$}
\EEN
the subscheme of nilpotent elements  of 
$\mbox{\BM$\GS^{\prime}$}_{\GA}$ and call it {\it nilpotent cone of} $\mbox{\BM$\GS^{\prime}$}_{\GA}$.
Then we know that the morphisms ${}^{\prime}d^{\pm}$ take their values in 
$\mbox{\BM${\cal N^{\prime}}$}_{\GA}$. Hence the diagram in (\ref{prd+-}) takes the following form
\BEN\label{prd+-1}
\xymatrix{
{\TPI} \ar[rr]^{{}^{\prime}d^{\pm}} \ar[dr]_{\tilde{\pi}} & &{ \mbox{\BM$\NI_{\GA}$}} \ar[dl]^{\gamma^{\prime}} \\
           &{\GAB} }
\EEN
The adjoint action\footnote{ in this case it is the usual conjugation.} of 
$\mbox{\BM$G^{\prime}$}_{\GA} ={\bf SL}(\FF^{\prime})$ on 
$\mbox{\BM$\NI_{\GA}$}$ divides it into $\mbox{\BM$G^{\prime}$}_{\GA}$-orbits. Our first task will be to clarify this orbit structure.

\subsection{The orbit structure of $\mbox{\BM$\NI_{\GA}$}$}
In this subsection we describe some basic properties of the orbits of the
$\mbox{\BM$G^{\prime}$}_{\GA}$-action on 
$\mbox{\BM$\NI_{\GA}$}$.

Fix such an orbit and denote it $\mbox{\BM$O(\NI_{\mbox{\UB$\GA$}})$}$ and consider the restriction
\BEN\label{rest-gampr}
\gamma^{\prime}_{\scriptstyle{\mbox{\BM$O(\NI_{\mbox{\UB$\GA$}})$}} }: \mbox{\BM$O(\NI_{\GA})$} \longrightarrow \GAB
\EEN
of $\gamma^{\prime}$ in (\ref{prd+-1}) to 
$ \mbox{\BM$O(\NI_{\GA} )$}$.
\begin{lem}\label{orbit=fb}
An orbit 
$\mbox{\BM$O(\NI_{\GA} )$}$
is a fibre bundle over
$\GAB$,
whose fibre is modeled on a fixed nilpotent orbit of ${\bf sl}_{d^{\prime}_{\GA}} (\CC)$,
where
$d^{\prime}_{\GA} = rk (\FF^{\prime})$.
\end{lem}
\begin{pf}
Choose a covering 
$\{ U_i \}_{i\in I}$ of 
$\GAB$ trivializing $\FF^{\prime}$, i.e. the restriction
$\FF^{\prime} \mid_{U_i}$ is isomorphic to the trivial bundle
$U_i \times \CC^{d^{\prime}_{\GA}}$ and let
\BEN\label{triv}
\phi_i : \FF^{\prime} \mid_{U_i}\stackrel{\cong}{\longrightarrow} U_i \times \CC^{d^{\prime}_{\GA}}
\EEN
be such a trivialization. Set 
\BEN\label{trans}
\phi_{ij} : U_{ij} \longrightarrow {\bf GL}_{\DE} (\CC)
\EEN
 to be the corresponding transition functions of $\FF^{\prime}$.

The trivializations $\phi_i$'s induce the trivializations
\BEN\label{triv1}
\psi_i : {\bf sl}(\FF^{\prime}) \mid_{U_i} \stackrel{\cong}{\longrightarrow} U_i \times {\bf sl}_{\DE} (\CC)
\EEN
with the transition functions
\BEN\label{trans1}
\psi_{ij} : U_{ij} \longrightarrow Aut ({\bf sl}_{\DE} (\CC))
\EEN
given by conjugation by $\phi_{ij}$, i.e. we have
\BEN\label{trans2}
\psi_{ij} (u) (A) =\phi_{ij} (u) A (\phi_{ij} (u))^{-1}\,,
\EEN
for every $u\in U_{ij}$ and every $A \in {\bf sl}_{\DE} (\CC)$.

Let 
{\BM
${\cal N} = {\cal N}(sl_{\DE} (C))$}
be the nilpotent cone of 
${\bf sl}_{\DE} (\CC)$, i.e.
{\BM
${\cal N}$} is the subvariety of nilpotent endomorphisms of 
$\CC^{\DE}$.
The trivializations $\psi_i$'s in (\ref{triv1}) induce the trivializations
\BEN\label{triv3}
\psi_i : \left.{\mbox{\BM$\NI_{\GA}$}}\right|_{U_i} \longrightarrow U_i \times \mbox{\BM${\cal N}$}\,,
\EEN
where
$ {\mbox{\BM$\NI$}_{\GA}}\mid_{U_i} =\gamma^{\prime -1} (U_i)$
is the inverse image of $U_i$ under the morphism $\gamma^{\prime}$ in (\ref{prd+-1}).

Set
\BEN\label{orbit-Ui}
\left.{\mbox{\BM$O(\NI_{\GA} )$}}\right|_{U_i} =\mbox{\BM$O(\NI_{\GA} )$} \bigcap (\left.{\mbox{\BM$\NI_{\GA}$}}\right|_{U_i} )\,.
\EEN
Then $\psi_i$ in (\ref{triv3}) maps 
$\left.{\mbox{\BM$O(\NI_{\GA} )$} }\right|_{U_i}$ onto a subset of 
$U_i \times \mbox{\BM${\cal N}$}$ of the form
$U_i \times \mbox{\BM$O_i$}$, for some nilpotent orbit
$\mbox{\BM$O_i$}$ of 
$\mbox{\BM${\cal N}$}$. On the intersection $U_{ij}$ the two isomorphisms $\psi_i$ and $\psi_j$
are related by the transition function $\psi_{ij}$
\BEN\label{diag-trans}
\xymatrix{
&U_{ij} \times \mbox{\BM$O_i$}\\
{\left. \mbox{\BM$O(\NI_{\GA})$} \right|_{U_{ij}}} \ar[ru]^{\psi_i} \ar[rd]_{\psi_j}& \\
&U_{ij} \times \mbox{\BM$O_j$} \ar[uu]_{id_{U_{ij}} \times \psi_{ij}} }
\EEN
Explicitly, for a section $s$ of 
$\mbox{\BM$O(\NI_{\GA})$}$ over $U_{ij}$, we have
\begin{eqnarray*}
\psi_i (s(u)) =&(u, A_i (u)) \in U_{ij} \times \mbox{\BM$O$}_i \subset U_{ij} \times {\bf sl}_{\DE} (\CC) \\ 
\psi_j (s(u)) =&(u, A_j (u)) \in U_{ij} \times \mbox{\BM$O$}_j \subset U_{ij} \times
 {\bf sl}_{\DE} (\CC)\,,
\end{eqnarray*}
for every $u \in U_{ij}$. These are related by conjugation (\ref{trans2})
\BEN\label{trans3}
A_i (u) = \phi_{ij} A_j (u) (\phi_{ij} (u) )^{-1}, \,\,\forall u\in U_{ij}\,.
\EEN
This implies that $A_i (u)$ and $A_j (u)$ are in the same nilpotent orbit of 
${\bf sl}_{\DE} (\CC)$. Hence
{\BM
$O_i =O_j$},
for all $i,j \in I$ with $U_{ij} =U_i \bigcap U_j \neq \emptyset$, and to
$\mbox{\BM$O(\NI_{\GA} )$}$ we can associate a unique nilpotent orbit
$\mbox{\BM$O$}$ 
of ${\bf sl}_{\DE} (\CC)$
such that 
$\mbox{\BM$O(\NI_{\GA} )$}$
is the fibre bundle over $\GAB$ with fibres isomorphic to $\mbox{\BM$O$}$ and the transition functions
$\psi_{ij}$'s in (\ref{trans1})
\end{pf}
\begin{cor}\label{orbits}
The orbits of
$\mbox{\BM$\NI_{\GA}$}$ under the adjoint (=conjugation) action of
${\bf SL}(\FF^{\prime})$ are in bijective correspondence with nilpotent orbits of 
${\bf sl}_{\DE} (\CC)$. This correspondence will be denoted as follows
\BEN\label{orbits1}
\mbox{\BM$O \longleftrightarrow O(\NI_{\mbox{\UB$\GA$}} )$}\,,
\EEN
for every nilpotent orbit 
{\BM
$O$}
of ${\bf sl}_{\DE} (\CC)$.
\end{cor}
\begin{pf}
Follows immediately from Lemma \ref{orbit=fb}
\end{pf}
\begin{rem}\label{or-part}
It is well known that nilpotent orbits of
${\bf sl}_{\DE} (\CC)$
are in bijective correspondence with the set of partitions of
$\DE$ (see e.g. \cite{[C-Gi]}). For a partition $\mu$ of
$\DE$, denote by
$\mbox{\BM$O_{\mu}$}$ the corresponding nilpotent orbit of 
${\bf sl}_{\DE} (\CC)$. 
Then 
$\mbox{\BM$O_{\mu} (\NI_{\mbox{\UB$\GA$}} )$}$
will denote
the orbit associated to
$\mbox{\BM$O_{\mu}$}$
by the correspondence in Corollary \ref{orbits}.
\end{rem}

Let us return to the diagram in (\ref{prd+-1}). For a nilpotent orbit
$\mbox{\BM$O_{\mu} (\NI_{\mbox{\UB$\GA$}} )$}$
define
\BEN\label{T-orbits}
\mbox{\BM$O_{\mu} $}(\TPI) := ({}^{\prime}d^{+})^{-1} (\mbox{\BM$O_{\mu} (\NI_{\mbox{\UB$\GA$}} )$})\,.
\EEN
This gives a partition of $\TPI$ into the disjoint union of locally closed strata
$\mbox{\BM$O_{\mu} $}(\TPI)$.
The above considerations imply the following.
\begin{pro}\label{pro-T-orbits}
Let $\GA$ be an admissible component in $\CS$ and assume it to be simple in the sense of Definition \ref{s-c}.
Then $\GA$ determines a finite collection of partitions
$$
P(\GA) =\left\{\left.\mu \in P_{\DE} \right| \mbox{\BM$O_{\mu} $}(\TPI) \neq \emptyset \right\}\,,
$$
where $P_{\DE}$ denotes the set of partitions of $\DE$. Equivalently, every admissible, simple component
$\GA$ in $\CS$ determines the finite collection of nilpotent orbits in ${\bf sl}_{\DE} (\CC)$
$$
O(\GA) =\left\{\left.\mbox{\BM$O_{\mu} $} \right| \mu \in P(\GA) \right\}\,.
$$
Furthermore, there exists a unique partition in 
$P(\GA)$ denoted $\mu_{\GA}$ such that the corresponding stratum
$\mbox{\BM$O$}_{\mu_{\GA}} (\TPI)$ is a dense Zariski open subset of 
$\TPI$.
\end{pro}
\begin{pf}
All but the last assertion is a combination of 
Corollary \ref{orbits}, Remark \ref{or-part} and (\ref{T-orbits}).

The last assertion follows from the fact that $\TPI$ is irreducible.
\end{pf} 

This result yields Theorem \ref{nilo} stated in the introduction. More precisely, we have the following
\begin{thm}\label{nilo1}
Let 
${\cal V}^r (L,d)$ denotes the set of admissible, simple components of $\CS$ and let
$$
{\cal V} (X;L,d) = \bigcup_{r\geq 1} {\cal V}^r (L,d)\,.
$$
Then the set ${\cal V} (X;L,d)$ is finite and every $\GA$ in it
determines a distinguished collection of nilpotent orbits 
$O(\GA)$ of ${\bf sl}_{\DE} (\CC)$ as in Proposition \ref{pro-T-orbits}.
\end{thm}
\begin{pf}
For every $r \geq 1$, the set 
${\cal V}^r (L,d)$ is finite, since there are finitely many components in $\CS$.
On the other hand, by definition of the index of $L$-speciality (see (\ref{is})) $r\leq  d$ . This yields
the finiteness of 
${\cal V} (L,d)$. 
The second assertion is the content of Proposition \ref{pro-T-orbits}.
\end{pf}

The stratification of $\TPI$ defined by the strata in (\ref{T-orbits}) should be compared 
to the one in Proposition \ref{pro-strat1}. The latter stratification is given by the strata
$T^{\lambda}$'s indexed by the set of admissible 
$\overrightarrow{h^{\prime}}_{\GA}$-graded partitions $\lambda$ in 
$P^a_{\DE} (\overrightarrow{h^{\prime}}_{\GA}) $ (see Remark \ref{part-lam=gr}). The partitions in 
$P(\GA)$ of Proposition \ref{pro-T-orbits} do not have the finer structure of 
$\overrightarrow{h^{\prime}}_{\GA}$-grading.
More precisely, we have the forgetful map
\BEN\label{map-forget}
F_{\GA} :  P^a_{\DE} (\overrightarrow{h^{\prime}}_{\GA}) \longrightarrow P(\GA)
\EEN
which sends an $\overrightarrow{h^{\prime}}_{\GA}$-graded partition
$$
\lambda =\bigcup^{\LG-1}_{p=0} \lambda^{(p)}
$$
with $\lambda^{(p)} =(1^{\mu_{0p}} 2^{\mu_{1p}} \ldots (p+1)^{\mu_{pp}})\,(p=0,\ldots,\LG-1)$, to the partition
$F_{\GA} (\lambda)$, where one forgets the grading. This partition is given in terms of multiplicities by the following formula
\BEN\label{map-forget-formula}
F_{\GA} (\lambda) =(1^{M_{0} (\lambda)} 2^{M_{1} (\lambda)} \ldots (\LG)^{M_{\LG-1} (\lambda)})\,,
\EEN
where 
\BEN\label{row}
M_s (\lambda) =\sum^{\LG-1}_{p=s} \mu_{sp},\,\,for \,\,s=0,\ldots,\LG-1\,.
\EEN
It is clear, that given a partition $\mu \in P(\GA)$, there might be several ways to define
$\overrightarrow{h^{\prime}}_{\GA}$-grading on it. The inverse image
$F^{-1}_{\GA}  (\mu)$ gives all such structures on $\mu$, determined by the morphisms $d^{\pm}$.

Using the above notation, we can express a relation between the stratification in Proposition \ref{pro-strat1} and the one given by the strata
$\mbox{\BM$O$}_{\mu} (\TPI)$ in (\ref{T-orbits}).
\begin{pro}\label{two-strat}
Let $\mu$ be a partition in the set $P(\GA)$ of Proposition \ref{pro-T-orbits}.
Then
$$
\mbox{\BM$O$}_{\mu} (\TPI) =\bigcup_{\lambda \in  F^{-1}_{\GA}  (\mu)} T^{\lambda}\,,
$$
where $T^{\lambda}$ are strata of Proposition \ref{pro-strat1}. Furthermore, the partition $\mu$ is related to the partitions
$\lambda$ in $F^{-1}_{\GA}  (\mu)$ by the formula in (\ref{row}).
\end{pro}

This relation and the results of \S\ref{sec-equations} show that the partitions
 in Proposition \ref{pro-T-orbits} distinguished by the nonabelian Jacobian
$\JA$ are closely related to various algebro-geometric properties of configurations of points on $X$ as well as curves
in the linear system $\left| L \right|$.

\subsection{Perverse sheaves and $\JA$}\label{sec-ps}

In the previous subsection we have seen how $\JA$ distinguishes a finite collection
${\cal V} (X;L,d)$ of subvarieties $\GA$ of the Hilbert scheme $\XD$ with the property that to each $\GA$
in ${\cal V} (X;L,d)$ one can attach a finite collection of partitions
$P(\GA)$ of $\DE$ as described in Proposition \ref{pro-T-orbits}. Recalling that the partitions of 
$\DE$ also parametrize irreducible representations of the symmetric group
$S_{\DE}$ we obtain an equivalent version of Theorem \ref{nilo1} formulated as Theorem \ref{irrSn} in the Introduction.
This theorem attaches to each $\GA$ in ${\cal V} (X;L,d)$ a finite collection 
\BEN\label{R-collec}
R_{\DE} (\GA) = \{ {\bf S}_{\mu} \mid \mu \in P(\GA) \}
\EEN
of irreducible $S_{\DE}$-modules ${\bf S}_{\mu}$ (up to an isomorphism) indexed by the set of partitions $P(\GA)$ in
 Proposition \ref{pro-T-orbits}.

Thinking of $\GA$ as a variety parametrizing geometric representatives of the second Chern class of rank 2 bundles on $X$,
one can view Theorem \ref{irrSn} as a way of elevating the topological invariant - the degree of the second Chern class - to the
category of modules of symmetric groups.

In this subsection we go further: we attach to each $\GA$ in ${\cal V} (X;L,d)$ a collection of the Intersection cohomology complexes
on the Hilbert scheme $\XD$, thus elevating the degree $d$ of the second Chern class to the category of perverse sheaves on $\XD$.
This will prove Theorem \ref{ps} stated in the Introduction.

Our construction is based on the Springer resolution
\BEN\label{Springer-res}
\sigma : \mbox{\BM$\tilde{\cal N}$} \longrightarrow \mbox{\BM${\cal N}$}
\EEN
of the nilpotent cone $ \mbox{\BM${\cal N}$}$ of 
${\bf sl}_{\DE} (\CC)$
and a well-known fact in the geometric representation theory\footnote{for this and other basic facts of the geometric representation theory
our reference is \cite{[C-Gi]}.} which realizes the cohomology groups of fibers of $\sigma$ (Springer fibres)
as modules of the Weyl group $W$ of ${\bf sl}_{\DE} (\CC)$. 

In our case we deal with the variety 
{\BM
$\NI_{\mbox{\UB$\GA$}}$}
which is fibred over $\GAB$ by the nilpotent cones isomorphic to 
{\BM
${\cal N}$}. So we need a relative version of Springer resolution.
For this we continue to view 
{\BM
$\GS^{\prime}_{\mbox{\UB$\GA$}}$}
as a vector bundle over $\GAB$ and define the {\it relative flag variety }
{\BM
${\cal B}_{\mbox{\UB$\GA$}}$}
of 
{\BM
$\GS^{\prime}_{\mbox{\UB$\GA$}}$}. 

By definition 
this variety comes with the natural projection
\BEN\label{rel-fl-B}
\beta_{\GA} : \mbox{\BM${\cal B}_{\mbox{\UB$\GA$}}$} \longrightarrow \GAB
\EEN
such that the fibre
{\BM
${\cal B}_{\mbox{\UB$Z$}}$}
of $\beta_{\GA}$ over a point $[Z] \in \GA$ is the variety of Borel subalgebras of
$\mbox{\BM$\GS^{\prime}_{\mbox{\UB$\GA$}}$} ([Z]) = {\bf sl} (\FF^{\prime} ([Z]))$, where $\FF^{\prime} ([Z])$ is the fibre of
$\FF^{\prime}$ at $[Z]$.
Set
\BEN\label{rel-Sp-res}
\mbox{\BM$\tilde{\cal N}^{\prime}_{\mbox{\UB$\GA$}}$} := {\bf T}^{\ast}_{\beta_{\GA}} =
\mbox{\BM$ T^{\ast}_{{\cal B}_{\mbox{\UB$\GA$}} /{\mbox{\UB$\GAB$}} } $}
\EEN
to be the relative cotangent bundle of $\beta_{\GA}$. Its closed points can be described as the following incidence
correspondence
$$
\mbox{\BM$\tilde{\NI}_{\mbox{\UB$\GA$}}$} = \left\{ \left.([Z],x,{\bf b}) \in
  \mbox{\BM$\NI_{\mbox{\UB$\GA$}}$} \times \mbox{\BM${\cal B}_{\mbox{\UB$\GA$}}$} \right|
{\bf b}\,\, \mbox{is a Borel subalgebra in}\,\, \mbox{\BM$\GS^{\prime}_{\mbox{\UB$\GA$}}$} ([Z]),\,\,
x\in {\bf b}\cap \mbox{\BM${\cal N}^{\prime}_{\mbox{\UB$\GA$}}$} ([Z]) \right\}\,.
$$
This gives the following commutative diagram
\BEN\label{rel-Sp-res1}
\xymatrix{
&{\mbox{\BM$\tilde{\cal N}^{\prime}_{\mbox{\UB$\GA$}}$}} \ar[dl]_{\sigma_{\GA}} \ar[dr]^{\sigma^{\prime}_{\GA}}&  \\
 {\mbox{\BM$\NI_{\mbox{\UB$\GA$}}$}} \ar[dr]_{\gamma^{\prime}_{\GA}} & &
{\mbox{\BM${\cal B}_{\mbox{\UB$\GA$}}$}} \ar[dl]^{\beta_{\GA}} \\
& {\GAB} & }
\EEN
The morphism
\BEN\label{morph-Sp-res}
\sigma_{\GA} : \mbox{\BM$\tilde{\cal N}^{\prime}_{\mbox{\UB$\GA$}}$} \longrightarrow \mbox{\BM$\NI_{\mbox{\UB$\GA$}}$} 
\EEN
in the above diagram is the relative Springer resolution of 
$\mbox{\BM$\NI_{\mbox{\UB$\GA$}}$} $, i.e. for every
$[Z] \in \GAB$, the restriction
$\sigma_{\GA, [Z]}$ of $\sigma_{\GA}$ to the fibre
$\mbox{\BM$\tilde{\cal N}^{\prime}_{\mbox{\UB$\GA$}}$} ([Z]) =(\beta_{\GA} \circ \sigma^{\prime}_{\GA} )^{-1} ([Z])$
of $\mbox{\BM$\tilde{\NI}_{\mbox{\UB$\GA$}}$}$ over $[Z]$ is the Springer resolution
\BEN\label{Sp-res-Z}
\sigma_{\GA, [Z]} : \mbox{\BM$\tilde{\cal N}^{\prime}_{\mbox{\UB$\GA$}}$} ([Z]) \longrightarrow 
\mbox{\BM$\NI_{\mbox{\UB$\GA$}}$} ([Z])
\EEN
of the nilpotent cone $\mbox{\BM$\NI_{\mbox{\UB$\GA$}}$} ([Z])$ in 
$\mbox{\BM$\GS^{\prime}_{\mbox{\UB$\GA$}}$} ([Z]) = {\bf sl} (\FF^{\prime} ([Z]))$.

Let
$\mbox{\BM$O_{\mu} (\NI_{\mbox{\UB$\GA$}})$}$ be a nilpotent orbit in 
$\mbox{\BM$\NI_{\mbox{\UB$\GA$}}$}$
and let
\BEN\label{Sp-res-mu}
\sigma_{\GA} : 
\mbox{\BM$\tilde{O}_{\mu} (\NI_{\mbox{\UB$\GA$}})$} = \sigma^{-1}_{\GA} (\mbox{\BM$O_{\mu} (\NI_{\mbox{\UB$\GA$}})$})
\longrightarrow \mbox{\BM$O_{\mu} (\NI_{\mbox{\UB$\GA$}})$}
\EEN
be the Springer resolution over the nilpotent orbit
$\mbox{\BM$O_{\mu} (\NI_{\mbox{\UB$\GA$}})$}$. From Lemma \ref{orbit=fb} and the properties of the Springer resolution it follows
that this is a fibre bundle over
$\mbox{\BM$O_{\mu} (\NI_{\mbox{\UB$\GA$}})$}$ 
with fibres modeled on a Springer fibre of $\sigma$ in (\ref{Springer-res}) over the nilpotent orbit $\mbox{\BM$O_{\mu}$}$
in ${\bf sl}_{\DE} (\CC)$. In particular, taking the $i$-th direct image of the constant sheaf $\underline{\CC}$ on 
$\mbox{\BM$\tilde{O}_{\mu} (\NI_{\mbox{\UB$\GA$}})$}$ we obtain
\BEN\label{dir-im}
\tilde{\cal L}^i_{\GA, \mu} = {\cal R}^i \sigma_{\GA \ast} \underline{\CC}
\EEN
local systems on $\mbox{\BM$O_{\mu} (\NI_{\mbox{\UB$\GA$}})$}$, for $i =0,\ldots, 2b_{\mu}$,
where $b_{\mu}$ is the complex dimension of a Springer fibre of $\sigma$ in (\ref{Springer-res}) 
over the nilpotent orbit $\mbox{\BM$O_{\mu}$}$.
\begin{rem}\label{dim-Sp-fib}
Let ${\bf B}$ be the flag variety of ${\bf sl}_{\DE} (\CC)$, i.e. the variety parametrizing Borel subalgebras of ${\bf sl}_{\DE} (\CC)$.
Recall that for a nilpotent element $x \in \mbox{\BM${O}_{\mu}$}$,
the Springer fibre
${\bf B}_x$ over $x$ of the Springer resolution $\sigma$ in (\ref{Springer-res}) is naturally identified with the subvariety of  
${\bf B}$ parametrizing Borel subalgebras of ${\bf sl}_{\DE} (\CC)$ containing $x$.
It is known that ${\bf B}_x$ is an equidimensional variety of complex dimension
$$
b_{\mu} =dim({\bf B}) - \frac{1}{2} dim \mbox{\BM$O_{\mu}$}\,.
$$
In our case $dim_{\CC} {\bf B} = \frac{1}{2} \DE (\DE -1)$. Substituting into the above formula we obtain
\BEN\label{dim-fib1}
b_{\mu} = \frac{1}{2}[ \DE (\DE -1) - dim \mbox{\BM$O_{\mu}$}]\,.
\EEN
Furthermore, if $\mu = (\mu_1 \geq \mu_2 \geq \cdots \geq \mu_s)$, then one has the following formula
(\cite{[C-Gi]}, Lemma 4.4.2)
$$
dim \mbox{\BM$O_{\mu}$} =(\DE)^2 - \sum^s_{k=1} (2k-1) \mu_k\,.
$$
Substituting into (\ref{dim-fib1}) yields
\BEN\label{dim-fib2}
b_{\mu} = \frac{1}{2} (\sum^s_{k=1} (2k-1) \mu_k -\DE) =\sum^s_{k=1} k \mu_k   - \DE =\sum^s_{k=1} (k-1) \mu_k\,.
\EEN
\end{rem}
\begin{lem}\label{loc-sys}
The local systems 
$\tilde{\cal L}^i_{\GA, \mu} $'s in (\ref{dir-im}) are the pullback under $\gamma^{\prime}$ (see (\ref{rel-Sp-res1}))
of the local systems on $\GAB$, i.e. for every $i \in \{0,\ldots, 2b_{\mu} \}$ there exists a unique, up to an isomorphism, local system
${\cal L}^i_{\GA, \mu} $ on $\GAB$ such that
$$
\tilde{\cal L}^i_{\GA,\mu}  = \gamma^{\prime \ast} {\cal L}^i_{\GA, \mu}\,.
$$
\end{lem}
\begin{pf}
Fix a base point $([Z_0],x_0) \in \mbox{\BM$O_{\mu} (\NI_{\mbox{\UB$\GA$}})$}$.
Then the local systems 
$\tilde{\cal L}^i_{\GA, \mu} $
correspond to representations
\BEN\label{rep-fg}
\tilde{\rho}_i : \pi_1 ( \mbox{\BM$O_{\mu} (\NI_{\mbox{\UB$\GA$}})$},([Z_0],x_0) ) \longrightarrow 
Aut (H^i ( {\bf B}_{([Z_0],x_0)}, \CC))
\EEN
of the fundamental group 
$\pi_1 ( \mbox{\BM$O_{\mu} (\NI_{\mbox{\UB$\GA$}})$}, ([Z_0],x_0) )$ of 
$\mbox{\BM$O_{\mu} (\NI_{\mbox{\UB$\GA$}})$}$ based at $([Z_0],x_0)$,
and where 
$Aut (H^i ( {\bf B}_{([Z_0],x_0)}, \CC))$ is the group of automorphisms of the 
$i$-th cohomology group of the Springer fibre 
${\bf B}_{([Z_0],x_0)}$, the fibre of $\sigma_{\GA}$ over $([Z_0],x_0)$.

From Lemma \ref{orbit=fb} it follows that
$\pi_1 ( \mbox{\BM$O_{\mu} (\NI_{\mbox{\UB$\GA$}})$}, ([Z_0],x_0) )$
fits into the following long exact sequence of groups
\BEN\label{exs-fg}
\xymatrix@1{
{\pi_1 ( \mbox{\BM$O_{\mu}$}, x_0)} \ar[r] &{\pi_1 ( \mbox{\BM$O_{\mu} (\NI_{\mbox{\UB$\GA$}})$}, ([Z_0],x_0) ) }\ar[r] &
{\pi_1 ( \GAB, [Z_0])} \ar[r]& {\pi_0 ( \mbox{\BM$O_{\mu}$}, x_0)}\,, }
\EEN
where we identified the fibre 
$\gamma^{\prime -1}_{\mbox{\BM$O_{\mu}(\NI_{\GA})$}} ([Z]) $
of $\gamma^{\prime}_{\mbox{\BM$O_{\mu} (\NI_{\GA})$}}$ in (\ref{rest-gampr})
with the orbit $\mbox{\BM$O_{\mu} $}$ in 
${\bf sl} (\FF^{\prime} ([Z]))$.
It is well-known that $\mbox{\BM$O_{\mu} $}$ is connected and simply connected (see \cite{[C-Gi]}), i.e. one has
$$
\pi_0 ( \mbox{\BM$O_{\mu}$}, x_0) = \pi_1 ( \mbox{\BM$O_{\mu}$}, x_0) =\{1\}\,.
$$
This together with (\ref{exs-fg}) yield an isomorphism
\BEN\label{iso-fg}
\pi_1 ( \mbox{\BM$O_{\mu} (\NI_{\mbox{\UB$\GA$}})$}, ([Z_0],x_0) ) \cong \pi_1 (  \GAB, [Z_0])\,.
\EEN
This isomorphism combined with (\ref{rep-fg}) gives  representations
\BEN\label{rep-fg1}
\rho^i_{\GA, \mu} : \pi_1 ( \GAB,[Z_0] ) \longrightarrow 
Aut (H^i ( {\bf B}_{([Z_0],x_0)}, \CC))\,,
\EEN
for $i=0,\ldots, 2b_{\mu}$.

Let ${\cal L}^i_{\GA,\mu}\, (i=0,\ldots, 2b_{\mu})$ be the local systems on $\GAB$
corresponding to the representations
$\rho^i_{\GA, \mu}$ in (\ref{rep-fg1}). Then by definition we have
$$
{\tilde{\cal L}}^i_{\GA,\mu} = \gamma^{\prime \ast}_{\GA} {\cal L}^i_{\GA,\mu}\,,
$$
for every $i\in \{0,\ldots ,2b_{\mu} \}$.
\end{pf}

Denote by
${\cal P}^i_{\GA, \mu}$
the Intersection cohomology complex 
$IC(\GAB, {\cal L}^i_{\GA,\mu} )$
of Deligne-Goresky-\linebreak MacPherson extended by zero to the entire Hilbert scheme $\XD$. This is an object of the bounded derived category of constructible sheaves 
$D^b_c (\XD)$ on $\XD$ which is characterized by the following properties.
\begin{eqnarray}\label{IC}
a)&{\cal P}^i_{\GA, \mu}\,\,\mbox{ is supported on the closure }\,\,
\overline{\GAB} =\overline{\GA}\,\,\mbox{of}\,\,\GA\,\,\mbox{in}\,\, \XD\,, \\ \nonumber
b)& \left. {\cal P}^i_{\GA, \mu} \right|_{\GAB} ={\cal L}^i_{\GA,\mu} [dim \GA ]\,, \\ \nonumber
c)& {\cal H}^k ( {\cal P}^i_{\GA, \mu} )=0,\,\,\mbox{if}\,\, k<-dim(\GA)\,, \\ \nonumber
d)& dim\left(supp ({\cal H}^k ( {\cal P}^i_{\GA, \mu} ))\right) <-k,\,\,\mbox{if}\,\, k> -dim(\GA)\,, \\ \nonumber
e)& dim\left(supp ({\cal H}^k ( ({\cal P}^i_{\GA, \mu})^{VD}) )\right) <-k,\,\,\mbox{if}\,\, k> -dim(\GA)\,,
\end{eqnarray}
where
$(\cdot)^{VD}$ stands for the Verdier dual complex.

Putting the complexes ${\cal P}^i_{\GA, \mu}$ together, we obtain the graded perverse sheaf
\BEN\label{IC-gr}
 {\cal P}^{\bullet}_{\GA, \mu} = \bigoplus^{2b_{\mu}}_{i=0} {\cal P}^i_{\GA, \mu} 
\EEN
which is the extension by zero to $\XD$ of the Intersection cohomology complex 
$IC(\GAB, {\cal L}^{\bullet}_{\GA,\mu})$, where
$ {\cal L}^{\bullet} =  \bigoplus^{2b_{\mu}}_{i=0} {\cal L}^i_{\GA, \mu}$. Thus we obtain the following
\begin{thm}\label{thm-IC}
Let
${\cal V}(X;L,d)$ be the collection of admissible, simple components as in Theorem \ref{nilo1}.
Then every $\GA$ in 
${\cal V}(X;L,d)$ 
determines the finite collection
$$
\mbox{\BM${\cal P}$} (\GA) = \left\{  \left.{\cal P}^{\bullet}_{\GA, \mu} \right| \mu\in P(\GA) \right\}
$$
of graded perverse sheaves 
$ {\cal P}^{\bullet}_{\GA, \mu}$ on $\XD$ indexed by the set of partitions $P(\GA)$ as in Proposition \ref{pro-T-orbits}.
\end{thm}

Taking the union of the collections $\mbox{\BM${\cal P}$} (\GA)$, as $\GA$ runs through the set ${\cal V}(X;L,d)$,
we obtain the finite collection
\BEN\label{P-collec}
\mbox{\BM${\cal P}$} (X;L,d) =\left\{ \left.{\cal P}^{\bullet}_{\GA, \mu} \right| \mu\in P(\GA) ,\,\, \GA\in {\cal V}(X;L,d) \right\}
\EEN 
of perverse sheaves on $\XD$ intrinsically associated to $(X,L,d)$. This is Theorem \ref{ps} of the Introduction. 
\begin{rem}\label{P-collec-geom}
It is clear that one can construct complexes
${\cal P}^{\bullet}_{\GA, \mu}$ for {\it any} partition $\mu$ of $\DE$. The main point of distinguishing the collection
$\mbox{\BM${\cal P}$} (\GA)$ is that the complexes of this collection pick out partitions of $\DE$ which are relevant
to the geometry of configurations of $X$ parametrized by $\GAB$. Indeed, in \S\ref{sec-equations} we have seen how partitions
in $P(\GA)$ are related to the equations defining configurations parametrized by $\GAB$. So, heuristically, one could say that
perverse sheaves ${\cal P}^{\bullet}_{\GA, \mu}\,(\mu \in P(\GA))$ condense in them those equations: the equations themselves
might be quite complicated (see e.g. (\ref{hom-eqns})) and one might want, for various purposes, to ``package" them neatly in the form
of perverse sheaves.
\end{rem} 

It was pointed out in the Introduction that the complexes in
$\mbox{\BM${\cal P}$} (\GA)$
also contain information about irreducible representations in  the collection $R_{\DE} (\GA)$ in (\ref{R-collec}).
To see this recall that one of the fundamental properties of the Springer resolution in (\ref{Springer-res}) is that the cohomology ring  of
the Springer fibres, the fibres of $\sigma$ in (\ref{Springer-res}), supports representations of the Weyl group $W$ of the Lie algebra in question.
In our situation the Lie algebra is ${\bf sl}_{\DE} (\CC)$. Hence the Weyl group
$$
W=S_{\DE}
$$
is the symmetric group $ S_{\DE}$ and Springer theory yields representations
\BEN\label{Sp-rep}
sp^{\bullet}_{\mu} : S_{\DE} \longrightarrow Aut (H^{\bullet} ( \sigma^{-1} (x), \CC))\,,
\EEN
where $x$ is a point of the nilpotent orbit
$\mbox{\BM$O_{\mu}$}$ of 
${\bf sl}_{\DE} (\CC)$ corresponding to a partition $\mu$ of $\DE$.
In particular, one knows (see \cite{[C-Gi]}) that the top degree cohomology group\footnote{recall: $b_{\mu} = dim_{\CC} (\sigma^{-1} (x))$.}
$H^{2b_{\mu}} ( \sigma^{-1} (x), \CC)$
is an irreducible 
$ S_{\DE} $-module corresponding to $\mu$.

The above discussion shows that the fibres of the cohomology sheaves of complexes
${\cal P}^{\bullet}_{\GA, \mu}$
are $S_{\DE} $-modules. In fact, this $S_{\DE} $-module structure is compatible with the action of
the fundamental group $\pi_1 ( \GAB, [Z_0])$ given by the representations in (\ref{rep-fg1}).
This is the meaning of the following statement.
\begin{pro}\label{Sp-rep-fg}
Let $\mu$ be a partition in $P(\GA)$ as in Proposition \ref{pro-T-orbits}. Set
\BEN\label{Sp-rep-fg1}
\rho^{\bullet}_{\GA,\mu} = \bigoplus^{2b_{\mu}}_{i=0} \rho^i_{\GA,\mu}\,,
\EEN
where $ \rho^i_{\GA,\mu}$'s are as in (\ref{rep-fg1}) and 
$b_{\mu}$ is the complex dimension  of Springer fibers over the nilpotent orbit $\mbox{\BM$O_{\mu}$}$.
Then the representation
$$
\rho^{\bullet}_{\GA,\mu} :\pi_1 ( \GAB,[Z_0]) \longrightarrow Aut (H^{\bullet} ({\bf B}_{([Z_0],x_0)}, \CC))
$$
factors through the Springer representation $sp^{\bullet}_{\mu}$ in (\ref{Sp-rep}).
\end{pro}
\begin{pf}
Recall the universal scheme $\ZD$ in (\ref{uc}) and consider its subscheme 
$\ZD_{\GAB} = p^{-1}_2 (\GAB)$  lying over $\GAB$ together with the unramified covering
$$
p_2 : \ZD_{\GAB} \longrightarrow \GAB\,.
$$
From Corollary \ref{cor-Fpr}, (\ref{Z-Zpr}), this factors as follows
\BEN\label{Z-fact}
\xymatrix{
{\ZD_{\GAB}} \ar[rr]^f \ar[rd]_{p_2} & &{\ZD^{\prime}_{\GAB}} \ar[dl]^{p^{\prime}_2} \\
&\GAB& }
\EEN
The morphism $p^{\prime}_2$, by Remark \ref{et}, is an unramified covering of degree $\DE$,
so its fibre over a point $[Z] \in \GAB$ is the set $Z^{\prime}$ of $\DE$ distinct points. Hence we obtain a group homomorphism
\BEN\label{fg-AutZpr}
\rho_{p^{\prime}_2} : \pi_1 ( \GAB, [Z]) \longrightarrow Aut (Z^{\prime})\,,
\EEN
where $Aut (Z^{\prime})$ is viewed as the set of all possible orderings of $Z^{\prime}$. 

We also recall that the fibre 
$\mbox{\BM$\GS^{\prime}_{\mbox{\UB$\GA$}}$} ([Z])$
of 
$\mbox{\BM$\GS^{\prime}_{\mbox{\UB$\GA$}}$}$
over $[Z]$ can be identified with
${\bf sl} (\HO{Z^{\prime}}))$ and $\HO{Z^{\prime}})$ is identified with a Cartan subalgebra of
${\bf gl} (\HO{Z^{\prime}}))$ via the multiplicative action of $\HO{Z^{\prime}})$ on itself (see Remark \ref{act-m}). Thus the subspace 
\BEN\label{hZ}
\goth{h}_{Z^{\prime}} =\left\{ \left.f\in \HO{Z^{\prime}}) \right| Tr(f) =\sum_{z^{\prime} \in Z^{\prime}} f(z^{\prime}) =0 \right\}
\EEN
can be identified with a Cartan subalgebra of 
${\bf sl} (\HO{Z^{\prime}}))$. In particular, permutations of  $Z^{\prime}$ act on the space
$\HO{Z^{\prime}})$ (resp. $\goth{h}_{Z^{\prime}}$) by the rule
\BEN\label{AutZpr-func}
(w^{\ast}f) (z) = f(w^{-1}(z)),\,\forall z \in Z^{\prime}\,.
\EEN
Hence the representation
\BEN\label{rep-stand}
Aut (Z^{\prime}) \longrightarrow Aut (\goth{h}_{Z^{\prime}})\,.
\EEN
This representation determines and is determined by the permutation action of $Aut (Z^{\prime})$  
on the set\footnote{the set ${\bf B}(\goth{h}_{Z^{\prime}})$ can be identified with the set of orderings of eigen spaces of
$\goth{h}_{Z^{\prime}}$-action on $\HO{Z^{\prime}})$. Those eigen spaces are generated by the delta-functions $\delta_{z^{\prime}}$, for
$z^{\prime} \in Z^{\prime}$. Thus ${\bf B}(\goth{h}_{Z^{\prime}})$ has a natural identification with the set of orderings of $Z^{\prime}$.} 
${\bf B}(\goth{h}_{Z^{\prime}})$ of Borel subalgebras of 
${\bf sl} (\HO{Z^{\prime}}))$ containing $\goth{h}_{Z^{\prime}}$. %(see \cite{[C-Gi]}, Claim 3.1.19)
 Thus we obtain the representation
\BEN\label{Sp-rep-sr}
sp_{rs}: Aut (Z^{\prime}) \longrightarrow Aut ({\bf B}(\goth{h}_{Z^{\prime}}))\,.
\EEN
Identifying $Aut (Z^{\prime})$ with the Weyl group of ${\bf sl}(\HO{Z^{\prime}})$ and ${\bf B}(\goth{h}_{Z^{\prime}})$
with a Springer fibre over a regular semisimple element of ${\bf sl}(\HO{Z^{\prime}})$, we see that $sp_{rs}$ in (\ref{Sp-rep-sr})
is the Springer representation over the set of regular semisimple
elements ${\bf sl}^{rs} (\HO{Z^{\prime}}))$ of ${\bf sl} (\HO{Z^{\prime}}))$.
Combining this with (\ref{fg-AutZpr}) yields
\BEN\label{fg-Sp-rep-sr}
\rho_{\GA,rs} = sp_{rs} \circ \rho_{p^{\prime}_2}: \pi_1 ( \GAB, [Z]) \longrightarrow Aut ({\bf B}(\goth{h}_{Z^{\prime}}))
 \EEN
which is the assertion of the proposition over the orbit of regular
 semisimple elements of  ${\bf sl} (\HO{Z^{\prime}}))$. 

Next recall that the Springer resolution in (\ref{Springer-res}) is a part of the universal resolution (\cite{[C-Gi]}, 3.1.31)
\BEN\label{univ-res}
\xymatrix{
&{\bf \tilde{sl}} (\HO{Z^{\prime}})) \ar[dl]_{pr_1} \ar[dr]^{pr_2}& \\
{\bf sl} (\HO{Z^{\prime}})) &  &{\bf B}_{[Z]} }
\EEN
where ${\bf B}_{[Z]} $ is the variety of Borel subalgebras of ${\bf sl} (\HO{Z^{\prime}}))$ and 
${\bf \tilde{sl}} (\HO{Z^{\prime}}))$ is the incidence correspondence
\BEN
{\bf \tilde{sl}} (\HO{Z^{\prime}})) =\left\{ \left.(x,{\bf b}) \in {\bf sl} (\HO{Z^{\prime}})) \times {\bf B}_{[Z]} \right| x \in {\bf b} \right\}\,.
\EEN
The set of regular semisimple elements 
${\bf sl}^{rs} (\HO{Z^{\prime}}))$ is a dense Zariski open subset of
${\bf sl} (\HO{Z^{\prime}}))$ and the Springer representation 
$sp^{\bullet}_{\mu}$ in (\ref{Sp-rep}) over the nilpotent orbit $\mbox{\BM$O_{\mu}$}$ can be obtained as the limit of Springer 
representations 
$$
sp_{rs, h}: S_{\DE} \longrightarrow Aut( pr^{-1}_1 (h)) 
$$ 
on the fibres  of $pr_1$ in (\ref{univ-res}) over $h $, regular semisimple, converging to 
$\mbox{\BM$O_{\mu}$}$ (see \cite{[C-Gi]}, 3.4). Hence the representation $\rho^{\bullet}_{\GA,\mu}$ in (\ref{Sp-rep-fg1}),
 which is the limit of
 the representations
$\rho_{\GA,rs, h}$, for $h \in \goth{h}_{Z^{\prime}}$ converging to 
$\mbox{\BM$O_{\mu}$}$, can now be expressed as follows
$$
\rho^{\bullet}_{\GA,\mu} =\lim_{h\to \mbox{\BM$O_{\mu}$}} \rho_{\GA,rs, h} =
\lim_{h\to \mbox{\BM$O_{\mu}$}} sp_{rs,h} \circ \rho_{p^{\prime}_2} =
\left(\lim_{h\to \mbox{\BM$O_{\mu}$}} sp_{rs,h} \right) \circ\rho_{p^{\prime}_2} = 
sp^{\bullet}_{\mu} \circ\rho_{p^{\prime}_2}\,,
$$
where the second equality is the validity of the proposition over regular semisimple elements obtained in (\ref{fg-Sp-rep-sr}).
\end{pf}
\subsection{Abelian category ${\cal A}(X;L,d)$}\label{A(XLd)}
The collection of perverse sheaves 
$\mbox{\BM${\cal P}$}(X;L,d)$ introduced in (\ref{P-collec}) gives rise to a distinguished full abelian subcategory of the bounded derived category ${\cal D}^b_c (\XD)$ of constructible sheaves on 
$\XD$.
This is defined as follows.

For every  
${\cal P}^{\bullet}_{\GA,\mu}$ in 
$\mbox{\BM${\cal P}$}(X;L,d)$ consider its components ${\cal P}^i_{\GA,\mu}$ as in (\ref{IC-gr}) and take its irreducible constituents in the category of perverse sheaves on $\XD$. More precisely, since ${\cal P}^i_{\GA,\mu}$ is essentially the intersection cohomology complex
$IC(\GAB, {\cal L}^i_{\GA,\mu} )$, the irreducible constituents of ${\cal P}^i_{\GA,\mu}$ are given by the decomposition of
${\cal L}^i_{\GA,\mu}$ into the direct sum of irreducible local systems, i.e. we go back to the representation
$\rho^i_{\GA,\mu}$ in (\ref{rep-fg1}) and decompose $H^i({\bf B}_{([Z_0],x_0)}, \CC)$ into the direct sum of irreducible
$\pi_1 (\GAB,[Z_0])$-modules
\BEN\label{Hi-fg-irr}
H^i({\bf B}_{([Z_0],x_0)}, \CC) =\bigoplus_{\chi} M^i_{\chi,\GA,\mu} \otimes V^i_{\chi,\GA,\mu}\,,
\EEN
where 
$V^i_{\chi,\GA,\mu}$'s are irreducible $\pi_1 (\GAB,[Z_0])$-modules occurring in 
$H^i({\bf B}_{([Z_0],x_0)}, \CC)$ and 
$M^i_{\chi,\GA,\mu}$'s are their respective multiplicity modules
\BEN\label{mult-mod-fg}
M^i_{\chi,\GA,\mu} = Hom_{\pi_1 (\GAB,[Z_0])}(V^i_{\chi,\GA,\mu},H^i({\bf B}_{([Z_0],x_0)}, \CC))\,.
\EEN

From Proposition \ref{Sp-rep-fg} we know that the representation $\rho^i_{\GA,\mu}$ factors through the representation
$\rho_{p^{\prime}_2}$. Hence the irreducible modules
$V^i_{\chi,\GA,\mu}$ are parametrized by irreducible characters of the finite subgroup
$$
Im(\rho_{p^{\prime}_2}) \subset Aut( Z^{\prime}) \cong S_{\DE}\,.
$$

Set $\Xi^i_{\GA,\mu}$ to be the collection of irreducible characters occurring in the decomposition (\ref{Hi-fg-irr}). For every
$\chi \in \Xi^i_{\GA,\mu}$, we have the representation
\BEN\label{chi-rep-fg}
 \rho^i_{\chi,\GA,\mu} : \pi_1 (\GAB,[Z_0]) \longrightarrow Aut( V^i_{\chi,\GA,\mu})\,.
\EEN
The corresponding local system on $\GAB$ will be denoted by 
${\cal L}^i_{\GA,\mu,\chi}$. This defines the Intersection cohomology complex $IC(\GAB,{\cal L}^i_{\GA,\mu,\chi})$. 
Its extension by zero to the whole of $\XD$ will be denoted by ${\cal C}^i_{\GA,\mu,\chi}$. This now is an irreducible perverse sheaf on $\XD$ and we have
\BEN\label{dec-P-irr}
{\cal P}^i_{\GA,\mu} = \bigoplus_{\chi \in \Xi^i_{\GA,\mu}} M^i_{\chi,\GA,\mu} \otimes {\cal C}^i_{\GA,\mu,\chi}\,,
\EEN
the decomposition of ${\cal P}^i_{\GA,\mu}$ into the direct sum of its irreducible constituents. Thus for every pair
$(\GA, \mu)$ we have the finite collection
\BEN\label{i-G-mu-collec}
\mbox{\BM$C$}_{\GA,\mu} =\left\{\left.{\cal C}^i_{\GA,\mu,\chi} \right| i=0,\ldots,2b_{\mu}, \,\,\chi \in \Xi^i_{\GA,\mu} \right\}
\EEN
of irreducible perverse sheaves on $\XD$. Taking the union of these collections, as $\mu$ runs through the set $P(\GA)$ defined in 
Proposition \ref{pro-T-orbits}, and $\GA$ runs through ${\cal V} (X;L,d)$ in Theorem \ref{nilo1}, we obtain the collection
\BEN\label{ps-irr-collec}
\mbox{\BM$C$} (X;L,d) = \bigcup_{\mu \in P(\GA), \GA \in {\cal V} (X;L,d)} \mbox{\BM$C$}_{\GA,\mu}
\EEN
of irreducible perverse sheaves on $\XD$ intrinsically associated to $(X,L,d)$.

We now define the abelian category
${\cal A}(X;L,d)$ as the full subcategory of the derived category of constructible sheaves on $\XD$ generated by the finite collection
$\mbox{\BM$C$} (X;L,d)$, i.e. objects of ${\cal A}(X;L,d)$ are isomorphic to finite direct sums
of objects in $\mbox{\BM$C$} (X;L,d)$ and their various translations.
\begin{defi}\label{non-ab-char}
\begin{enumerate}
\item[1)]
Elements of the collection $\mbox{\BM$C$} (X;L,d)$ will be called irreducible non-abelian characters of $\JA$.
\item[2)]
The abelian category ${\cal A}(X;L,d)$ will be called the category of non-abelian characters of $\JA$ and objects of 
${\cal A}(X;L,d)$ will be called non-abelian characters of $\JA$.
\end{enumerate}
\end{defi}
\begin{rem}\label{analogy}
 Recall that for a smooth irreducible curve $C$ its Jacobian
$J(C)$ is an abelian variety and its fundamental group
$\pi_1 (J(C)) =H_1 (J(C), {\bf Z})$ is just the first homology group of $J(C)$. Thus irreducible local systems on
$J(C)$ are given by homomorphisms or characters
$$ 
H_1 (J(C), {\bf Z}) \longrightarrow {\CC}^{\times}
$$
and the group of characters $Hom (H_1 (J(C), {\bf Z}), {\CC}^{\times})$ parametrizes isomorphism classes of local systems on $J(C)$.

With the above in mind, the collection $\mbox{\BM$C$} (X;L,d)$ defined in (\ref{ps-irr-collec}) can be envisaged as a non-abelian analogue of 
$Hom (H_1 (J(C), {\bf Z}), {\CC}^{\times})$, while
the abelian category ${\cal A}(X;L,d)$ can be viewed as an analogue of the group-ring of $Hom (H_1 (J(C), {\bf Z}), {\CC}^{\times})$.
This, hopefully, justifies the terminology in Definition \ref{non-ab-char}.
\end{rem}

Though the objects of the abelian category ${\cal A}(X;L,d)$ are complexes of sheaves on the Hilbert scheme $\XD$, one should really keep in mind that they descended from the Jacobian $\JA$. There is even more subtle connection which relates the sections of the relative 
tangent/cotangent sheaf of 
$\JA \setminus \mbox{\BM$\Theta$}(X;L,d)$ with the abelian category ${\cal A} (X;L,d)$.
\begin{pro}\label{tan-ps}
Let $\GA$ be a component in ${\cal V} (X;L,d)$ and let ${\cal T}_{\pi,\GA}$ be the relative tangent sheaf of the projection
 $\pi: \JABG \longrightarrow \GAB$.
Then there is a natural map
\BEN\label{tan-ps-map}
exp\left(\int_{\GA}\right) : H^0 (\JABG, {\cal T}_{\pi,\GA}) \longrightarrow {\cal A} (X;L,d)\,.
\EEN
\end{pro}
\begin{pf}
Let $\theta$ be a section of ${\cal T}_{\pi, \GA}$. Interpreting it geometrically, we view $\theta$ as the corresponding morphism
$$
\theta: \JABG \longrightarrow T_{\pi,\GA}
$$
of $\GAB$-schemes. We now consider the intersection of the image of $\theta$ with the strata $\mbox{\BM$O_{\mu}$}(T_{\pi,\GA})$ defined in
(\ref{T-orbits}). Denote by $P(\GA,\theta)$ the subset of partitions $\mu$ in $P(\GA)$ (see Proposition \ref{pro-T-orbits} for notation) such that 
$\theta^{-1} (\mbox{\BM$O_{\mu}$}(T_{\pi,\GA}))$ is non-empty. The map $exp(\int_{\GA})$ we are after can now be defined by sending $\theta$ to the direct sum of perverse sheaves ${\cal P}^{\bullet}_{\GA,\mu}$ as in (\ref{IC-gr}), where $\mu \in P(\GA,\theta)$:
 \BEN\label{tan-ps-map1}
exp\left(\int_{\GA}\right)(\theta) = \bigoplus_{\mu \in P(\GA,\theta)} {\cal P}^{\bullet}_{\GA,\mu} = \bigoplus_{\mu \in P(\GA,\theta)} 
\left( \bigoplus^{2b_{\mu}}_{i=0} \left(\bigoplus_{\chi \in \Xi^i_{\GA,\mu}} M^i_{\chi,\GA,\mu} \otimes {\cal C}^i_{\GA,\mu,\chi}\right) \right)\,. 
\EEN
\end{pf}
\begin{rem}\label{analogy1}
The notation $exp(\int_{\GA})$ is an allusion to the operation of integration followed by the exponential. This comes from pursuing the
 analogy with the classical Jacobian made in Remark \ref{analogy}. Namely, in the classical case one has a map
\BEN\label{classic-map}
H^0 (J(C), \Omega_{J(C)}) \longrightarrow Hom (H_1 (J(C), {\bf Z}), {\CC}^{\times})
\EEN
which sends a holomorphic $1$-form $\omega$ on $J(C)$ to the exponential of the linear functional
$$
\int \omega : H_1 (J(C), {\bf Z}) \longrightarrow \CC
$$
given by integrating $\omega$ over $1$-cycles on $J(C)$. 

The map $exp(\int_{\GA})$ in (\ref{tan-ps-map}) clearly has the same flavor, except that in our story we are ``integrating" (relative) vector fields on
$\JA$. However, the relative tangent sheaf ${\cal T}_{\pi,\GA}$ is self-dual\footnote{the self-duality of ${\cal T}_{\pi,\GA}$ comes from the isomorphism $M$ in (\ref{M}) and the identification of the quotient-sheaf 
$\HT /{\OO_{\JABG}}$ with the orthogonal complement $\HH =(\OO_{\JABG})^{\perp}$ of $\OO_{\JABG}$ in $\HT$. Now the quadratic form
$\QT$ in (\ref{q}) restricts to a non-degenerate quadratic form on $\HH$, thus making it self-dual.} and naturally isomorphic to the relative cotangent bundle
${\cal T}^{\ast}_{\pi,\GA}$. This gives an identification
\BEN\label{tan-cotan}
H^0 (\JABG,{\cal T}^{\ast}_{\pi,\GA}) \cong H^0 (\JABG, {\cal T}_{\pi,\GA})\,. 
\EEN
With this identification in mind, we can say that in the map (\ref{tan-ps-map}) we are ``integrating" the (relative) $1$-forms after all 
\BEN\label{cotan-ps-map}
exp\left(\int_{\GA}\right) : H^0 (\JABG,{\cal T}^{\ast}_{\pi,\GA}) \longrightarrow {\cal A} (X;L,d)\,,
\EEN
thus making Proposition \ref{tan-ps} conceptually analogous to the classical map (\ref{classic-map}).
\end{rem}

One can put the maps in (\ref{cotan-ps-map}) together as $\GA$ varies in ${\cal V} (X;L,d)$ to obtain the following.
\begin{thm}\label{thm-cot-ps-map}
Let $\stackrel{\circ}{\JAA}(X;L,d) = \JA \setminus \mbox{\BM$\Theta$}(X;L,d)$ be the complement of the theta-divisor 
$\mbox{\BM$\Theta$}(X;L,d)$ in $\JA$ and
let $ {\cal A} (X;L,d)$ be the category of non-abelian characters of $\JA$ (see Definition \ref{non-ab-char}). Then there is a natural map
\BEN\label{cot-ps-map}
exp\left(\int\right): H^0 ({\cal T}^{\ast}_{\stackrel{\circ}{\JAA}(X;L,d)/ {\XD}} ) \longrightarrow {\cal A} (X;L,d)\,,
\EEN
where ${\cal T}^{\ast}_{\stackrel{\circ}{\JAA}(X;L,d) / {\XD}}$ is the relative cotangent sheaf of $\stackrel{\circ}{\JAA}(X;L,d)$ over $\XD$.
\end{thm}
\begin{pf}
Let $\omega$ be a global section of ${\cal T}^{\ast}_{\stackrel{\circ}{\JAA}(X;L,d) / {\XD}}$. For every $\GA$ in ${\cal V} (X;L,d)$ denote by 
$\omega_{\GA}$ the restriction of $\omega$ to $\JABG$. This is a section of 
${\cal T}^{\ast}_{\JABG/{\GAB}}$ and 
$exp(\int_{\GA} ) (\omega_{\GA})$ has been defined in (\ref{cotan-ps-map}). One can now define
\BEN\label{cot-ps-map1}
exp\left(\int\right) (\omega) = \bigoplus_{\GA \in {\cal V} (X;L,d)} exp\left(\int_{\GA} \right) (\omega_{\GA}) =\bigoplus_{\GA \in {\cal V} (X;L,d)}
\left( \bigoplus_{\mu \in P(\GA, \theta_{\GA})} {\cal P}^{\bullet}_{\GA,\mu} \right)\,,
\EEN
where $\theta_{\GA}$ is the section of the relative tangent sheaf ${\cal T}_{\JABG / {\GAB}}$ corresponding to $\omega_{\GA}$ under the isomorphism in (\ref{tan-cotan}) and $P(\GA, \theta_{\GA})$ is the subset of $P(\GA)$ defined in the proof of Proposition \ref{tan-ps}.
\end{pf}
\subsection{Generalized Macdonald functions on $\XD$}\label{sec-M-fun}
Here we define several invariants inspired by the well-known constructions in the theory of symmetric functions and the theory of representations
of symmetric groups.

By construction, the Intersection cohomology complexes
${\cal P}^{\bullet}_{\GA,\mu}$
are endowed with an action of the symmetric group $S_{\DE}$.
Hence the cohomology sheaves 
\BEN\label{coh-sh}
{\cal H}^{j,i}_{\GA,\mu} := {\cal H}^j ({\cal P}^i_{\GA,\mu})
\EEN
are $S_{\DE}$-sheaves, for every $i,j$.

Let $P_{\DE}$ be the set of partitions of $\DE$ and let $\lambda$ be a partition in $P_{\DE}$.
Denote by ${\bf S}_{\lambda}$ an irreducible 
$S_{\DE}$-module corresponding to $\lambda$ and, for every closed point $[Z] \in \XD$, define
\BEN\label{mult-ji}
m^{(i,j)}_{\lambda,\mu,\GA} ([Z])= dim_{\CC} Hom_{S_{\DE}} ({\bf S}_{\lambda}, {\cal H}^{j,i}_{\GA,\mu} ([Z]))\,,
\EEN
where
$ {\cal H}^{j,i}_{\GA,\mu} ([Z])$ is the fibre of $ {\cal H}^{j,i}_{\GA,\mu}$ in (\ref{coh-sh}) at $[Z]$.
This defines constructible functions on $\XD$
\BEN\label{func-m-ij}
m^{(i,j)}_{\lambda,\mu,\GA}: \XD \longrightarrow {\bf Z_{+}}
\EEN
which assigns to $[Z] \in \XD$ the value 
$m^{(i,j)}_{\lambda,\mu,\GA} ([Z])$ in (\ref{mult-ji}).

By definition of ${\cal P}^{i}_{\GA,\mu}$, the functions 
$m^{(i,j)}_{\lambda,\mu,\GA}$ are identically zero unless $i\in \{0,\ldots, 2b_{\mu} \}$ and
$j\in \{-dim_{\CC} \GA, \ldots, -1 \}$ (see (\ref{IC}), a)-d)), where $b_{\mu}$ is as in Remark \ref{dim-Sp-fib}.
Furthermore, from (\ref{IC}), a), it follows that 
the support of $m^{(i,j)}_{\lambda,\mu,\GA}$ is contained in the closure
$\overline{\GA}$ of $\GA$ in $\XD$.

Using the functions $m^{(i,j)}_{\lambda,\mu,\GA}$ as coefficients, we define polynomials 
$P_{\lambda,\mu,\GA}$ in two variables $q$ and $t$
\BEN\label{qt-pol}
P_{\lambda,\mu,\GA} = \sum_{i,k} n^{(i,k)}_{\lambda,\mu,\GA} q^i t^k\,,
\EEN
where $i =0,\ldots, 2b_{\mu}$ and $k=0,\ldots, dim_{\CC} \GA -1$, and the coefficients 
$n^{(i,k)}_{\lambda,\mu,\GA}$ are constructible functions on $\XD$ defined by the following identity
\BEN\label{coef-n}
n^{(i,k)}_{\lambda,\mu,\GA} = m^{(i,k-dim_{\CC} \GA)}_{\lambda,\mu,\GA} \,\,(k=0,\ldots, dim_{\CC} \GA -1)\,.
\EEN
These polynomials will be called {\it generalized Kostka-Macdonald coefficients of } $\GA$.

The final step of our construction is to put together the generalized Kostka-Macdonald coefficients for various $\lambda$.
Namely, consider the graded ring
$$
{\bf \Lambda} =\bigoplus^{\infty}_{n=0} {\bf \Lambda}^n
$$
of symmetric functions in infinitely many formal variables $ {\bf x} =(x_k)_{k\in {\bf N}}$, where $ {\bf \Lambda}^n$ is the subspace of ${\bf \Lambda}$ of homogeneous symmetric 
functions of degree $n$. Let
\BEN\label{shur}
\left\{\left.s_{\lambda} \right| \lambda \in P_{\DE} \right\}
\EEN
be the basis of $ {\bf \Lambda}^{\DE}$ formed by Schur functions\footnote{for basic facts and terminology concerning symmetric functions our reference is \cite{[Mac]}.} $s_{\lambda} \,(\lambda \in P_{\DE}  )$.

Define the function
\BEN\label{M-fun}
M_{\GA,\mu} : \XD \longrightarrow {\bf \Lambda}^{\DE} [q,t]
\EEN
by the following identity
\BEN\label{M-fun1}
M_{\GA,\mu} = \sum_{\lambda \in P_{\DE}} P_{\lambda,\mu,\GA} s_{\lambda}\,.
\EEN
This way we obtain a collection of functions 
\BEN\label{collec-M}
{\bf M} (X;L,d) =\left\{ \left.M_{\GA,\mu} \right| \GA \in {\cal V} (X;L,d), \,\,\mu \in P(\GA) \right\}
\EEN
intrinsically associated to $(X,L,d)$, where ${\cal V} (X;L,d)$ is as in Theorem \ref{nilo1} and $ P(\GA) $ is as in Proposition \ref{pro-T-orbits}).
\begin{rem}
Our definitions (\ref{mult-ji}), (\ref{qt-pol}) and (\ref{M-fun1}), as well as the terminology, are modeled on the ones in the work
of Haiman, \cite{[Hai]}, where he introduced Kostka-Macdonald coefficients as graded character multiplicities of certain
doubly graded $S_n$-modules, for every partition $\mu$ of $n$ (see \cite{[Hai]}, 2.2, for more details). In our situation it is the 
cohomology sheaf
${\cal H}^{\bullet} ({\cal P}^{\bullet}_{\GA,\mu})$
on $\XD$ which is double graded and endowed with $S_{\DE}$-action. So our constructions can be viewed as a natural generalization of
Haiman's bigraded modules and  the functions $M_{\GA,\mu}$ in (\ref{collec-M}) can be considered as
an analogue of Macdonald functions. 
\end{rem}

Below we summarize properties of functions $M_{\GA,\mu}$ defined in (\ref{M-fun1}).
\begin{pro-defi}
Let $\GA$ be a component of the set 
${\cal V} (X;L,d)$ in Theorem \ref{nilo1} and let $\mu$ be a partition in $ P(\GA) $, the set defined in Proposition \ref{pro-T-orbits}.
Then the function
$M_{\GA,\mu}$ has the following properties
\begin{enumerate}
\item[1)]
$M_{\GA,\mu}$ is a constructible function on $\XD$ with values in 
${\bf \Lambda}^{\DE} [q,t]$.
\item[2)]
The support of $M_{\GA,\mu}$ is contained in the closure $\overline{\GA}$ of $\GA$ in $\XD$.
\item[3)]
For every $[Z] \in \GAB$, the value
$M_{\GA,\mu} ([Z])$ is the $q$-polynomial with coefficients in ${\bf \Lambda}^{\DE}$
$$
M_{\GA,\mu} ([Z]) ({\bf x};q,t) =M_{\GA,\mu} ([Z])( {\bf x};q,0)
$$
which computes the graded character of the cohomology ring of a Springer fibre over the nilpotent orbit
$\mbox{\BM$O_{\mu}$}$ of ${\bf sl}_{\DE} (\CC)$ (see (\ref{M-fun-Z}) below for precise expression).
\end{enumerate}

The function $M_{\GA,\mu}$ will be called \rm{Macdonald function of $\XD$ of type $(\GA,\mu)$}.
\end{pro-defi}
\begin{pf}
All the properties of $M_{\GA,\mu}$, but the last one, are immediate from the definition of $M_{\GA,\mu}$.
For the last assertion use (\ref{IC}), b), to deduce
\BEN\label{coh-sh1}
\left. {\cal H}^{j,i} \right|_{\GAB} =\left\{\begin{array}{cc}
{\cal L}^i_{\GA,\mu},&if \,\, j=-dim_{\CC} \GA,\\
0,&otherwise.
\end{array} \right.
\EEN
This implies that the coefficients $n^{(i,k)}_{\lambda,\mu,\GA} ([Z]) =0$, unless $k=0$, and the generalized
Kostka-Macdonald coefficients
$P_{\lambda,\mu, \GA}$ evaluated at $[Z] \in \GAB$ have the following form
$$
P_{\lambda,\mu, \GA} (q,t) ([Z]) = P_{\lambda,\mu, \GA} (q,0) ([Z]) =\sum^{2b_{\mu}}_{i=0} n^{(i,0)}_{\lambda,\mu,\GA} ([Z]) q^i = 
\sum^{2b_{\mu}}_{i=0} m^{(i,-dim_{\CC} \GA)}_ {\lambda,\mu,\GA} ([Z]) q^i\,, 
$$
where, by (\ref{mult-ji}), the coefficient
$ m^{(i,-dim_{\CC} \GA)}_ {\lambda,\mu,\GA} ([Z])$ is the multiplicity of the irreducible 
$S_{\DE}$-module ${\bf S}_{\lambda}$ in the $S_{\DE}$-module
$H^i ({\bf B} ([Z], x), \CC)$ and
where ${\bf B} ([Z], x)$ is the Springer fibre of the morphism $\sigma_{\GA}$ in (\ref{morph-Sp-res}) over a point
$([Z],x)$ in the nilpotent orbit
$\mbox{\BM$O_{\mu} (\NI_{\mbox{\UB$\GA$}})$}$ (see (\ref{Sp-res-mu})). 
Hence $P_{\lambda,\mu, \GA} (q,0) ([Z])$ is the graded character multiplicity of the irreducible module 
 ${\bf S}_{\lambda}$ in the 
graded $S_{\DE}$-module 
$$
H^{\bullet} ({\bf B} ([Z], x), \CC) = \bigoplus^{2b_{\mu}}_{i=0} H^i ({\bf B} ([Z], x), \CC)
$$
and it can be written as follows
\BEN\label{P-Z}
P_{\lambda,\mu, \GA} (q,0) ([Z]) = \sum^{2b_{\mu}}_{i=0} \langle \chi^{\lambda}, ch(H^i ({\bf B} ([Z], x), \CC)) \rangle q^i\,, 
\EEN
where $\chi^{\lambda}$ is the character of ${\bf S}_{\lambda}$, $ ch(H^i ({\bf B} ([Z], x), \CC))$ is the character
of $ H^i ({\bf B} ([Z], x), \CC)$ and $\langle \cdot, \cdot \rangle$ stands for the standard pairing in the ring of characters 
$R(S_{\DE})$ of the symmetric group $S_{\DE}$. Setting
$$
ch(H^{\bullet} ({\bf B} ([Z], x), \CC),q) = \sum^{2b_{\mu}}_{i=0} ch(H^i ({\bf B} ([Z], x), \CC)) q^i
$$
to be the graded character of the cohomology ring $H^{\bullet} ({\bf B} ([Z], x), \CC)$, we can rewrite (\ref{P-Z}) in the following way 
\BEN\label{P-Z1}
P_{\lambda,\mu, \GA} (q,0) ([Z]) = \langle \chi^{\lambda}, ch(H^{\bullet} ({\bf B} ([Z], x), \CC),q) \rangle\,.
\EEN
Hence the generalized Macdonald function evaluated at $[Z] \in \GAB$
gives the graded character of the cohomology ring of the Springer fibre ${\bf B} ([Z], x)$
\BEN\label{M-fun-Z}
M_{\GA,\mu} ([Z]) ({\bf x}; q,t) =M_{\GA,\mu} ([Z]) ({\bf x}; q,0)=
\sum_{\lambda \in P_{\DE} } \langle \chi^{\lambda}, ch(H^{\bullet} ({\bf B} ([Z], x), \CC),q) \rangle s_{\lambda}\,.
\EEN
\end{pf}
\begin{rem}\label{chi-map}
The set of ${\bf M} (X;L,d)$ of generalized Macdonald functions defined in (\ref{collec-M}), can be viewed as a functional
counterpart of the collection of perverse sheaves
$\mbox{\BM${\cal P}$} (X;L,d)$ in (\ref{P-collec}). They are related by a sort of `character' map
\BEN\label{char-map}
\mbox{\BM$\chi$}_{ (X;L,d)} : \mbox{\BM${\cal P}$} (X;L,d) \longrightarrow {\bf M} (X;L,d)
\EEN
sending every Intersection cohomology complex ${\cal P}^{\bullet}_{\GA,\mu}$ in 
$\mbox{\BM${\cal P}$} (X;L,d)$ to its generalized 
Macdonald function $M_{\GA, \mu}$ defined in (\ref{M-fun1}), i.e.
\BEN\label{char-map1}
\mbox{\BM$\chi$}_{ (X;L,d)} ({\cal P}^{\bullet}_{\GA,\mu} ) =M_{\GA, \mu}\,,
\EEN
for every $\GA \in {\cal V} (L,d)$ and every $\mu \in P(\GA)$.
\end{rem}

\section{$\JA$ and the Langlands Duality}\label{sec-LD}
This section corresponds to the discussion in \S\ref{AffLie} of the Introduction. 
For $\GA$ an admissible, simple component of $\CS$, we consider the relative Infinite Grassmannian over $\GAB$ 
associated to the sheaf of Lie
algebras 
$\mbox{\BM$\GS^{\prime}_{\mbox{\UB$\GA$}}$}$ (see (\ref{LAGpr}) for notation). Then we show that there is a natural map
of a certain Zariski open subset of the relative tangent bundle 
$\TPI$ of $\JABG$ over $\GAB$ to this relative Infinite Grassmannian. This establishes, via the geometric Satake isomorphism, a link of our nonabelian Jacobian with
 the Langlands
Duality.
\subsection{Some preliminaries}\label{sec-IG}
In this subsection we fix notation and recall some known facts about Infinite Grassmannians. Our references on the subject are
\cite{[Gi]}, \cite{[Lu]}, \cite{[P-S]}.

We fix $\GA$ in the collection 
${\cal V} (X;L,d)$ defined in Theorem \ref{nilo1} and consider the locally free sheaf
\BEN\label{sh-Fpr}
\FF^{\prime} =p^{\prime}_{2\ast} (\OO_{\ZD^{\prime}})\,,
\EEN
where $\ZD^{\prime}$ and $p^{\prime}_2$ are as in (\ref{Z-fact}). 

Set
\BEN\label{LAGpr1}
\mbox{\BM$\GS^{\prime}_{\mbox{\UB$\GA$}}$} ={\bf sl} (\FF^{\prime})\,,
\EEN
 the sheaf of germs of traceless endomorphisms of $\FF^{\prime}$ and let
\BEN\label{LGpr}
\mbox{\BM$G^{\prime}_{\mbox{\UB$\GA$}}$} = {\bf SL}  (\FF^{\prime})
\EEN
be the corresponding sheaf of Lie groups. This will be regarded as a fibre bundle over $\GAB$
with the natural projection
\BEN\label{LG-proj}
\varpi: \mbox{\BM$G^{\prime}_{\mbox{\UB$\GA$}}$} \longrightarrow \GAB\,.
\EEN

Let ${\bf k} =\CC[t^{-1},t]$ be the ring of complex valued Laurent polynomials in a variable $t$ and let
${\bf o} =\CC[t]$ be its subring of polynomials in $t$. Consider the scheme
$\mbox{\BM$G^{\prime\mbox{rel}}_{\mbox{\UB$\GA$}} (k)$}$ of `vertical'
${\bf k}$-valued points of 
$ \mbox{\BM$G^{\prime}_{\mbox{\UB$\GA$}}$}$, where by {\it vertical} ${\bf k}$-valued point we mean a morphism
$$
Spec ({\bf k}) \longrightarrow  \mbox{\BM$G^{\prime}_{\mbox{\UB$\GA$}}$} 
$$
for which the diagram
$$
\xymatrix{
Spec ({\bf k}) \ar[r] \ar[d]& {\mbox{\BM$G^{\prime}_{\mbox{\UB$\GA$}}$}} \ar[d]^{\varpi} \\
Spec (\CC) \ar[r]& {\GAB} }
$$
commutes. Thus $\mbox{\BM$G^{\prime\mbox{rel}}_{\mbox{\UB$\GA$}} (k)$}$ comes with natural projection
\BEN\label{LG(k)-proj}
\varpi(\BK) : \mbox{\BM$G^{\prime\mbox{rel}}_{\mbox{\UB$\GA$}} (k)$} \longrightarrow \GAB\,,
\EEN
whose fibre $\mbox{\BM$G^{\prime\mbox{rel}}_{\mbox{\UB$\GA$}} (k)$} ([Z])$ over a point 
$[Z] \in \GAB$ is the group of $\BK$-valued points of the group
$ \mbox{\BM$G^{\prime}_{\mbox{\UB$\GA, [Z]$}}$}  ={\bf SL}(\HO{Z^{\prime}}))$, the fibre of 
$\varpi$ in (\ref{LG-proj}) over $[Z]$ and where $Z^{\prime}$ is the fibre over $[Z]$ of $p^{\prime}_2$ in (\ref{Z-fact}).

Replacing $\BK$ by $\BO$ gives the scheme
$\mbox{\BM$G^{\prime\mbox{rel}}_{\mbox{\UB$\GA$}} (o)$} $ and the coset space
\BEN\label{IGr}
\mbox{\BM$Gr_{\mbox{\UB$\GA$}} = 
G^{\prime\mbox{rel}}_{\mbox{\UB$\GA$}} (k) / G^{\prime\mbox{rel}}_{\mbox{\UB$\GA$}} (o)$}
\EEN
which we call {\it the Infinite Grassmannian of $\LIG$}. Thus $\IG$ is a scheme over $\GAB$ with the natural projection
\BEN\label{IG-proj}
\omega : \IG \longrightarrow \GAB\,,
\EEN
whose fibres are modeled on the Infinite Grassmannian of ${\bf SL}_{\DE}$.

We now give a description of $\IG$ which will reveal its ind-scheme structure.
For this 
consider the sheaf $\mbox{\BM$\GS^{\prime}_{\mbox{\UB$\GA$}} (k)$}$
(resp. $\mbox{\BM$\GS^{\prime}_{\mbox{\UB$\GA$}} (o)$}$) of Lie algebras over $\GAB$ with the bracket operation
defined in the usual way, i.e.
\BEN\label{bracket}
[x\otimes P, y\otimes Q ] =[x,y] \otimes PQ\,,
\EEN
for any local sections $x,y$ of $\mbox{\BM$\GS^{\prime}_{\mbox{\UB$\GA$}}$}$ and any
$P,Q \in \BK$.

Let $[Z]$ be a closed point of $\GAB$. Following Lusztig in \cite{[Lu]}, we consider
$ \BO$-submodules 
${\cal L}$ of $\mbox{\BM$\GS^{\prime}_{\mbox{\UB$\GA$}}$} ([Z]) (\BK) ={\bf sl}(H^0(\OO_{Z^{\prime}}))(\BK) $
 of maximal rank and which are closed under the Lie bracket
in (\ref{bracket}). For such  an ${\cal L}$ one defines its `dual'
${\cal L}^{\vee}$ whose elements $x$ are characterized by the condition
\BEN\label{L-dual}
\langle x,y \rangle \in \BO\,,
\EEN
for all elements $y$ of ${\cal L}$, and where 
$\langle \cdot, \cdot \rangle$ denotes the $\BK$-valued Killing form of 
 $\mbox{\BM$\GS^{\prime}_{\mbox{\UB$\GA$}} (k)$}$.

With these notions in mind one has the following description of 
$\IG$. Let 
$\mbox{\BM$Grass (\GS^{\prime}_{\mbox{\UB$\GA$}} (k))$}$
be the set whose points are pairs $([Z], {\cal L})$, where $[Z]$ is a closed point of $\GAB$ and ${\cal L}$
is  an $\BO$-submodule of 
$\mbox{\BM$\GS^{\prime}_{\mbox{\UB$\GA$}}$} ([Z]) (\BK)$. Then we have a map
\BEN\label{LG-Grass}
\LGR (\BK) \longrightarrow  \mbox{\BM$Grass (\GS^{\prime}_{\mbox{\UB$\GA$}} (k))$}
\EEN
which takes a local section $a$ of the fibration $\varpi(\BK)$ in (\ref{LG(k)-proj}) to the local section
of $ \mbox{\BM$Grass (\GS^{\prime}_{\mbox{\UB$\GA$}} (k))$}$
defined by the $\BO$-submodule
$Ad(a) ( \mbox{\BM$\GS^{\prime}_{\mbox{\UB$\GA$}} (o)$})$ of
$\mbox{\BM$\GS^{\prime}_{\mbox{\UB$\GA$}} (k)$}$, where
$Ad(a)$ denotes the adjoint action of
$\LGR (\BK)$ on $\mbox{\BM$\GS^{\prime}_{\mbox{\UB$\GA$}} (k)$}$.
According to \cite{[Lu]}, this establishes a bijection of 
$\IG$ in (\ref{IGr}) with the set
$\mbox{\BM$Grass^{0} (\GS^{\prime}_{\mbox{\UB$\GA$}} (k))$}$
of pairs $([Z], {\cal L})$, where $[Z]$ is a closed point of $\GAB$ and 
${\cal L}$ is an $\BO$-submodule of maximal rank in  
$\mbox{\BM$\GS^{\prime}_{\mbox{\UB$\GA$}}$} ([Z]) (\BK)$
which is closed under the Lie bracket and self-dual, i.e.
${\cal L} ={\cal L}^{\vee}$.

Using this description of $\IG$, one obtains the following stratification of $\IG$
by $\LGR(\BO)$-stable subschemes
\BEN\label{strat-IG}
\IG(1) \subset \IG(2) \subset \ldots \subset \IG(i) \subset \ldots
\EEN
where each stratum $\IG(i)$ is a finite dimensional algebraic variety defined as follows
\BEN\label{IG-i}
\IG(i) =\left\{ \left.([Z], {\cal L}) \in \mbox{\BM$Grass^{0} (\GS^{\prime}_{\mbox{\UB$\GA$}} (k))$} \right|
t^i \mbox{\BM$\GS^{\prime}_{\mbox{\UB$\GA$}}$} ([Z]) (\BO) \subset  {\cal L}\subset
t^{-i} \mbox{\BM$\GS^{\prime}_{\mbox{\UB$\GA$}}$} ([Z]) (\BO) \right\}\,.
\EEN
\subsection{$\LGR(\BO)$-orbits of $\IG$}\label{sec-IG-orbits}
By construction in \S\ref{sec-IG}, the ind-scheme $\IG$ is a fibre space over $\GAB$ with the natural projection
$$
\omega : \IG \longrightarrow \GAB
$$
 in (\ref{IG-proj}).  The fibres of $\omega$ are modeled on the Infinite Grassmannian ${\bf Gr}$ of 
${\bf SL}_{\DE} (\CC)$.

From the description in (\ref{IGr}) of $\IG$ as a quotient, it follows that $\IG$ admits a natural left action of 
$\LGR(\BO)$ and our goal is to understand the orbits of this action. 

By definition $\LGR(\BO)$ acts fibrewise on the fibration $\omega$ and its action on the fibres amounts to
the action of ${\bf SL}_{\DE} (\BO)$ on ${\bf Gr}$. It is well-known that ${\bf SL}_{\DE} (\BO)$-orbits of ${\bf Gr}$
are indexed by coweights (up to the action of the Weyl group) of 
${\bf SL}_{\DE}$ with respect to a maximal torus of ${\bf SL}_{\DE}$ (see \cite{[P-S]}).
It turns out that the same holds in our relative situation. To be more precise, we need to describe the relative, fibre version,
of the coweight lattice and the Weyl group. This in turn will clarify the structure of the bundles 
$\mbox{\BM$\GS^{\prime}_{\GA}$}$ and $\mbox{\BM$G^{\prime}_{\GA}$}$ as well as the structure of the Infinite Grassmannian
$\IG$.

Let us return to the morphism 
$$
p^{\prime}_2 : \ZD^{\prime}_{\GA} \longrightarrow \GAB
$$
 in (\ref{Z-fact}) and recall that 
the direct image 
$\FF^{\prime} =p^{\prime}_{2 \ast} ( \OO_{ \ZD^{\prime}_{\GA}} ) $ can be identified, via the multiplicative action of 
$\FF^{\prime}  $ on itself,  with the subsheaf of Cartan
subalgebras of
${\bf gl} (\FF^{\prime})$.  So we think of $\FF^{\prime}  $ as a subsheaf of ${\bf gl} (\FF^{\prime})$ and define
\BEN\label{Cartan1}
\mbox{\BM${\cal H}_{\GA}$} = \FF^{\prime} \bigcap \mbox{\BM$\GS^{\prime}_{\mbox{\UB$\GA$}}$}\,,
\EEN
 where 
$\mbox{\BM$\GS^{\prime}_{\mbox{\UB$\GA$}}$} ={\bf sl} (\FF^{\prime})$ is as in  (\ref{LAGpr1}).
This is a subsheaf of Cartan subalgebras in $\mbox{\BM$\GS^{\prime}_{\mbox{\UB$\GA$}}$}$ or, more geometrically,
it is a vector bundle over $\GAB$ which fits into the following commutative diagram
\BEN\label{Cartan-sc}
\xymatrix{
{\mbox{\BM${\cal H}_{\GA}$}} \ar@{^{(}{-}{>}}[rr] \ar[dr]_{\gamma^{\prime 0}} & & 
{\mbox{\BM$\GS^{\prime}_{\mbox{\UB$\GA$}}$}}\ar[dl]^{\gamma^{\prime}} \\
& {\GAB}& }
\EEN
where $\gamma^{\prime 0}$ is the restriction to
$\mbox{\BM${\cal H}_{\GA}$}$ of $\gamma^{\prime}$, the projection of
$\mbox{\BM$\GS^{\prime}_{\mbox{\UB$\GA$}}$}$ onto $\GAB$ in (\ref{prd+-}).
Thus the fibres of 
$\gamma^{\prime 0}$
are Cartan subalgebras in the fibres of $\gamma^{\prime}$, i.e. for a point $[Z] \in \GAB$, the fibre
$\mbox{\BM${\cal H}$}_{\GA, [Z]}$ of $\mbox{\BM${\cal H}$}_{\GA}$ over $[Z]$ is a Cartan subalgebra in
${\bf sl} (\HO{Z^{\prime}}))$, the fibre of $\gamma^{\prime}$ over $[Z]$.
\begin{rem}\label{HZ=funZ}
The fibre $\mbox{\BM${\cal H}$}_{\GA, [Z]}$, as a subspace of $\HO{Z^{\prime}})$, is identified 
with the subspace of functions on $Z^{\prime}$
whose trace is $0$, i.e.
\BEN\label{HZ1}
\mbox{\BM${\cal H}$}_{\GA, [Z]} = \left\{\left.f\in \HO{Z^{\prime}}) \right| \sum_{z^{\prime} \in Z^{\prime}} f(z^{\prime}) =0 \right\}\,.
\EEN
\end{rem}

We will now describe the scheme of roots and coroots of 
$\mbox{\BM$\GS^{\prime}_{\mbox{\UB$\GA$}}$}$ with respect to the subsheaf of Cartan subalgebras
$\mbox{\BM${\cal H}$}_{\GA}$.
 
First remark, that viewing $\FF^{\prime \ast}$, the dual of $\FF^{\prime}$, as a scheme over $\GAB$,
we have the natural inclusion
\BEN\label{Z-dFpr}
ev: \ZD^{\prime} \hookrightarrow  \FF^{\prime \ast}
\EEN
of $\GAB$-schemes given by evaluating functions on the underlying points, i. e. on the fibre $Z^{\prime}$
of $\ZD^{\prime}$ over a point $[Z] \in \GAB$, the above inclusion 
$$
ev_Z : Z^{\prime}  \hookrightarrow  (\FF^{\prime} ([Z])) ^{\ast}
$$
sends a point $z^{\prime} \in  Z^{\prime} $ to the linear function $ev_Z (z^{\prime})$ on $\FF^{\prime} ([Z])=\HO{Z^{\prime}}) $,
 the fibre of 
$\FF^{\prime}$ at $[Z]$, defined by evaluation at $z^{\prime}$:
$$
ev_Z (z^{\prime})(f) = f(z^{\prime}), \,\forall f\in \HO{Z^{\prime}})\,.
$$
Viewing the dual $\mbox{\BM${\cal H}^{\ast}$}_{\GA}$ of $\mbox{\BM${\cal H}$}_{\GA}$ as a quotient bundle
of $\FF^{\prime \ast}$ induces the morphism
\BEN\label{Z-Hg}
ev^{\prime} : \ZD^{\prime} \longrightarrow  \mbox{\BM${\cal H}^{\ast}$}_{\GA}\,.
\EEN
Denote its image by $\mbox{\BM$\Phi$}_{\GA}$ and observe that the fibres of 
$\mbox{\BM$\Phi$}_{\GA}$ over $\GAB$ are weights of fibres\footnote{these are the Lie algebras
${\bf sl} (\HO{Z^{\prime}}))$, for $[Z] \in \GAB$.} 
$\mbox{\BM$\GS^{\prime}_{\mbox{\UB$\GA$}}$}$ with respect to the Cartan subalgebras given by the fibres of the sheaf
$\mbox{\BM${\cal H}$}_{\GA}$. So the subvariety 
$\mbox{\BM$\Phi$}_{\GA}$  is {\it the scheme of weights of 
$\mbox{\BM$\GS^{\prime}_{\mbox{\UB$\GA$}}$}$ with respect to the subsheaf of Cartan subalgebras
$\mbox{\BM${\cal H}$}_{\GA}$}. Observe that $ev^{\prime}$ is still an inclusion, so 
$\mbox{\BM$\Phi$}_{\GA}$ is isomorphic to $\ZD^{\prime}$, thus giving a representation theoretic meaning of the scheme
$\ZD^{\prime}$. 

It should also be observed that the weights parametrized by the scheme $\mbox{\BM$\Phi$}_{\GA}$ span the sheaf of weight lattices of
$\mbox{\BM$\GS^{\prime}_{\mbox{\UB$\GA$}}$}$ with respect to the sheaf of Cartan subalgebras $\mbox{\BM${\cal H}$}_{\GA}$.
This is the content of the following result.
\begin{pro}\label{weight-lat}
Let $\mbox{\BM$\Lambda$}_{\GA}$ be ${\bf Z}$-span of $\mbox{\BM$\Phi$}_{\GA}$ in
$\mbox{\BM${\cal H}^{\ast}$}_{\GA}$. Then 
$\mbox{\BM$\Lambda$}_{\GA}$ is the sheaf of weight lattices of 
$\mbox{\BM$\GS^{\prime}_{\mbox{\UB$\GA$}}$}$ with respect to the subsheaf of Cartan subalgebras
$\mbox{\BM${\cal H}$}_{\GA}$. Furthermore, the isomorphism
$$
 ev^{\prime} : \ZD^{\prime} \longrightarrow \mbox{\BM$\Phi$}_{\GA}
$$
induces the isomorphism of sheaves
\BEN\label{wei-lat=lsys}
p^{\prime}_{2 \ast} (\underline{\bf Z}_{\ZD^{\prime}}) / {\underline{\bf Z}_{\GAB}} \cong \mbox{\BM$\Lambda$}_{\GA}\,,
\EEN
where $\underline{\bf Z}_{\ZD^{\prime}}$ (resp. $\underline{\bf Z}_{\GAB}$) denotes the constant sheaf ${\bf Z}$
on $\ZD^{\prime}$ (resp. $\GAB$).
\end{pro}
\begin{rem}\label{coweights}
Set $\mbox{\BM$\check{\Lambda}$}_{\GA}$ to be the subset of integer valued functions in $ \mbox{\BM${\cal H}$}_{\GA}$. This
is a fibre space over $\GAB$ whose fibre over a point $[Z] \in \GAB$ is the lattice 
$$
\mbox{\BM$\check{\Lambda}$}_{\GA, [Z]} =\left\{\left. f \in \mbox{\BM${\cal H}$}_{\GA, [Z]} \right| 
 f(z^{\prime}) \in {\bf Z},\,\forall z^{\prime} \in Z^{\prime} \right\}\,.
$$
Thus the fibre space $\mbox{\BM$\check{\Lambda}$}_{\GA}$ is the sheaf of the coweight lattices of 
$\mbox{\BM$\GS^{\prime}_{\mbox{\UB$\GA$}}$}$ with respect to the subsheaf of Cartan subalgebras
$\mbox{\BM${\cal H}$}_{\GA}$.
\end{rem}

Set $(\ZD^{\prime}) ^2  = \ZD^{\prime} \times_{\GAB} \ZD^{\prime}$ to be the fibre product of 
$\ZD^{\prime}$ with itself over $\GAB$ and let
$\Delta_{\ZD^{\prime}}$ be the diagonal in $(\ZD^{\prime})^ 2$. Define
\BEN\label{diag-comp}
\ZD^{\prime (2)} =(\ZD^{\prime})^2 \setminus \Delta_{\ZD^{\prime}}
\EEN
to be the complement of the diagonal in $(\ZD^{\prime})^2$. The following statement gives the representation theoretic meaning of this scheme.
\begin{pro-defi}\label{Z2-roots}
There is a natural morphism
$$
ev^{(2)}: \ZD^{\prime (2)} \longrightarrow \mbox{\BM${\cal H}^{\ast}$}_{\GA}\,,
$$
whose image $\mbox{\BM${\cal R}$}_{\GA}$ parametrizes the roots of 
$\mbox{\BM$\GS^{\prime}_{\mbox{\UB$\GA$}}$}$ with respect to the subsheaf of Cartan subalgebras
$\mbox{\BM${\cal H}$}_{\GA}$.

The scheme $\mbox{\BM${\cal R}$}_{\GA}$ will be called the root variety of 
$\mbox{\BM$\GS^{\prime}_{\mbox{\UB$\GA$}}$}$ with respect to the subsheaf of Cartan subalgebras
$\mbox{\BM${\cal H}$}_{\GA}$.
\end{pro-defi}
\begin{pf}
Consider the fibre product of the morphism $ev^{\prime}$ in (\ref{Z-Hg})
$$
(ev^{\prime})^2 : (\ZD^{\prime}) ^2 \longrightarrow \mbox{\BM${\cal H}^{\ast}$}_{\GA} \times_{\GAB} \mbox{\BM${\cal H}^{\ast}$}_{\GA}
$$
and compose it with the ``difference'' morphism
$$
\delta :\mbox{\BM${\cal H}^{\ast}$}_{\GA} \times_{\GAB} \mbox{\BM${\cal H}^{\ast}$}_{\GA}  \longrightarrow
 \mbox{\BM${\cal H}^{\ast}$}_{\GA}
$$
defined by $\delta (x,y) =x-y$, for any local sections $x$ and $y$ of $\mbox{\BM${\cal H}^{\ast}$}_{\GA}$.
The resulting morphism restricted to
$\ZD^{\prime (2)}$ will be denoted by $ev^{(2)}$.

To see the second assertion it is enough to consider the restriction $ev^{(2)}_Z$ of 
$ev^{(2)}$ to the fibre of $\ZD^{\prime (2)} $ over a point $[Z] \in \GAB$
$$
ev^{(2)}_Z : Z^{\prime (2)} \longrightarrow \mbox{\BM${\cal H}^{\ast}$}_{\GA, [Z]}\,.
$$
This map sends a pair of distinct points $(p^{\prime}, q^{\prime}) \in Z^{\prime (2)}$ to
the functional
$$
r_{p^{\prime}, q^{\prime}} = ev_Z (p^{\prime}) -  ev_Z (q^{\prime})
$$
which acts on $\HO{Z^{\prime}})$ as follows
\BEN\label{root-fun}
r_{p^{\prime}, q^{\prime}} (f) =f(p^{\prime}) - f(q^{\prime})\,.
\EEN
We claim that $r_{p^{\prime}, q^{\prime}}$, for $(p^{\prime}, q^{\prime}) \in Z^{\prime (2)} $, are roots of 
${\bf sl} (\HO{Z^{\prime}}))$ with respect to the Cartan subalgebra
$\mbox{\BM${\cal H}$}_{\GA, [Z]}$. Indeed, the action (by multiplication) of 
$\mbox{\BM${\cal H}$}_{\GA, [Z]}$ on 
$\HO{Z^{\prime}})$ determines the basis of $\HO{Z^{\prime}})$ consisting of the delta-functions
$\delta_{p^{\prime}}$, for $p^{\prime} \in Z^{\prime}$. Then for every pair
  $(p^{\prime}, q^{\prime})  \in Z^{\prime (2)}$ the endomorphism
$E_{p^{\prime}, q^{\prime}} \in {\bf sl} (\HO{Z^{\prime}}))$ which sends $\delta_{q^{\prime}}$ to 
$\delta_{p^{\prime}}$ and annihilates all other basis vectors, is a root vector of ${\bf sl} (\HO{Z^{\prime}}))$ with respect 
to the Cartan subalgebra
$\mbox{\BM${\cal H}$}_{\GA, [Z]}$. The root vectors $E_{p^{\prime}, q^{\prime}}$, as 
$(p^{\prime}, q^{\prime}) $ runs trough $Z^{\prime (2)}$, are linearly independent and give a complete set of representatives of 
the root spaces of ${\bf sl} (\HO{Z^{\prime}}))$. Furthermore, the action of 
$\mbox{\BM${\cal H}$}_{\GA, [Z]}$ on $E_{p^{\prime}, q^{\prime}}$ is as described in (\ref{root-fun}), i.e.
$r_{p^{\prime}, q^{\prime}}$ is a root of ${\bf sl} (\HO{Z^{\prime}}))$ with respect to the Cartan subalgebra
$\mbox{\BM${\cal H}$}_{\GA, [Z]}$, for every $(p^{\prime}, q^{\prime}) \in Z^{\prime (2)}$.
\end{pf}

From the above proof it also follows that the groups of automorphisms of fibres of the covering
$$
p^{\prime}_2 : \ZD^{\prime}_{\GA} \longrightarrow \GAB
$$
are naturally identified with the Weyl groups of the fibres of $\mbox{\BM$\GS^{\prime}_{\mbox{\UB$\GA$}}$}$.
Furthermore, the way these fibre groups fit together to form a fibre bundle of groups over $\GAB$ is given by the
representation 
$\rho_{p^{\prime}_2}$ which was already used in the proof of Proposition \ref{Sp-rep-fg}, (\ref{fg-AutZpr}). 
More precisely,
fix a point $[Z]$ in $\GAB$ and let
\BEN\label{univ-cov}
\tilde{\GAB} \longrightarrow \GAB
\EEN
be the universal covering of $\GAB$. Then the fibre bundle $\mbox{\BM$ W$}_{\GA}$ of the Weyl groups of
$\mbox{\BM$\GS^{\prime}_{\mbox{\UB$\GA$}}$}$ with respect to the subsheaf of Cartan subalgebras
$\mbox{\BM${\cal H}$}_{\GA}$ can be described  as follows
\BEN\label{Weyl-fg}
\mbox{\BM$ W$}_{\GA} = \tilde{\GAB} \times_{(\pi_1 (\GAB, [Z]),\rho_{p^{\prime}_2})} Aut (Z^{\prime})\,,
\EEN
where the action of $\pi_1 (\GAB, [Z])$ on the first factor is by the deck transformations and on the second via 
the representation $\rho_{p^{\prime}_2}$ in (\ref{fg-AutZpr}). 
\begin{rem}\label{fg-desc}
The fibre bundle $\mbox{\BM$ W$}_{\GA}$ can be viewed as the {\rm{abstract}} Weyl group of 
$\mbox{\BM$\GS^{\prime}_{\mbox{\UB$\GA$}}$}$ with respect to the subsheaf of Cartan subalgebras
$\mbox{\BM${\cal H}$}_{\GA}$. The actual action of $\mbox{\BM$ W$}_{\GA}$ on the variety of roots
$\mbox{\BM${\cal R}$}_{\GA}$ is given via the identification
$\ZD^{\prime (2)}$ with $\mbox{\BM${\cal R}$}_{\GA}$ provided by Proposition-Definition \ref{Z2-roots}.
In particular, this identification over the fixed point $[Z] \in \GAB$ induces the representation
\BEN\label{fg-rep-roots}
\rho^{roots}_{p^{\prime}_2} : \pi_1 (\GAB, [Z]) \longrightarrow Aut (\mbox{\BM${\cal R}$}_{\GA, [Z]} )\,,
\EEN
where $\mbox{\BM${\cal R}$}_{\GA, [Z]}$ is the fibre of $\mbox{\BM${\cal R}$}_{\GA}$ over $[Z]$.
This gives the following description of the variety of roots
\BEN\label{rootvar-fg}
\mbox{\BM${\cal R}$}_{\GA} = \tilde{\GAB} \times_{(\pi_1 (\GAB, [Z]), \rho^{roots}_{p^{\prime}_2})} \mbox{\BM${\cal R}$}_{\GA, [Z]}\,.
\EEN

The analogous description holds for all relative objects related to the bundle of Lie algebras
$\mbox{\BM$\GS^{\prime}_{\mbox{\UB$\GA$}}$}$ discussed so far. For example, we have the induced representation
\BEN\label{fg-rep-Cartan}
\rho^{\scriptscriptstyle{Cartan}}_{p^{\prime}_2} : \pi_1 (\GAB, [Z]) \longrightarrow Aut (\mbox{\BM${\cal H}$}_{\GA, [Z]})
\EEN
giving the identification
\BEN\label{Cartan-fg}
\mbox{\BM${\cal H}$}_{\GA} = 
\tilde{\GAB} \times_{(\pi_1 (\GAB, [Z]), \rho^{\scriptscriptstyle{Cartan}}_{p^{\prime}_2})} \mbox{\BM${\cal H}$}_{\GA, [Z]}\,.
\EEN
The representation (\ref{fg-rep-Cartan}) gives rise to the representation on the coweight lattice 
\BEN\label{fg-rep-cowei}
\rho^{\scriptscriptstyle{coweights}}_{p^{\prime}_2} : \pi_1 (\GAB, [Z]) \longrightarrow Aut (\mbox{\BM$\check{\Lambda}$}_{\GA, [Z]})\,,
\EEN
where $\mbox{\BM$\check{\Lambda}$}_{\GA, [Z]}$ is the fibre of the sheaf of coweight lattices 
$\mbox{\BM$\check{\Lambda}$}_{\GA}$ introduced in Remark \ref{coweights}. This gives the following identification
\BEN\label{cowei-fg}
\mbox{\BM$\check{\Lambda}$}_{\GA} =
\tilde{\GAB} \times_{(\pi_1 (\GAB, [Z]), \rho^{\scriptscriptstyle{coweights}}_{p^{\prime}_2})} \mbox{\BM$\check{\Lambda}$}_{\GA, [Z]}\,.
\EEN
\end{rem}

It is well-known that the representation of the Weyl group on a Cartan subalgebra of a semisimple Lie algebra lifts to
a representation on the whole Lie algebra. Thus the  representation $\rho^{\scriptscriptstyle{Cartan}}_{p^{\prime}_2}$
in (\ref{fg-rep-Cartan}) lifts to the representation
\BEN\label{fg-rep-La}
\rho^{\scriptscriptstyle{la}}_{p^{\prime}_2} : \pi_1 (\GAB, [Z]) \longrightarrow Aut \left({\bf sl} (\HO{Z^{\prime}}))\right)\,.
\EEN
This yields the following identification of 
$\mbox{\BM$\GS^{\prime}_{\mbox{\UB$\GA$}}$}$ as a fibre bundle over $\GAB$:
\BEN\label{LA-fg}
\mbox{\BM$\GS^{\prime}_{\mbox{\UB$\GA$}}$} =
  \tilde{\GAB} \times_{(\pi_1 (\GAB, [Z]),\rho^{\scriptscriptstyle{la}}_{p^{\prime}_2})} {\bf sl} (\HO{Z^{\prime}}))\,.
\EEN
The same representation will give the description of the fibre bundle of groups
$\LIG$ in (\ref{LGpr}):
\BEN\label{Lg-fg}
\LIG = 
\tilde{\GAB} \times_{(\pi_1 (\GAB, [Z]),\rho^{\scriptscriptstyle{la}}_{p^{\prime}_2})} {\bf SL} (\HO{Z^{\prime}}))\,.
\EEN

We can now identify  the Infinite Grassmannian $\IG$ as a fibre bundle induced by representation of the fundamental group of $\GAB$ as well.
\begin{thm}\label{th-IG-fg}
Let ${\bf Gr}_{[Z]}$ be the Infinite Grassmannian of 
${\bf SL} (\HO{Z^{\prime}}))$. Then the representation 
$\rho^{\scriptscriptstyle{la}}_{p^{\prime}_2}$ in (\ref{fg-rep-La}) determines the representation
$$
\rho^{\scriptscriptstyle{loop}}_{p^{\prime}_2} : \pi_1 (\GAB, [Z]) \longrightarrow Aut ( {\bf Gr}_{[Z]})
$$
which gives rise to the following identification
\BEN\label{IG-fg}
\IG = \tilde{\GAB} \times_{(\pi_1 (\GAB, [Z]),\rho^{\scriptscriptstyle{loop}}_{p^{\prime}_2})} {\bf Gr}_{[Z]}\,.
\EEN
\end{thm}
\begin{pf}
The representation $\rho^{\scriptscriptstyle{la}}_{p^{\prime}_2}$ in (\ref{LA-fg})
gives representation
\BEN\label{k-rep-fg}
\rho^{\scriptscriptstyle{\BK}}_{p^{\prime}_2} : \pi_1 (\GAB, [Z]) \longrightarrow Aut \left({\bf SL} (\HO{Z^{\prime}})) (\BK)\right)\,.
\EEN
Replacing $\BK$ by $\BO$ yields the representation 
\BEN\label{o-rep-fg}
\rho^{\scriptscriptstyle{\BO}}_{p^{\prime}_2} : \pi_1 (\GAB, [Z]) \longrightarrow Aut \left({\bf SL} (\HO{Z^{\prime}})) (\BO)\right)\,.
\EEN
These representation together with the description of $\LIG$ in (\ref{Lg-fg}) give the following
\begin{eqnarray}\label{loopLg-fg}
\LGR (\BK) &= &
\tilde{\GAB} \times_{(\pi_1 (\GAB, [Z]),\rho^{\scriptscriptstyle{\BK}}_{p^{\prime}_2})} {\bf SL} (\HO{Z^{\prime}}))(\BK)\,, \\ \nonumber
\LGR (\BO) &= &
\tilde{\GAB} \times_{(\pi_1 (\GAB, [Z]),\rho^{\scriptscriptstyle{\BO}}_{p^{\prime}_2})} {\bf SL} (\HO{Z^{\prime}}))(\BO)\,. \\ \nonumber
\end{eqnarray}
Viewing the Infinite Grassmannian $\IG$ (resp. ${\bf Gr}_{[Z]}$) as the quotient in (\ref{IGr}) \linebreak
(resp.
${\bf Gr}_{[Z]} = {\bf SL} (\HO{Z^{\prime}})) (\BK) / {{\bf SL} (\HO{Z^{\prime}})) (\BO)}$ ), we obtain 
$$
\IG = \tilde{\GAB} \times_{(\pi_1 (\GAB, [Z]),\rho^{\scriptscriptstyle{loop}}_{p^{\prime}_2})} {\bf Gr}_{[Z]}\,,
$$
where the representation
$$
\rho^{\scriptscriptstyle{loop}}_{p^{\prime}_2} : \pi_1 (\GAB, [Z]) \longrightarrow Aut ( {\bf Gr}_{[Z]})
$$
is the representation induced by representations in 
(\ref{k-rep-fg}) and (\ref{o-rep-fg}) on the quotient
$ {\bf Gr}_{[Z]} ={\bf SL} (\HO{Z^{\prime}})) (\BK) / {{\bf SL} (\HO{Z^{\prime}})) (\BO)}$.
\end{pf}

The description of $\IG$ given in Theorem \ref{th-IG-fg} completely determines the structure of 
$\LGR(\BO)$-orbits of $\IG$.
\begin{cor}\label{orbits-IG}
Let $[Z]$ be a point in $\GAB$ as in Theorem \ref{th-IG-fg} and let 
$\mbox{\BM$ W$}_{\GA,[Z]}$ and $\mbox{\BM$\check{\Lambda}$}_{\GA,[Z]}$ be 
the Weyl group and the coweight lattice of 
${\bf sl} (\HO{Z^{\prime}}) )$ respectively.

Let $\check{\lambda} \in \mbox{\BM$\check{\Lambda}$}_{\GA,[Z]}$ be a coweight
of 
${\bf sl} (\HO{Z^{\prime}}) )$ and let
$[\check{\lambda} ]$ be its $\mbox{\BM$ W$}_{\GA,[Z]}$-orbit in 
$\mbox{\BM$\check{\Lambda}$}_{\GA,[Z]}$.
Set 
$\mbox{\BM$O$}_{[\check{\lambda} ], [Z]}$ to be the
${\bf SL} (\HO{Z^{\prime}}) )(\BO)$-orbit in $Gr_{[Z]}$ corresponding to $[\check{\lambda} ]$. Then
\BEN\label{orbit-IG-fg}
\mbox{\BM$O$}_{[\check{\lambda} ], \GA} =
\tilde{\GAB} \times_{(\pi_1 (\GAB, [Z]),\rho^{\scriptscriptstyle{loop}}_{p^{\prime}_2})} \mbox{\BM$O$}_{[\check{\lambda} ],[Z]}
\EEN
is an $\LGR(\BO)$-orbit of $\IG$. Furthermore, every 
$\LGR(\BO)$-orbit of $\IG$ arises in this way. In particular, one has a bijective correspondence between
the set of $\LGR(\BO)$-orbits of $\IG$ and the cosets of the quotient
 $\mbox{\BM$\check{\Lambda}$}_{\GA,[Z]} / {\mbox{\BM$ W$}_{\GA,[Z]}}$.
\end{cor}
\begin{pf}
The representation $\rho^{\scriptscriptstyle{loop}}_{p^{\prime}_2}$ in Theorem \ref{th-IG-fg} induces the representation
$$
\rho^{\scriptscriptstyle{loop}}_{p^{\prime}_2} : \pi_1 (\GAB, [Z]) \longrightarrow Aut ( \mbox{\BM$O$}_{[\check{\lambda} ], [Z]})
$$
thus giving sense to the right hand side in (\ref{orbit-IG-fg}). The fact that it is 
an $\LGR(\BO)$-orbit of $\IG$ follows from the description of $\LGR(\BO)$ in (\ref{loopLg-fg}).

It is  well-known (see \cite{[P-S]}, Ch8) that $\mbox{\BM$O$}_{[\check{\lambda} ], [Z]}$,
as $[\check{\lambda} ]$ runs through the set  $\mbox{\BM$\check{\Lambda}$}_{\GA,[Z]} / {\mbox{\BM$ W$}_{\GA,[Z]}}$,
form a complete set of ${\bf SL} (\HO{Z^{\prime}}) )(\BO)$-orbits in $Gr_{[Z]}$. Hence the last two assertions of the corollary.
\end{pf}

\subsection{Relating $\JABG$ and $\IG$}
In this subsection we construct a `loop' analogue of the morphism $d^{+}$ in (\ref{d+}). For this we recall the linear stratification
of $\TPI$ in (\ref{linear-strat}) and set
\BEN\label{To}
\stackrel{\circ}{T}_{\pi} = \TPI \setminus T^{(\LG-2)}_{\pi}\,. 
\EEN

Next we show how to go from $\stackrel{\circ}{T}_{\pi}$ to the points of the Infinite Grassmannian $\IG$.
\begin{pro}\label{d+loop}
There is a natural map of $\GAB$-schemes
$$
Ld^{+}\! : \, \stackrel{\circ}{T}_{\pi} \longrightarrow \IG
$$
such that the following holds.
\begin{enumerate}
\item[1)]
On every stratum $\stackrel{\circ}{T^{\lambda}_{\pi}} = \stackrel{\circ}{T}_{\pi} \bigcap T^{\lambda}_{\pi}$, where
$T^{\lambda}_{\pi}$ is as in Proposition \ref{pro-strat1}, the map $Ld^{+}$ is a morphism and it is constant along the fibres of the projection
$\tau: \TPI \longrightarrow \JABG$.
\item[2)]
 There exists a positive integer $n$ such that the image of 
$Ld^{+} $ is contained in the stratum
$\IG(n\LG)$ of the stratification of $\IG$ in (\ref{strat-IG}).
\end{enumerate}
\end{pro}
\begin{pf}
Both $\TPI$ and $\IG$ are fibre bundles over $\GAB$ so it will be enough to show that there is a natural map on the fibres
over $\GAB$.

Let $[Z] \in \GAB$ and let $\JAB_Z$ be the fibre of $\JABG$ over $[Z]$.
Denote by
$T_{\JAB_Z}$ the tangent bundle of $\JAB_Z$. This is the fibre over $[Z]$ of the natural projection
$$
\tilde{\pi} : \TPI \longrightarrow \GAB
$$
(see (\ref{prd+-}) for notation). Denote by $\stackrel{\circ}{T}_{\JAB_Z}$ the intersection of $T_{\JAB_Z}$ with $\stackrel{\circ}{T}_{\pi}$. Our objective is to construct a map
$$
Ld^{+}([Z]) :\,  \stackrel{\circ}{T}_{\JAB_Z} \longrightarrow \IG ([Z])\,,
$$
where $\IG([Z])$ is the Infinite Grassmannian of 
${\bf SL} (\HO{Z^{\prime}}))$, the fibre of $\IG$ over $[Z]$.
Using the description of 
$\IG([Z])$ as a subset of 
$\mbox{\BM$Grass$}\left({\bf sl} (\HO{Z^{\prime}}))(\BK)\right)$ (see \S\ref{sec-IG}), our strategy is to define a map of 
$\stackrel{\circ}{T}_{\JAB_Z}$
into 
$\mbox{\BM$Grass$}\left({\bf sl} (\HO{Z^{\prime}}))(\BK)\right)$
and then to check that it lands into the subset
$\mbox{\BM$Grass^{0}$}\left({\bf sl} (\HO{Z^{\prime}}))(\BK)\right)$.

Let $v$ be a tangent vector of $\JAB_Z$ at a point $[\alpha] \in \JAB_Z$ such that $d^{+}_p (v) \neq 0$, for all $p\in [0,\LG-2]$, i.e.
we are at the point $([Z],[\alpha], v)$ of $\stackrel{\circ}{T}_{\pi} $.
Identifying
${\bf sl} (\HO{Z^{\prime}})) (\BK) $ with\footnote{recall from  (\ref{LAG=sl}):
{\BM
$\LAG = \mbox{\UB$\pi^{\ast} {\bf sl}(\FF^{\prime})$}$}.} 
{\BM
$\LAG$}$\ZA (\BK)$, 
 we define 
$Ld^{+} ([Z]) ([\alpha], v): ={\cal L} ( [Z],[\alpha], v)$, the value of $Ld^{+}$ at $([Z],[\alpha], v)$,
 as the $\BO$-submodule of {\BM
$\LAG$}$\ZA (\BK)$ generated, as a Lie algebra,   
by
$\mbox{\BM$\GS^0_{\GA}$} \ZA$ 
and 
$D^{\pm} (v,t)$, where
\BEN\label{dv-t}
D^{+} (v,t) = \sum^{\LG-2}_{p=0} t^{ a_p} d^{+}_p (v)\,\,\, and\,\,D^{-} (v,t) = \sum^{\LG-2}_{p=0} t^{- a_p} d^{-}_{p+1} (v)\,,
\EEN
with the exponents $a_p$ defined by
\BEN\label{exp}
a_p =tr(h_{p+1} (v)) -tr(h_{p} (v))\,,
\EEN
where $h(v)$ is a semisimple element of $\mbox{\BM$\GS^0_{\GA}$} \ZA$ which comes along with
$d^{+} (v)$ in an ${\bf sl_2}$-triple\footnote{see \S\ref{sec-sl2} for details about ${\bf sl_2}$-triples associated to
$d^{+} (v)$.}  in 
{\BM
$\LAG$}$\ZA$.
 Observe that the (graded) trace of $h(v)$ depends on $d^{+} (v)$ only.

We need to verify now that ${\cal L} ( [Z],[\alpha], v)$ is in 
$\mbox{\BM$Grass^{0}$}\left({\bf sl} (\HO{Z^{\prime}}))(\BK)\right)$. This is done by giving an explicit description
of this submodule. The essential ingredient of this description is the decomposition of 
{\BM
$\LAG$}
obtained in \S\ref{sec-gr}, (\ref{LAG-roots}), which says that 
{\BM
$\LAG$}$\ZA$ can be written as the following direct sum of 
$\mbox{\BM$\GS^0_{\GA}$} \ZA$-submodules
\BEN\label{g0-subm}
\mbox{\BM$\GS_{\GA}$} \ZA =\mbox{\BM$\GS^0_{\GA}$} \ZA \oplus \left( \bigoplus_{\nu \in R_{\LG}} \mbox{\BM${\cal W}$}_{\nu} \ZA \right)\,,
\EEN
where $R_{\LG}$ is the set of roots of ${\bf sl}_{\LG}$ as in (\ref{sl-lg-roots}), and the summands 
$\mbox{\BM${\cal W}$}_{\nu} \ZA)$ are the fibres at $\ZA$ of the sheaves $\mbox{\BM${\cal W}$}_{\nu}$'s in (\ref{W-nu}).

Recall that the set of roots $R_{\LG}$ comes with a preferred choice of simple roots
$\nu_p\, (p=0,\ldots ,\LG-2)$ (see Remark \ref{quiver-sl}). Now we use the exponents $a_p$ defined in (\ref{exp}) to associate to
our tangent vector $v$ 
the coroot $\check{v}$ of ${\bf sl}_{\LG}$ defined as follows
\BEN\label{cowei-v}
\nu_p (\check{v}) =a_p\,,
\EEN
for $p=0,\ldots ,\LG-2$. This and Proposition \ref{pro-Liei-s} imply the following description of 
${\cal L} ( [Z],[\alpha], v)$:
\BEN\label{Ld+Z}
{\cal L} ( [Z],[\alpha], v) =\mbox{\BM$\GS^0_{\GA}$} \ZA (\BO) \oplus
 \left(\bigoplus_{\nu \in R_{\LG}} t^{\nu(\check{v})} \mbox{\BM${\cal W}$}_{\nu} \ZA (\BO) \right)\,.
\EEN
From this identity it follows immediately that ${\cal L} ( [Z],[\alpha], v)$ is closed under the bracket operation and is self-dual.
Hence ${\cal L} ( [Z],[\alpha], v) \in \IG([Z])$ and the assignment
\BEN\label{Zav-IG}
 \stackrel{\circ}{T}_{\pi} \ni ( [Z],[\alpha], v) \longmapsto {\cal L} ( [Z],[\alpha], v) \in \IG
\EEN
gives a well-defined map of $\stackrel{\circ}{T}_{\pi}$ into $\IG$. Furthermore, ${\cal L} ( [Z],[\alpha], v)$ depends on $v$ only through
the exponents $a_p$ in (\ref{exp}) and these are constant, for all $v$ in a stratum $T^{\lambda}_{\pi}$ as in 
Proposition \ref{pro-strat1}. Hence the map in (\ref{Zav-IG}) is constant along the fibres of the projection $\tau$ on every stratum
$\stackrel{\circ}{T^{\lambda}_{\pi}}$. Since the assignment in (\ref{Zav-IG}) varies holomorphically with respect to the parameters
$\ZA$, it follows that $Ld^{+}$ is a morphism of $\GAB$-schemes on every stratum $\stackrel{\circ}{T^{\lambda}_{\pi}}$. This completes the proof of the first part of the proposition.

Turning to the part 2), recall from Lemma \ref{Hp-dgr} that the weights 
$w(h(v))$ of $h(v)$ on $\HH^p \ZA$ are bounded as follows
$$
p-\LG+1 \leq w(h(v)) \leq p\,.
$$
Hence the trace of $h(v)$ on $\HH^p \ZA$ is subject to the following inequalities
$$
(p-\LG+1)h^p_{\GA} \leq tr(h_p (v)) \leq ph^p_{\GA}\,,
$$
where $h^p_{\GA} = rk (\HH^p)$. This implies
\begin{eqnarray*}
a_p &=tr(h_{p+1} (v)) - tr(h_p (v))  \leq  (p+1)h^{p+1}_{\GA} - (p-\LG+1)h^p_{\GA} =
\LG h^p_{\GA} + (p+1)(h^{p+1}_{\GA} - h^{p}_{\GA})\,, \\
a_p &=tr(h_{p+1} (v)) - tr(h_p (v)) \geq (p+2-\LG)h^{p+1}_{\GA} - ph^{p}_{\GA} = -(\LG-2)h^{p+1}_{\GA} + p(h^{p+1}_{\GA} - h^{p}_{\GA})\,.
\end{eqnarray*}
From these inequalities it follows 
\BEN\label{exp1}
-2h^{max}_{\GA} \LG \leq a_p \leq 2 h^{max}_{\GA} \LG\,,
\EEN
for $p=0,\ldots, \LG-2$, where
$$
h^{max}_{\GA}  =\max \left\{\left.h^{p}_{\GA} \right| p=0,\ldots, \LG-1 \right\}\,.
$$
This and the grading in (\ref{sLie-grad}) imply 
\BEN\label{IG(n)}
{\cal L} ( [Z],[\alpha], v)  \in \IG ([Z]) (2 h^{max}_{\GA} \LG^2)\,,
\EEN
for all $([Z],[\alpha], v) \in \stackrel{\circ}{T}_{\pi}$.
\end{pf}

The map $Ld^{+}$ together with the orbit structure of the Infinite Grassmannian imply the following.
\begin{pro}\label{LO(Gam)}
Let $\GA$ be an admissible, simple component of $\CS$. Then it determines a distinguished collection
$LO(\GA)$ of 
$\LGR(\BO)$-orbits of the Infinite Grassmannian $\IG$. These are $\LGR(\BO)$-orbits intersecting the image of the map
 $Ld^{+}$ in Proposition \ref{d+loop}. 
\end{pro}
\begin{pf}
By Proposition \ref{d+loop} the image of $Ld^{+}$ is contained in the stratum
$\IG(n\LG)$ which is a $\LGR(\BO)$-stable, finite dimensional algebraic subvariety of $\IG$.
If we choose a base point $[Z] \in \GAB$, then, using the results and notation of \S\ref{sec-IG-orbits},
this stratum can be written as follows
$$
 \IG(n\LG) = \tilde{\GAB} \times_{(\pi_1 (\GAB, [Z]),\rho^{\scriptscriptstyle{loop}}_{p^{\prime}_2})} {\bf Gr}_{[Z]} (n\LG)\,,
$$
where ${\bf Gr}_{[Z]} (n\LG)$ is the corresponding stratum of ${\bf Gr}_{[Z]}$, the Infinite Grassmannian
of \linebreak
${\bf SL} (\HO{Z^{\prime}}))$, and $\rho^{\scriptscriptstyle{loop}}_{p^{\prime}_2}$ is the representation of the fundamental group
 in Theorem \ref{th-IG-fg}. This implies that  $\IG(n\LG)$ consists of $\LGR(\BO)$-orbits induced by the representation
$\rho^{\scriptscriptstyle{loop}}_{p^{\prime}_2}$ from the \linebreak 
${\bf SL} (\HO{Z^{\prime}}))(\BO)$-orbits of 
${\bf Gr}_{[Z]} (n\LG)$. It is known that the latter set of orbits is finite.\footnote{see \cite{[P-S]}, Ch8, or \cite{[Gi]}, Proposition 1.2.2.}
 Hence $\IG(n\LG)$ consists of finitely many
$\LGR(\BO)$-orbits of $\IG$. The finite collection $LO(\GA)$ of the corollary consists of the orbits of 
$\IG(n\LG)$ having nonempty intersection with the image of the map $Ld^{+}$. In view of Corollary \ref{orbits-IG}
 this collection 
 can be described as follows:
\BEN\label{LOG}
LO(\GA) =\left\{ \left.[\check{\lambda}] \in \mbox{\BM$\check{\Lambda}$}_{\GA,[Z]} / {\mbox{\BM$ W$}_{\GA,[Z]}} \right|
  \mbox{\BM$O$}_{[\check{\lambda}] , \GA} \bigcap Ld^{+} (\stackrel{\circ}{T}_{\pi}) \neq \emptyset \right\}\,.
\EEN
\end{pf}

 We appeal now to the profound and beautiful geometric version of the Satake isomorphism proved by 
Ginzburg in \cite{[Gi]} and Mirkovi\v{c} and Vilonen in \cite{[M-V]}. In our case their result says that
the tensor category $P({\bf Gr})$ of 
${\bf SL}_{\DE} (\BO)$-equivariant perverse sheaves on the Infinite Grassmannian 
${\bf Gr}$ of ${\bf SL}_{\DE} (\CC)$ is equivalent to the tensor category
$Rep_{{\bf {}^L SL}_{\DE} (\CC)}$ of finite dimensional representations of the Langlands dual group
${\bf {}^L SL}_{\DE} (\CC)( = {\bf PGL}_{\DE} (\CC) )$ of ${\bf SL}_{\DE} (\CC)$. 
 Combining this with Proposition \ref{LO(Gam)}
we obtain the following result stated as Theorem \ref{LD} in the Introduction.
\begin{pro}\label{LD1}
Let ${\cal V}(X;L,d)$ be the collection of admissible simple components as in Theorem \ref{nilo1}.
Then every $\GA$ in ${\cal V}(X;L,d)$
determines the finite collection
${}^L R (\GA)$ of irreducible representations of the Langlands dual group
${\bf {}^L SL}_{\DE} (\CC) = {\bf PGL}_{\DE} (\CC) $. Furthermore, upon a choice of a base point
$[Z]$ in $\GA$ the set
${}^L R (\GA)$ is identified, via geometric Satake isomorphism, with the set
$LO(\GA)$ in Proposition \ref{LO(Gam)} (see (\ref{Satake-cor}) below for the precise relation).
\end{pro}
\begin{pf}
Once a base point $[Z] \in \GA$ is fixed, Corollary \ref{orbits-IG} says that the lattice of coweights
$\mbox{\BM$\check{\Lambda}$}_{\GA, [Z]}$ of 
${\bf sl} (\HO{Z^{\prime}}))$ factored out by the action of its Weyl group
$\mbox{\BM$ W$}_{\GA,[Z]}$ parametrizes 
the ${\bf SL} (\HO{Z^{\prime}}))(\BO)$-orbits of the Infinite Grassmannian
${\bf Gr}_{[Z]}$ of ${\bf SL} (\HO{Z^{\prime}}))$ as well as $\LGR(\BO)$-orbits of the Infinite Grassmannian $\IG$.

Let $\mbox{\BM$O$}_{[\check{\lambda}]}$ be the 
${\bf SL} (\HO{Z^{\prime}}))(\BO)$-orbit of ${\bf Gr}_{[Z]}$ corresponding
to the coset \linebreak
$[\check{\lambda}] \in   \mbox{\BM$\check{\Lambda}$}_{\GA, [Z]} /{ \mbox{\BM$ W$}_{\GA,[Z]}}$
of a coweight $\check{\lambda}$ in $\mbox{\BM$\check{\Lambda}$}_{\GA, [Z]}$.
Let  $IC(\mbox{\BM$O$}_{[\check{\lambda}]}, \underline{\CC})$ be the Intersection cohomology complex
 of Deligne-Goresky-MacPherson corresponding to the trivial
local system $\underline{\CC}$ on $\mbox{\BM$O$}_{[\check{\lambda}]}$. This is a complex supported on the closure
$\mbox{\BM$\overline{O}$}_{[\check{\lambda}]}$ of $\mbox{\BM$O$}_{[\check{\lambda}]}$ in 
${\bf Gr}_{[Z]}$. Denote by
$IC(\mbox{\BM$O$}_{[\check{\lambda}]})$ its extension by zero to the whole of ${\bf Gr}_{[Z]}$.
This is an ${\bf SL} (\HO{Z^{\prime}}))(\BO)$-equivariant perverse sheaf on ${\bf Gr}_{[Z]}$ and the geometric Satake isomorphism
sends it to the object denoted $V_{[\check{\lambda}]}$ of the category
$Rep_{{}^L \! {\bf SL} (\HO{Z^{\prime}}))}$ of finite dimensional representations of the Langlands dual group
${}^L \! {\bf SL} (\HO{Z^{\prime}}))$
of ${\bf SL} (\HO{Z^{\prime}}))$. Thus setting
\BEN\label{Satake-cor}
{}^L R (\GA) =\left\{\left. V_{[\check{\lambda}]} \right| [\check{\lambda}]\in LO(\GA) \right\}\,,
\EEN
where $LO(\GA)$ is as in Proposition \ref{LO(Gam)}, yields the assertion.
\end{pf}

The geometric Satake isomorphism of Ginzburg, Mirkovi\v{c} and Vilonen, also tells us that the functor sending the perverse sheaf
$IC(\mbox{\BM$O$}_{[\check{\lambda}]} )$ to the object $V_{[\check{\lambda}]}$ of the category $Rep_{{}^L \! {\bf SL} (\HO{Z^{\prime}}))}$
is the hypercohomology functor, i.e. $V_{[\check{\lambda}]}$ is isomorphic to
$\mbox{\BM$H^{\bullet}$} (IC(\mbox{\BM$O$}_{[\check{\lambda}]}) )$, the hypercohomology of the bounded constructible complex of sheaves 
$IC(\mbox{\BM$O$}_{[\check{\lambda}]})$. We will show now that in our relative version the functor produces local systems on $\GAB$, for 
every $\GA \in {\cal V}(X;L,d)$.

Let $\GA$ be a component of ${\cal V}(X;L,d)$. Fix a base point $[Z] \in \GAB$ and let 
$\check{\lambda}$ be a coweight in $\mbox{\BM$\check{\Lambda}$}_{\GA, [Z]}$. Then the description of the orbit
$\mbox{\BM$O$}_{[\check{\lambda}],\GA}$ in Corollary \ref{orbits-IG}, (\ref{orbit-IG-fg}), implies that 
$\mbox{\BM$H^{\bullet}$} (IC(\mbox{\BM$O$}_{[\check{\lambda}],[Z]}) )$ is a $\pi_1 (\GAB, [Z])$-module. Namely, the representation
$\rho^{\scriptscriptstyle{loop}}_{p^{\prime}_2}$ in (\ref{orbit-IG-fg}) induces the representation
\BEN\label{H-orb-fg}
\mbox{\BM$H^{\bullet}$} (\rho^{\scriptscriptstyle{loop}}_{p^{\prime}_2}) : \pi_1 (\GAB, [Z]) \longrightarrow 
Aut \left(\mbox{\BM$H^{\bullet}$} (IC(\mbox{\BM$O$}_{[\check{\lambda}],[Z]}) )\right)\,.
\EEN
This representation gives rise to the local system on $\GAB$ which will be denoted ${\cal L}^{\bullet}_{\GA,[\check{\lambda}]}$.

Let $IC (\GAB,{\cal L}^{\bullet}_{\GA,[\check{\lambda}]})$ be the Intersection cohomology complex of Deligne-Goresky-MacPherson and let
${\cal P}^{\bullet}_{\GA,[\check{\lambda}]}$ be its extension by zero to the entire Hilbert scheme. Thus for every $\GA$ in 
${\cal V}(X;L,d)$
we obtained a distinguished collection
\BEN\label{cP-G-collec} 
\mbox{\BM$\check{\cal P}$}(\GA) =\left\{\left.{\cal P}^{\bullet}_{\GA,[\check{\lambda}]} \right| [\check{\lambda}] \in LO(\GA) \right\}
\EEN
of perverse sheaves on $\XD$ parametrized by the set $LO(\GA)$ in (\ref{LOG}).

Taking the union of these collections as $\GA$ runs through ${\cal V}(X;L,d)$, we obtain the collection
\BEN\label{cP-collec}
\mbox{\BM$\check{\cal P}$}(X;L,d) = \bigcup_{\GA \in {\cal V}(X;L,d)} \mbox{\BM$\check{\cal P}$}(\GA)\,.
\EEN

Repeating the same construction as in \S\ref{A(XLd)}, yields the collection
$\mbox{\BM$\check{C}$} (X;L,d)$ of irreducible perverse sheaves $\check{\cal C}_{\GA,[\check{\lambda}],\zeta}$ on $\XD$, indexed by
$\GA \in {\cal V}(X;L,d),\,[\check{\lambda}] \in LO(\GA)$ and irreducible $\pi_1 (\GAB, [Z])$-submodules $\zeta$ of
$\mbox{\BM$H^{\bullet}$} (IC(\mbox{\BM$O$}_{[\check{\lambda}],[Z]}) )$ in the representation 
$\mbox{\BM$H^{\bullet}$} (\rho^{\scriptscriptstyle{loop}}_{p^{\prime}_2})$ in (\ref{H-orb-fg}).

Similar to \S\ref{A(XLd)}, we use the collection $\mbox{\BM$\check{C}$} (X;L,d)$ to define the abelian category
$\check{\cal A} (X;L,d)$. This is the full abelian subcategory of 
${\cal D}^b_c (\XD)$ whose objects are isomorphic to finite direct sums of complexes of the form
$\check{\cal C}[n]$, where $\check{\cal C} \in \mbox{\BM$\check{C}$} (X;L,d)$ and $n \in {\bf Z}$.
\begin{defi}\label{non-ab-cowei}
\begin{enumerate}
\item[1)]
The perverse sheaves in $\mbox{\BM$\check{C}$} (X;L,d)$ are called irreducible non-abelian coweights of $\JA$.
\item[2)]
The abelian category $\check{\cal A} (X;L,d)$ is called the category of non-abelian coweights of $\JA$.
\end{enumerate}
\end{defi}

As in Theorem \ref{thm-cot-ps-map} we can relate $\JA$ with $\check{\cal A} (X;L,d)$
via the ``exponential integration" of relative $1$-forms on $\stackrel{\circ}{\JAA}(X;L,d) = \JA \setminus \mbox{\BM$\Theta$}(X;L,d)$.
\begin{thm}\label{thm-cot-cps-map}
Let $ \check{\cal A} (X;L,d)$ be the category of non-abelian coweights of $\JA$ (see Definition \ref{non-ab-cowei}). Then there is a natural map
\BEN\label{cot-cps-map}
\check{exp}\left(\int\right): H^0 ({\cal T}^{\ast}_{\stackrel{\circ}{\JAA}(X;L,d) / {\XD}} ) \longrightarrow \check{\cal A} (X;L,d)\,,
\EEN
where ${\cal T}^{\ast}_{\stackrel{\circ}{\JAA}(X;L,d) / {\XD}}$ is the relative cotangent sheaf of $\stackrel{\circ}{\JAA}(X;L,d)$ over $\XD$.
\end{thm}
\begin{pf}
Let $\omega$ be a global section of ${\cal T}^{\ast}_{\stackrel{\circ}{\JAA}(X;L,d) / {\XD}}$. For every $\GA$ in 
${\cal V}(X;L,d)$, let $\omega_{\GA}$ be the restriction of $\omega$ to $\JABG$. By Remark \ref{analogy1}, (\ref{tan-cotan}), we have a natural
isomorphism
$$
H^0 ({\cal T}^{\ast}_{\JABG /{\GAB}} ) \cong H^0 ({\cal T}_{\JABG /{\GAB}} )\,.
$$
Let $\theta_{\GA}$ be the section of ${\cal T}_{\JABG /{\GAB}}$ corresponding to $\omega_{\GA}$ under the above isomorphism.
We view $\theta_{\GA}$ as the corresponding morphism
\BEN\label{theta-mor}
\theta_{\GA} : \JABG \longrightarrow T_{\JABG /{\GAB}}
\EEN
of $\GAB$-schemes. At this stage we recall that the map $Ld^{+}$ in Proposition \ref{d+loop} is defined only on the complement
of the linear stratum $T^{(\LG-2)}_{\JABG /{\GAB}}$ in (\ref{linear-strat}). So we distinguish two cases.

{\bf Case 1:} $\theta_{\GA} ( \JABG ) \subset T^{(\LG-2)}_{\JABG /{\GAB}}$. 
\\
\noindent
In this case we send $\omega_{\GA}$ to the zero object of
$\check{\cal A} (X;L,d)$, i.e. we define
\BEN\label{c-exp0}
\check{exp}\left(\int_{\GA}\right) (\omega_{\GA}) = 0\,.
\EEN

{\bf Case 2:} $\theta_{\GA} ( \JABG ) \not\subset T^{(\LG-2)}_{\JABG /{\GAB}}$. 
\\
\noindent
In this case there is non-empty Zariski open subset
${\bf U}$ of $\JABG$ such that $ \theta_{\GA} ({\bf U})$ is contained in $\stackrel{\circ}{T}_{\JABG /{\GAB}}$, the complement
of $T^{(\LG-2)}_{\JABG /{\GAB}}$ in $T_{\JABG /{\GAB}}$. Composing $\theta_{\GA}$ with $Ld^{+}$ yields
the map 
$$
Ld^{+} \circ \theta_{\GA} : {\bf U}\longrightarrow \IG\,.
$$

Set $LO (\GA,\theta_{\GA})$ to be the collection of orbits of $\IG$ intersecting 
$Ld^{+} \circ \theta_{\GA} ( {\bf U})$. This is a subset of $LO (\GA)$ in Proposition \ref{LO(Gam)}. Using the identification in
(\ref{LOG}) we set
\BEN\label{c-exp1}
\check{exp}\left(\int_{\GA}\right) (\omega_{\GA}) = \bigoplus_{[\check{\lambda}] \in LO (\GA,\theta_{\GA})} {\cal P}^{\bullet}_{\GA,[\check{\lambda}]}\,,
\EEN
where ${\cal P}^{\bullet}_{\GA,[\check{\lambda}]}$ are perverse sheaves from the the collection $\mbox{\BM$\check{\cal P}$}(\GA)$ in
(\ref{cP-G-collec}).
Combining the definitions  (\ref{c-exp0}) and (\ref{c-exp1}) we define
$$
\check{exp}\left(\int\right) (\omega) = \bigoplus_{\GA \in {\cal V}(X;L,d)} \check{exp}\left(\int_{\GA}\right) (\omega_{\GA})
$$
\end{pf}

It is plausible that the categories ${\cal A} (X;L,d)$ and $\check{\cal A} (X;L,d)$ are functorially related.
Such hypothetical functors could be viewed as a manifestation of the Langlands duality for surfaces. This and other
issues concerning these categories will be treated elsewhere.

\vspace{1cm}
\begin{flushright}    
Universit\'e d'Angers\\
D\'epartement de Math\'ematiques
\\
2, boulevard Lavoisier\\
49045 ANGERS Cedex 01 \\
FRANCE\\
{\em{E-mail addres:}} reider@univ-angers.fr
\end{flushright}

 \end{document}